%% file: main.tex
\documentclass[a4paper, 10pt, english, reqno, dvipsnames]{amsart}
\usepackage[utf8]{inputenc}

\usepackage{amsmath}
\usepackage{amsthm}
\usepackage{amssymb}
\usepackage{mathtools}  
\usepackage{mathrsfs}   
\usepackage{amsfonts}
\usepackage{mathdots}   
\usepackage[bb=libus]{mathalpha} 
\usepackage{amsbsy}
\usepackage{tikz}
\usepackage{tikz-cd}
\usepackage{quiver}

\usepackage{geometry}
\usepackage{lmodern}

\theoremstyle{plain}
\newtheorem{theorem}{Theorem}[section]

\newtheorem{proposition}[theorem]{Proposition}
\newtheorem{lemma}[theorem]{Lemma}
\newtheorem{corollary}[theorem]{Corollary}

\newtheorem{definition-theorem}[theorem]{Definition-Theorem}
\newtheorem{definition-proposition}[theorem]{Definition-Proposition}

\theoremstyle{definition}

\newtheorem{definition}[theorem]{Definition}

\theoremstyle{remark}
\newtheorem{remark}[theorem]{Remark}

\usepackage{paralist}   
\usepackage{enumerate}   
\usepackage{enumitem}
\usepackage[figurewithin=none]{caption}
\usepackage{graphicx}
\usepackage{booktabs} 
\usepackage{tabularx} 
\usepackage{array}    
\usepackage{float}
\usepackage{xstring}
\usepackage{array}
\usepackage{longtable}

\usepackage{xcolor}

\usepackage{mycom}
\usepackage{myarrows}

\makeatletter
\def\l@subsection{\@tocline{2}{0pt}{2pc}{5pc}{}}
\makeatother

\usepackage[style=alphabetic,sorting=nyt, backend=biber, backref=true]{biblatex}
\usepackage{fancyvrb}

\DefineVerbatimEnvironment{codebox}{Verbatim}{
    frame=single,
    numbers=left,
    numbersep=5pt,
    fontsize=\small
}

\usepackage{url}

\usepackage[unicode, colorlinks, linktocpage=true]{hyperref}
\hypersetup{
	linkcolor = blue,
	citecolor = Red,
	urlcolor = PineGreen,
	bookmarksnumbered = true
}

\title[The BSD formula for rank-one CM elliptic curves]{The Birch and Swinnerton-Dyer formula for CM elliptic curves of analytic rank one}
\date{\today}
\usepackage{orcidlink}
\author{Xin Wan, Zerui Xiang, Ruichen Xu \,\orcidlink{0009-0003-4555-2116}}

\address[Xin Wan]{Academy of Mathematics and Systems Science, Chinese Academy of Sciences and University of Chinese Academy of Sciences, Haidian District, Beijing, China}
\email{xwan@math.ac.cn}

\address[Zerui Xiang]{Academy of Mathematics and Systems Science, Chinese Academy of Sciences, No. 55 Zhongguancun East Road, Beijing, 100190, China.}
\email{xiangzerui@amss.ac.cn}

\address[Ruichen Xu]{Faculty of Mathematics, University of Vienna. Oskar-Morgenstern-Platz 1, 1090 Vienna, Austria.}
\email{ruichen.xu@univie.ac.at}

\begin{document}

\begin{abstract}
In this article, we prove the Birch and Swinnerton-Dyer formula for elliptic curves over the field of rational numbers with complex multiplication and analytic rank one. 
\end{abstract}

\maketitle
\tableofcontents

\input{101-Introduction}
\input{102-Setup}

\part{Transfer mechanism: small primes, split and ramified cases} \label{part:1}

In Part \ref{part:1}, we assume that $p$ does not split in $\calK$, i.e., $p$ is inert or ramified in $\calK$. We will treat the case where $p$ splits in $\calK$ in Section \ref{sec:split}.

\input{103-Plus-minus}
\input{104-GZW-comparison}
\input{105-magic}

\input{106-BSD-inert-2-3}
\input{107-BSD-ramified}

\input{108-split}

\part{Bertolini--Darmon--Prasanna main conjecture: general inert cases} \label{part:BDP}

\input{201-Setup}
\input{202-Selmer}
\input{203-analytic-comparison}
\input{204-main-conjecture}
\input{205-Waldspurger}
\input{206-final}
\input{207-Kolyvagin}

\appendix
\part{Appendices}
\input{AppA-Elliptic-units}
\input{AppBC-calculations}

\printbibliography
\end{document}

%% file: 101-Introduction.tex
\section{Introduction}

Let $E/\QQ$ be an elliptic curve, and let $r$ be the order of vanishing of $L(E,s)$ at $s=1$. The Birch and Swinnerton-Dyer (BSD) conjecture predicts that $r$ is the rank of $E(\QQ)$, that the Tate--Shafarevich group $\Sha(E/\QQ)$ is finite, and that
\begin{equation}\label{eq:intro-bsd}
\frac{L^{(r)}(E,1)}{r!\,\Omega_E\Reg(E/\QQ)}=\frac{\#\Sha(E/\QQ)\prod_q\Tam_q(E/\QQ)}{(\#E(\QQ)_{\mathrm{tor}})^2}.
\end{equation}
Here $\Omega_E$ is the integral of the absolute value of a N\'eron differential over $E(\mathbb{R})$, $\Reg(E/\QQ)$ is the N\'eron--Tate regulator, and $\Tam_q(E/\QQ)$ is the Tamagawa number at $q$. When $r=0$ or $1$, the work of Gross--Zagier and Kolyvagin gives the rank equality and the finiteness of $\Sha(E/\QQ)$ \cite{GrossZagier,Kolyvagin}. In rank one, the Gross--Zagier formula relates the derivative of a Rankin--Selberg $L$-function to the height of a Heegner point. In this paper, we prove the following result for elliptic curves $E$ with complex multiplication given by an imaginary quadratic field $\calK$.

\begin{theorem}\label{thm:intro-cm-bsd}
Let $E/\QQ$ be an elliptic curve with complex multiplication. If $E/\QQ$ has analytic rank one, i.e., $\ord_{s=1} L(E,s) = 1$, then the full Birch and Swinnerton-Dyer formula \eqref{eq:intro-bsd} holds for $E/\QQ$.
\end{theorem}

For CM elliptic curves over $\QQ$ of analytic rank zero, the full formula follows from the work of Burungale and Flach \cite{BurungaleFlach}, building on the work of Rubin and of Johnson-Leung and Kings on the main conjecture for elliptic units \cite{Rubin,jlk}. In analytic rank one, Li--Tian--Yan--Zhu prove the $p$-part of the formula for primes split in the CM field, assuming good reduction when $p=2$ \cite[Theorem~1.1]{LiTianYanZhu}. Their argument extends Perrin-Riou's use of the $p$-adic Gross--Zagier formula and draws on work of Coates at $2$. Our proof is different. We compare the complex Gross--Zagier formula in \cite{YZZ, CST} with a $p$-adic Waldspurger formula in \cite{FanWan}. One consequence of this comparison is a transfer principle stating that, at a prime $p$ such that $p$ does not split in $\calK$, the defect in the $p$-part of BSD depends only on the local CM character.

\subsection{Transfer principle}

We explain this transfer principle first. Write $\calK$ for the CM field of $E$, and let $\psi_{E}$ be the Hecke character over $\calK$ attached to $E$ with $L(E/\QQ,s) = L(\psi_E, s)$. By isogeny invariance of the BSD formula, we may work with a curve with complex multiplication by the maximal order $\mathcal{O}_{\calK}$. Let $\Delta_p(E/\QQ)$ be the $p$-adic valuation of the ratio of the two sides of \eqref{eq:intro-bsd}. Thus the $p$-part of BSD is the assertion that $\Delta_p(E/\QQ)=0$.

\begin{theorem}[Transfer principle]\label{thm:intro-cm-transfer}
For $j=1,2$, let $E_j/\QQ$ have complex multiplication by $\mathcal{O}_{\calK_j}$ and analytic rank one. Suppose that $p$ is nonsplit in both CM fields and that their completions at $p$ are isomorphic over $\QQ_p$. Under this identification, suppose that the local CM characters $\psi_{E_1,p}$ and $\psi_{E_2,p}$ are equal. Then $\Delta_p(E_1/\QQ)=\Delta_p(E_2/\QQ)$. In particular, the $p$-part of BSD holds for $E_1/\QQ$ if and only if it holds for $E_2/\QQ$.
\end{theorem}

To see how the condition that $\psi_{E_1,p} = \psi_{E_2,p}$ enters the proof, suppose that $p$ is nonsplit in $\calK$, put $F=\calK_p$, and fix a coefficient embedding $\sigma$ of $\mathcal{O}_F$. Let $\calO_{L}$ be the ring of integers of a sufficiently large finite extension $L/\QQ_p$ and use it as the coefficient ring. Put $R_p=\calO_{\calK}\otimes\ZZ_p$ and $T_E=T_pE$ as a rank-two $G_{\QQ}$-representation over $\ZZ_p$. The module $T_E$ is free of rank one over $R_p$. For a coefficient embedding $\sigma:R_p\to\calO_L$, define
\[
T_\sigma=T_E\otimes_{R_p,\sigma}\calO_L,\qquad
 V_\sigma=T_\sigma\otimes_{\calO_L}L,\qquad W_\sigma=V_\sigma/T_\sigma.
\]

\subsubsection{Arithmetic side}

On a small disc about the trivial character in anticyclotomic weight space over $\calK$, Fan and the first-named author \cite{FanWan} constructed a local plus-minus theory, decomposing the local cohomology $\rmH^{1}(F, V_{\sigma})$ into two rank-one summands, namely the plus part $\rmH^{1}_{+}(F, V_{\sigma})$ and the minus part $\rmH^{1}_{-}(F,V_{\sigma})$, such that the plus summand $\rmH^{1}_{+}(F, V_{\sigma})$ is the local Bloch--Kato cohomology group $\rmH^1_{\rmf}(F,V_\sigma)$. Projection along the minus summand identifies this line with $\rmH^1(F,V_\sigma)/\rmH^1_-(F,V_\sigma)$. The quotient carries an integral lattice, namely the image of $\rmH^1(F,T_\sigma)$, which we compare with the Bloch--Kato lattice $\rmH^{1}_{\rmf}(F, T_{\sigma})$. 

On the arithmetic side, the Iwasawa main conjecture for elliptic units by Rubin and Johnson-Leung-Kings, formulated by Burungale and Flach \cite{BurungaleFlach} in the form here, gives a generator of the global cohomology determinant $\Det^{-1}_{\calO_{L}} \RG(\calO_{\calK,S}, T_{\sigma})$. The index of its localization in the quotient by the minus line gives the order of the minus Selmer group together with the local and global torsion factors. Poitou--Tate duality then compares the minus Selmer condition and Bloch--Kato conditions at $p$, as in the rank-one argument of Jetchev--Skinner--Wan \cite{JSW}. Taking the product over the two conjugate coefficient embeddings recovers the $p$-part of the order of $\Sha(E/\calK)$ and of the Tamagawa factors $\Tam_{v}(E/\calK)$ for $v$ away from $p$. These arguments are carried out in Section \ref{sec:plus-minus}.

\subsubsection{Analytic side}
For the analytic comparison, we choose the following Rankin--Selberg data: Let \(f=f_\xi\) be the weight-two CM eigenform associated with a Hecke character \(\xi\) of \(\calK\). As in Section \ref{sec:GZW-comparison}, choose a finite-order Hecke character \(\chi_0\) such that \(\xi\chi_0=\psi_E\) and \(L(\xi\chi_0^c,1)\ne0\), and put \(\psi_0=\xi^c\chi_0\). The Rankin–Selberg \(L\)-function appearing in the Gross–Zagier formula then factors as
\[
L(f,\chi_0,s)=L(\xi\chi_0,s)L(\xi^c\chi_0,s)
=L(\psi_E,s)L(\psi_0,s).
\]
Conjugating the Hecke characters shows that \(L(\psi_0,1)=L(\xi\chi_0^c,1)\ne0\). Since \(L(\psi_E,1)=0\), we obtain $L^{\prime}(f,\chi_0,1)=L^{\prime}(\psi_E,1)L(\psi_0,1)$.

For two choices of test vectors, the following formulas apply:
\begin{itemize}
    \item the $p$-adic Waldspurger formula of \cite{FanWan} gives the squares of the logarithms of the associated Heegner points, and
    \item the Gross--Zagier formula of Yuan--Zhang--Zhang \cite{YZZ, CST} gives the heights obtained by pairing each point with the dual of the other.
\end{itemize}
The reason for taking two test vectors is that the local toric pairing has rank at most one. If $\calA_{ij,v}$ is its value on the $i$-th vector and the $j$-th dual vector at the trivial character, then
\[
\calA_{11,v}\calA_{22,v}=\calA_{12,v}\calA_{21,v}.
\]
The two Waldspurger formulas supply the diagonal terms, and the two Gross--Zagier formulas supply the crossed terms. Thus the finite toric pairings away from $p$ cancel without having to evaluate them individually. It turns out that the auxiliary rank-zero value cancels as well, leaving the toric factor at $p$ and the period coming from the N\'eron differential. Multiplying over the two embeddings cancels the ratios of conjugate Heegner coefficients (over a Mordell-Weil group generator). The localization index now occurs with the same exponent on both sides and drops out of the BSD defect. These arguments are carried out in Section \ref{sec:GZW-comparison}. It then remains only to compare a certain factor for the two elliptic curves, as we now explain.  

\subsubsection{Comparison of periods}
There remains the question of how to compare the periods for two curves. The trivial character lies outside the interpolation range for the rank-one branch, although the Waldspurger comparison gives a nonzero value there. For $E_j$ ($j=1,2$), divide the localized elliptic-unit class by the renormalized $p$-adic $L$-function, and write $w_{j,\sigma}(k)$ for its coordinate in a common local basis of the plus line at weight $k$. For matched curves we can choose the same test data at $p$. At positive interpolation weights the explicit reciprocity law for elliptic units and the interpolation formula give
\[
v_p\!\left(\frac{w_{1,\sigma}(k)}{w_{2,\sigma}(k)}\right)=(2k-1)v_p\!\left(\frac{\delta'_{2,\sigma}}{\delta'_{1,\sigma}}\right),
\]
where $\delta'_{j,\sigma}$ is the period factor obtained by comparing the (Archimedean and $p$-adic) CM period with the logarithm of an integral generator of the Bloch--Kato line. Take $k=Np^a$ on the interpolation branch. These weights tend $p$-adically to zero, and the quotient on the left is analytic and nonzero at zero. Its valuation is therefore eventually constant. Since the $k$'s are unbounded as integers, the valuation on the right must vanish. It follows that the two coordinates have the same valuation at the origin as well. This proves equality of the defects of $E_1$ and $E_2$ without requiring an absolute evaluation of either period factor. This argument is carried out in Section \ref{sec:magic}. 

\subsubsection{Transfer principle in practice: numerical experiments}
We apply the transfer principle at the inert primes $2$ and $3$ and at every ramified prime. Since $E$ is defined over $\QQ$, its CM field has class number one. In these cases, only finitely many local CM characters occur, and for each of them we exhibit a curve of analytic rank one. Sections \ref{sec:inert BSD 2,3} and \ref{sec:ramified BSD} give the representatives, and Appendix \ref{sec:special values calculation} computes their local characters. The analytic BSD quotients for these curves are taken from the LMFDB, with a separate numerical computation in Appendix \ref{sec:p-Sha calculation} for the representative outside the database. On the arithmetic side we use $2$- and $3$-descent, executed by Magma, at the small primes and the $p$-isogeny descent of Miller--Stoll \cite{miller2013explicit} at ramified primes $p>3$.

\subsubsection{The case $p$ splits in $\calK$} 

When $p$ splits in $\calK$, the local plus-minus theory of \cite{FanWan} directly gives the two local summands associated with the two places above $p$. We apply the two-variable regulator to the elliptic-unit class and compare the resulting $p$-adic $L$-function with the renormalized anticyclotomic $p$-adic $L$-function. In this case the toric factor at $p$ can be computed explicitly. For a ramified local character, its Gauss sum pairs with the one in the regulator interpolation formula. We use the interpolation formulas of Loeffler--Zerbes \cite{LZ-ijnt} at odd primes and of Sano--Venjakob \cite{Venjakob-regulator} at $2$. Comparing at arithmetic weights tending to zero then gives $\Delta_p(E/\calK)=0$ directly. This includes the case where $2$ splits and $E$ has bad reduction at $2$.

\subsection{The Bertolini--Darmon--Prasanna main conjecture for $p$ inert in $\calK$}
Finally, suppose that $p\geq 5$ is inert in $\calK$. We obtain the lower and upper bounds for the order of $\Sha(E/\QQ)[p^\infty]$ separately. 

For the lower bound, following an idea of Ye Tian, we choose an auxiliary imaginary quadratic field $\calL$ in which $p$ and all primes dividing the conductor of $E$ split, and for which the quadratic twist of $E$ has analytic rank zero. Put $\calM=\calK\calL$, and let $\calF$ be its real quadratic subfield. Although $p$ is inert in $\calK$, the quartic CM field $\calM$ admits a $p$-ordinary CM type induced from $\calL$ since $p$ is inert in $\calF$ and the unique prime of $\calF$ lying over $p$ splits in $\calM$. We can therefore use the divisibility in the CM main conjecture for $\calM/\calF$ proved by the third-named author in his thesis using Eisenstein congruences \cite{Xu} and independently in \cite{ABL}, extending work of Hsieh \cite{hsieh2014eisenstein}. We have the following descent theorems:
\begin{itemize}
    \item On the arithmetic side, essentially by Shapiro's lemma, the dual Greenberg Selmer group over the $\mathbb{Z}_p^3$-extension of $\calM$ specializes to the dual Bertolini--Darmon--Prasanna (BDP) Selmer group over the anticyclotomic $\mathbb{Z}_p$-extension of $\calL$.
    \item On the analytic side, the same specialization takes the Katz $p$-adic $L$-function to the BDP $p$-adic $L$-function up to a $\Lambda$-adic unit. 
\end{itemize}
It follows that the characteristic ideal of the dual BDP Selmer group is contained in the ideal generated by the BDP $p$-adic $L$-function.

Let $P\in E(\calL)$ be the Heegner point attached to an optimal modular parametrization $\phi: X_{0}(N) \rightarrow E$, let $m$ be its index in $E(\calL)$ modulo torsion, and let $c_\phi$ be the Manin constant. Specializing the divisibility of the BDP main conjecture at the trivial character and applying the $p$-adic Waldspurger formula gives the lower bound in
\begin{equation} \tag{lower-bound} \label{eq:lower-bound}
2v_p(m/c_\phi)\ \leq\ v_p\bigl(\#\Sha(E/\calL)[p^\infty]\bigr)
\end{equation}
For good reduction the required Waldspurger formula is that of Bertolini--Darmon--Prasanna \cite{BDP}. For bad reduction we apply the formula of Liu--Zhang--Zhang \cite{LZZ} to the test vector defining our BDP $p$-adic $L$-function. Note that the local Tamagawa numbers are $p$-adic units since we are in the CM case and $p\geq 5$. 

The upper bound is Kolyvagin's \cite[Corollary~13]{kolyvagin2007euler}, in the formulation of Matar and Nekov\'a\v{r} in \cite{matar2019kolyvagin}, which gives
\begin{equation} \tag{upper-bound} \label{eq:upper-bound}
v_p\bigl(\#\Sha(E/\calL)[p^\infty]\bigr)\ \leq\ 2v_p(m).
\end{equation}

For $p\geq 11$, Edixhoven's theorem \cite[Theorem~3]{edixhoven1991manin} implies that the Manin constant $c_\phi$ is a $p$-adic unit, since $E$ has potentially supersingular reduction at an inert prime. The bounds in \eqref{eq:lower-bound} and \eqref{eq:upper-bound} therefore agree. For $p=5$ and $7$, we exhibit optimal curves of analytic rank one with Manin constant one representing each of the six local characters at $5$ and the four at $7$. The bounds \eqref{eq:lower-bound} and \eqref{eq:upper-bound} agree for these representatives, and the transfer principle gives the $p$-part of BSD for every curve with the same local character. This shows that
\begin{equation} \label{eq:equality-for-sha}
v_{p}(\# \Sha(E/\calL)[p^{\infty}]) = 2 v_p(m/c_\phi)
\end{equation}
for $p \geq 5$ when $p$ remains inert in $\calK$. It then follows from the complex Gross--Zagier formula that equation \eqref{eq:equality-for-sha} is equivalent to the BSD conjecture for $E/\calL$. Together with the BSD conjecture in the rank-zero case for the quadratic twist $E^{\calL}$, proved in \cite{BurungaleFlach}, this completes the proof of Theorem~\ref{thm:intro-cm-bsd}.

\subsubsection*{Acknowledgements}
The work of X.W. and Z.X. was supported by \emph{NSFC grants 12288201, 11621061}, \emph{CAS Project for Young Scientists in Basic Research grant no. YSBR-033}, the \emph{Strategic Priority Research Program of the Chinese Academy of Sciences (Grant No. XDA0480503)} and the \emph{National Key R\&D Program of China (Grant No. 2023YFA1009701)}.

R.X. carried out much of this work during his postdoctoral appointment at the University of Vienna under the supervision of Professor Harald Grobner and the co-supervision of Professor Adel Betina, whom he sincerely thanks for their support and guidance. Part of R.X.'s contribution to this work was completed during his doctoral studies at the Academy of Mathematics and Systems Science, Chinese Academy of Sciences, under the supervision of Professor Xin Wan, and during his visit to the Laboratoire de mathématiques de Besançon (LmB) in November 2026, hosted by Professor Adel Betina. R.X. thanks both institutions for their stimulating research environments. R.X. would also like to thank Professors Adel Betina and Olivier Fouquet for their interest in this work and for helpful discussions during his visit to LmB. R.X.'s postdoctoral research is supported by the \emph{FWF Principal Investigator Project (PAT-2584625): Eigenvarieties and $p$-adic $L$-functions}.

\subsubsection*{AI Disclosure}
This project was begun about a year ago as a continuation of the work of Fan--Wan \cite{FanWan} and as an application of the thesis of the third-named author (Xu). In particular, the main ideas underlying the proofs were developed independently and were not influenced by recent releases of OpenAI models, which use completely different strategies. OpenAI's GPT-6 Astra was used to find the numerical examples and check the writing.

%% file: 102-Setup.tex
\section{CM representations and the \texorpdfstring{$p$}{p}-part of the BSD formula}\label{sec:setup}

\subsection{Coefficient rings and Galois representations}

Let $E/\QQ$ be an elliptic curve with complex multiplication by an order in the imaginary quadratic field $\calK$, and assume $\ord_{s=1}L(E/\QQ,s)=1$. Without loss of generality, we take $\End_{\overline\QQ}(E)=\calO_{\calK}$. Let $\psi := \psi_{E}$ be the associated Hecke character, with $\psi_{E,\infty}(z)=z^{-1}$. Fix a rational prime $p$ (which is allowed to be $2$) and embeddings
\[
 \overline{\QQ}\hookrightarrow\CC, \quad \iota:\overline{\QQ}\hookrightarrow\CC_p.
\]
Let $M$ be a CM field containing the values of all Hecke characters used below, let $L/\QQ_p$ be a finite extension containing the image of $M$ under $\iota$, and write $\calO_{L}$ for its ring of integers. Let $\varpi_{L}$ be a fixed uniformizer of $\calO_{L}$, and normalize the valuation $v_L$ on $L$ by $v_L(\varpi_L)=1$.

\subsubsection{Galois representations}
Put $R_p=\calO_{\calK}\otimes\ZZ_p$ and $T_E=T_pE$ as a rank-two $G_{\QQ}$-representation over $\ZZ_p$. The module $T_E$ is free of rank one over $R_p$. For a coefficient embedding $\sigma:R_p\to\calO_L$, define
\begin{equation}\label{eq:native-CM-lattice}
 T_\sigma=T_E\otimes_{R_p,\sigma}\calO_L,\qquad
 V_\sigma=T_\sigma\otimes_{\calO_L}L,\qquad W_\sigma=V_\sigma/T_\sigma.
\end{equation}
We denote by $\rho_{\sigma}: G_{\calK} \rightarrow \calO_{L}^{\times}$ the Galois character on $T_{\sigma}$. After fixing one $\sigma$, we shall write $T,V,W$ when no confusion is possible.  Let $S$ be a finite set of places of $\calK$ that contains the archimedean places, the places above $p$, and every ramified place of $T$.

We explain the duality notation:
\begin{itemize}
\item We put $T^\vee=\Hom_{\calO_L}(T,\calO_L)$ and $T^\dagger=T^\vee(1)$.
\item We put $V^\dagger=T^\dagger\otimes_{\calO_L}L$ and $W^\dagger=V^\dagger/T^\dagger$.
\item The Weil pairing on $E/\QQ$ gives $\rho_\sigma\rho_\sigma^c=\chi_{\cyc}$. Here $c$ on a representation means precomposition by conjugation on $G_{\calK}$. Conjugation on the base field identifies this representation with the other coefficient realization of $T_E$.
\end{itemize}
Local Tate duality pairs $\rmH^i(\calK_v,T)$ with $\rmH^{2-i}(\calK_v,W^\dagger)$ in $L/\calO_L$.

\subsubsection{Order ideal}
Suppose that $p$ does not split in $\calK$, and denote by $\frp$ the unique prime ideal of $\calO_{\calK}$ lying over $p$. In this case, we write $F=\calK_{\frp}$ and $R=\calO_F=R_p$.

For every finitely generated $\calO_{L}$-module $M$, set
\[
 M_{\tf}:=M/M_{\tor},\quad
 M_{\tor}:=\{m\in M: \varpi_L^n \cdot m=0\text{ for some }n\geq 0\}.
\]
For a cofinitely generated discrete $\calO_L$-module $D$, let $D_{\divi}$ be its maximal divisible submodule and put $D_{\fin}=D/D_{\divi}$. A divisible module attached to a free $\calO_L$-module $M$ is $M\otimes_{\calO_L}L/\calO_L$. For a fractional ideal $(a)$, put $v_L((a))=v_L(a)$ and $v_p((a))=v_L(a)/e(L/\QQ_p)$.

If $N$ is a finite $\calO_{L}$-module, its \emph{order ideal} is
\[
 \abs{N}_L:=\left(\varpi_L^{\length_{\calO_{L}}N}\right).
\]
Equalities of fractional $\calO_{L}$-ideals will always be interpreted after extension to a common coefficient field. We record the following lemma, which is easily verified.

\begin{lemma} \label{lem:order-norm}
For a finite $R$-module $N$, put $N_{\tau} := N \otimes_{R, \tau} \calO_{L}$ for $\tau \in \{\sigma, \sigma^{c}\}$. Then we have $\Fitt_{\calO_{L}}(N_\tau) =\tau \left(\Fitt_R(N)\right) \calO_{L}$ and $\Fitt_{\calO_{L}}(N_\sigma) \cdot \Fitt_{\calO_{L}}(N_{\sigma^c}) = (\#N) \calO_{L}$.
\end{lemma}

\subsubsection{The Bloch--Kato Selmer group}
For a place $v$ of $\calK$, $\rmH^{1}_{\rmf}(\calK_v,V)$ denotes the Bloch--Kato finite subspace.  At $v\nmid p$ this is the unramified subspace; at $v\mid p$ it is the kernel of the usual comparison map to $\rmH^{1}(\calK_v,V\otimes \BB_{\cris})$. The \emph{Bloch--Kato Selmer group} is
\[
 \rmH^{1}_{\rmf}(\calK,V)=\ker \left( \rmH^{1}(\calK_S/\calK,V)\rightarrow \prod_{v\in S}\frac{\rmH^{1}(\calK_v,V)}{\rmH^{1}_{\rmf}(\calK_v,V)} \right).
\]
We write
\[
 \calK_p:=\calK\otimes_{\QQ}\QQ_p=\prod_{\frP\mid p}\calK_{\frP},
 \quad
 \rmH^{i}(\calK_p,-):=\bigoplus_{\frP\mid p}\rmH^{i}(\calK_{\frP},-),
\]
and use the same direct-sum convention for the local Bloch--Kato subspaces. The integral finite condition is the inverse image
\[
 \rmH_{\rmf}^1(\calK_v,T)=\{x\in \rmH^1(\calK_v,T):x_L\in \rmH_{\rmf}^1(\calK_v,V)\}.
\]
The discrete finite condition $\rmH^{1}_{\rmf}(\calK_{v}, W)$ is the image of $\rmH_{\rmf}^1(\calK_v,V)$ in $\rmH^1(\calK_v,W)$.

\subsection{The Birch and Swinnerton-Dyer formula}

\subsubsection{Heights, regulators and periods}
We shall fix the height and period conventions. Let $\widehat h$ be the absolute Néron--Tate canonical height. We put
\begin{equation}\label{eq:height-convention}
 b_{\mathrm{abs}}(P,Q)=\frac{\widehat h(P+Q)-\widehat h(P)-\widehat h(Q)}2,
 \qquad b_{\calF}=[\calF:\QQ]b_{\mathrm{abs}}.
\end{equation}
The height pairing $b_{\mathrm{abs}}$ is called the \emph{absolute height} and the height pairing $b_{\calF}$ is called the \emph{height pairing over $\calF$}. On points over a finite extension $\calH/\calF$, $b_{\calF}$ means $[\calH:\calF]^{-1}b_{\calH}$. The regulator $\Reg(E/\calF)$ is then defined as the Gram determinant of $b_{\calF}$ in the free Mordell--Weil group $E(\calF)/E(\calF)_{\tor}$. 

We put $D_{\calK}=|\disc(\calK)|$. For a nonzero invariant differential $\omega$ for $E/\calK$ and local invariant differentials $\omega_v$ given by the Néron module of $E/\calK_v$, define the period of $E$ over $\calK$, following \cite[I, \S 7]{MilneADT}, by
\begin{equation}\label{eq:period-convention}
 \Omega_{E/\calK} :=D_{\calK}^{-1/2} \int_{E(\CC)}|\omega\wedge\bar\omega| \prod_{v<\infty}|\omega/\omega_v|_v, \qquad |\varpi_v|_v=(Nv)^{-1}.
\end{equation}

\begin{remark}
The definition of $\Omega_{E/\calK}$ above is independent of the choice of $\omega$. In fact, if $\omega$ is replaced by $a\omega$ for some $a\in \calK^\times$, the complex integral is multiplied by $|\sigma(a)|^2$, while the finite product is multiplied by $\prod_{v<\infty}|a|_v$. For an imaginary quadratic field, its sole complex-place absolute value is $|a|_\infty=|\sigma(a)|^2$. Hence
\[
\frac{\Omega_{E/K}(a\omega)}{\Omega_{E/K}(\omega)}
=|\sigma(a)|^2\prod_{v<\infty}|a|_v=1
\]
by the product formula.
\end{remark}

\subsubsection{The BSD formula}
Let $\calF$ be either $\QQ$ or $\calK$ and suppose that the analytic and algebraic ranks of $E/\calF$ agree and that $\Sha(E/\calF)$ is finite. Let $r_\calF$ be the common rank. Define the BSD quotient, which is often called the \emph{analytic Tate--Shafarevich invariant}:
\begin{equation}\label{eq:bsd-quotient}
 \calQ(E/\calF):= \frac{L^{(r_\calF)}(E/\calF,1)}{r_\calF!\,\Omega_{E/\calF}\Reg(E/\calF)} \frac{\#E(\calF)_{\tor}^{2}}{\prod_v \Tam_v(E/\calF)},
\end{equation}
while $\# \Sha(E/\calF)$ is often called the \emph{arithmetic Tate--Shafarevich invariant}. The $p$-primary part of the Birch and Swinnerton-Dyer (BSD) formula for $E/\calF$ means
\[
v_p\bigl(\calQ(E/\calF)\bigr)=v_p\bigl(\#\Sha(E/\calF)\bigr).
\]
Equivalently, define the error exponent
\begin{equation}\label{eq:error}
 \Delta_p(E/\calF):= v_p\bigl(\calQ(E/\calF)\bigr) -v_p\bigl(\#\Sha(E/\calF)\bigr).
\end{equation}
The $p$-primary part of the BSD formula predicts that $\Delta_p(E/\calF)=0$.

\begin{lemma}\label{lem:base-change}
Let $E/\QQ$ have CM by $\calK$ and analytic rank one. Then
\begin{enumerate}[label = \rm (\arabic*)]
    \item The $\ZZ$-rank of $E(\QQ)$ is $1$, and $\rank_{\ZZ} E(\calK) = \ord_{s=1} L(E/\calK,s) = 2$.
    \item The Tate--Shafarevich groups $\Sha(E/\QQ)$ and $\Sha(E/\calK)$ have finite cardinality.
    \item $ \Delta_p(E/\calK)=2\Delta_p(E/\QQ)$.
\end{enumerate}
In particular, the $p$-primary part of the BSD formula holds over $\calK$ if and only if it holds over $\QQ$.
\end{lemma}

\begin{proof}
Let $\eta=\eta_{\calK/\QQ}$ be the quadratic character attached to $\calK/\QQ$ by class field theory.

\underline{Step 1: the CM twist is isogenous to $E$.}
For every prime $\ell$ (possibly $\ell=p$), after extending coefficients so that the CM character is one-dimensional, complex multiplication gives
\[
 V_\ell(E)\simeq \Ind_{G_\calK}^{G_\QQ}V_{\psi,\ell}.
\]
Since $\eta|_{G_\calK}=1$, the projection formula $(\Ind_{G_\calK}^{G_\QQ} V_{\psi,\ell})\otimes\eta\simeq \Ind_{G_\calK}^{G_\QQ}(V_{\psi,\ell} \otimes \Res_{G_{\calK}}^{G_{\QQ}}\eta)$ gives
\[
 V_\ell(E) \simeq V_\ell(E)\otimes\eta\simeq V_\ell(E^{\eta}).
\]
Faltings's isogeny theorem \cite{Faltings} therefore gives a $\QQ$-isogeny $E^{\eta}\sim_{\QQ}E$. In particular
\begin{equation}\label{eq:cm-twist-L}
 L(E^{\eta}/\QQ,s)=L(E/\QQ,s).
\end{equation}

\underline{Step 2: restriction of scalars.} Let $B=\Res_{\calK/\QQ}(E_\calK)$. Then there is a $\QQ$-isogeny
\begin{equation}\label{eq:weil-restriction-isogeny}
 B\sim_{\QQ}E\times E^{\eta}.
\end{equation}
Indeed, over $\calK$ the sum-and-difference map is $(x,y)\mapsto(x+y,x-y)$ after identifying the twist with $E$. Conjugation swaps the two factors on the source and acts by $+1,-1$ on the two target factors. The map therefore descends to $\QQ$. Its composite with the analogous reverse map is multiplication by $2$, so its kernel is finite. 

Artin--Brauer formalism applied to $B$ (i.e. comparing the Euler factors attached to the Galois representations) consequently gives
\begin{equation}\label{eq:artin-base-change}
 L(E/\calK,s)=L(B/\QQ,s) =L(E/\QQ,s)L(E^{\eta}/\QQ,s)=L(E/\QQ,s)^2,
\end{equation}
where the last equality follows from \eqref{eq:cm-twist-L}. In particular, this shows that $\ord_{s=1} L(E/\calK, s) = 2$. 

Moreover, the Gross--Zagier--Kolyvagin theorem \cite{GrossZagier, Kolyvagin} says that $\rank_{\ZZ} E(\QQ) = 1$, and since $E^{\eta}$ is $\QQ$-isogenous to $E$, we see that $\rank_{\ZZ} E^{\eta}(\QQ) = 1$. This implies
\[
\rank_{\ZZ} E(\calK) = \rank_{\ZZ} B(\QQ) = \rank_{\ZZ} E(\QQ) + \rank_{\ZZ} E^{\eta}(\QQ) = 2.
\]
This proves (1). 

For (2), we first note that finiteness of the Tate--Shafarevich groups is invariant under isogeny: for an isogeny $u:C\to D$, choose an isogeny $v:D\to C$ such that $u\circ v=[m]_D$. The image of $v_{\ast}:\Sha(D/\QQ)\to\Sha(C/\QQ)$ is finite if $\Sha(C/\QQ)$ is finite, while its kernel is contained in the finite group $\Sha(D/\QQ)[m]$ (note that $\Sha(D/\QQ)[m]$ is a finite group since it is a quotient of the finite $m$-Selmer group $\Sel_{m}(D/\QQ)$). The Gross--Zagier--Kolyvagin theorem \cite{GrossZagier, Kolyvagin} gives the finiteness of $\Sha(E/\QQ)$. Since $E^\eta$ is $\QQ$-isogenous to $E$ by Step 1, $\Sha(E^\eta/\QQ)$ is finite as well. Put $A=E\times E^\eta$. Then
$\Sha(A/\QQ)=\Sha(E/\QQ)\times\Sha(E^\eta/\QQ)$ is finite.
The $\QQ$-isogeny \eqref{eq:weil-restriction-isogeny} then shows that $\Sha(B/\QQ)$ is finite. Finally, Shapiro's lemma identifies $\Sha(B/\QQ)\simeq\Sha(E/\calK)$. Thus $\Sha(E/\calK)$ is finite. This proves (2).

\underline{Step 3: the full BSD quotient respects the isogenies.}
For an abelian variety $A/\calF$, the analytic invariant $\calQ(E/\calF)$ and the arithmetic invariant $\#\Sha(A/\calF)$ are invariant under isogeny, are multiplicative under products, and are compatible with restriction of scalars: see \cite[I, \S7]{MilneADT} and \cite[\S1]{MilneRes}. Thus
\[
 \Delta_p(E/\calK)=\Delta_p(B/\QQ) =\Delta_p(E/\QQ)+\Delta_p(E^{\eta}/\QQ)=2\Delta_p(E/\QQ).
\]
This proves the lemma.
\end{proof}

\begin{remark} \label{rem:L-function-over-K}
Differentiating $L(E/\QQ, s)$ at its simple zero at $s=1$, \eqref{eq:artin-base-change} also proves that $L^{\prime\prime}(E/\calK, 1)/2 = L^{\prime}(E/\QQ,1)^{2}$.
\end{remark}

\subsubsection{A regulator identity}
We have seen in Lemma \ref{lem:base-change} that $E(\calK)$ has $\ZZ$-rank $2$, and hence it has $\calO_{\calK}$-rank $1$. We choose $P_{0}$ to be a generator of $E(\calK)_{\tf}$, and denote $m := [E(\calK)_{\tf}:\calO_{\calK} P_0]$. Also set $R_{0} := b_{\calK,M}(P_0, P_0)$. We have the following lemma relating the regulator $\Reg(E/\calK)$ to the height $R_{0}$.

\begin{lemma}\label{lem:CM-regulator}
With the notation above,
\begin{equation}\label{eq:CM-regulator}
 \Reg(E/\calK)=\frac{D_{\calK}}{4m^2} R_0^2.
\end{equation}
\end{lemma}
\begin{proof}
For any $a,b\in\calK^{\times}$, one checks that $b_{\calK}(aP_0,bP_0)=\Re(\bara b) R_0$. Choose an integral basis $\{1,\theta\}$ of $\calO_{\calK}$, then $P_{0}, \theta P_0$ is a $\ZZ$-basis of the sublattice $\calO_{K} P_0$. The Gram matrix is therefore $R_0 \cdot \begin{bmatrix}
       1 & (\theta + \overline{\theta})/2 \\ (\theta + \overline{\theta})/2 & \theta \overline{\theta} 
\end{bmatrix}$, and its determinant is
\[
R_{0}^{2} \left(\theta \bar{\theta}- \frac{(\theta+\bar\theta)^2}{4}\right) =-\frac{R_0^2}{4}(\theta-\bar\theta)^2 =\frac{D_{\calK}}4R_{0}^2,
\]
because the discriminant of $\{1,\theta\}$ is $-D_{\calK}$. An index-$m$ sublattice $\calO_{K} P_0$ has regulator $m^2$ times that of the full lattice $E(\calK)_{\tf}$. The lemma is thus proved.
\end{proof}

\subsection{Classification of CM elliptic curves}
Classically, in Iwasawa theory, we study elliptic curves according to their reduction type at $p$. For CM elliptic curves, the situation is different. We classify the cases by the decomposition type of $p$ in the CM field $\calK$. The two classifications are different, and we compare them in Table \ref{table:CM_classification}.

Recall that, in our setup, $E/\QQ$ is an elliptic curve with complex multiplication, and let $\calK$ be its imaginary quadratic field. It is attached to an algebraic Hecke character $\psi: \calK^{\times} \bs \AA_{\calK}^{\times} \rightarrow \CC^{\times}$ of infinite type $(-1,0)$, which means $\psi_{\infty}(z)=z^{-1}$ in our paper. It is known that the conductor of $E$ can be written as 
\[
N_{E} = \abs{D_{\calK}} \Norm_{\calK/\QQ}(\frf_{\psi})
\]
where $\frf_{\psi}$ is the conductor of $\psi$. For any $\frp \mid p$ of $\calK$, we set 
\[
a_{\frp}(\psi) := v_{\frp}(\frf_{\psi}), \quad n_{p} := v_{p}(N_{E}), \quad d_{\calK,p} := v_{p}(D_{\calK}),
\]
and write $f(\frp/p)$ for the residue degree. Taking the $p$-adic valuation of $N_{E}$ gives
\[
n_{p} = d_{\calK,p} + \sum_{\frp \mid p} f(\frp/p) a_{\frp}(\psi).
\]
There are three cases for the decomposition type of $p$ in $\calK$:
\begin{itemize}
  \item If $p$ splits in $\calK$, then $p \calO_{\calK} = \frp \bar{\frp}$ and hence $n_{p} = a_{\frp}(\psi)+ a_{\bar{\frp}}(\psi)$. One checks that conjugate conductor exponents agree, and hence $n_{p} = 2 a_{\frp}(\psi)$.
  \item If $p$ is inert in $\calK$, then $p \calO_{\calK} = \frp$ with $f(\frp/p) = 2$. This shows that $n_{p} = 2 a_{\frp}(\psi)$.
  \item If $p$ is odd and ramified in $\calK$, then $p \calO_{\calK} = \frp^{2}$, $f(\frp/p) = 1$ and $d_{\calK,p} = 1$. Thus $n_{p} = 1+a_{\frp}(\psi)$.
\end{itemize}

Moreover, in \cite[Proposition 6.1]{Schuett}, all possible local exponents $a_{\frp}(\psi)$ are classified. We list his result here.
\begin{proposition}
Suppose $p$ is an odd prime. Then we have 
\begin{enumerate}[label = \rm (\arabic*)]
  \item If $p$ is unramified in $\calK$, then $a_{\frp}(\psi) \in \{0,1\}$,
  \item if $p \geq 5$ and ramified in $\calK$, then $a_{\frp}(\psi) = 1$,
  \item if $p = 3$ and $\calK = \QQ(\sqrt{-3})$, then $a_{\frp}(\psi) \in \{1,2,4\}$.
\end{enumerate} 
\end{proposition}

Combining these considerations with Deuring's criterion, namely, that 
\begin{itemize}
  \item if $p$ splits in $\calK$, then $E$ has potentially ordinary reduction over $\QQ_p$, and 
  \item if $p$ does not split in $\calK$, then $E$ has potentially supersingular reduction over $\QQ_p$,
\end{itemize}
we obtain Table \ref{table:CM_classification}.

\begin{table}[H]
\centering
\small
\caption{Decomposition type and the reduction type}
\label{table:CM_classification}
\begin{tabularx}{\textwidth}{
>{\raggedright\arraybackslash}X
>{\raggedright\arraybackslash}X
>{\centering\arraybackslash}p{0.10\textwidth}
>{\raggedright\arraybackslash}X}
\toprule
Decomposition of \(p > 2\) in \(\calK\) & 
ramification of \(\psi_{\frp}\) and $a_{\frp} := v_{\frp}(\frf_{\psi})$ & 
\(v_{p}(N_{E})\) & 
reduction type over \(\QQ_p\) \\
\midrule

\(p=\mathfrak p \overline{\mathfrak p}\) split & 
\textbf{unramified}, \(a_{\frp}= a_{\bar{\frp}} = 0\) & 
\(0\) & 
\textbf{good} ordinary \\

\(p=\mathfrak p \overline{\mathfrak p}\) split & 
\emph{ramified}, \(a_{\frp}= a_{\bar{\frp}} = 1\) & 
\(2\) & 
\emph{additive}, potentially good ordinary \\

\midrule
\(p\) inert & 
\textbf{unramified}, \(a_{\frp}=0\) & 
\(0\) & 
\textbf{good} supersingular \\

\(p\) inert & 
\emph{ramified}, \(a_{\frp}=1\) & 
\(2\) & 
\emph{additive}, potentially good supersingular \\

\midrule
\(p\geq5\) ramified in \(\calK\) & 
necessarily \emph{ramified}, \(a_{\frp}=1\) & 
\(2\) & 
\emph{additive}, potentially good supersingular \\

\(p=3\), ramified in \(\calK=\QQ(\sqrt{-3})\) & 
necessarily \emph{ramified}, \(a_{\frp}\in\{1,2,4\}\) & 
\(2,3,5\) & 
\emph{additive}, potentially good supersingular\\
\bottomrule
\end{tabularx}
\end{table}
The $p=2$ case is more complicated. Since this table is used only to illustrate the overall framework, we do not list the case $p=2$ in detail.

%% file: 103-Plus-minus.tex
\section{Arithmetic side: A local plus-minus theory} \label{sec:plus-minus}

Throughout this section $p$ is nonsplit in $\calK$ (and $p=2$ is allowed). We use the lattice in \eqref{eq:native-CM-lattice}, put $F=\calK_{\frp}$, $R=\calO_F$, and fix one coefficient embedding $\sigma: R_{p} \rightarrow \calO_{L}$.

\subsection{The small disc}
Let $\calK_\infty^-/\calK$ be the anticyclotomic $\ZZ_p$-extension and $\Gamma_\calK^-:=\Gal(\calK_\infty^-/\calK)$. Choose a topological generator $\gamma$ and identify
\[
 \calO_{L}\lrbracket{\Gamma_\calK^-}\simeq \calO_{L}\lrbracket{X},\quad \gamma\longmapsto 1+X.
\]
For $m \geq 0$, let\footnote{Note that the notation here is a little bit different from \cite{FanWan}.} $\Lambda_{m} := \calO_{L} \lrbracket{p^{-m}X} \otimes_{\calO_{L}} L$. This is the algebra of analytic functions on a sufficiently small closed disc around the trivial character. It lies between the Tate algebras
\[
L \lrangle{p^{-m} X} \subseteq \Lambda_{m} \subseteq L \lrangle{p^{-m-1}X}
\]
as described in \cite[\S 5.3]{FanWan}. Let $\bpsi$ denote the anticyclotomic deformation of $\psi$ with coefficients in $\Lambda_{m}$. The specialization at the origin (i.e. at the trivial character) is the character $\psi$.

\subsection{The plus-minus subspaces and Selmer groups}
An arithmetic point $\phi:\Lambda_m\to\overline{\QQ}_p$, as in \cite[\S 5.3]{FanWan}, has an integer weight $\kappa_\phi$.  Following the notation of \cite[Definition 5.5]{FanWan}, put
\[
 \frX_{1}=\{\phi:\kappa_\phi\geq 0, \, 2(p^2 - 1) \mid \kappa_{\phi} \},\quad \frX_2=\{\phi:\kappa_\phi\leq -2, \, 2(p^2 - 1) \mid \kappa_{\phi}\}.
\]
Let $\phi_0$ be the trivial specialization (also called the origin of the family $\boldsymbol{\psi}$).

\subsubsection{Plus-minus Selmer groups}
The plus-minus Selmer condition defined in \cite{FanWan} is based on the following proposition.

\begin{proposition}\label{prop:fw-splitting}
After increasing $m$ if necessary, there are free rank-one direct summands
\[
 \rmH^{1}_\pm(\calK_p, V_{\bpsi}) \subset \rmH^{1}(\calK_p,V_{\bpsi})
\]
such that
\begin{equation}\label{eq:pm-splitting}
 \rmH^{1}(\calK_p,V_{\bpsi})=\rmH^{1}_+(\calK_p,V_{\bpsi})\oplus\rmH^{1}_-(\calK_p,V_{\bpsi}).
\end{equation}
Moreover, $\rmH^{1}_+$ is characterized by specialization to $\rmH^{1}_{\rmf}$ at every point of $\frX_{1}$, and $\rmH^{1}_-$ is characterized by specialization to $\rmH^{1}_{\rmf}$ at every point of $\frX_2$.
\end{proposition}

\begin{proof}
This is \cite[Definition 5.5, Definition~5.24, Proposition~5.25 and Lemma~5.26]{FanWan}, with the two summands constructed in \cite[Definition 5.24]{FanWan}.
\end{proof}

The decomposition \eqref{eq:pm-splitting} is over the field-valued ring $\Lambda_{m}$. We define
\begin{align}
 H^1_\pm(\calK_p,T)&=\{x\in H^1(\calK_p,T):x_L\in H^1_\pm(\calK_p,V)\},\notag\\
 H^1_\pm(\calK_p,W^\dagger)&=H^1_\pm(\calK_p,T)^\perp,\label{eq:pm-W-defn}
\end{align}
where the orthogonal complement in \eqref{eq:pm-W-defn} is taken under the local Tate duality between $T$ and $W^{\dagger}$. We further define compact and discrete Selmer groups by
\begin{align}
 \rmH^{1}_\pm(\calK,T) &:=\ker \left( \rmH^{1}(\calK_S/\calK,T)\rightarrow \frac{\rmH^{1}(\calK_p,T)}{\rmH^{1}_\pm(\calK_p,T)}
 \right),\label{eq:compact-selmer}\\
 \Sel_\pm(\calK,W) &:=\ker \left( \rmH^{1}(\calK_S/\calK,W)\rightarrow \prod_{v\in S,\,v\nmid p}\rmH^{1}(\calK_v,W)\times\frac{\rmH^{1}(\calK_p,W)}{\rmH^{1}_\pm(\calK_p,W)} \right). \label{eq:discrete-selmer}
\end{align}

\subsubsection{Two results on Bloch--Kato Selmer groups}
A direct consequence of the defining Proposition \ref{prop:fw-splitting} is that, for the specialization at the origin, $\rmH^{1}_{+}(\calK_p, T) = \rmH^{1}_{\rmf}(\calK_p, T)$, and globally $\rmH^{1}_{+}(\calK, T)$ is the Selmer group with the Bloch--Kato Selmer condition at $p$ and the relaxed condition at $v \nmid p$. We record the following proposition. Although it is stated in terms of the Bloch--Kato Selmer group, as remarked above, it also describes the plus Selmer group at the origin.

\begin{proposition}\label{prop:finite-lines}
With the notation above, in particular, $E/\QQ$ is an elliptic curve with CM given by $\calK$, and $\ord_{s=1} L(E/\QQ,s) = 1$. Then we have
\begin{equation}\label{eq:finite-line-dimensions}
 \dim_L\rmH^{1}_{\rmf}(\calK,V_{\sigma})=\dim_L\rmH^{1}_{\rmf}(\calK_p,V_{\sigma})=1,
\end{equation}
and the localization map
\begin{equation}\label{eq:finite-line-localization}
 \loc_p:\rmH^{1}_{\rmf}(\calK,V_{\sigma})\rightarrow \rmH^{1}_{\rmf}(\calK_p,V_{\sigma})
\end{equation}
is nonzero, hence is an isomorphism.
\end{proposition}

\begin{proof}
We have seen in Lemma \ref{lem:base-change} that $\rank_{\ZZ} E(\calK) = 2$ and $\Sha(E/\calK)$ is finite. Hence $E(\calK)_{\tf}$ is a free $\calO_{\calK}$-module of rank one, and consequently $E(\calK)_p^{\wedge}$ has rank one over $R$.

The inverse limit of the Kummer exact sequence gives (for the identification of the rational Bloch--Kato condition with the Kummer local condition, see, for example, \cite[\S 1.6]{rubin2000euler}):
\[
0 \rightarrow E(\calK)_p^{\wedge} \rightarrow \rmH^{1}_{\rmf}(\calK, T_p E) \rightarrow T_p \Sha(E/\calK) \rightarrow 0.
\]
The last term vanishes since $\Sha(E/\calK)$ is finite. After the flat scalar extension $R \rightarrow \calO_{L}$, it gives
\[
\rmH^{1}_{\rmf}(\calK, V_{\sigma}) \simeq E(\calK)_{p}^{\wedge} \otimes_{R,\sigma} L.
\]
The right-hand side has $L$-dimension one, and hence so does $\rmH^{1}_{\rmf}(\calK, V_{\sigma})$.

For the local statement, by the local Kummer description of the Bloch–Kato finite condition for an abelian variety \cite[\S 1.6]{rubin2000euler}, we have
\[
\rmH_{\rmf}^1(\calK_p,V_p E) \simeq E(F)^\wedge_p\otimes_{\ZZ_p} \QQ_p.
\]
Cohomology and the Bloch--Kato finite condition commute with the flat coefficient extension $\sigma: R \rightarrow \calO_{L}$. This shows that
\[
\rmH_{\rmf}^1(\calK_p,V_\sigma) \simeq (E(\calK_p)^\wedge_p\otimes_{\ZZ_p} \QQ_p) \otimes_{\calK_p,\sigma} L.
\]
Choose a sufficiently deep formal subgroup $E_m(F)\subset E(F)$. It has finite index, and the formal logarithm identifies it with a full $\calO_{\calK_p}$-lattice in $\Lie(E/\calK_p)$. Consequently it induces an isomorphism
\[
\log_E: E(\calK_p)^\wedge_p \otimes_{\ZZ_p} \QQ_p \xrightarrow{\sim} \Lie(E/\calK_p).
\]
The right-hand side is one-dimensional over $\calK_p$. Combining the above two isomorphisms gives 
\[
\rmH_{\rmf}^1( \calK_p ,V_\sigma) \simeq \Lie(E/\calK_p) \otimes_{\calK_p,\sigma}L, \quad \dim_L \rmH_{\rmf}^1(\calK_p,V_\sigma)=1.
\]

It remains to show that the localization map $\loc_p$ is nonzero. Choose a nontorsion point $P_0 \in E(\calK)$. By the arguments above, it suffices to show that $\log_E(P_0)\neq 0$ in $\Lie(E/\calK_p)$. Otherwise, because $E_m(\calK_p)$ has finite index in $E(\calK_p)$, some positive multiple $nP_0$ lies in $E_m(\calK_p)$. Then $\log_E(nP_0)=n\log_E(P_0)=0$. The logarithm is injective on $E_m(\calK_p)$, so $nP_0 = 0$, contradicting the choice of $P_0$. Thus $\log_E(P_0)\neq0$. This shows that the localization map $\loc_p$ is nonzero, and a nonzero map between two one-dimensional $L$-vector spaces is an isomorphism.
\end{proof}

\begin{lemma}\label{lem:global-vanishing}
At every finite place $v$ of $\calK$, we have $\rmH^0(\calK_v,V) = \rmH^2(\calK_v,V)=0$. At $v \nmid p$, we also have $\rmH^1(\calK_v,V)=0$, so
\[
 \rmH_{\rmf}^1(\calK_v,T)= \rmH^1(\calK_v,T), \quad \rmH_{\rmf}^1(\calK_v,W)=0.
\]
Globally, we have
\begin{equation}\label{eq:global-rational}
\rmH^1(G_{\calK,S},V)=\rmH_{\rmf}^1(\calK,V),\qquad
 \rmH^i(G_{\calK,S},V)=0\quad(i\ne1).
\end{equation}
\end{lemma}
\begin{proof}
Put $F = \calK_{p}$. We use Proposition \ref{prop:finite-lines} both for $V$ and $V^{\dagger}$. In particular, the localization map $\loc_p$ is an isomorphism of $L$-vector spaces.

First consider a finite place $v$ of $\calK$. The torsion subgroup of $E(\calK_v)$ is finite, so
\[
\rmH^0(\calK_v,T_p E) =T_p (E(\calK_v)[p^\infty])=0.
\]
This gives $\rmH^0(\calK_v,T)=0$. If $x\in \rmH^0(\calK_v,V)$, a sufficiently large power of $p$ carries $x$ into $T$; that multiple is still invariant and must therefore be zero. Hence $\rmH^{0}(\calK_v, V) = 0$. The same conclusion holds for $V^{\dagger}$. Local Tate duality now yields
\[
\rmH^2(\calK_v,V) \simeq \rmH^0(\calK_v,V^\dagger)^{\ast}=0.
\]
When \(v\nmid p\), the local Euler–Poincaré formula has Euler characteristic zero:
\[
\dim_L \rmH^0(\calK_v,V) -\dim_L \rmH^1(\calK_v,V) +\dim_L \rmH^2(\calK_v,V)=0.
\]
We have seen that the outer cohomology groups vanish, and hence $\rmH^1(\calK_v,V)=0$. In particular $\rmH_{\rmf}^1(\calK_v,V)=0$. It then follows from the definition that
\begin{equation} \label{eq:local-vanishing-BK-condition}
\rmH_{\rmf}^1(\calK_v,T)=\rmH^1(\calK_v,T),
\quad \rmH_{\rmf}^1(\calK_v,W)=0.
\end{equation}

We next prove the global assertion in degree one. At places outside $S$, classes in $\rmH^1(G_{\calK,S},V)$ are unramified and hence satisfy the finite condition. At places of $S$ away from $p$, the condition is automatic by \eqref{eq:local-vanishing-BK-condition}. Thus only localization at $F$ can distinguish $\rmH^1(G_{\calK,S},V)$ from $\rmH_{\rmf}^1(\calK,V)$.

Poitou–Tate duality (see Section \ref{sec:recap-Poitou-Tate}) for the finite and $p$-relaxed conditions gives an exact sequence
\begin{equation} \label{eq:poitou-tate-relaxed-finite}
0\rightarrow \rmH_{\rmf}^1(\calK,V) \rightarrow \rmH^1(G_{\calK,S},V) \rightarrow \frac{H^1(F,V)}{\rmH_{\rmf}^1(F,V)} \xrightarrow{\Theta}\rmH_{\rmf}^1(\calK,V^\dagger)^{\ast}.
\end{equation}
The local Bloch–Kato finite subspaces for \(V\) and \(V^\dagger\) are orthogonal complements under local Tate duality. Consequently local duality identifies the source of \(\Theta\) with \(\rmH_{\rmf}^1(F,V^\dagger)^*\). Under this identification, \(\Theta\) is the dual of $\loc_{p}: \rmH_{\rmf}^1(\calK,V^\dagger)\rightarrow \rmH_{\rmf}^1(F,V^\dagger)$. More explicitly, if $\barx$ is represented by $x \in \rmH^{1}(F,V)$, then $\Theta(\bar x)(y) =\lrangle{x, \loc_p(y)}_{F}$. This is independent of the choice of $x$. Proposition \ref{prop:finite-lines} applied to $V^{\dagger}$ says that the localization map is an isomorphism. Therefore $\Theta$ is an isomorphism. Exactness of \eqref{eq:poitou-tate-relaxed-finite} proves
\[
\rmH^1(G_{\calK,S},V)=\rmH_{\rmf}^1(\calK,V).
\]
This space has dimension one over $L$ by Proposition \ref{prop:finite-lines}.

Finally, $\rmH^0(G_{\calK,S},V)=0$, since global invariants are local invariants at $F$. Because $\calK$ has no real places and $S$ contains the places above $p$, the $p$-cohomological dimension of $G_{\calK,S}$ is two. Thus $\rmH^i(G_{\calK,S},V)=0$ for $i \geq 3$. The global Euler–Poincaré formula gives
\[
\sum_{i=0}^{2}(-1)^i\dim_L \rmH^i(G_{\calK,S},V) =-\dim_L V=-1.
\]
The dimensions in degrees zero and one are $0$ and $1$, so the formula forces $\dim_{L} \rmH^{2}(G_{\calK,S}, V) = 0$. All the assertions are therefore proved.
\end{proof}

\subsection{The classical CM main conjecture}

Let $z_{\bpsi} \in \rmH^{1}(G_{\calK,S},V_\bpsi)$ be the family of elliptic units in \cite[\S 5.4]{FanWan}, and let $z_{\psi} \in \rmH^{1}(G_{\calK,S},V)$ be its specialization at the origin. By Lemma~\ref{lem:global-vanishing}, the elliptic unit $z_{\psi}$ gives an $L$-basis of the fundamental line $\det_{L}^{-1} \RG(\calO_{\calK,S}, V)$. The classical CM main conjecture is then the following theorem, proved in \cite[Theorem 4.1]{BurungaleFlach}.

\begin{theorem}[Classical CM main conjecture]\label{thm:classical_CM_main_conjecture}
We have the following trivialization 
\begin{equation}\label{eq:integral-zeta-generator}
 \calO_{L} z_\psi = \det_{\calO_{L}}^{-1} \RG(\calO_{\calK,S},T)
\end{equation}
inside the rational fundamental line $\det_L^{-1}\RG(\calO_{\calK,S},V)$, induced by the elliptic unit zeta morphism $z_{\psi}$.
\end{theorem}

\subsection{A plus-minus control theorem}
The main goal of this subsection is the integral determinant comparison in Theorem~\ref{thm:fan_wan_control}.

\begin{remark}[Determinant functors]
We fix the convention for the determinant functors. If $R$ is a discrete valuation ring and $C$ is a perfect $R$-complex, then
\begin{equation}\label{eq:det-convention}
 \det_R C:=\bigotimes_i(\det_R \rmH^{i}(C))^{(-1)^i}, \quad \det_R^{-1}C:=\bigotimes_i(\det_R \rmH^{i}(C))^{(-1)^{i+1}}.
\end{equation}
\end{remark} 

Proposition~\ref{prop:finite-lines} and Lemma~\ref{lem:global-vanishing} give the following isomorphism:
\begin{equation} \label{eq:q-loc-p}
q_{-} \circ \loc_p: \rmH^{1}(G_{\calK,S}, V) = \rmH^{1}_{\rmf}(\calK,V) \xrightarrow[\sim]{\loc_p} \rmH^{1}_{\rmf}(\calK_p, V) = \rmH^{1}_{+}(\calK_p, V) \xrightarrow[\sim]{q_{-}} \dfrac{\rmH^1(\calK_p,V)}{\rmH^1_-(\calK_p,V)}.
\end{equation}
Note that the kernel of $\rmH^{1}(\calK_p, T) \rightarrow \rmH^{1}(\calK_p,V)$ is exactly the $\calO_{L}$-torsion of $\rmH^{1}(\calK_p, T)$, and hence $\rmH^{1}(\calK_p, T)_{\tf}$ embeds as an $\calO_{L}$-lattice inside $\rmH^{1}(\calK_p, V)$. We therefore have a short exact sequence
\[
0 \rightarrow \rmH^{1}(\calK_p, T)_{\tf} \cap \rmH^{1}_{-}(\calK_p,V) \rightarrow \rmH^{1}(\calK_p, T)_{\tf} \rightarrow \dfrac{\rmH^1(\calK_p,V)}{\rmH^1_-(\calK_p,V)}.
\]
As a result, $Q := \frac{\rmH^{1}(\calK_p, T)_{\tf}}{\rmH^{1}(\calK_p, T)_{\tf} \cap \rmH^{1}_{-}(\calK_p, V)}$ is a $\calO_{L}$-lattice inside $\frac{\rmH^1(\calK_p,V)}{\rmH^1_-(\calK_p,V)}$, and hence we have $Q \otimes_{\calO_{L}} L \xrightarrow{\sim} \frac{\rmH^1(\calK_p,V)}{\rmH^1_-(\calK_p,V)}$. Composing with the map \eqref{eq:q-loc-p}, we have
\begin{equation} \label{eq:q-loc-p-2}
q_{-} \circ \loc_p: \rmH^{1}(G_{\calK,S}, V) \xrightarrow{\sim} Q \otimes_{\calO_{L}} V.    
\end{equation}

Lemma~\ref{lem:global-vanishing} canonically identifies $\det_L^{-1}\RG(\calO_{\calK,S},V)$ with $\det_L\rmH^{1}(\calO_{\calK,S},V)$. Taking the determinant of \eqref{eq:q-loc-p-2} therefore gives the following trivialization map, denoted by $\iota_{-}$:
\begin{equation}\label{eq:iota-minus}
 \iota_-: \det_L^{-1}\RG(\calO_{\calK,S},V) \simeq \det_{L} \rmH^{1}(G_{\calK,S}, V) \xrightarrow[\sim]{\det_{L}(q_{-} \circ \loc_p)} \det_L \left(Q \otimes_{\calO_L} L\right).
\end{equation}

\begin{theorem}[Control theorem]\label{thm:fan_wan_control}
There is an equality of fractional $\calO_{L}$-lattices inside $Q\otimes_{\calO_{L}} L$:
\begin{equation}\label{eq:det-proposed}
 \iota_-\bigl(\det_{\calO_{L}}^{-1}\RG(\calO_{\calK,S},T)\bigr)
 =\frac{|\Sel_-(\calK,W^\dagger)|_L}{|\rmH^0(\calK,W)|_L |\rmH^0(\calK,W^\dagger)|_L}
 \left(\prod_{v\in S_0} |\rmH^0(\calK_v,W^\dagger)|_L \right)Q,
\end{equation}
where $S_0$ is the set of finite places in $S$. In particular $\Sel_-(\calK,W^\dagger)$ is finite.
\end{theorem}
We prove this by the Selmer complex method underlying \cite[Proposition~2.1]{BurungaleFlach}, replacing $\iota$ in \textit{op.cit.} by our $\iota_{-}$.

\subsubsection{Selmer complex computation}
Let $\RG(\calO_{\calK,S},T)$ denote the \'etale cohomology complex.  Since $\calK$ has no real embeddings, \cite[Proposition 2.1(b)]{BurungaleFlach} proves that it is a perfect $\calO_{L}$-complex for every $p$, including $p =2$. 

We now reproduce the Selmer complex calculation in the framework of the plus-minus theory in \cite{FanWan}. For $v\nmid p$, let
\[
 U_v^+(T):=E(\calK_v)^{\wedge}_p[-1], \quad i_v^{+}: E(\calK_v)^{\wedge}_p[-1] \rightarrow R\Gamma(\calK_v,T),
\]
where $i_{v}^{+}$ is induced by the Kummer maps, i.e., the inverse limit of the connecting homomorphisms
\[
E(\calK_v) / p^{\nu} \rightarrow \rmH^{1}(\calK_v, A[p^{\nu}])
\]
induced by the Kummer sequence. Here for any abelian group $M$, we denote by $M^{\wedge}_p := \varprojlim_{\nu} M/p^{\nu}$ its underived $p$-adic completion. At $v \mid p$, put
\[
 U_{v}^-(T):=\rmH^{1}_-(\calK_{v},T)[-1], \quad i_p^{-}: \rmH^{1}_-(\calK_{v},T)[-1] \rightarrow \RG(\calK_v,T),
\]
with the map induced by the inclusion into $H^1(\calK_v,T)$ and the preceding truncation map. We define the quotient local complex as
\[
\RG_{/\rmf}(\calK_v,T) := \Cone \left(i_{v}^{\circ}: U_{v}^{\circ}(T) \rightarrow \RG(\calK_v, T) \right)
\]
for $\circ = +$ when $v \nmid p$ and $\circ = -$ when $v \mid p$. We are now ready to define the \emph{signed Selmer complex} by \footnote{Recall that for any morphism \(f:A^\bullet\to B^\bullet\), the mapping fiber is defined by
\[
\operatorname{Fib}(f):=\operatorname{Cone}(f)[-1].
\]
Thus the fiber is nothing but the mapping cone shifted by \(-1\).}
\[
 C_-(T):=\operatorname{Fib} \left(\RG(\calO_{\calK,S},T)\rightarrow \bigoplus_{v\in S_0}\RG_{/\rmf}(\calK_v,T) \right).
\]
Equivalently, the complex $C_{-}(T)$ fits into the following distinguished triangle
\begin{equation}\label{eq:minus-triangle}
 C_-(T)\rightarrow \RG(\calO_{\calK,S},T) \rightarrow\bigoplus_{v\in S_0}\RG_{/\rmf}(\calK_v,T) \rightarrow +1.
\end{equation}

\begin{lemma}[Local cohomology computations]
\label{lem:minus-local-factors}
Under the map $\iota_{-}$ in \eqref{eq:iota-minus}, we have
\[
 \det_{\calO_{L}}^{-1}\RG_{/\rmf}(\calK_v,T)=
 \begin{cases}
 |\rmH^0(\calK_v,W^\dagger)|_L ,&v\nmid p,\\
 |\rmH^0(\calK_\frp,W^\dagger)|_L \cdot Q,&v=\frp.
 \end{cases}
\]
\end{lemma}

\begin{proof}
(1) At $v\nmid p$, taking inverse limits in the Kummer sequences gives
\[
0 \rightarrow E(\calK_v)_p^{\wedge} \xrightarrow{\Kum_{v}} \rmH^{1}(\calK_v, T_p E) \rightarrow T_{p} \rmH^{1}(\calK_v, E) \rightarrow 0.
\]
Local Tate duality identifies the last term with $\Hom_{\cts}(E(\calK_{v}), \ZZ_p)$. Since $\ell \neq p$, this is zero. This implies that $\Kum_{v}$ is an isomorphism, i.e., $\rmH^{1}(U_{v}^{+}(T)) \simeq \rmH^{1}(\calK_v,T)$. Moreover, $\rmH^{0}(\calK_v, T) = T_p(E(\calK_v)) \otimes_{R_p,\sigma} \calO_{L} = 0$. It follows from the definition of the cone $\RG_{/\rmf}(\calK_v,T) = \Cone(i_{v}^{+})$ that there is a long exact sequence
\begin{multline*}
0 \rightarrow \rmH^{0}_{/\rmf}(\calK_v,T) \rightarrow \rmH^{1}(U_{v}^{+}(T)) \xrightarrow[\sim]{\Kum_{v}} \rmH^{1}(\calK_v,T) \rightarrow \rmH^{1}_{/\rmf}(\calK_v, T) \\ 
\rightarrow 0 \rightarrow \rmH^{2}(\calK_v,T) \rightarrow \rmH^{2}_{/\rmf}(\calK_v,T) \rightarrow 0.
\end{multline*}
This shows that $\rmH^{0}_{/\rmf}(\calK_v,T) = \rmH^{1}_{/\rmf}(\calK_v,T)$ and $\rmH^{2}_{/\rmf}(\calK_v,T) \simeq \rmH^{2}(\calK_v,T)$. Local Tate duality gives $\rmH^{2}(\calK_v, T) \simeq \rmH^{0}(\calK_v, W^{\dagger})^{\vee}$. Taking the determinant gives the result for $v \nmid p$. 

(2) We next focus on the computation at $p$. It follows from the definition of the cone $R\Gamma_{/f}(\calK_v,T) = \Cone(i_v^{-})$ that there is a long exact sequence
\begin{multline*}
    0 \rightarrow \rmH^{0}_{/f}(\calK_v,T) \rightarrow \rmH^{1}_{-}(\calK_v,T) \xrightarrow{i_v^{-}} \rmH^{1}(\calK_v,T) \rightarrow \rmH^{1}_{/f}(\calK_v,T) \\ \rightarrow 0 \rightarrow \rmH^{2}(\calK_v,T) \rightarrow \rmH^{2}_{/f}(\calK_v,T) \rightarrow 0.
\end{multline*}
Hence we read from the sequence that $\rmH^{0}_{/f}(\calK_v,T)=0$ since $i_{v}^{-}$ is injective, \[
 \rmH^{1}_{/f}(\calK_v,T)=\rmH^{1}(\calK_v,T)/\rmH^{1}_-(\calK_v,T),\quad \rmH^2_{/f}(\calK_v,T)= \rmH^2(\calK_v,T),
\]
and other cohomologies vanish. Local Tate duality gives $\rmH^2(\calK_v,T)\simeq \rmH^0(\calK_v,W^{\dagger})^\vee$. Taking the determinant gives the result for $v = \frp$.
\end{proof}

\begin{lemma}[Cohomology of the signed Selmer complex]
\label{lem:minus-selmer-cohomology}
The complex $C_-(T)\otimes L$ is acyclic and has the following cohomology:
\[
 \rmH^i(C_-(T))=\begin{cases}
 \rmH^0(\calK,W),&i=1,\\
 \Sel_-(\calK,W^\dagger)^\vee,&i=2,\\
 \rmH^0(\calK,W^\dagger)^\vee,&i=3,\\
 0,&\text{otherwise}.
 \end{cases}
\]
As a result,
\[
 \det_{\calO_{L}}^{-1}C_-(T)= \frac{|\Sel_-(\calK,W^\dagger)|_L}{|\rmH^0(\calK,W)|_L |\rmH^0(\calK,W^\dagger)|_L}.
\]
\end{lemma}

\begin{proof}
  We adapt the proof of \cite[Lemma~2]{BurungaleFlach}.

  \underline{(1) Calculate $\rmH^{1}(C_{-}(T))$}. We have seen in the proof of Lemma \ref{lem:minus-local-factors} that $\rmH^{0}_{/\rmf}(\calK_v,T) = 0$; the long exact sequence of the defining triangle of $C_{-}(T)$ begins with
\[
0\rightarrow \rmH^1(C_-(T)) \rightarrow \rmH^1(\calO_{\calK,S},T) \xrightarrow{\lambda_-^1} \bigoplus_{v \in S_0}\rmH^1_{/\rmf}(\calK_v,T),
\]
and thus $\rmH^{1}(C_{-}(T)) = \ker(\lambda_-^1)$. If $x\in \rmH^1(C_-(T))$, then after extension of scalars its image $x_L$ belongs to the kernel of
\[
q_-\circ\loc_p: \rmH^1(\calO_{\calK,S},V) \rightarrow
\rmH^1(\calK_{p},V)/\rmH^1_-(\calK_p,V).
\]
This map is an isomorphism. Hence $x_L=0$, so $x$ is torsion. Therefore $\rmH^{1}(C_{-}(T)) \subseteq \rmH^{1}(\calO_{\calK,S}, T)_{\tor}$. Conversely, local conditions are full away from $p$, and the inverse-image minus condition at $p$ contains all torsion. Every global torsion class therefore satisfies them. We conclude that $\rmH^1(C_-(T)) = \rmH^1(\calO_{\calK,S},T)_{\tor} \simeq \rmH^0(\calK,W)$.

  \underline{(2) Auxiliary compact-support cohomology}. Here we use the auxiliary distinguished triangle involving compact-support cohomology
\[
\RG_c(\calO_{\calK,S},T) \rightarrow \RG(\calO_{\calK,S},T) \rightarrow \bigoplus_{v\in S_0}\RG(\calK_v,T) \rightarrow +1
\]
which fits into the following diagram whose rows and columns are distinguished triangles
\[\begin{tikzcd}
	{\RG_{c}(\calO_{\calK,S},T)} && {C_{-}(T)} && {\bigoplus_{v \in S_0} U_{v}^{\circ}(T)} \\
	{\RG_{c}(\calO_{\calK,S},T)} && {\RG(\calO_{\calK,S}, T)} && {\bigoplus_{v \in S_0}\RG(\calK_v,T)} \\
	&& {\bigoplus_{v \in S_0} \RG_{\/\rmf}(\calK_v,T)} && {\bigoplus_{v \in S_0} \RG_{\/\rmf}(\calK_v,T)}
	\arrow[from=1-1, to=1-3]
	\arrow["{=}"{description}, no head, from=1-1, to=2-1]
	\arrow[from=1-3, to=1-5]
	\arrow[from=1-3, to=2-3]
	\arrow[from=1-5, to=2-5]
	\arrow[from=2-1, to=2-3]
	\arrow[from=2-3, to=2-5]
	\arrow[from=2-3, to=3-3]
	\arrow[from=2-5, to=3-5]
	\arrow["{=}"{description}, no head, from=3-3, to=3-5]
\end{tikzcd}\]
This diagram is parallel to \cite[Diagram (4)]{BurungaleFlach} in our setting. In particular, the first row is a distinguished triangle, and since $U^{\circ}_{v}(T)$ has cohomology only in degree $1$, its long exact sequence gives
\begin{equation} \label{eq:compact-support}
\bigoplus_{v \in S_0} H^1(U_{v}^{\circ}(T)) \rightarrow \rmH^{2}_{c}(\calO_{\calK,S}, T) \rightarrow \rmH^{2}(C_{-}(T)) \rightarrow 0 \rightarrow \rmH^{3}_{c}(\calO_{\calK,S}, T) \rightarrow \rmH^{3}(C_{-}(T)) \rightarrow 0.
\end{equation}
Continuing this long exact sequence and noting that $\rmH^{q}(\calO_{\calK,S}, T)$ vanishes for $q \geq 4$ shows that $\rmH^{q}(C_{-}(T))$ vanishes unless $q = 1, 2, 3$. Equation \eqref{eq:compact-support} enables us to compute $\rmH^{2}(C_{-}(T))$ and $\rmH^{3}(C_{-}(T))$.

\underline{(3) Calculate $\rmH^{3}(C_{-}(T))$}. 
The exact sequence \eqref{eq:compact-support} in particular gives an isomorphism $\rmH^{3}(C_{-}(T)) \simeq \rmH^{3}_{c}(\calO_{\calK,S}, T)$. Moreover, Artin--Verdier duality gives
\[
\rmH_c^3(\calO_{\calK,S},T)^\vee \simeq \rmH^0(\calO_{\calK,S},W^\dagger) = \rmH^0(\calK,W^\dagger).
\]
Therefore,
\[
\rmH^3(C_-(T)) \simeq \rmH^0(\calK,W^\dagger)^\vee.
\]

\underline{(4) Compute $\rmH^{2}(C_{-}(T))$}. 
We write $L_{v}^{\circ}(T) := \rmH^1(U_{v}^{\circ}(T)) = U_{v}^{\circ}(T)[1]$ for simplicity. When $v \nmid p$, it is $\rmH^1(\calK_v,T)$ by the proof of Lemma \ref{lem:minus-local-factors} and, when $v \mid p$, it is $\rmH^{1}_{-}(\calK_v,T)$. For every $v$, define
\[
(L_v^{\circ}(T))^\perp := \left\{ y\in H^1(\calK_v,W^\dagger): \langle x,y\rangle_v=0 \text{ for all }x\in L_v^{\circ}(T)
\right\}.
\]
The perfectness of local Tate duality implies $(L_v^\circ(T))^\vee\simeq H^1(\calK_v,W^\dagger)/(L_v^\circ(T))^\perp$.

Dualizing the preceding sequence \eqref{eq:compact-support} therefore gives
\[
0\rightarrow \rmH^2(C_-(T))^\vee \rightarrow \rmH^1(\calO_{\calK,S},W^\dagger) \rightarrow \bigoplus_{v\in S_0} \dfrac{\rmH^1(\calK_v,W^\dagger)}{(L_v^{\circ}(T))^\perp}.
\]
Here for the second term, we have invoked the Artin--Verdier duality $\rmH_c^2(\calO_{\calK,S},T)^\vee \simeq \rmH^1(\calO_{\calK,S},W^\dagger)$. From this we can already see that $\rmH^{2}(C_{-}(T))^{\vee}$ is a global Selmer group with local conditions $(L_v^{\circ}(T))^{\perp}$ at $v \in S_0$. We unpackage these dual local conditions explicitly.

For $v\nmid p$, $L_{v}^{\circ}(T) = \rmH^{1}(\calK_v,T)$, and hence its orthogonal complement is zero. At \(v \mid p\),
\[
(L_v^{-})^\perp = (\rmH^1_-(\calK_v,T))^\perp =: \rmH^1_-(\calK_v,W^\dagger)
\]
where the last equality is precisely the definition in \eqref{eq:pm-W-defn}. Therefore we have proved that
\[
\rmH^{2}(C_{-}(T))^{\vee} \simeq \Sel_{-}(\calK, W^{\dagger}),
\]
i.e. it is precisely the strict-away-from-$p$ and minus-at-$p$ Selmer group. Taking the dual gives the result.
\end{proof}

Now we are ready to prove Theorem \ref{thm:fan_wan_control}.

\begin{proof}[Proof of Theorem \ref{thm:fan_wan_control}]
Apply the determinant functor to \eqref{eq:minus-triangle}. The result is then a combination of Lemmas \ref{lem:minus-local-factors} and \ref{lem:minus-selmer-cohomology}.
\end{proof}

\subsection{Canonical vector and the analytic defect}
Let $\mathscr L_\sigma\in L^\times$ be a nonzero scalar (which shall be specified later). For now, we define
\begin{equation}\label{eq:vcan}
 v_{E,\sigma}=\frac{\loc_p(z_{\psi_\sigma})}{\mathscr L_\sigma} \in H_{\rmf}^1(\calK_p,V_\sigma).
\end{equation}
We define the \emph{analytic defect} $\delta(E,\scL_{\sigma})$ as a fractional ideal by
\begin{equation}\label{eq:delta}
 \calO_{L} q_-(v_{E,\sigma})= \delta(E,\scL_\sigma) \cdot Q.
\end{equation}

We then have the following theorem.

\begin{theorem}["Tamagawa number conjecture"]\label{thm:fan_wan_tnc}
We have
\begin{equation}\label{eq:lvalue-selmer}
 (\mathscr L_\sigma)=\frac{|\Sel_-(\calK,W^\dagger)|_L}{|\rmH^0(\calK,W)|_L |\rmH^0(\calK,W^\dagger)|_L}
 \left(\prod_{v\in S_0} |\rmH^0(\calK_v,W^\dagger)|_L \right) \delta(E,\scL_\sigma)^{-1}.
\end{equation}
\end{theorem}
\begin{proof}
    This is a combination of the classical CM main conjecture (Theorem~\ref{thm:classical_CM_main_conjecture}) of Burungale--Flach and the control theorem (Theorem~\ref{thm:fan_wan_control}) in the plus-minus Iwasawa theory. The theorem then follows from the definition of the analytic defect $\delta(E, \scL_{\sigma})$. 
\end{proof}

\subsection{Appearance of torsion points and Tamagawa numbers} \label{sec:appear-tam}
Up to this point, the "Tamagawa number conjecture" (Theorem \ref{thm:fan_wan_tnc}) is purely cohomological. In this section, we shall see that the global torsion points and the Tamagawa numbers away from $p$ appear.

\begin{lemma} \label{lem:torsion-appear}
    We have
    \[
    |\rmH^0(\calK,W)|_L = \sigma(\Fitt_{R} E(\calK)[p^{\infty}]) \calO_{L}, \quad |\rmH^0(\calK,W^{\dagger})|_L = \sigma^{c}(\Fitt_{R} E(\calK)[p^{\infty}]) \calO_{L}.
    \]
    Consequently, $|\rmH^0(\calK,W)|_L \cdot |\rmH^0(\calK,W^{\dagger})|_L = (\# E(\calK)[p^{\infty}]) \calO_{L}$.
\end{lemma}
\begin{proof}
    By definition, we have $W_{\sigma} \simeq E[p^{\infty}] \otimes_{R,\sigma} \calO_{L}$, and likewise $W_{\sigma}^{\dagger} \simeq E[p^{\infty}] \otimes_{R,\sigma^{c}} \calO_{L}$. Taking the $G_{\calK}$-invariants and noting that $\calO_{L}$ is flat over $R$ via both $\sigma$ and $\sigma^{c}$, we have
    \[
    \rmH^{0}(\calK, W_{\sigma}) \simeq E(\calK)[p^{\infty}] \otimes_{R,\sigma} \calO_{L}, \quad \rmH^{0}(\calK, W_{\sigma}^{\dagger}) \simeq E(\calK)[p^{\infty}] \otimes_{R,\sigma^{c}} \calO_{L}.
    \]
    Fitting ideals commute with this extension, proving the first two equalities. The last consequence follows from Lemma \ref{lem:order-norm}.
\end{proof}

\begin{lemma} \label{lem:tamagawa-appear}
    For $v \nmid p$ such that $E/\calK$ has bad reduction at $v$, let $\Phi_{v}$ be the component group of the Néron model. Then we have
    \[
    |\rmH^0(\calK_v,W^\dagger)|_L = \sigma^{c}(\Fitt_{R}(\Phi_v(\kappa_v)[p^{\infty}])) \calO_{L}.
    \]
    For $v \mid p$, we have
    \[
    |\rmH^0(\calK_\frp,W^\dagger)|_L = \abs{E(\calK_{\frp}[p^{\infty}]) \otimes_{R,\sigma^{c}} \calO_{L}}_{L}.
    \]
\end{lemma}
\begin{proof}
    The same argument as in Lemma \ref{lem:torsion-appear} gives
    \[
    \rmH^{0}(\calK_v, W_{\sigma}^{\dagger}) \simeq E(\calK_{v})[p^{\infty}] \otimes_{R,\sigma^{c}} \calO_{L}.
    \]
    This proves the result for $v = \frp$. For $v \nmid p$, we further have the following short exact sequences
    \[
    0 \rightarrow E_1(\calK_v) \rightarrow E_{0}(\calK_{v}) \rightarrow \tilE_{\ns}(\kappa_v) \rightarrow 0, \quad 0 \rightarrow E_{0}(\calK_v) \rightarrow E(\calK_v) \rightarrow \Phi_v(\kappa_v) \rightarrow 0.
    \]
    Because $E$ has CM, its bad reduction at $v$ is therefore additive, and hence $\tilE_{\ns}(\kappa_v)$ has $\# \kappa_v$ points. Moreover, $E_{1}(\calK_{v})$ is pro-$\ell$. Since $\ell \neq p$, multiplication by $p$ is consequently an automorphism of $E_{0}(\calK_v)$. The second short exact sequence then gives
    \[
    E(\calK_v)[p^{\infty}] \simeq \Phi_v(\kappa_v)[p^{\infty}].
    \]
    This shows the first equality, as desired.
\end{proof}

\subsection{Appearance of the Tate--Shafarevich group: Poitou--Tate bookkeeping} \label{sec:appear-sha}

In this section, the crucial invariant $\Sha(E/\calK)$ will finally appear, supplying the last missing piece in the BSD formula. We first recall the setup. In the plus-minus control theorem (Theorem \ref{thm:fan_wan_control}), the term $\Sel_-(\calK,W^{\dagger})$ appears. Our goal in this section is to show that, up to a few $p$-adic local defects, it is the Tate--Shafarevich group. The main goal of this subsection is to establish Proposition \ref{prop:length}.

\subsubsection{A recap of the Poitou--Tate exact sequence} \label{sec:recap-Poitou-Tate}
Fix a finite set $S$ containing the places above $p$, the archimedean places, and all places at which $T$ is ramified. Write $\rmH^1(\calK,T):=\rmH^1(G_{\calK,S},T)$. A \emph{compact Selmer structure} $\calF$ on $T$ consists of submodules
\[
\rmH^1_{\calF}(\calK_v,T) \subseteq H^1(\calK_v,T)
\]
for every finite $v$. Its global Selmer group is
\[
\rmH^1_{\calF}(\calK,T) := \ker\left( \rmH^1(G_{\calK,S},T) \rightarrow \bigoplus_{v \in S_0} \dfrac{\rmH^1(\calK_v,T)}{\rmH^1_{\calF}(\calK_v,T)}\right).
\]
Local Tate duality gives a perfect pairing $\rmH^1(\calK_v,T)\times \rmH^1(\calK_v,W^\dagger) \rightarrow L/\calO_L$. The \emph{discrete dual Selmer structure} is 
\[
\rmH^1_{\calF^*}(\calK_v,W^\dagger) :=\rmH^1_{\calF}(\calK_v,T)^\perp.
\]
The corresponding discrete dual Selmer group is
\[
\Sel_{\calF^*}(\calK,W^\dagger) := \ker\left(\rmH^1(G_{\calK,S},W^\dagger) \rightarrow \bigoplus_{v\in S_0} \dfrac{\rmH^1(\calK_v,W^\dagger)}{\rmH^1_{\calF^*}(\calK_v,W^\dagger)}\right).
\]
We fix some compact local conditions $\calF_v$ for $v \nmid p$, and vary only the condition at $p$. The compact \emph{$p$-relaxed} Selmer group is
\[
\rmH^1_{\rel}(\calK,T) := \ker\left( \rmH^1(G_{\calK,S},T) \rightarrow
\bigoplus_{v \in S, \, v \nmid p} \dfrac{\rmH^1(\calK_v,T)}{\rmH^1_{\calF}(\calK_v,T)}
\right).
\]
In other words, there is no Selmer condition at $p$. The compact \emph{$p$-strict Selmer group} is
\[
\rmH^1_{\str}(\calK,T) := \ker\left( H^1_{\rel}(\calK,T) \xrightarrow{\loc_p} \rmH^1(\calK_p,T)\right).
\]
In other words, its local condition at $p$ is zero. For the discrete $p$-strict group used below, the local condition is zero at $p$ and also zero away from $p$; it is dual to the full compact local conditions.

Let $\calF \subseteq \calG$ be two compact Selmer structures on $T$, namely $\rmH^1_{\calF}(\calK_v,T)\subseteq \rmH^1_{\calG}(\calK_v,T)$ at every $v$. The Poitou--Tate exact sequence compares the two Selmer conditions as follows:
\begin{equation} \label{eq:original-PT-sequence}
0 \rightarrow \rmH^{1}_{\calF}(\calK,T) \rightarrow \rmH^{1}_{\calG}(\calK,T) \rightarrow \bigoplus_{v \in S_0} \dfrac{\rmH^{1}_{\calG}(\calK_v,T)}{\rmH^{1}_{\calF}(\calK_v,T)} \rightarrow \Sel_{\calF^*}(\calK,W^\dagger)^\vee \rightarrow
\Sel_{\calG^*}(\calK,W^\dagger)^\vee \rightarrow 0.
\end{equation}

In this article, we take $\calF$ to be the plus-minus Selmer condition
\[
\rmH^{1}_{\calF_{\pm}}(\calK_v,T) := \begin{cases}
  \rmH^{1}_{\pm}(\calK_v, T), &\quad v \mid p, \\
  \rmH^1(\calK_v,T), &\quad v \nmid p.
\end{cases}
\]
and take $\calG = \calF_{\rel}$ to be the relaxed Selmer condition of $\calF$. Applying \eqref{eq:original-PT-sequence}, we have 
\begin{equation} \label{eq:plus-minus-PT-sequence}
0 \rightarrow \rmH^{1}_{\pm}(\calK,T) \rightarrow \rmH^{1}_{\rel}(\calK,T) \xrightarrow{q_{\pm} \circ \loc_p} \dfrac{\rmH^{1}(\calK_p,T)}{\rmH^{1}_{\pm}(\calK_p, T)} \rightarrow \Sel_{\pm}(\calK, W^{\dagger})^{\vee} \rightarrow \Sel_{\str}(\calK, W^{\dagger})^{\vee} \rightarrow 0,
\end{equation}
where $\Sel_{\str}(\calK, W^{\dagger})$ denotes the strict Selmer group of $W^{\dagger}$. We now compute the finite kernels and cokernels in this sequence.

\subsubsection{The invariant $\locind$}

Since $\ord_{s=1} L(E/\QQ,s) =1$, Proposition \ref{prop:finite-lines} guarantees that the following invariant is well-defined.

\begin{definition} \label{def:locind}
We define
\begin{equation}\label{eq:locind}
 \locind_p(T):=  \length_{\calO_{L}} \coker \left( \rmH^{1}_{\rmf}(\calK,T)_{\tf} \xrightarrow{\loc_p} \rmH^{1}_{\rmf}(\calK_p,T)_{\tf} \right)
\end{equation}
Here the subscript $\tf$ means \emph{torsion-free quotient}. 
\end{definition}

In Proposition~\ref{prop:finite-lines}, we have seen that both local and global Bloch--Kato cohomology groups $\rmH^{1}_{\rmf}(\calK,T)$ and $\rmH^{1}_{\rmf}(\calK_p, T)$ have torsion-free quotients of rank one and that the localization  map $\loc_p$ is an isomorphism over $L$. Consequently localization induces an injective map of rank-one torsion-free $\calO_{L}$-modules with finite cokernel. 

\begin{remark}
Explicitly, choose generators $x$ and $u$ of the two rank-one torsion-free lattices.  There is a unique $a\in\calO_{L}\setminus\{0\}$, up to multiplication by $\calO_{L}^\times$, such that $\loc_p(x)=au$, and
\begin{equation}\label{eq:locind-scalar}
 \coker(\loc_p)\simeq\calO_{L}/a\calO_{L},\quad \locind_p(T)=v_L(a). 
\end{equation}
We shall also need the multiplicative version of this invariant:
\begin{equation}\label{eq:local-index-ideal}
 \frl_p(T):=\operatorname{Fitt}_{\calO_{L}}
 \operatorname{coker}(\loc_p)=(a),
 \quad v_L(\frl_p(T))=\locind_p(T).
\end{equation}
\end{remark}

\subsubsection{The plus/minus quotient map and the localization map}
We first note that there exists a short exact sequence
\[
 0 \rightarrow \rmH^{1}_{+}(\calK,T) \xrightarrow{q_{-}} \left(\rmH^{1}(\calK_p,T)/\rmH^{1}_-(\calK_p,T)\right)_{\tf} \rightarrow \frac{\rmH^{1}(\calK_p,T)_{\tf}}
 {\rmH^{1}_+(\calK_p,T)_{\tf}+\rmH^{1}_-(\calK_p,T)_{\tf}} \rightarrow 0.
\]
We define the \emph{second $p$-adic defect} by
\begin{equation}\label{eq:delta2}
 \delta^{\prime\prime}_\sigma= \Fitt_{\calO_{L}}\left( \frac{\rmH^{1}(\calK_p,T)_{\tf}}
 {\rmH^{1}_+(\calK_p,T)_{\tf}+\rmH^{1}_-(\calK_p,T)_{\tf}} \right).
\end{equation}
It follows that
\begin{equation}\label{eq:coker-q-minus}
 \length_{\calO_{L}}\coker(q_{-})=v_L(\delta''_\sigma).
\end{equation}
By Proposition \ref{prop:finite-lines}, the source and the target of the localization map
\[
\loc_{p}: \rmH^{1}_{+}(\calK,T)_{\tf} \rightarrow \rmH^{1}_{+}(\calK_p,T)_{\tf}
\]
are both finitely generated $\calO_L$-modules with free rank-one torsion-free quotients, and it has finite cokernel. By definition, 
\begin{equation} \label{eq:review-locind}
\locind_{p}(T) = \length_{\calO_{L}}(\coker(\loc_p)).
\end{equation}
Consider the composition map
\[
\rmH^{1}_{+}(\calK,T)_{\tf} \xrightarrow{\loc_p} \rmH^{1}_{+}(\calK_p,T)_{\tf} \xrightarrow{q_{-}} \left(\rmH^{1}(\calK_p,T)/\rmH^{1}_-(\calK_p,T)\right)_{\tf},
\]
By injectivity, we have the kernel-cokernel exact sequence (see \cite[Lemma A.2 in Chapter II]{MilneCFT})
\[
0 \rightarrow \coker(\loc_p) \rightarrow \coker(q_{-} \circ \loc_p) \rightarrow \coker(q_{-}) \rightarrow 0,
\]
and therefore, using \eqref{eq:coker-q-minus} and \eqref{eq:review-locind}:
\begin{align} 
\length_{\calO_L} \coker(q_{-} \circ \loc_p) 
&= \length_{\calO_{L}} \coker(\loc_p) + \length_{\calO_{L}} \coker(q_{-}) \notag \\ 
&= \locind_{p}(T) + v_{L}(\delta''_\sigma). \label{eq:coker-q-minus-loc-p}
\end{align}

\subsubsection{The plus part of the Poitou--Tate sequence} The main goal of this subsubsection is to establish \eqref{eq:plus-part-conclusion2}.

Integrally, one checks directly that $H^1_{\rel}(\calK,T)=H^1_+(\calK,T)$ by Lemma~\ref{lem:global-vanishing}. It then follows from the plus part of the Poitou--Tate sequence (Equation \eqref{eq:plus-minus-PT-sequence}) that
\[
0 \rightarrow \dfrac{\rmH^{1}(\calK_p,T)}{\rmH^{1}_{+}(\calK_p, T)} \rightarrow \Sel_{+}(\calK, W^{\dagger})^{\vee} \rightarrow \Sel_{\str}(\calK, W^{\dagger})^{\vee} \rightarrow 0.
\]
Taking Pontryagin duality gives
\[
0 \rightarrow \Sel_{\str}(\calK, W^{\dagger}) \xrightarrow{\iota_{+}} \Sel_{+}(\calK, W^{\dagger}) \xrightarrow{\lambda_{+}} \left(\dfrac{\rmH^{1}(\calK_p,T)}{\rmH^{1}_{+}(\calK_p, T)}\right)^{\vee} \rightarrow 0.
\]
Consider the $\dagger$-dual case: now we have the localization map
\[
\loc_p^{\dagger}: \rmH^{1}_{+}(\calK,T^{\dagger})_{\tf} \hookrightarrow \rmH^{1}_{+}(\calK_p, T^{\dagger})_{\tf},
\]
and the maximal divisible subgroup 
\[ 
\Sel_{+}(\calK, W^{\dagger})_{\ddiv} \simeq \dfrac{\rmH_{\rmf}^1(\calK,V^\dagger)}{\rmH_{\rmf}^1(\calK,T^\dagger)_{\tf}} \simeq \rmH^{1}_{+}(\calK,T^{\dagger})_{\tf} \otimes_{\calO_L} L/\calO_L.
\]
Similarly, the local Tate duality and the self-dual property of the Bloch--Kato Selmer condition gives 
\[
\left(\dfrac{\rmH^{1}(\calK_p,T)}{\rmH^{1}_{+}(\calK_p, T)}\right)^{\vee} \simeq \dfrac{\rmH_{\rmf}^1(F,V^\dagger)}{\rmH_{\rmf}^1(F,T^\dagger)_{\tf}} \simeq \rmH^{1}_{+}(\calK_p,T^{\dagger})_{\tf} \otimes_{\calO_L} L/\calO_L.
\]
Under these identifications, the restriction of the localization map $\lambda_{+}$ to the maximal divisible subgroup is exactly
\[
\loc_p^{\dagger} \otimes 1: \rmH^{1}_{+}(\calK,T^{\dagger})_{\tf} \otimes_{\calO_L} L/\calO_L \rightarrow \rmH^{1}_{+}(\calK_p, T^{\dagger})_{\tf} \otimes_{\calO_L} L/\calO_L.
\]

For the map on divisible groups we use the following lemma \footnote{It may also be proved quite directly using snake lemma, or elementary verifications without using $\Tor_1$-functors.}. Readers can compare the statement with "$\loc_{p}^{\dagger} \otimes 1$".

\begin{lemma}\label{lem:lattice-map}
Let $f:A\hookrightarrow B$ be an injection of free rank-one $\calO_{L}$-modules with finite cokernel $C$. Then there is a canonical connecting homomorphism $\partial_{f}: C \rightarrow A \otimes_{\calO_L} L/\calO_L$ such that
\begin{equation}\label{eq:divisible-exact}
 0\rightarrow C \xrightarrow{\partial_{f}} A\otimes_{\calO_{L}}L/\calO_{L}\xrightarrow{f\otimes1}B\otimes_{\calO_{L}}L/\calO_{L}\rightarrow0
\end{equation}
is exact. In particular, $\ker(f \otimes 1) \simeq \coker(f)$, and hence
\[
\length_{\calO_{L}} (\coker(f)) = \length_{\calO_L} \ker(f \otimes 1).
\]
\end{lemma}

\begin{proof}
Tensoring the short exact sequence $0 \rightarrow A \rightarrow B \rightarrow C \rightarrow 0$ with $L/\calO_{L}$ gives a long exact Tor sequence 
\[
0 \rightarrow \Tor_{1}^{\calO_L}(C, L/\calO_L) \rightarrow A\otimes_{\calO_{L}}L/\calO_{L}\xrightarrow{f\otimes1}B\otimes_{\calO_{L}}L/\calO_{L}\rightarrow C\otimes_{\calO_{L}}L/\calO_{L} \rightarrow 0.
\]
Since $C$ is finite, we have that $C \otimes_{\calO_L} L/\calO_L$ is zero. It remains to identify $\Tor_{1}^{\calO_L}(C, L/\calO_L)$. To compute this Tor group, tensor the short exact sequence 
\[
0 \rightarrow \calO_L \rightarrow L \rightarrow L/\calO_L \rightarrow 0
\]
with $C$. Since $C \otimes_{\calO_L} L = 0$ and $\Tor_{1}^{\calO_{L}}(C, L) = 0$ because $L$ is flat over $\calO_L$, the resulting sequence gives a canonical isomorphism $\Tor_{1}^{\calO_{L}}(C, L/\calO_L) \simeq C$. The result follows.
\end{proof}

We apply Lemma \ref{lem:lattice-map} to the map $\loc_{p}^{\dagger}$; it gives a short exact sequence
\[
0 \rightarrow \coker(\loc_p^{\dagger}) \xrightarrow{\partial_{\loc_p^{\dagger}}} \rmH^{1}_{+}(\calK,T^{\dagger})_{\tf} \otimes_{\calO_L} L/\calO_L \xrightarrow{\loc_{p}^{\dagger} \otimes 1} \rmH^{1}_{+}(\calK_p, T^{\dagger})_{\tf} \otimes_{\calO_L} L/\calO_L \rightarrow 0.
\]
We have therefore shown that $\ker(\lambda_{+} |_{\Sel_{+}(\calK,W^{\dagger})_{\ddiv}}) = \coker(\loc_p^{\dagger})$, and moreover the latter is contained in $\Sel_{\str}(\calK, W^{\dagger})$. The map from the strict group to the finite quotient is surjective: lift a class to the plus Selmer group, and subtract an element of its divisible subgroup with the same image under $\lambda_+$. The kernel is the kernel on the divisible subgroup already computed. Thus we obtain the short exact sequence: 
\[
0 \rightarrow \coker(\loc_p^{\dagger}) \rightarrow \iota_{+}(\Sel_{\str}(\calK, W^{\dagger})) \rightarrow \Sel_{+}(\calK,W^{\dagger})_{\fin} \rightarrow 0
\]
where $\Sel_{+}(\calK,W)_{\fin} := \Sel_{+}(\calK,W)/\Sel_{+}(\calK,W)_{\ddiv}$. Taking lengths gives 
\begin{equation} \label{eq:plus-part-conclusion}
\length_{\calO_L} \Sel_{\str}(\calK, W^{\dagger}) = \length_{\calO_L} \coker(\loc_{p}^{\dagger}) + \length_{\calO_L} \Sel_{+}(\calK,W^{\dagger})_{\fin}.
\end{equation}
Note that
\begin{equation}\label{eq:locp-dagger}
\length_O\coker(\loc_p^\dagger)=\length_O\coker(\loc_p)=\locind_p(T).
\end{equation}
Indeed, in our case, $T^{\dagger} = T^{c}$, and this identification is compatible with global and local Bloch--Kato conditions and with the localization maps; we have 
\[
\rmH^{1}_{+}(\calK,T^{\dagger}) \simeq \rmH^{1}_{+}(\calK,T)^{c}, \quad \rmH^{1}_{+}(\calK_p, T^{\dagger}) \simeq \rmH^{1}_{+}(\calK,T)^{c}, \quad \loc_{p}^{\dagger} = \loc_{p}^{c}.
\]
This implies that, if we choose a generator so that $\loc_p$ is multiplication by a nonzero $a \in \calO_{L}$, then $\loc_p^{\dagger}$ is multiplication by $c(a)$ up to a unit in $\calO_L^{\times}$. Since conjugation preserves the normalized valuation, we have $v_{L}(c(a)) = v_L(a)$. This proves \eqref{eq:locp-dagger}.

Substituting \eqref{eq:locp-dagger} into \eqref{eq:plus-part-conclusion}, we obtain 
\begin{equation} \label{eq:plus-part-conclusion2}
\length_{\calO_L} \Sel_{\str}(\calK, W^{\dagger}) = \locind_{p}(T) + \length_{\calO_L} \Sel_{+}(\calK,W^{\dagger})_{\fin}.
\end{equation}

\subsubsection{The minus part of the Poitou--Tate sequence} 
The composition
\[
\rmH^{1}_{\rel}(\calK,V) \xrightarrow{\loc_p} \rmH^{1}_{+}(\calK_p, V) \rightarrow \rmH^{1}(\calK_p,V)/\rmH^{1}_{-}(\calK_p,V)
\]
is an isomorphism. Indeed, the first map is an isomorphism by Proposition \ref{prop:finite-lines}, and the second is an isomorphism because
\[
\rmH^{1}(\calK_p,V) = \rmH^{1}_{+}(\calK_p,V) \oplus \rmH^{1}_{-}(\calK_p,V).
\]
Consequently, the kernel of the map $\rmH^{1}_{\rel}(\calK,T) \rightarrow \frac{\rmH^{1}(\calK_p,T)}{\rmH^{1}_{-}(\calK_p, T)}$ is precisely the torsion module of $\rmH^{1}_{\rel}(\calK,T)$. The minus part of the Poitou--Tate sequence (equation \eqref{eq:original-PT-sequence}) then reduces to 
\[
0 \rightarrow \rmH^{1}_{+}(\calK,T)_{\tf} \xrightarrow{q_{-} \circ \loc_p} \dfrac{\rmH^{1}(\calK_p,T)}{\rmH^{1}_{-}(\calK_p, T)} \rightarrow \Sel_{-}(\calK, W^{\dagger})^{\vee} \rightarrow \Sel_{\str}(\calK, W^{\dagger})^{\vee} \rightarrow 0.
\]
Equivalently, splitting it into a short exact sequence gives 
\[
0 \rightarrow \coker(q_{-} \circ \loc_p) \rightarrow \Sel_{-}(\calK, W^{\dagger})^{\vee} \rightarrow \Sel_{\str}(\calK, W^{\dagger})^{\vee} \rightarrow 0,
\]
and dualizing gives
\[
0 \rightarrow \Sel_{\str}(\calK, W^{\dagger}) \rightarrow \Sel_{-}(\calK, W^{\dagger}) \rightarrow \coker(q_{-} \circ \loc_p)^{\vee} \rightarrow 0. 
\]
For a finite $R$-module, its Pontryagin dual and the module itself have the same length; hence
\[
\length_{\calO_{L}} \Sel_{-}(\calK, W^{\dagger}) = \length_{\calO_L} \coker(q_{-} \circ \loc_p) + \length_{\calO_L} \Sel_{\str}(\calK, W^{\dagger}).
\]
Using \eqref{eq:coker-q-minus-loc-p}, it becomes
\begin{equation} \label{eq:minus-part-conclusion}
\length_{\calO_{L}} \Sel_{-}(\calK, W^{\dagger}) = \locind_{p}(T) + v_{L}(\delta''_\sigma) + \length_{\calO_L} \Sel_{\str}(\calK, W^{\dagger})
\end{equation}

\subsubsection{Conclusion for the bookkeeping, and the Tate--Shafarevich group} \label{sec:Sha-appear}
From equations \eqref{eq:plus-part-conclusion2} and \eqref{eq:minus-part-conclusion}, we obtain the following length equality.

\begin{proposition}\label{prop:length}
  With the notation above, we have 
\begin{equation}\label{eq:length-identity}
\length_{\calO_{L}} \Sel_{-}(\calK, W^{\dagger}) = 2 \locind_{p}(T) + v_{L}(\delta''_\sigma) + \length_{\calO_L} \Sel_{+}(\calK,W^{\dagger})_{\fin}
\end{equation}
\end{proposition}

The following proposition identifies the finite quotient with the required coefficient component of the Tate--Shafarevich group.

\begin{proposition} \label{prop:find_sha}
With the notation above, we have
\[
\Sel_{+}(\calK, W^{\dagger})_{\fin} \simeq \Sha(E/\calK)[p^{\infty}] \otimes_{R_{p}, \sigma^c} \calO_{L}.
\]
In particular, 
\[
    \abs{\Sel_{+}(\calK, W^{\dagger})_{\fin}}_{L} = \sigma^c\!\left(\operatorname{Fitt}_{R_p}\Sha(E/\calK)[p^\infty]\right)\calO_L
\]
\end{proposition}
\begin{proof}
The classical Kummer sequence over $\calK$ gives a short exact sequence of $R_{p}$-modules
    \[
    0 \rightarrow E(\calK) \otimes_{\ZZ} \QQ_p/\ZZ_p \rightarrow \Sel_{p^{\infty}}(E/\calK) \rightarrow \Sha(E/\calK)[p^{\infty}] \rightarrow 0.
    \]
Since $\calO_{L}$ is flat over the discrete valuation ring $R_{p}$, tensoring this sequence via $\sigma^c$ gives
    \[
    0 \rightarrow (E(\calK) \otimes_{\ZZ} \QQ_p/\ZZ_p) \otimes_{R_{p}, \sigma^c} \calO_{L} \rightarrow \Sel_{p^{\infty}}(E/\calK) \otimes_{R_p,\sigma^c} \calO_{L} \rightarrow \Sha(E/\calK)[p^{\infty}] \otimes_{R_{p}, \sigma^c} \calO_{L} \rightarrow 0.
    \]
Under a compatible integral $\Gal_{\calK}$-equivariant identification
    \[
    W^{\dagger} \simeq E[p^{\infty}] \otimes_{R_{p}, \sigma^c} \calO_{L},
    \]
the middle term $\Sel_{p^{\infty}}(E/\calK) \otimes_{R_p,\sigma^c} \calO_{L}$ can be identified with $\Sel_{+}(\calK, W^{\dagger})$: at $p$, the $+$-condition is the Kummer condition, while at $v \nmid p$, both discrete local conditions are zero. Thus, the short exact sequence can be written as
    \[
    0 \rightarrow (E(\calK) \otimes_{\ZZ} \QQ_p/\ZZ_p) \otimes_{R_{p}, \sigma^c} \calO_{L} \rightarrow \Sel_{+}(\calK, W^{\dagger}) \rightarrow \Sha(E/\calK)[p^{\infty}] \otimes_{R_{p}, \sigma^c} \calO_{L} \rightarrow 0.
    \]
    The first term is divisible, and the last term is finite by Lemma \ref{lem:base-change}(2). Every divisible subgroup of the middle term has zero image in that finite quotient $\Sha(E/\calK)[p^{\infty}] \otimes_{R_{p}, \sigma^c} \calO_{L}$, and hence the first term is precisely its maximal divisible subgroup. This proves that
    \[
    \Sel_{+}(\calK, W^{\dagger})_{\fin} \simeq \Sha(E/\calK)[p^{\infty}] \otimes_{R_{p}, \sigma^c} \calO_{L},
    \]
    as desired. Taking order ideals gives \[
    \abs{\Sel_{+}(\calK, W^{\dagger})_{\fin}}_{L} = \sigma^c\!\left(\operatorname{Fitt}_{R_p}\Sha(E/\calK)[p^\infty]\right)\calO_L,
    \]
    as desired.
    \end{proof}

%% file: 104-GZW-comparison.tex
\section{Analytic side: The Gross--Zagier--Waldspurger comparison} \label{sec:GZW-comparison}

\subsection{Rankin--Selberg data setup} \label{sec:Rankin-Selberg-setup}
We use the setup in \cite[\S5.6.1]{FanWan}. Choose a weight-two CM eigenform $f=f_\xi$ attached to a Hecke character $\xi$ such that $\xi$ is unramified at $p$, and a finite-order character $\chi_0$ such that
\begin{equation}\label{eq:factor-data}
 \xi\chi_0=\psi, \quad L(\xi\chi_0^c,1)\ne0.
\end{equation}
We put $\psi = \psi_E = \xi \chi_0$ and $\psi_0 := \xi^{c}\chi_0$, so $L(\psi_0,1) \neq 0$. Let $A=A_f$ be the modular abelian variety attached to $f$.

\subsection{$p$-adic Waldspurger formula}
We recall the $p$-adic Waldspurger formula in \cite{FanWan}. 

\subsubsection{Setup: cuspidal toric period integral}
We start with the setup in \cite[Section 3.3]{FanWan}. In this subsection only, $\sigma$ denotes an automorphic representation rather than a coefficient embedding. Let $\sigma$ be a cuspidal automorphic representation with central character $\omega \abs{-}^{-r}$ and lowest weight $r \geq 2$. Let $B$ be the unique indefinite quaternion algebra $B$ over $\QQ$ such that $B_{p}$ is split and $\epsilon(B_{v})\mu_{v}(-1)\eta(-1) = \epsilon(1/2, \sigma_v \times \mu_v)$ for $v \nmid p\infty$. Let $\sigma := \sigma_{A}$ be the automorphic representation associated with $A$. Then $\sigma_p$ is a principal series, so such an indefinite quaternion algebra $B$ exists.

Let $\pi$ be the automorphic $B^{\times}$-representation whose Jacquet--Langlands transfer is $\sigma$. For any $\varphi \in \pi$ and $\mu \in \Sigma^{+}(r)$, we set
\[
P_{\mu}(\varphi) = \int_{\AA_{\QQ}^{\infty, \times} \calK^{\times} \bs \AA_{\calK}^{\infty, \times}} \varphi(t) \mu(t) \dif^{\times} t.
\]
The classical Waldspurger formula states that
\[
P_{\mu}(\varphi) P_{\mu^{-1}}(\varphi^{\vee}) = \dfrac{\zeta^{*}(2) \Lambda(1/2, \sigma \times \mu)}{8 \Lambda(1,\eta)^{2} \Lambda(1,\sigma, \ad)} \prod_{v} \alpha_{v}(\varphi_v, \varphi_v^{\vee}; \mu_v),
\]
where
\begin{equation} \label{eq:define-alpha-v}
\alpha_v(\varphi_v, \varphi_v^{\vee}; \mu_v) := \dfrac{L_v(1,\eta) L_v(1,\sigma,\ad)}{\zeta_v(2) L_v(1/2, \sigma_{v}\times\mu_{v})} \int_{\QQ_v^{\times} \bs \calK_{v}^{\times}} (\pi_v(t) \varphi_v, \varphi_v^{\vee}) \mu_v(t) \dif^{\times}t_v.
\end{equation}
Here the global automorphic $L$-values in the quoted formula are completed; we write them as $\Lambda$ below, and $\zeta^{\ast}(s)$ denotes the completed Riemann zeta function. Moreover, we use the following conventions:
\begin{itemize}
    \item The measure $\dif^{\times}t$ is the Haar measure on $\AA_{\QQ}^{\infty, \times} \calK^{\times} \bs \AA_{\calK}^{\infty, \times}$ with total volume $1$,
    \item The local measure $\dif^{\times} t_{v}$ is the Haar measure on $\QQ_{v}^{\times} \bs \calK_{v}^{\times}$ such that $\prod_{v} \dif^{\times} t_v = 2\Lambda(1,\eta) \dif^{\times} t$, and the volume of $\ZZ_{v}^{\times} \bs \calO_{\calK_v}^{\times}$ is $1$ for almost all $v$.
    \item The pairing $(-,-)$ is the Petersson inner product on $\pi \times \pi^{\vee}$ given by 
    \[
    (\varphi, \varphi^{\vee}) = \int_{B^{\times} \bs \BB^{\times} / \AA_{\QQ}^{\times}} \varphi(g) \varphi^{\vee}(g) \dif g,
    \]
    where $\dif g$ is the Haar measure on $B^{\times} \bs \BB^{\times} / \AA_{\QQ}^{\times}$ with total volume $2$. Moreover, $(-,-)_{v}$ is an invariant pairing on $\pi_v \times \pi_v^{\vee}$ such that $(-,-) = \prod_{v} (-,-)_{v}$.
\end{itemize}

We take an element $J \in B$ such that $J t^{c} = tJ$ for all $t \in \calK$ and $J^{2} \in \QQ^{\times}$. Since $\mu(t^{c}) = \mu^{-1}(t)\omega^{-1}\abs{tt^{c}}^{r}$, a change of variables $t \mapsto JtJ^{-1} = t^{c}$ and then $t \mapsto t^{-1}$ gives
\[
P_{\mu}(\varphi) = P_{\mu^{-1}}(\pi(J)\varphi \otimes \omega^{-1} \abs{-}^{r}).
\]

We now fix test vectors
\[
\varphi_1 = \otimes_{v} \varphi_{1,v} \in \pi, \quad \varphi_2 = \otimes_{v} \varphi_{2,v} \in \pi,
\]
where $\varphi_2 \in \pi$ is uniquely defined by $\overline{\varphi_2} = \pi(J)\varphi_1 \otimes \omega^{-1} \abs{-}^{r}$, with $(\varphi_1, \overline{\varphi_2}) \neq 0$, and we fix test vectors
\[
\varphi_3 = \otimes_{v} \varphi_{3,v} \in \sigma, \quad \varphi_4 = \otimes_{v} \varphi_{4,v} \in \sigma^{\vee}
\]
of weight $r$ and $-r$, respectively, such that $(\varphi_3, \varphi_4) \neq 0$. Assume that the $\psi_{\QQ}$-Whittaker functional of $\varphi_3$ is $W_{3} = \otimes_v W_{3,v}$, and that the $\barpsi_{\QQ}$-Whittaker functional of $\varphi_4$ is $W_{4} = \otimes_v W_{4,v}$. Here $\psi_{\QQ}$ is an additive character on $\QQ \bs \AA$. We put
\[
\lrangle{W_{3,v}, W_{4,v}}_{v} := \dfrac{\int_{\QQ_{v}^{\times}} W_{3,v}(a(y)) W_{4,v}(a(y)) \dif^{\times} y}{L_{v}(1,\sigma,\ad) \zeta_{v}(1) \zeta_{v}(2)^{-1}},
\]
where $a(y) = \diag[y,1]$. Here we fix the Haar measure $\dif^{\times} y$ as in \cite{CST}. For the test vectors $\phi_{i}$ above (for $i=1,\ldots, 4$), let
\begin{equation} \label{eq:define-beta-tilde}
\tbeta(\varphi_v, \mu_v; \varphi_{1,v}, \overline{\varphi_{2,v}}, W_{3,v}, W_{4,v}) := \dfrac{\alpha_v(\varphi_{v}, \pi(J)\varphi_v \otimes \omega^{-1} \abs{-}^{r}; \mu_v) \lrangle{W_{3,v}, W_{4,v}}_v}{(\varphi_{1,v}, \overline{\varphi_{2,v}})_{v}}.
\end{equation}

\subsubsection{Setup: Eisenstein toric integral}
For the toric integral for Eisenstein series, we follow the setup in \cite[page 26]{FanWan}. For $j \in \{E, 0\}$, we have the normalized toric functional
\[
\lambda_{j,v}(h) := \left.\dfrac{L_v(2z, \omega_{j}^{-1})}{L_v(z-\frac{1}{2}, \mu_v)} \int_{\QQ_{v}^{\times} \bs \calK_{v}^{\times}} h_{j,v,z}(t) \lambda_{j}(t) \abs{t}_{\calK_v}^{-1/2} \dif^{\times} t_{v} \right\vert_{z = 1/2}
\]
with the test vector chosen following \cite[Appendix A]{FanWan}.

\subsubsection{Anticyclotomic $p$-adic $L$-functions}
We choose parameters $\ell$ and $n$ as in \cite[Appendix A]{FanWan}, and, following \textit{op.cit.}, we put
\[
s_p = \begin{cases}
    1, &\text{ if $p$ remains inert in $\calK$}, \\ 
    2, &\text{ if $p$ ramifies in $\calK$}
\end{cases}, \quad d_{p} := \ell - 2n + s_p - 1,
\]
and 
\[
g_{\ell, n} := \begin{bmatrix}
    p^{2n-\ell+1-s_p} & p^{-\ell+n+1-s_p} \\ 0 & 1
\end{bmatrix}.
\]
We abbreviate the \emph{nonarchimedean cuspidal and Eisenstein integrals} as
\[
\calB_{i}(k) := \prod_{v < \infty} \tbeta_{v}(\pi_v(g_{\ell,n}) \varphi_{i,v}, \mu_v), \quad \calE_{j}(k) := \prod_{v < \infty} \lambda_{j,v}(h_{j,v}(-g_{\ell, n})),
\]
and we define the \emph{archimedean cuspidal integral} $\calC_{f}(k)$ and the \emph{archimedean Eisenstein integrals} $\calC_{E}(k)$ and $\calC_{0}(k)$ by 
\[
\calC_{f}(k) := \tbeta_{\infty}(\Delta^{k-1} \varphi_\infty, \mu_\infty(k)), \quad \calC_{E}(k) := \lambda_{E,\infty}(\Delta^{k-1} h_{E,\infty}), \quad \calC_{0}(k) := \lambda_{0,\infty}(\Delta^{k} h_{0,\infty}),
\]
where $\Delta := -(8 \pi i)^{-1} \begin{bmatrix} 1 & i \\ i & -1 \end{bmatrix}$ is the level raising operator. Then we have the following three $p$-adic $L$-functions constructed in \cite[Proposition 3.10]{FanWan} \footnote{For the reader's convenience, we note that the parameter $k$ in \cite{FanWan} is the parameter $r$ in this section. To obtain Theorem \ref{thm:fan-wan-p-adic-l-function}, we take $r=2$ and the parameter $\nu$ in \cite{FanWan} to be $k-1$.}.

\begin{theorem}[Anticyclotomic $p$-adic $L$-functions] \label{thm:fan-wan-p-adic-l-function}
    There are analytic functions $\calL_{i}$, $\calL_{E}$ and $\calL_{0}$ with the following property. Let $\phi$ be an anticyclotomic Hecke character over $\calK$ with infinite type $(k,-k)$ such that $k \geq 1$ on a fixed branch $N\mid k$, where $N$ is divisible by $2(p^2-1)$, and put $\rho_E=\psi_E\phi$, $\rho_0=\psi_0\phi$. Then 
    \begin{subequations}
    \begin{align}
        \calL_{i}(\phi) &= r_{f}^{4k}p^{2d_p(k-1)} \dfrac{\Lambda(\rho_E, 1) \Lambda(\rho_0, 1)}{4 \Lambda(1, \eta)^{2}} \dfrac{(\varphi_{1}, \overline{\varphi_{2}})}{(\varphi_3, \varphi_4)} \calC_{f}(k) \calB_{i}(k). \\ 
        \calL_{E}(\phi) &= r_{E}^{2k-1}p^{d_p(k-1)} \dfrac{\Lambda(\rho_E, 1) }{2\Lambda(1, \eta) \Lambda(1,\omega_{E}^{-1})} \calC_{E}(k) \calE_{E}(k). \\ 
        \calL_{0}(\phi) &= r_{0}^{2k+1}p^{d_p k} \dfrac{\Lambda(\rho_0, 1)}{2\Lambda(1, \eta) \Lambda(1,\omega_{0}^{-1})} \calC_{0}(k) \calE_{0}(k),
    \end{align}
    \end{subequations}
where $r_{f} := \Omega_{p, f, \chi_0}/\Omega_{\infty,f, \chi_0}$ and $r_{j} := \Omega_{p,j}/\Omega_{\infty,j}$ are quotients of periods. Consequently, the quotient $\calF_{i} := \calL_{i}/(\calL_E \calL_{0})$ satisfies
    \[
    \calF_{i}(\phi) = p^{-d_p} \Lambda(1,\omega_{E}^{-1}) \Lambda(1,\omega_{0}^{-1}) \dfrac{(\varphi_{1}, \overline{\varphi_{2}})}{(\varphi_3, \varphi_4)} \dfrac{\calC_{f}(k)}{\calC_{E}(k)\calC_{0}(k)} \dfrac{\calB_{i}(k)}{\calE_{E}(k) \calE_{0}(k)} \calR(\phi),
    \]
where $\calR(\phi) := \dfrac{r_{f}^{4k}}{r_{E}^{2k-1} r_{0}^{2k+1}}$.
\end{theorem}

We now consider the Heegner points attached to test vectors and their dual analogues. For $a \in \End^{0}(A)$, we write $[a]$ for the action of $a$ on $A$, and let $a$ act on $A^{\vee}$ by the dual endomorphism $[a]^{\vee}$. Fix an isogeny $\beta: A \rightarrow A^{\vee}$ satisfying
\[
\beta\circ[a]=[a^c]^\vee\circ\beta, \quad\text{ for $a \in \End^{0}(A)$},
\]
where $c$ is complex conjugation on $\End^{0}(A)$. At a common level group $U$, we represent $\varphi_i$ by $F_i: J_U \rightarrow A$ where $J_{U}$ is the Jacobian of the Shimura curve $X_{U}$ at level $U$. Then we set
\[
\varphi_i^\top:=\beta\circ F_i:J_U\rightarrow A^\vee,
\]
and in particular, $\varphi_i^{\top} \in \pi_{A^{\vee}} \simeq \pi_{A}^{\vee}$. Moreover, for any test vector $\varphi_{i}$, we see that for any finite place $v$, there exists a fixed scalar $c_v$ such that 
\begin{equation} \label{eq:compare-top-conjugation}
\overline{\varphi_v} = c_{v} \varphi_v^{\top},
\end{equation}
by Schur's lemma, for example, since $\pi_{v}$ and $\pi_{v}^{\vee}$ are irreducible. We put $P_i:=P_{\chi_0}(\varphi_i)$ and $P_i^\top:=P_{\chi_0^{-1}}(\varphi_i^\top)$ for \emph{Heegner points} attached to the test vector $\varphi_{i}$. 

Then we have the $p$-adic Waldspurger formula given by \cite[Theorem 4.2]{FanWan}:
\begin{theorem}[$p$-adic Waldspurger formula] \label{thm:fan-wan-p-adic-waldspurger}
    $\calL_i(0)=\log_{\omega_E}(P_i)^2$.\label{eq:two-padic-waldspurger}
\end{theorem}

\subsection{Gross--Zagier formula}
Next we recall the Gross--Zagier formula following \cite{YZZ}. Here we take the convention and formulation in \cite[Section 2C, Equation (2-4)]{CST}.

\begin{theorem}[Gross--Zagier formula] \label{thm:fan-wan-gross-zagier}
We have
\begin{subequations}
    \begin{align}
        \lrangle{P_1, P_2^{\top}}_{\calK,M} &= \dfrac{\zeta^{*}(2)}{4 \Lambda(1,\eta)^{2} \Lambda(1,\pi,\ad)} \Lambda^{\prime}(1,f,\chi_0) \dfrac{(\varphi_1, \varphi_2^{\top})_{U}}{\vol(X_{U})} \calB_{1,2}, \\ 
        \lrangle{P_2, P_1^{\top}}_{\calK,M} &= \dfrac{\zeta^{*}(2)}{4 \Lambda(1,\eta)^{2} \Lambda(1,\pi,\ad)} \Lambda^{\prime}(1,f,\chi_0) \dfrac{(\varphi_2, \varphi_1^{\top})_{U}}{\vol(X_{U})} \calB_{2,1}, 
    \end{align}
\end{subequations}
where the local integral $\calB_{i,j} := \prod_{v} \beta_{i,j,v}^{\GZ}$ with the local integral defined as
\[
\beta_{i,j,v}^{\GZ} := \dfrac{L(1,\eta_v) L(1,\pi_v, \ad)}{L(\frac{1}{2}, \pi_v, \chi_v)\zeta_v(2)} \int_{\calK_{v}^{\times}/\QQ_{v}^{\times}} \dfrac{(\pi(t_v)\varphi_{i,v}, \varphi_{j,v}^{\top})_v}{(\varphi_{i,v}, \varphi_{j,v}^{\top})_v} \chi(t_v) \dif t_{v}.
\]
Here the $L$-function $\Lambda(1,f,\chi_0)$ is the completed $L$-function.
\end{theorem}

\begin{proof}
This is precisely \cite[Section 2C, Equation (2-4)]{CST} with the local integral $\beta_{i,j,v}^{\GZ}$ being the local integral $\beta(\varphi_{i,v}, \varphi_{j,v}^{\top})$ using the notation of \textit{op.cit.}. There are two essential shifts of notation between \cite{FanWan} and \cite{YZZ}.  If $\pi_f$ is the cohomological representation generated by $f$, put $\sigma=\pi_f\otimes|\det|^{-1/2}$.  Then 
\begin{equation}\label{eq:unitary-shift}
 L(s,\sigma\times\chi_0)=L(f,\chi_0,s+1/2), \quad L(1,\sigma,\ad)=L(1,\pi_f,\ad).
\end{equation}
Thus the automorphic central point $s=1/2$ is the classical point $s=1$. We are using this shift to replace the term $L^{\prime} \left(\dfrac{1}{2}, \pi_{A}, \chi \right)$ in \textit{op.cit.} by $\Lambda'(1,f,\chi_0)$ after also restoring the archimedean factors.
\end{proof}

\subsection{Comparison of Gross--Zagier--Waldspurger}

\subsubsection{Comparison of Gross--Zagier--Waldspurger} \label{sec:compare-gzw}
We are now finally ready to state the desired comparison formula.

Let $P_0$ be a generator of the free rank-one $\calO_{K}$-module $E(K)_{\tf}$. Also write $P_{0}$ for its image in $(E \otimes_{\calO_{K}} \calO_{M})(K)$ induced by the inclusion $\ZZ \hookrightarrow \calO_{M}$. Let $k_1, k_2 \in M$ be such that $P_i=k_i P_0$ for $i = 1,2$.  

On the {logarithmic side}, since the logarithm is linear, we see that 
\[
\log_{\omega_{E}}(P_1) = k_1 \log_{\omega_{E}}(P_0), \quad \log_{\omega_{E}}(P_2) = k_2 \log_{\omega_{E}}(P_0).
\]
Substituting this into Theorem~\ref{thm:fan-wan-p-adic-waldspurger}, we obtain
\[
\calL_{1}(0) = k_1^{2} \log_{\omega_{E}}(P_0)^{2}, \quad \calL_{2}(0) = k_2^{2} \log_{\omega_{E}}(P_0)^{2}.
\]
On the {Néron--Tate pairing side}, note that $P_i^{\top} = k_i^{c} P_0^{\top}$ for $i = 1,2$. Therefore we see that
\[
 \lrangle{P_1,P_2^{\top}}_{\calK,M} =k_1k_2^c \cdot \lrangle{P_0,P_0^{\top}}_{\calK,M}, \quad \lrangle{P_2,P_1^{\top}}_{\calK,M} = k_2k_1^c \cdot {\lrangle{P_0,P_0^{\top}}_{\calK,M}}.
\]
Combining the two sides, we have 
\begin{equation} \label{eq:comparison_start}
\dfrac{\calL_1(0) \calL_2(0)}{\lrangle{P_1,P_2^{\top}}_{\calK,M} \cdot \lrangle{P_2,P_1^{\top}}_{\calK,M}} = \dfrac{k_1 k_2}{k_1^{c} k_2^{c}} \dfrac{L_0^{4}}{R_{0}^{2}}
\end{equation}
where $L_{0} := \log_{\omega_{E}}(P_0)$ and $R_0 := \lrangle{P_0, P_0^{\top}}_{\calK,M}$.

We substitute $\calL_i(0)=\calF_i(0)\calL_E(0)\calL_0(0)$ and apply the interpolation formula for $\calL_{0}(0)$:
\[
\calL_0(0)= \frac{r_0\calC_0(0)\calE_0(0)}{2\Lambda(1,\eta)\Lambda(1,\omega_0^{-1})} \Lambda(\psi_0,1).
\]
Moreover, since $L(\psi_E,1) = 0$, we have $\Lambda'(1,f,\chi_0)=\Lambda'(\psi_E,1)\Lambda(\psi_0,1)$. With these observations, applying Theorem \ref{thm:fan-wan-gross-zagier} yields
\begin{equation}\label{eq:semi-comparision-with-Theta}
 \left(\frac{\calL_E(0)}{\Lambda'(\psi_E, 1)}\right)^2
 =\Theta\,\frac{k_1k_2}{k_1^ck_2^c}\frac{L_0^4}{R_0^2},
\end{equation}
where
\[
\Theta:= \left( \frac{\zeta^{*}(2)}{4\Lambda(1,\eta)^2\Lambda(1,\pi,\ad)} \right)^{2}\frac{\calG_{12} \calG_{21}}
 {\calF_1(0)\calF_2(0)} \left( \frac{r_0\calC_0(0)\calE_0(0)}{2\Lambda(1,\eta)\Lambda(1,\omega_0^{-1})}  \right)^{-2}, \quad \calG_{ij}:=\frac{(\varphi_i,\varphi_j^\top)_U}{\vol(X_U)}\calB_{i,j}.
\]

\subsubsection{Cancellation of the finite toric pairings} \label{sec:compare-finite-toric}
Recall we have defined the local toric integral $\alpha_{v}$ in \eqref{eq:define-alpha-v}. We put
\[
\calA_{i,j, v}(0) := \alpha_{v}(\varphi_{i,v}, \varphi_{j,v}^{\top}; \mu_{v}(0)).
\]

\begin{lemma}\label{lem:rank-one-toric-Nov}
    For any finite rational place $v$ of $\calK$, we have the equality
    \[
    \calA_{1,1,v}(0) \calA_{2,2,v}(0) = \calA_{1,2,v}(0) \calA_{2,1,v}(0).
    \]
\end{lemma}
\begin{proof}
We recall the definition of $\alpha_v$ in \eqref{eq:define-alpha-v}:
\[
\alpha_v(\varphi_v, \varphi_v^{\vee}; \mu_v) := \dfrac{L_v(1,\eta) L_v(1,\sigma,\ad)}{\zeta_v(2) L_v(1/2, \sigma_{v}\times\mu_{v})} \int_{\QQ_v^{\times} \bs \calK_{v}^{\times}} (\pi_v(t) \varphi_v, \varphi_v^{\vee}) \mu_v(t) \dif^{\times}t_v.
\]
The constant on the right-hand side is the same number for all four pairs $(i,j)$. It therefore remains to deal with the integral
\[
\talpha_v(\varphi_v, \varphi_v^{\vee}; \mu_v) := \int_{\QQ_v^{\times} \bs \calK_{v}^{\times}} (\pi_v(t) \varphi_v, \varphi_v^{\vee}) \mu_v(t) \dif^{\times}t_v.
\]
For any fixed $s_v \in \calK_v^{\times}$, changing the variable from $t_v$ to $t_v s_v$ gives
\[
\talpha_v(\pi_{v}(s_v)\varphi_v, \varphi_v^{\vee}; \mu_v) = \mu_v(s_v)^{-1} \talpha_v(\varphi_v, \varphi_v^{\vee}; \mu_v).
\]
Therefore, for any fixed $\varphi_v^{\vee}$, the map $\varphi_v \mapsto \talpha_v(\varphi_v, \varphi_v^{\vee}; \mu_v)$ belongs to $\Hom_{\calK_v^{\times}}(\pi_v, \mu_v^{-1})$. The local toric multiplicity-one theorem (by a theorem of Tunnell and Saito \cite{Tunnell,Saito}) says that this space has dimension at most one.\footnote{This method is also applied in \cite[\S 3A]{CST}.} For $\varphi_v^{0}, \varphi_v^{\vee,0}$ such that $\talpha_v(\varphi_v^{0}, \varphi_v^{\vee,0}; \mu_v)$ is nonzero, any two functionals
\[
\varphi_{v} \mapsto \talpha_v(\varphi_v, \varphi_v^{\vee}; \mu_v), \quad \varphi_v \mapsto \talpha_v(\varphi_v, \varphi_v^{\vee,0}; \mu_v)
\]
are proportional. Evaluating at $\varphi_{v}^{0}$ determines the constant, which gives
\[
\talpha_v(\varphi_v, \varphi_v^{\vee}; \mu_v) = \dfrac{\talpha_v(\varphi_v, \varphi_v^{\vee,0}; \mu_v) \talpha_v(\varphi_{v}^{0}, \varphi_v^{\vee}; \mu_v)}{\talpha_v(\varphi_{v}^{0}, \varphi_v^{\vee,0}; \mu_v)}.
\]
As a result, we see that
\begin{align*}
\talpha_v(\varphi_{1,v}, \varphi_{1,v}^{\top}; \mu_v) &\talpha_v(\varphi_{2,v}, \varphi_{2,v}^{\top}; \mu_v) \\ 
&= \dfrac{\talpha_v(\varphi_{1,v}, \varphi_v^{\vee,0}; \mu_v) \talpha_v(\varphi_{v}^{0}, \varphi_{1,v}^{\top}; \mu_v)}{\talpha_v(\varphi_{v}^{0}, \varphi_v^{\vee,0}; \mu_v)} \dfrac{\talpha_v(\varphi_{2,v}, \varphi_v^{\vee,0}; \mu_v) \talpha_v(\varphi_{v}^{0}, \varphi_{2,v}^{\top}; \mu_v)}{\talpha_v(\varphi_{v}^{0}, \varphi_v^{\vee,0}; \mu_v)} \\
&= \dfrac{\talpha_v(\varphi_{1,v}, \varphi_v^{\vee,0}; \mu_v) \talpha_v(\varphi_{2,v}, \varphi_v^{\vee,0}; \mu_v)}{\talpha_v(\varphi_{v}^{0}, \varphi_v^{\vee,0}; \mu_v)} \dfrac{ \talpha_v(\varphi_{v}^{0}, \varphi_{2,v}^{\top}; \mu_v) \talpha_v(\varphi_{v}^{0}, \varphi_{1,v}^{\top}; \mu_v)}{\talpha_v(\varphi_{v}^{0}, \varphi_v^{\vee,0}; \mu_v)} \\
&= \talpha_v(\varphi_{1,v}, \varphi_{2,v}^{\top}; \mu_v) \talpha_v(\varphi_{2,v}, \varphi_{1,v}^{\top}; \mu_v).
\end{align*}
This proves the lemma.
\end{proof}

Multiplying \eqref{eq:define-beta-tilde} over all places, the global Whittaker function formula in \cite[Proposition 2.1]{CST} applied to the pairing $(\varphi_3, \varphi_4)$ gives
\begin{equation} \label{eq:waldspurger-CST}
\dfrac{(\varphi_{1}, \overline{\varphi_{2}})}{(\varphi_3, \varphi_4)} \calC_{f}(k) \calB_{i}(k) = \dfrac{\zeta^{*}(2)}{2\Lambda(1,\pi,\ad)} \calA_{\infty}(k) \prod_{v < \infty}\calA_{i,i,v}(k).
\end{equation}

We put
\[
\calJ_{\geom}:=
\frac{(\varphi_i,\varphi_j^\top)_U \beta_{i,j,\infty}^{\GZ}}
{\vol(X_U)\prod_{v<\infty} (\varphi_{i,v},\varphi_{j,v}^\top)}.
\]
The next lemma computes this term. We can see from the definition of $\calJ_{\geom}$ that it is actually the factor comparing two ways of measuring the same pair of test vectors $\varphi_{i}$ and $\varphi_{j}^{\top}$, one by the geometric pairing on the Shimura curve $X_{U}$ and the other by their Petersson pairing as automorphic forms.

To state the result, we define 
\[
I_{E} := \int_{E(\CC)} i \omega_{E} \wedge \overline{\omega_{E}} = \int_{E(\CC)} \abs{\omega_{E} \wedge \overline{\omega_{E}}},
\]
where $\omega_E$ is an invariant differential on $E$.

\begin{lemma}\label{lem:compute-J-geom}
We have $\calJ_{\geom} = 16\pi^{2}/I_{E}$.
\end{lemma}

\begin{proof}
    Let $X_{U}$ be the Shimura curve of level $U$, and let $J_{U} := \Jac(X_{U})$ be its Jacobian, with $\lambda_{U}: J_U \rightarrow J_{U}^{\vee}$ its principal polarization. Then the two test vectors $\varphi_i$ and $\varphi_j^{\top}$ give two maps
    \[
    F_i: J_U \rightarrow A, \quad G_{j}: J_U \rightarrow A^{\vee},
    \]
    and their geometric pairing $(\varphi_i, \varphi_j^{\top})_{U}$ is given by the composition
    \[
    q_{ij, U}:= F_{i} \circ \lambda_{U}^{-1} \circ G_{j}^{\vee} \in \End^{0}(A)
    \]
    Here the latter $\End^{0}(A)$ is a field and by fixing an embedding $\tau: \End^{0}(A) \rightarrow \CC$, the element $q_{ij, U}$ is identified with $(\varphi_i, \varphi_j^{\top})_{U}$. This follows the convention in \cite[\S 1B]{CST}.

    For the fixed invariant differentials $\omega_{A}$ and $\omega_{A^{\vee}}$ for $A$ and $A^{\vee}$, respectively, we have $\lrangle{\omega_{A}, \omega_{A^{\vee}}}_{\infty} = I_{E}$, and the adjunction of pullback with respect to $\lambda_{U}$ then gives 
    \begin{equation} \label{eq:differential-formula}
    \lrangle{F_{i}^{\ast}\omega_{A}, G_{j}^{\ast}\omega_{A^{\vee}}}_{X_{U}} = \lrangle{q_{ij,U}^{\ast}\omega_{A}, \omega_{A^{\vee}}}_{\infty} = (\varphi_i, \varphi_j^{\top})_{U} I_{E}.
    \end{equation}
    We write
    \[
    F_{i}^{\ast} \omega_{A} = 2 \pi i g_i(z) \dif z, \quad G_{j}^{\ast} \omega_{A^{\vee}}^{\top} = 2 \pi i g_j(z) \dif z.
    \]
    Then the equation \eqref{eq:differential-formula} becomes
    \begin{equation} \label{eq:geometric-pairing}
    8 \pi^{2} \int_{X_U} g_i(z) \overline{g_j(z)} \dif x \dif y = (\varphi_i, \varphi_j^{\top})_{U} I_{E},
    \end{equation}
    noting that $i \dif z \wedge \dif \barz = 2 \dif x \dif y$. Meanwhile, we emphasize that, in the convention of \cite{CST}, the volume $\vol(X_U)$ is computed using the measure $\dif x \dif y /(4 \pi y^{2})$; this gives the Petersson pairing
    \begin{equation} \label{eq:petersson-pairing}
       (\varphi_{i}, \varphi_{j}^{\top})_{\Pet} = \dfrac{2}{\vol(X_U)} \int_{X_U} g_i(z) \overline{g_j(z)} y^{2} \dfrac{\dif x \dif y}{4 \pi y^{2}} = \dfrac{1}{2 \pi \vol(X_U)} \int_{X_U} g_i(z) \overline{g_j(z)} \dif x \dif y. 
    \end{equation}
    Comparing \eqref{eq:geometric-pairing} and \eqref{eq:petersson-pairing}, we have
    \begin{equation} \label{eq:quotient-pairing}
    \dfrac{(\varphi_i, \varphi_j)}{(\varphi_i, \varphi_j^{\top})_{U}} = \dfrac{I_{E}}{16 \pi^{3} \vol(X_U)}.
    \end{equation}
    The remaining term is $\beta_{i,j,\infty}^{\GZ}$. With the convention in \cite{CST} (see, in particular, \cite[Lemma 3.14]{CST}), we see that the Euler factor appearing in the definition of $\beta_{i,j,\infty}^{\GZ}$ is
    \[
     \frac{L_\infty(1,\eta)L_\infty(1,\pi,\ad)}{\zeta_\infty(2)L_\infty(1/2,\pi\otimes\chi_0)}=\frac{\pi^{-1}(2\pi^3)^{-1}}{\pi^{-1}\Gamma_\CC(1)^2}=\frac1{2\pi},
    \]
    while the toric integral equals $2$ (see the proof of \cite[Lemma 3.14]{CST}). This altogether gives $\beta_{i,j,\infty}^{\GZ} = \pi^{-1}$. Combining with \eqref{eq:quotient-pairing}, the lemma is therefore proved.
\end{proof}

From Lemma \ref{lem:compute-J-geom}, we see that $\calJ_{\geom}$ is independent of the pair $(\varphi_{i}, \varphi_j)$, and hence we can compute $\calG_{i,j}$ as
\[
\calG_{i,j} = \dfrac{(\varphi_{i}, \varphi_j^{\top})}{\vol(X_U)} \calB_{i,j} = \calJ_{\geom} \prod_{v < \infty} \calA_{i,j,v}(0),
\]
and hence
\begin{align} 
\calG_{1,2} \calG_{2,1} = 
\dfrac{(\varphi_{1}, \varphi_{2}^{\top})(\varphi_{2}, \varphi_{1}^{\top})}{\vol(X_U)^{2}} \calB_{1,2} \calB_{2,1} &= \calJ_{\geom}^{2} \prod_{v < \infty} \calA_{1,2,v}(0) \calA_{2,1,v}(0) \notag \\&= \calJ_{\geom}^{2} \prod_{v < \infty} \calA_{1,1,v}(0) \calA_{2,2,v}(0). \label{eq:after-J-geom} 
\end{align}
where the last equality uses Lemma~\ref{lem:rank-one-toric-Nov}. Equation \eqref{eq:after-J-geom} is the most crucial identity in this part: it transfers the local integrals $\calA_{1,2}$ and $\calA_{2,1}$ appearing in the Gross--Zagier formula (Theorem \ref{thm:fan-wan-gross-zagier}) to the local integrals appearing in the $p$-adic Waldspurger formula (Theorems \ref{thm:fan-wan-p-adic-l-function}, \ref{thm:fan-wan-p-adic-waldspurger}) involving $\calA_{1,1}$ and $\calA_{2,2}$.

\subsubsection{Computation of the archimedean quotient and the $p$-adic quotient} \label{sec:compare-archimedean}
In this section, we compute the archimedean part of the quotient $\calF_{i}(k)$ in Theorem \ref{thm:fan-wan-p-adic-l-function}. We fix the common test vector $\varphi_{\infty} := \varphi_{1,\infty} = \varphi_{2,\infty}$ and also its dual vector $\varphi^{\vee}$. Then we put
\[
\calA_{\infty}(k) := \alpha_\infty \left( \Delta^{k-1} \varphi_\infty, \mu_\infty(k)\right).
\]

\begin{lemma} \label{lem:archimedean-computation-comparison}
    When $k \geq 1$, we have
    \[
    \calA_{\infty}(k) = (-1)^{k-1} 2^{3-2k}, \quad \calC_{E}(k) = i^{k-1}2^{2-k}, \quad \calC_{0}(k) = i^{k}2^{1-k}.
    \]
In particular, we have
    \begin{equation} \label{eq:archimedean-quotient}
    \dfrac{\calA_{\infty}(k)}{\calC_{E}(k)\calC_{0}(k)} = -i, \quad \calC_{0}(0) = 2.    
    \end{equation}
\end{lemma}

\begin{proof}
We use the same notation as in the proof of Lemma \ref{lem:rank-one-toric-Nov}. We start with the vector $\varphi_{\infty}$ with $\talpha_{\infty}(\varphi_{\infty}, \mu_{\infty}(k)) = 2$. Then by the definition of the raising operator $\Delta$, we have
\[
\talpha_{\infty}(\Delta^{\nu+1}\varphi_{\infty}, \mu_{\infty}(k)) = - \dfrac{(\nu+1)(\nu+2)}{(4\pi)^{2}} \talpha_{\infty}(\Delta^{\nu}\varphi_{\infty}, \mu_{\infty}(k)).
\]
Iterating gives
\[
\talpha_{\infty}(\Delta^{\nu}\varphi_{\infty}, \mu_{\infty}(k)) = 2 (-1)^{\nu} \dfrac{\nu! (\nu+1)!}{(4 \pi)^{2\nu}}.
\]
Inserting this into the definition \eqref{eq:define-alpha-v}, we obtain
\[
\calA_\infty(k) =\frac{1}{2\pi^3} \frac{\talpha_{\infty}(\Delta^{k-1}\varphi_{\infty}, \mu_{\infty}(k))}{\Gamma_{\CC}(k)\Gamma_{\CC}(k+1)} = (-1)^{k-1} 2^{3-2k}.
\]
For the two Eisenstein factors $\calC_{E}(k)$ and $\calC_{0}(k)$, we use the same normalized weight one Eisenstein section $h_{1}$ at $\infty$ that satisfies $h_1(1) = 1$. By the definition of the raising operator, we have
\begin{equation} \label{eq:raising-Eisenstein-section}
\Delta h_{2m+1} = \dfrac{i(m+1)}{4\pi} h_{2m+3}, \quad\text{ and hence }\Delta^{\nu} h_1 = \dfrac{i^{\nu} \nu!}{(4 \pi)^{\nu}} h_{2\nu+1}.
\end{equation}
Moreover, note that
\[
h_{2\nu+1}(r_{\theta})\mu_{\nu,\infty}(r_{\theta}) = e^{i(2\nu+1)\theta} e^{-i(2\nu+1)\theta} = 1, \quad\text{ with } r_{\theta} = \begin{bmatrix}
    \cos \theta & \sin \theta \\ -\sin \theta & \cos \theta
\end{bmatrix}.
\]
The toric integral in the definition of $\lambda_{\infty}(h_{2\nu+1})$ is the volume of $\RR^{\times} \bs \CC^{\times}$, which is $2$ in our convention following \cite{CST}. We then have
\[
\lambda_{\infty}(h_{2\nu+1}) = \dfrac{2\pi^{-1}}{\Gamma_{\CC}(\nu+1)}.
\]
It follows from \eqref{eq:raising-Eisenstein-section} that $\lambda_{\infty}(\Delta^{\nu}h_1) = i^{\nu} 2^{1-\nu}$. Taking $\nu = k-1$ and $\nu = k$ respectively, we have
\[
\calC_{E}(k) = i^{k-1}2^{2-k}, \quad \calC_{0}(k) = i^{k} 2^{1-k}.
\]
For $k = 0$, the raising index $\nu = 0$, so we have
\[
\calC_{0}(0) = \lambda_{0,\infty}(h_1) = 2 \dfrac{\pi^{-1}}{\Gamma_{\CC}(1)} = 2.
\]
The value of $\calA_{\infty}(k)/(\calC_{E}(k)\calC_{0}(k))$ follows directly from the computation above.
\end{proof}

Next, we put $\lambda_{j,v}(k) := \lambda_{j,v}(h_{j,v}(-g_{\ell,n}))$, and we write
\[
\calE_{j}^{(p)}(k) = \prod_{v \nmid p\infty} \lambda_{j,v}(k).
\]
As at the archimedean places, at the $p$-adic place, define $\calA_{p}(k) := \calA_{1,1,p}(k) = \calA_{2,2,p}(k)$ as the common $p$-vector. Meanwhile, put
\[
\calA_p^0 :=\alpha_p\!\left(\pi_p(g_{\ell,n})\varphi_{1,p}, \mu_p(0)\right),
\]
and define the quotient
\[
\calV_{p}(k) := \dfrac{\lambda_{E,p}(k) \lambda_{0,p}(k)}{\calA_{p}(k)}
\]
for $k \geq 1$ and $\calV_{p}(0)$ is the continuation of this quotient. This is well-defined since one can read off the term $\calV_{p}(k)$ in the formula of $\calF_{i}(k)$ in Theorem \ref{thm:fan-wan-p-adic-l-function}. We put
\[
r_{p}^{\tor} = p^{d_p} \dfrac{\calA_{p}^{0}}{\lambda_{0,p}(0)} \calV_{p}(0).
\]
Then we prove the following proposition.

\begin{proposition} \label{prop:simplify-Theta}
    Let $r := r_E$. The error term $\Theta$ can be simplified as
    \[
    \Theta = - \dfrac{\calJ_{\geom}^{2}}{4r^{2} \Lambda(1,\eta)^{2} \Lambda(1,\omega_{E}^{-1})^{2}} \cdot \calE_{E}^{(p)}(0)^{2} \cdot (r_p^{\tor})^{2}, 
    \]
i.e., the local integral factors away from $p$ and the rank-zero counterpart cancel in the defining equation of the error term $\Theta$. As a result,
\[
\left(\dfrac{\calL_{E}(0)}{\Lambda^{\prime}(\psi_E, 1)} \right)^{2}  = - \dfrac{\calJ_{\geom}^{2}}{4r^{2} \Lambda(1,\eta)^{2} \Lambda(1,\omega_{E}^{-1})^{2}} \cdot \calE_{E}^{(p)}(0)^{2} \cdot (r_p^{\tor})^{2} \cdot \dfrac{k_1 k_2}{k_1^{c} k_2^{c}} \dfrac{L_0^{4}}{R_{0}^{2}}.
\]
\end{proposition}

\begin{proof}
    Using \eqref{eq:waldspurger-CST} and \eqref{eq:archimedean-quotient} in the factorization formula in Theorem \ref{thm:fan-wan-p-adic-l-function}, we have
    \[
    \calF_{i}(k) = -i \dfrac{\zeta^{*}(2) \Lambda(1,\omega_{E}^{-1})\Lambda(1,\omega_{0}^{-1})}{2 \Lambda(1,\pi,\ad)} p^{-d_p} \dfrac{\prod_{v \nmid p\infty} \calA_{i,i,v}(k)}{\calE_{E}^{(p)}(k) \calE_{0}^{(p)}(k)} \calV_p(k)^{-1}\calR(k).
    \]
    Multiplying this equality for $i=1, 2$, we have
    \begin{equation} \label{eq:F1-times-F2}
    \calF_{1}(0)\calF_2(0) = -\left( \dfrac{\zeta^{*}(2) \Lambda(1,\omega_{E}^{-1})\Lambda(1,\omega_{0}^{-1})}{2 \Lambda(1,\pi,\ad)} p^{-d_p} \right)^{2} \dfrac{\prod_{v \nmid p\infty} \calA_{1,1,v}(0)\calA_{2,2,v}(0)}{\calE_{E}^{(p)}(0)^{2} \calE_{0}^{(p)}(0)^{2}} \calV_p(0)^{-2}\calR(0)^2.
    \end{equation}
    Next, the following term appearing in the interpolation formula of $\calL_0$ at $0$ can be simplified:
    \[
    \dfrac{r_0}{2\Lambda(1,\eta)\Lambda(1,\omega_0^{-1})} \calC_0(0)\calE_0(0) = \dfrac{r_0}{\Lambda(1,\eta)\Lambda(1,\omega_0^{-1})} \calE_{0}^{(p)}(0) \lambda_{0,p}(0)
    \]
    using $\calC_{0}(0) = 2$.

    Substituting \eqref{eq:after-J-geom} into the definition of $\Theta$ gives
    \[
    \Theta = \dfrac{\zeta^{*}(2)^{2} \Lambda(1, \omega_{0}^{-1})^{2} }{r_{0}^{2} 4\Lambda(1,\eta)^{2} \Lambda(1,\pi,\ad)^{2}} \dfrac{1}{\calF_1(0) \calF_2(0) \calC_{0}(0)^{2} \calE_{0}(0)^{2}} \calJ_{\geom}^{2} \prod_{v < \infty} \calA_{1,1,v}(0) \calA_{2,2,v}(0),
    \]
    Using \eqref{eq:F1-times-F2}, we further obtain
    \[
    \Theta = \dfrac{1}{-4r_{0}^{2}\calR(0)^2 \Lambda(1,\eta)^{2} \Lambda(1,\omega_{E}^{-1})^{2}p^{-2d_p}} \dfrac{\calE_{E}^{(p)}(0)^{2} \calV_p(0)^{2}}{\lambda_{0,p}(0)^{2}} \calJ_{\geom}^{2} (\calA_{p}^{0})^{2}
    \]
    Here we use $\calC_{0}(0) = 2$ and note that $\calE_{0,p}(0) = \lambda_{0,p}(0)$, together with $\calA_{p}^{0} = \calA_{1,1,p}(0) = \calA_{2,2,p}(0)$ for the fixed common local vector. Rearranging the terms gives
    \[
    \Theta = - \dfrac{\calJ_{\geom}^{2}}{4r^{2} \Lambda(1,\eta)^{2} \Lambda(1,\omega_{E}^{-1})^{2}} \cdot \calE_{E}^{(p)}(0)^{2} \cdot (r_p^{\tor})^{2}, \quad r_{p}^{\tor} := \dfrac{p^{d_p} \calA_{p}^{0}}{\lambda_{0,p}(0)} \calV_{p}(0).
    \]
    Here $r=r_E$ and $\calR(0)=r_E/r_0$, so $r_0^2\calR(0)^2=r^2$. The factor $r_p^{\tor}$ gathers local terms at $p$.
\end{proof}

We renormalize the $p$-adic $L$-function $\calL_{E}$ by
\begin{equation} \label{eq:renormalize-calL_E}
\scL_{E}(\phi) := \dfrac{\Lambda(1,\eta) \Lambda(1,\omega_{E}^{-1})}{4\calE_{E}^{(p)}(\phi)} \cdot \calL_{E}(\phi),
\end{equation}
then we have
\[
\left(\dfrac{\scL_{E}(0)}{\Lambda^{\prime}(\psi_E, 1)} \right)^{2}  = - \dfrac{\calJ_{\geom}^{2}}{4 \cdot 16 r^{2} } \cdot (r_p^{\tor})^{2} \cdot \dfrac{k_1 k_2}{k_1^{c} k_2^{c}} \dfrac{L_0^{4}}{R_{0}^{2}} = - \dfrac{4 \pi^{4}}{r^{2} I_{E}^{2}} \cdot (r_p^{\tor})^{2} \cdot \dfrac{k_1 k_2}{k_1^{c} k_2^{c}} \dfrac{L_0^{4}}{R_{0}^{2}}.
\]
Here, for the second equality, we use Lemma \ref{lem:compute-J-geom}. Finally, since $\Lambda^{\prime}(\psi_E, 1) = \pi^{-1} L^{\prime}(\psi_E,1)$, we have
\[
\left(\dfrac{\scL_{E}(0)}{L'(\psi_E,1)} \right)^{2}  = \dfrac{-4 \pi^{2}}{r^{2} I_{E}^{2}} \cdot (r_p^{\tor})^{2} \cdot \dfrac{k_1 k_2}{k_1^{c} k_2^{c}} \dfrac{L_0^{4}}{R_{0}^{2}} = \left(\dfrac{r_p^{\tor}}{\hatr I_{E}}\right)^{2} \dfrac{k_1 k_2}{k_1^{c} k_2^{c}} \dfrac{L_0^{4}}{R_{0}^{2}},
\]
where $\hatr := r / (2\pi i)$. We have therefore proved the following result.

\begin{theorem}[Gross--Zagier--Waldspurger comparison] \label{thm:final_comparison}
    With the preceding conventions,
    \[
    \left(\dfrac{\scL_{E}(0)}{L'(\psi_E,1)} \right)^{2} = \left(\dfrac{r_p^{\tor}}{\hatr I_{E}}\right)^{2} \dfrac{k_1 k_2}{k_1^{c} k_2^{c}} \dfrac{L_0^{4}}{R_{0}^{2}}.
    \]
\end{theorem}

%% file: 105-magic.tex
\section{The $p$-adic BSD defect formula and the transfer principle} \label{sec:magic}
This section derives the $p$-adic BSD defect formula and the transfer argument from Theorem~\ref{thm:final_comparison} and the arithmetic identities in Section~\ref{sec:plus-minus}.

\subsection{$p$-adic BSD defect formula} \label{sec:BSD-defect-formula}
Use the coefficient embedding $\iota$ of Section~\ref{sec:setup} and put $\iota^c=\iota\circ c$. For each $\sigma \in \{\iota,\iota^{c} \}$, we carry out the previous construction and let $\scL_{E, \sigma}$ be the renormalized anticyclotomic $p$-adic $L$-function of $E$ in Theorem \ref{thm:fan-wan-p-adic-l-function}, and let $\delta_{\sigma} := \delta(E, \scL_{\sigma})$ be defined as in \eqref{eq:delta} for the renormalized $\scL_{E,\sigma}$. We write $\frl_\sigma := \Fitt_{\calO_L} \coker \left(\rmH_{\rmf}^1(\calK,T_\sigma)_{\tf} \rightarrow \rmH_{\rmf}^1(\calK_{\frp},T_\sigma)_{\tf}\right)$ for short.

On the arithmetic side, we have the following proposition.

\begin{proposition}[Arithmetic side identity] \label{prop:embedding-arithmetic}
We have the identity
    \[
    (\scL_{E,\iota}(0) \cdot \scL_{E,\iota^{c}}(0)) = \left( \dfrac{\# \Sha(E/\calK)}{\# E(\calK)[p^{\infty}]^{2}} \prod_{v \nmid p} \Tam_{v}(E/\calK) \cdot \# E(\calK_{\frp})[p^{\infty}] \right) (\frl_{\iota} \frl_{\iota^{c}})^{2} \dfrac{\delta_{\iota}^{\prime\prime} \delta_{\iota^{c}}^{\prime\prime}}{\delta_{\iota} \delta_{\iota^{c}}}.
    \]
\end{proposition}
\begin{proof}
    This is a consequence of Theorem \ref{thm:fan_wan_tnc} by taking the product of the equalities under two embeddings $\iota$ and $\iota^{c}$. The appearance of the term $E(\calK)[p^{\infty}]$ is then given by Lemma \ref{lem:torsion-appear}. The appearance of Tamagawa numbers is given by Lemma \ref{lem:tamagawa-appear} together with Lemma \ref{lem:order-norm}. The appearance of the Tate--Shafarevich group term is a combination of Propositions \ref{prop:length} and  \ref{prop:find_sha}, together with Lemma \ref{lem:order-norm}.
\end{proof}

Next, we focus on the analytic side. For $e_{k} := k_1 k_2/(k_1^{c} k_2^{c})$, there is an equality in the coefficient field: $e_{k}^{c} = e_{k}^{-1}$. Thus, after taking the embedding, we have $\iota(e_{k}) \iota^{c}(e_{k}) = 1$. Choose a generator $\omega$ of the module of invariant differentials of the N\'eron model $\calE/\calO_{\calK}$ of $E/\calK$. Let
\[
\int_{\rmH_{1}(E(\CC),\ZZ)} \omega = a_{E} \calO_{\calK},
\]
then we have
\begin{equation} \label{eq:Omega-E-over-K}
\Omega_{E/\calK} = a_{E} \bara_{E}, \quad I_{E} = \sqrt{\abs{D_{\calK}}} a_E \bara_{E}.
\end{equation}
Since $E$ is defined over $\QQ$, write $\omega^c=\alpha\omega$ for $\alpha\in\calK^\times$. Moreover, $\bar a_E\calO_\calK=\alpha a_E\calO_\calK$, so $\bar a_E/a_E\in\alpha\calO_\calK^\times$ is algebraic. Thus the quotient $a_{E}/\bara_{E}$ is algebraic, and $(\bara_{E}/a_{E})^{c} = (\bara_{E}/a_{E})^{-1}$. Choose $P=P_0$ to be a generator of $E(\calK)_{\tf}$ over $\calO_{\calK}$. For each $\sigma$, put $A_\sigma=H_f^1(\calK_{\frp},T_\sigma)_{\tf}$, a free rank-one $\calO_L$-module. Choose a generator $u_\sigma$ of $A_\sigma$. We write
\[
\loc_{p}(P) = a_{\sigma} u_{\sigma}, \quad \log_{\omega}(u_{\sigma}) = C_{\sigma}^{-1}.
\]
Thus $(a_{\sigma}) = \frl_{\sigma}$. We put $h := b_{\calK}(P,P)$. Then we have
\begin{equation} \label{eq:R0-L0}
L_{0,\sigma} = \ell_{\sigma}(P) = a_{\sigma}/C_{\sigma}, \quad R_{0,\sigma} = h/2.
\end{equation}
Additionally, we define
\begin{equation} \label{eq:pi_sigma-and-delta-prime-sigma}
    \pi_{\sigma} := \hatr_{\sigma} a_{E,\sigma}, \quad \delta^{\prime}_{\sigma} = C_{\sigma} \pi_{\sigma} = C_{\sigma} \hatr_{\sigma} a_{E,\sigma},
\end{equation}
and call $\delta_{\sigma}^{\prime}$ the \emph{first $p$-adic defect}.

We are then ready to state and prove the following analytic side identity.

\begin{proposition}[Analytic side identity] \label{prop:embedding-analytic}
We have the following identity
    \[
    \left( \dfrac{L^{\ast}(E/\calK, 1)}{\Omega_{E/\calK} \Reg(E/\calK)} \right) (\frl_{\iota}\frl_{\iota^{c}})^{2} = (\scL_{E,\iota}(0) \cdot \scL_{E,\iota^{c}}(0)) \prod_{\sigma = \iota, \iota^{c}} \left( \dfrac{C_{\sigma} \delta_{\sigma}^{\prime}}{r_{p,\sigma}^{\tor}} \right).
    \]
\end{proposition}

\begin{proof}
    Recall that by Theorem \ref{thm:final_comparison}, we have
    \[
    \left(\dfrac{\scL_{E, \sigma}(0)}{L^{\prime}(1,\psi_{E})} \right)^{2} = \left(\dfrac{r_{p,\sigma}^{\tor}}{\hatr_{\sigma} I_E}\right)^{2} \dfrac{k_1 k_2}{k_1^{c} k_2^{c}} \dfrac{L_{0,\sigma}^{4}}{R_{0,\sigma}^{2}}.
    \]
Using \eqref{eq:R0-L0} together with Lemma \ref{lem:CM-regulator}, we have 
\[
 \left(\frac{\scL_{E,\sigma}(0)}{L'(\psi_E,1)}\right)^2
 =\left(\frac{r_{p,\sigma}^{\tor}}{\hatr_\sigma I_E}\right)^2
 \sigma(e_k)\frac{4a_\sigma^4}{C_\sigma^4h^2}
 =\left(\frac{r_{p,\sigma}^{\tor}}{\hatr_\sigma}\right)^2
 \sigma(e_k)\frac{a_\sigma^4}{C_\sigma^4\Reg(E/\calK)\Omega_{E/\calK}^2}.
\]
Here, by \eqref{eq:Omega-E-over-K}, $I_E^2=D_\calK\Omega_{E/\calK}^2$. Rearranging and noting that $L^{\ast}(E/\calK) = L^{\prime}(\psi_E,1)^{2}$ gives
\[
 \dfrac{L^*(E/\calK,1)}{\Omega_{E/\calK}\Reg(E/\calK)}=\scL_{E,\sigma}(0)^2
 \left(\frac{\hatr_\sigma}{r_{p,\sigma}^{\tor}}\right)^2
 \sigma(e_k)^{-1}\frac{C_\sigma^4\Omega_{E/\calK}}{a_\sigma^4}.
\]
Applying \eqref{eq:pi_sigma-and-delta-prime-sigma}, we can write the above equation as
\[
\dfrac{L^{\ast}(E/\calK, 1)}{\Omega_{E/\calK} \Reg(E/\calK)} = \scL_{E,\sigma}(0)^{2} \left(\dfrac{C_{\sigma} \delta_{\sigma}^{\prime}}{r_{p,\sigma}^{\tor}}\right)^{2} \sigma(e_k)^{-1} \dfrac{ \sigma(a_{E}^{c}/a_{E})}{a_{\sigma}^{4}}.
\]
Multiplying the above identities for $\iota$ and $\iota^{c}$, we obtain
\[
\left( \dfrac{L^{\ast}(E/\calK, 1)}{\Omega_{E/\calK} \Reg(E/\calK)} \right)^{2} (a_{\iota} a_{\iota^{c}})^{4} = \scL_{E,\iota}(0)^{2}\scL_{E,\iota^{c}}(0)^{2} \prod_{\sigma = \iota, \iota^{c}} \left(\dfrac{C_{\sigma} \delta_{\sigma}^{\prime}}{r_{p,\sigma}^{\tor}}\right)^{2}.
\]
Pass to principal fractional $\calO_L$-ideals. Their group is torsion-free, so equality of squares implies equality of the ideals themselves. Thus
\[
\left( \dfrac{L^{\ast}(E/\calK, 1)}{\Omega_{E/\calK} \Reg(E/\calK)} \right) (a_{\iota} a_{\iota^{c}})^{2} = \scL_{E,\iota}(0) \scL_{E,\iota^{c}}(0) \prod_{\sigma = \iota, \iota^{c}} \left(\dfrac{C_{\sigma} \delta_{\sigma}^{\prime}}{r_{p,\sigma}^{\tor}}\right).
\]
As $(a_{\sigma}) = \frl_{\sigma}$ for $\sigma \in \{\iota, \iota^{c}\}$, this proves the proposition as desired.
\end{proof}

Combining Propositions \ref{prop:embedding-arithmetic} and \ref{prop:embedding-analytic}, we have (note that $(\frl_{\iota} \frl_{\iota^{c}})$ cancels out):
\begin{equation} \label{eq:all-error-appear}
    \dfrac{L^{\ast}(E/\calK, 1)}{\Omega_{E/\calK} \Reg(E/\calK)}= \left( \dfrac{\# \Sha(E/\calK)}{\# E(\calK)_{\tors}^{2}} \prod_{v \nmid p} \Tam_{v}(E/\calK) \cdot \# E(\calK_{\frp})[p^{\infty}] \right) \prod_{\sigma = \iota, \iota^{c}} \left( \dfrac{C_{\sigma} \delta_{\sigma}^{\prime} \delta_{\sigma}^{\prime\prime}}{r_{p,\sigma}^{\tor} \delta_{\sigma}} \right).
\end{equation}

To continue, we have the following lemma.
\begin{lemma} \label{lem:tam-at-p-appear}
    We have
    \[
    \prod_{\sigma = \iota, \iota^{c}} (C_{\sigma}) = \left( \dfrac{\Tam_{\frp}(E/\calK) \calE_{p}(E/\calK)}{\# E(\calK_{p})[p^{\infty}]} \right),
    \]
where $\calE_{p}(E/\calK) = P_{\frp}(q^{-1}) = \# \tilE_0(\kappa_{\frp})/q$.
\end{lemma}

\begin{proof}
Let $F := \calK_{\frp}$, let $\calE$ be the Néron model over $\calO_{F}$, and let $E_0(F)$ be the group of points in $E(F)$ reducing to the identity component of $\calE/\calO_{F}$, with $E_{1}(F) = \ker(E_0(F) \rightarrow \tilE_0(\kappa_{\frp}))$. For $n \geq 1$, let $E_{n}(F)$ consist of the points whose formal group parameter $X := -x/y$ belongs to $\frp^{n}$, where $x, y$ are the variables in the minimal Weierstrass equation for $E/F$. Then there are successive quotients
\[
\#(E(F)/E_0(F)) = \Tam_{\frp}(E/\calK), \, \# (E_0(F)/E_1(F)) = \# \tilE_{0}(\kappa_{\frp}), \, \#(E_{m}(F)/E_{m+1}(F)) = q_{\frp}.
\]
Multiplying these orders gives $[E(F):E_n(F)] = \Tam_{\frp}(E/\calK) \cdot \# \tilE_{0}(\kappa_{\frp}) \cdot q_{\frp}^{n-1}$.

We write $M = (E(F)_{p}^{\wedge})_{\tf}$. Choose $n$ sufficiently large so that the formal group logarithm is an isomorphism $\log_{\omega}: E_{n}(F) \simeq \frp^{n}$. In particular, $E_{n}(F)$ is torsion-free. Consequently, we have
\[
[\log_{\omega}(M):\frp^{n}] = [M : E_n(F)] = \dfrac{(\Tam_{\frp}(E/\calK) \cdot \# \tilE_{0}(\kappa_{\frp}) \cdot q_{\frp}^{n-1})_{p}}{\# E(F)[p^{\infty}]}. 
\]
Taking the norm, we have
\[
\Norm_{F/\QQ_p}(\log_{\omega}(M)) = \left( \dfrac{q_{\frp}^{n} \# E(F)[p^{\infty}]}{\Tam_{\frp}(E/\calK) \cdot \# \tilE_{0}(\kappa_{\frp}) \cdot q_{\frp}^{n-1}} \right) = \left( \dfrac{q_{\frp}\# E(F)[p^{\infty}]}{\Tam_{\frp}(E/\calK) \cdot \# \tilE_{0}(\kappa_{\frp})} \right)
\]
as fractional $\ZZ_p$-ideals. By definition of $C_{\sigma}$, we see that $(C_{\sigma}) = (\sigma(\log_{\omega}(M)) \calO_{L})^{-1}$. Multiplying the two embeddings $\sigma = \iota, \iota^{c}$ then gives
\[
\prod_{\sigma = \iota, \iota^{c}} (C_{\sigma}) = \Norm_{F/\QQ_p}(\log_{\omega}(M)^{-1})\calO_{L} = \left( \dfrac{\Tam_{\frp}(E/\calK) \cdot \# \tilE_{0}(\kappa_{\frp})}{q_{\frp}\# E(F)[p^{\infty}]} \right).
\]
Finally, the local Euler factor satisfies $\calE_{p}(E/\calK) = \# \tilE_{0}(\kappa_{\frp})/q_{\frp}$; this gives the lemma, as desired.
\end{proof}

Combining \eqref{eq:all-error-appear} and Lemma \ref{lem:tam-at-p-appear}, we obtain the following theorem, which we call the \emph{$p$-adic BSD defect formula}.

\begin{theorem}[$p$-adic BSD defect formula] \label{thm:p-adic-defect}
    We have
    \[
    \Delta_{p}(E/\calK) = v_p(\calE_{p}(E/\calK)) + \sum_{\sigma = \iota, \iota^{c}} v_p \left( \dfrac{\delta_{\sigma}^{\prime} \delta_{\sigma}^{\prime\prime}}{r_{p,\sigma}^{\tor} \delta_{\sigma}} \right)
    \]
\end{theorem}

\subsection{Transfer principle} \label{sec:transfer-principle}
We define the following relation between CM elliptic curves.

\begin{definition} \label{defn:matched}
Let $E_{1}/\QQ$ and $E_{2}/\QQ$ be two CM elliptic curves, with $\ord_{s=1} L(E_{j}/\QQ, s) = 1$ for $j=1, 2$. We say that $E_{1}$ and $E_{2}$ are \emph{matched at $p$} if
\begin{enumerate}[label = \rm (\arabic*)]
    \item For $j=1,2$, each $E_j$ has CM by the maximal order $\calO_{\calK_j}$ for some imaginary quadratic field $\calK_{j}$.
    \item The prime $p$ is nonsplit in $\calK_1$ and $\calK_2$, and there is a $\QQ_p$-field isomorphism $\calK_{1, \frp_1} \simeq \calK_{2, \frp_2} =: F$.
    \item After the identification in (2), their local Galois representations
    \[
    T_p E_1|_{\Gal_F}\simeq T_p E_2|_{\Gal_F}
    \]
    are isomorphic as $\calO_{F}[\Gal_{F}]$-modules. Equivalently, the local characters $\psi_{E_{1}, p} = \psi_{E_{2},p}$ are equal after the identification in (2).
\end{enumerate}
\end{definition}

Choose $u=u_\sigma$ generating the torsion-free lattice $A_\sigma=\rmH_{\rmf}^1(\calK_\frp,T_\sigma)_{\tf}$, and let $z_\sigma$ be the elliptic unit attached to $\psi$ (at the origin point of $\bpsi$). In $A_\sigma\otimes L$, write
\begin{equation}\label{eq:define-w-sigma}
 \frac{\loc_p(z_\sigma)}{\scL_{E,\sigma}(0)}=w_{E,\sigma}(0)u
\end{equation}
and moreover we have
\begin{equation}\label{eq:Nov-family-w}
 \dfrac{\loc_p(z_{k})}{\scL_{E,\sigma}(k)}= w_{E,\sigma}(k) u(k).
\end{equation}
Here, $u(k) \in \rmH^{1}_{+}(\calK_{\frp}, \bpsi)_{\tf}$ satisfies $u(0) = u$, and $z_{k}$ is given by the element $z_{\Psi_{k},\frf}(\gamma_{E,k})$ in Section \ref{sec:ERC-basis}.

To check the defects, choose a basis $q$ of $Q_\sigma$ and write $q_-(u)=d q$. Then $\delta''_\sigma=(d)$ by \eqref{eq:delta2}. Applying $q_-$ to \eqref{eq:define-w-sigma} gives $q_-\loc_p(z_\sigma)=\scL_{E,\sigma}(0)w_{E,\sigma}(0)dq$. The definition \eqref{eq:delta} therefore gives $\delta_\sigma=(w_{E,\sigma}(0))\delta''_\sigma$. It follows that
\[
\Delta_p(E/\calK)=v_p(\calE_p(E/\calK))+
\sum_{\sigma=\iota,\iota^c}v_p\!\left(
\frac{\delta'_\sigma}{r_{p,\sigma}^{\tor}w_{E,\sigma}(0)}\right).
\]

Now suppose that $E_1$ and $E_2$ are matched at $p$. We choose the same local test vectors at $p$, so $r_{E_1,p,\sigma}^{\tor}=r_{E_2,p,\sigma}^{\tor}$. The local Euler polynomials $\calE_{p}(E/\calK)$ for $E_1$ and $E_2$ agree by the identified local Galois representations at $p$. We therefore have two remaining factors, as stated in the lemma below.

\begin{lemma} \label{lem:transfer-prepare}
Let $E_{1}/\QQ$ and $E_{2}/\QQ$ be two CM elliptic curves, with $\ord_{s=1} L(E_{j}/\QQ, s) = 1$ for $j=1, 2$, such that $E_{1}$ and $E_{2}$ are {matched at $p$} in the sense of Definition \ref{defn:matched}. Then
\[
\Delta_{p}(E_1/\calK_1) - \Delta_{p}(E_{2}/\calK_2) = \sum_{\sigma = \iota, \iota^{c}} \left[ v_p \left( \dfrac{\delta^{\prime}_{E_1,\sigma}}{\delta^{\prime}_{E_2,\sigma}} \right) - v_p \left( \dfrac{w_{E_1, \sigma}(0)}{w_{E_2, \sigma}(0)} \right) \right]
\]
\end{lemma}

Next we deal with the remaining terms using a continuity trick. The following proposition gives the relation between the valuations of $\delta_{E,\sigma}^{\prime}$ and $w_{E,\sigma}$.

\begin{proposition} \label{prop:compare-two-L-functions}
With the notation above, we have
\[
    v_p \left( \dfrac{w_{1,\sigma}(k)}{w_{2,\sigma}(k)} \right) = (2k-1) v_p \left( \dfrac{\delta_{2,\sigma}^{\prime}}{\delta_{1,\sigma}^{\prime}} \right).
\]
\end{proposition}

\begin{proof}
    We use the notation and results in Appendix \ref{app:elliptic-units}. We fix an embedding $\sigma \in \{\iota, \iota^{c} \}$ and temporarily omit it from the notation for simplicity, and we write $j = 1,2$ for objects corresponding to the elliptic curve $E_{j}$.

Recall our normalization that
\[
\scL_{j}(\phi) = \dfrac{\Lambda(1,\eta) \Lambda(1,\omega_{j}^{-1})}{4 \calE_{j}^{(p)}(\phi)} \calL_{j}(\phi).
\]
Together with the interpolation formula in Theorem \ref{thm:fan-wan-p-adic-l-function}, we have
\begin{equation} \label{eq:scL-phi-k}
\scL_{j}(k) := \scL_{j}(\phi_{k}) = \dfrac{1}{8} r_j^{2k-1} p^{d_p(k-1)} \calC_{j}(k) \Gamma_{\CC}(k) L(\psi_{j}\phi_{k}, 1) \lambda_{j,p}(k)
\end{equation}
and we put $B(k) = \dfrac{1}{8} p^{d_p(k-1)} \calC_{j}(k) \Gamma_{\CC}(k)$. We choose the common local translation parameters $\ell,n$, hence the same $d_p$, for the two elliptic curves. Note that here $\calC_{j}(k) = i^{k-1} 2^{2-k}$ is computed in Lemma \ref{lem:archimedean-computation-comparison}, so it is independent of the elliptic curve $E_{j}$, and hence the term $B(k)$ is independent of the elliptic curve $E_{j}$ as well.

The explicit reciprocity law in Theorem \ref{thm:eu-article-ERC} gives
\begin{equation} \label{eq:ERC-elliptic-unit}
\exp^{\ast}(\loc_p z_{j,k}) = \calE_{j,p}(k) L(\psi_{j}\phi_{k}, 1) \dfrac{(2\pi i )^{k-1}}{(a_{E_j} d_{j})^{2k-1}} \ttb_{j,k},
\end{equation}
where $d_{j}$ generates the different of $\calK_j/\QQ$, and $\ttb_{j,k}$ denotes the de Rham tensor basis associated with $\omega_{j}$. Meanwhile, we apply the dual exponential map $\exp^{\ast}$ to \eqref{eq:Nov-family-w} to obtain
\begin{equation} \label{eq:exp-dual-L-function}
\exp^{\ast}(\loc_{p} z_{j,k}) = \scL_{j}(k) w_{j}(k) \exp^{\ast}(u(k)) = B(k) r_j^{2k-1} L(\psi_j\phi_k,1) \lambda_{j,p}(k) w_j(k) \exp^{\ast}(u(k)). 
\end{equation}
At critical weights where the common interpolated values are nonzero, for each individual $j$ we compare \eqref{eq:ERC-elliptic-unit} and \eqref{eq:exp-dual-L-function} to cancel the term $L(\psi_j\phi_k,1)$ and obtain
\[
\calE_{j,p}(k) \dfrac{(2\pi i)^{k-1}}{(a_{E_j}d_j)^{2k-1}} b_{j,k} = B(k) r_j^{2k-1} \lambda_{j,p}(k) w_{j}(k)\exp^*(u(k)).
\]
Taking the ratio for $j=1,2$ cancels this common nonzero vector, and using $\ttb_{1,k}/\ttb_{2,k} = (\omega_1/\omega_2)^{2k-1}$ gives the following cancellation:
\[
\dfrac{w_{1}(k)}{w_2(k)} = \left(\dfrac{d_2}{d_1}\right)^{2k-1} \dfrac{\calE_{1,p}(k) \lambda_{2,p}(k)}{\calE_{2,p}(k) \lambda_{1,p}(k)}\left(\dfrac{\omega_1}{\omega_2} \dfrac{r_2 a_{E_2}}{r_1 a_{E_1}}\right)^{2k-1}.
\]
The common local representations, test vectors and measures give $\calE_{1,p}(k)\lambda_{2,p}(k)=\calE_{2,p}(k)\lambda_{1,p}(k)$. Under the same de Rham identification, $\log_{\omega_j}(u)=C_j^{-1}$ gives $\omega_1/\omega_2=C_2/C_1$. Finally $\delta_j'=(2\pi i)^{-1}C_jr_ja_{E_j}$, so
\[
 \frac{\omega_1}{\omega_2}\frac{r_2a_{E_2}}{r_1a_{E_1}}
 =\frac{\delta_2'}{\delta_1'}.
\]
Therefore, we obtain
\[
\dfrac{w_{1}(k)}{w_2(k)} = \left(\dfrac{d_2}{d_1}\right)^{2k-1}\left(\dfrac{\delta_{2}^{\prime}}{\delta_1^{\prime}}\right)^{2k-1}.
\]
The identification of the local fields identifies their different ideals, so $d_2/d_1\in\calO_F^\times$. Taking the $p$-adic valuation on both sides and restoring the embedding $\sigma$ in the notation, we obtain
\begin{equation} \label{eq:comparison-w-delta-prime}
v_p \left( \dfrac{w_{1,\sigma}(k)}{w_{2,\sigma}(k)} \right) = (2k-1) v_p \left( \dfrac{\delta_{2,\sigma}^{\prime}}{\delta_{1,\sigma}^{\prime}} \right).
\end{equation}
This proves the proposition.
\end{proof}

We now apply a continuity trick to this exact valuation identity (Proposition \ref{prop:compare-two-L-functions}) to prove the following theorem.

\begin{theorem}[Transfer principle] \label{thm:main_transfer}
Let $E_{1}/\QQ$ and $E_{2}/\QQ$ be two CM elliptic curves, with $\ord_{s=1} L(E_{j}/\QQ, s) = 1$ for $j=1, 2$, such that $E_{1}$ and $E_{2}$ are {matched at $p$} in the sense of Definition \ref{defn:matched}. Then
\[
\Delta_{p}(E_1/\calK_1) = \Delta_{p}(E_2/\calK_2).
\]
In particular, the following statements are equivalent:
\begin{enumerate}[label= \rm (\roman*)]
\item the $p$-part of BSD holds for $E_1/\calK_1$;
\item the $p$-part of BSD holds for $E_2/\calK_2$;
\item the $p$-part of BSD holds for $E_1/\QQ$;
\item the $p$-part of BSD holds for $E_2/\QQ$.
\end{enumerate}
\end{theorem}

\begin{proof}
Put $f_\sigma(k)=w_{1,\sigma}(k)/w_{2,\sigma}(k)$ and $d_\sigma=v_p(\delta'_{2,\sigma}/\delta'_{1,\sigma})$. The function $f_\sigma$ is analytic and nonzero at zero. For $k_a=Np^a$, continuity and the ultrametric inequality give $v_p(f_\sigma(k_a))=v_p(f_\sigma(0))$ for all sufficiently large $a$. Equation~\eqref{eq:comparison-w-delta-prime} gives
\[
v_p \left( \dfrac{w_{1,\sigma}(k)}{w_{2,\sigma}(k)} \right) = (2k-1) v_p \left( \dfrac{\delta_{2,\sigma}^{\prime}}{\delta_{1,\sigma}^{\prime}} \right),
\]
and therefore gives
\[
 v_p(f_\sigma(0))=(2Np^a-1)d_\sigma.
\]
The left side is fixed, whereas the positive integers $Np^a$ are unbounded. Hence $d_\sigma=0$ and $v_p(f_\sigma(0))=0$. Substitution in Lemma~\ref{lem:transfer-prepare} proves equality of the defects, and therefore the equivalence of (i) and (ii). Lemma~\ref{lem:base-change} identifies (i) with (iii) and (ii) with (iv).
\end{proof}

\subsection{The transfer mechanism in practice}\label{sst:practice}

In this section, we discuss a general mechanism for using the transfer principle, Theorem \ref{thm:main_transfer}. Fix a prime $p$ and a quadratic extension $F$ of $\QQ_p$. Our main goal is to classify all local characters 
\[
\psi_p: F^{\times}\rightarrow\CC^{\times}
\]
that can occur as CM local characters. For each possibility, we then find an elliptic curve $E/\QQ$ of analytic rank $1$, with CM field $\calK$, such that $p$ is nonsplit, $\calK_p\cong F$, and the local character $\psi_{E,p}$ is isomorphic to $\psi_p$. Finally, we check the $p$-BSD formula for these representatives. Note that any elliptic curve $E/\QQ$ with CM and analytic rank $1$ has a $p$-local character isomorphic to one of the representatives above. Therefore, Theorem \ref{thm:main_transfer} shows that the $p$-BSD formula holds for $E$.

We fix the embedding $\calK\hookrightarrow\CC$ sending $\sqrt{D_\calK}$ to the upper half-plane. Requiring $\psi$ to be a CM local character of infinite type $(-1,0)$, one can show that $\psi_p$ takes values in $\calK^{\times}$ and satisfies the normalization
\[
\psi_p\vert_{\QQ_p^{\times}}=\eta_{F/\QQ_p}\lvert\cdot\rvert_p^{-1},
\]
where $\eta_{F/\QQ_p}$ is the unique non-trivial quadratic character on $\QQ_p^{\times}$, with kernel $\Norm_{F/\QQ_p}(F^{\times})$.
For the remainder of this part, we fix a uniformizer $\pi$ of $F$ and let $U^n$ denote $1+\pi^n\mathcal{O}_F$.

\subsubsection{The inert case}\label{sst:inert}

Assume $p$ is inert in $\calK$, and again denote $p\mathcal{O}_\calK$ by $(p)$. Then $F\cong \calK_p$ is the unique unramified quadratic extension of $\QQ_p$. In this situation, $\psi_p(p)=-p$ and $\psi_p$ is trivial on $\ZZ_p^{\times}$.

For any $x\in \calK^{\times}$, one can write $x$ as $p^{v_p(x)}u$ for some $u\in\mathcal{O}_F^{\times}$, since $F/\QQ_p$ is unramified. Thus $\psi_p(x)=(-p)^{v_p(x)}\psi_p(u)$. The question reduces to classifying
\[
\psi_p:\mathcal{O}_F^{\times}\rightarrow \calK^{\times}.
\]

Let $\mathfrak{f}_{\psi}$ denote the conductor of ${\psi}$, which is an ideal in $\mathcal{O}_\calK$. Recall that $v_p(\mathfrak{f}_{\psi})$ is the minimal $n$ making $\psi_p$ trivial on $U^n=1+\pi^n\mathcal{O}_F$.

Thus we consider the structure of $G_n:=\mathcal{O}_F^{\times}/(\ZZ_p^{\times}U^n)$, so any $\psi_p$ must factor through $G_{v_p(\mathfrak{f}_{\psi})}$. There is a full classification of $G_n$ for $F/\QQ_p$ unramified quadratic:
\begin{equation}\label{eq:G_n classification}
    G_n\cong
    \begin{cases*}
        C_{p+1}\oplus C_{p^{n-1}} & \text{$p\geq3$}\\
        C_3\oplus C_2\oplus C_{2^{n-2}} & \text{$p=2,~n\geq2$}\\
        C_3 & \text{$p=2,~n=1$}.\\
    \end{cases*}
\end{equation}
Therefore $\psi_p$ is reduced to a character on finite groups.

To determine all such $\psi_p$, we may choose $\calK$ to be either $\QQ(\sqrt{-1})$ or $\QQ(\sqrt{-3})$, among all the CM fields, since they contain the most roots of unity, namely $\mu_4$ or $\mu_6$. (Of course, one should first check that $p$ is inert in the chosen $\calK$.) Note that constructions of $\psi_p$ from other choices of $\calK$ are already duplicated by constructions from $\calK=\QQ(\sqrt{-1})$ or $\QQ(\sqrt{-3})$. View $\Hom(G_n,\calK^{\times})$ as a subgroup of $\Hom(G_{n+1},\calK^{\times})$ via inflation. One may check that
\[
\Hom(G_n,\mu_4)~\text{and}~\Hom(G_n,\mu_6)
\]
eventually stabilize as $n$ grows. Thus there are only finitely many choices for $v_p(\mathfrak{f}_{\psi})$ and hence for $\psi_p$.

Finally, one has to find an $E/\QQ$ with local CM character isomorphic to $\psi_p:F^{\times}\rightarrow\CC^{\times}$ and check the needed information, such as the analytic $\Sha$ or the Manin constant. Note that the CM field of such an $E$ need not be $\QQ(\sqrt{-1})$ or $\QQ(\sqrt{-3})$; for instance, the constructed $\psi_p$ takes values only in $\{\pm1\}$ on the unit group.

The method we use to recover such $E/\QQ$ is to search LMFDB directly for elliptic curves with CM by $\QQ(\sqrt{-1})$ or $\QQ(\sqrt{-3})$, analytic rank $1$, and a conductor with the expected $v_p(N_E)$. We then calculate their local $\psi_{E,p}$ until all possibilities are exhausted. Although, in principle, the possibilities need not have representatives with CM by $\QQ(\sqrt{-1})$ or $\QQ(\sqrt{-3})$, this route works in practice.

\subsubsection{The ramified case}\label{sst:ramified}

Now we assume $p$ is ramified in the CM field $\calK$ and write $\mathfrak{p}$ for the unique prime above $p$. So the only cases of $(p,\calK)$ are:
\[
(2,\QQ(\sqrt{-1})),(2,\QQ(\sqrt{-2}))~\text{and}~(p,\QQ(\sqrt{-p})),
\]
with $p\in\{3,7,11,19,43,67,163\}$. Note that $\calK_{\mathfrak{p}}$ are all distinct among those $(p,\calK)$.

Therefore $F=\calK_{\mathfrak{p}}$ is a totally ramified extension of degree $2$ over $\QQ_p$, with ring of integers $\mathcal{O}_F\cong\mathbb{Z}_p[\pi]$. Note that $\pi$ can be taken as $\sqrt{-p}$ for $p\geq2$ and $\calK=\QQ(\sqrt{-p})$. For $p=2$ and $\calK=\QQ(i)$, we take $\pi=1+i$. In all these cases one can check that $p=\Norm_{F/\QQ_p}(\pi)$ and that $-1$ is not a norm element, hence $\psi_p(p)=p$ and $\psi_p(-1)=-1$. In particular $\psi_p(\pi)\in\{\pm\pi\}$ for $\calK\neq\QQ(\sqrt{-1})$.

Any $x\in F^{\times}$ can be written as $x=\pi^{v_F(x)}u$ for some $u\in\mathcal{O}_F^{\times}$, so $\psi_p(x)=\psi_p(\pi)^{v_F(x)}\psi_p(u)$ and $\psi_p$ is determined by $\psi_p(\pi)$ and its restriction to $\mathcal{O}_F^{\times}$. Since $\psi_p$ is no longer trivial on $\ZZ_p^{\times}$, we consider $H_n:=\mathcal{O}_F^{\times}/U^n$. The following is a complete classification of $(p,\calK,H_n)$:

\begin{equation}\label{eq:H_n classification}
    H_n\cong
    \begin{cases*}
        C_{p-1}\oplus C_{p^{\lfloor n/2\rfloor}}\oplus C_{p^{\lfloor(n-1)/2\rfloor}} & \text{$p\geq7,~n\geq1$}\\
        C_2 & \text{$p=3,~n=1$}\\
        C_2\oplus C_3\oplus C_{3^{\lceil n/2\rceil-1}}\oplus C_{3^{\lfloor n/2\rfloor-1}} & \text{$p=3,~n\geq2$}\\
        C_1 & \text{$p=2,~n=1$}\\
        C_2 & \text{$p=2,~n=2$}\\
        C_4\oplus C_{2^{\lceil(n-3)/2\rceil}}\oplus C_{2^{\lfloor(n-3)/2\rfloor}} & \text{$p=2,~\calK=\QQ(\sqrt{-1}),~n\geq3$}\\
        C_4 & \text{$p=2,~\calK=\QQ(\sqrt{-2}),~n=3$}\\
        C_2\oplus C_{2^{\lfloor n/2\rfloor}}\oplus C_{2^{\lceil n/2\rceil-2}} & \text{$p=2,~\calK=\QQ(\sqrt{-2}),~n\geq4$}.
    \end{cases*}
\end{equation}
Here $n\geq1$, and $C_1$ denotes the trivial group.

Since $F$ corresponds to a unique $\calK$, we cannot arbitrarily assume that $\psi_p$ takes values in $\mu_4$ or $\mu_6$, unlike in the inert case, in which a fixed $F$ may come from any of these $\calK$. We will consider $\Hom(H_n,\calK^{\times})=\Hom(H_n,\mu(\calK))$ for each pair, where $\mu(\calK)$ denotes the roots of unity in $\calK$. It is easy to see that, for each pair $(p,\calK)$, $\Hom(H_n,\mu(\calK))$ eventually stabilizes as $n$ grows. Hence, $\psi_p$ has only finitely many choices. The last step is to find a representative for each possibility and verify its BSD invariants, as in the inert case.

\subsubsection{The computation of $\psi_{E,p}$ for special families}\label{sst:calculation}

Next, we consider $2$ special families of elliptic curves over $\QQ$, given by $X_a:y^2=x^3+a$ and $Y_b:y^2=x^3+bx$, for $a,b\in\QQ^{\times}$. We know that $X_a$ has CM by the maximal order in $\calK=\QQ(\sqrt{-3})$ and that $Y_b$ has CM by the maximal order in $\calK=\QQ(\sqrt{-1})$. We choose them as candidates for realizing $\psi_p$ since one can calculate $\psi_{X_a,p}$ and $\psi_{Y_b,p}$ efficiently.

In the following parts, we fix the embedding $\QQ(\sqrt{-p})\hookrightarrow\CC$ sending $\sqrt{-p}$ to the upper half plane, fix $\sqrt{-1}$ to be $i$, and fix $\zeta_n:=e^{2\pi i/n}$, $\omega:=\zeta_3=\zeta_6^2=\frac{-1+\sqrt{-3}}{2}$.

\begin{proposition}\label{prop:calculation 3,5,7 inert}
    If $p$ is inert in $\calK$, let $g_p$ denote the Teichmuller lift of a generator of $(\mathcal{O}_F/p\mathcal{O}_F)^{\times}\cong\FF_{p^2}^{\times}$ whose image generates $G_1$.
    \begin{enumerate}
        \item In the case $(p,\calK)=(5,\QQ(\sqrt{-3}))$, choose $g_5:=\Teich(2+\bar{\omega})$, then $\psi_{X_a,5}(g_5)=\zeta_6^{-v_5(a)}$.
        \item In the case $(p,\calK)=(7,\QQ(\sqrt{-1}))$, choose $g_7:=\Teich(2+\bar{i})$, then $\psi_{Y_b,7}(g_7)=i^{v_7(b)}$.
        \item In the case $(p,\calK)=(3,\QQ(\sqrt{-1}))$, choose $g_3:=\Teich(1-\bar{i})$, then $\psi_{Y_b,3}(g_3)=i^{v_3(b)}$.
    \end{enumerate}
\end{proposition}
\begin{proof}
    See Appendix \ref{sec:special values calculation}.
\end{proof}

%% file: 106-BSD-inert-2-3.tex
\section{Birch--Swinnerton-Dyer formula in the inert case for $p=2$ and $3$}\label{sec:inert BSD 2,3}

From now on, we apply the transfer mechanism to concrete cases that are difficult to handle using other existing tools. We start with small inert primes. Unlike in the ramified cases, there are infinitely many inert primes in a CM field $\calK$. Thus, in this section, we focus only on the inert primes $p=2,3$, while the rest will be handled by a BDP-type method later in Part \ref{part:BDP}. Recall that Section \ref{sst:inert} defines $F$ to be the unique unramified quadratic extension of $\QQ_p$, $U^n:=1+p^n\mathcal{O}_F$, and $G_n:=\mathcal{O}_F^{\times}/(\ZZ_p^{\times}U^n)$ for $n\geq1$.

\subsection{The case $p=3$}

For $p=3$ inert in $\calK$, we have $G_n\cong C_4\oplus C_{3^{n-1}}$ for $n\geq1$ by the classification \eqref{eq:G_n classification}. To exhaust all possibilities, we assume $\calK=\QQ(\sqrt{-1})$ and $\psi_3|_{\mathcal O_F^{\times}}$ takes values in $\mu(\calK)=\mu_4$. Note that $3$ is ramified in $\QQ(\sqrt{-3})$, so $\mu_6$ is excluded.

Since, via inflation, $\Hom(G_n,\mu_4)=\Hom(G_1,\mu_4)\cong C_4$, the unit restriction of $\psi_3$ is trivial on $U^1$. Hence $v_3(\mathfrak{f}_{\psi})\in\{0,1\}$, and the target curve $E$ has $v_3(N_E)\in\{0,2\}$. Let $g_3=\Teich(1-i)$ be as in Proposition \ref{prop:calculation 3,5,7 inert}. The following table, based on LMFDB, gives a realization for each possible $\psi_3$.

\begingroup\small\setlength{\tabcolsep}{7pt}\renewcommand{\arraystretch}{1.25}
\begin{longtable}{@{}cccccccc@{}}
\caption{Representatives for each $\psi_3$ with $3$ inert}\label{tab:p_3 inert database}\\
\hline
Equation of $E$ & LMFDB label & $N_E$ & $\calK$ & $r_{\rm an}$ & $S_{3}$ & $S_{\rm an}$ & $\psi_{E,3}(g_3)$\\\hline\endfirsthead
\multicolumn{8}{c}{\tablename\ \thetable\ (continued)}\\\hline
Equation of $E$ & LMFDB label & $N_E$ & $\calK$ & $r_{\rm an}$ & $S_{3}$ & $S_{\rm an}$ & $\psi_{E,3}(g_3)$\\\hline\endhead
\hline\endfoot
    $y^2=x^3-2x$ & {256.b1} & $2^{8}$ & $\QQ(i)$ & 1 & 1 & 1 & $1$ \\
    $y^2=x^3+3x$ & {288.a2} & $2^{5}\cdot 3^{2}$ & $\QQ(i)$ & 1 & 1 & 1 & $i$ \\
    $y^2=x^3+9x$ & {576.c4} & $2^{6}\cdot 3^{2}$ & $\QQ(i)$ & 1 & 1 & 1 & $-1$ \\
    $y^2=x^3-54x$ & {2304.m1} & $2^{8}\cdot 3^{2}$ & $\QQ(i)$ & 1 & 1 & 1 & $-i$ \\
\end{longtable}\endgroup

In the tables, $r_{\rm an}$ denotes the analytic rank and the two $\Sha$ columns mean
\[
 S_{p}=\#\Sha(E/\QQ)[p^{\infty}],\qquad
 S_{\rm an}=\frac{L'(E/\QQ,1)\,\#E(\QQ)_{\rm tors}^{2}}
 {\Omega_{E/\QQ} \Reg(E/\QQ) \prod_{\ell} \Tam_\ell(E/\QQ)}.
\]
The $S_3$ terms are calculated in Appendix \ref{sec:p-Sha calculation}, as are the $S_p$ terms for representatives in the following parts. In particular, the $3$-BSD formula holds for the curves above.

\subsection{The case $p=2$}

For $p=2$ inert in $\calK$, we have $G_1\cong C_3$ and $G_n\cong C_3\oplus C_2\oplus C_{2^{n-2}}$ for $n\geq2$. For similar reasons, we may assume $\calK=\QQ(\sqrt{-3})$ and $\psi_2|_{\mathcal O_F^{\times}}$ takes values in $\mu_6$. Since $\Hom(G_1,\mu_6)\cong C_3$, $\Hom(G_2,\mu_6)\cong C_6$, and $\Hom(G_n,\mu_6)\cong C_6\oplus C_2$ for $n\geq3$, it follows that the unit restriction of $\psi_2$ factors through $G_3$ with conductor exponent at most $3$, and thus $v_2(N_E)\in\{0,2,4,6\}$.

Thus, in this case, we need to find generators for $G_1$, $G_2$, and $G_3$, not just a Teichmuller lift generating $G_1$. Put
\[
f_2=-\omega-2,\qquad h_2=1+4\omega.
\]
One can show that the image of $f_2$ generates $G_2\cong C_6$, the image of $f_2^2$ generates $G_1\cong C_3$, and the images of $f_2,h_2$
generate $G_3\cong C_6\oplus C_2$.

Consequently we have the following calculation, proved in Appendix \ref{sec:special values calculation}.

\begin{proposition}\label{prop:calculation 2 inert}
    For $X_a:y^2=x^3+a$ with $a\in\QQ^{\times}$, write $a=2^m u$ with $u\in\ZZ_2^{\times}$ and
    $e=(u-1)/2\pmod2$. Under the identifications as in Section \ref{sst:calculation}, the special values are
    \[
     \psi_{X_a,2}(f_2)=\omega^{-m+1}(-1)^{e+m},\qquad
     \psi_{X_a,2}(h_2)=(-1)^m.
    \]
\end{proposition}

The following table gives a realization of each possible $\psi_2$.

\begingroup\small\setlength{\tabcolsep}{7pt}\renewcommand{\arraystretch}{1.25}
\begin{longtable}{@{}ccccccccc@{}}
\caption{Representatives for each $\psi_2$ with $2$ inert}\label{tab:p_2 inert database}\\
\hline
Equation of $E$ & LMFDB label & $N_E$ & $\calK$ & $r_{\rm an}$ & $S_{2}$ & $S_{\rm an}$ & $\psi_{E,2}(f_2)$ & $\psi_{E,2}(h_2)$\\\hline\endfirsthead
\multicolumn{9}{c}{\tablename\ \thetable\ (continued)}\\\hline
Equation of $E$ & LMFDB label & $N_E$ & $\calK$ & $r_{\rm an}$ & $S_{2}$ & $S_{\rm an}$ & $\psi_{E,2}(f_2)$ & $\psi_{E,2}(h_2)$\\\hline\endhead
\hline\endfoot
    $y^2+y=x^3+1$ & {225.c2} & $3^2\cdot 5^2$ & $\QQ(\sqrt{-3})$ & 1 & 1 & 1 & $1$ & $1$ \\
    $y^2=x^3+9$ & {972.b2} & $2^2\cdot 3^5$ & $\QQ(\sqrt{-3})$ & 1 & 1 & 1 & $\omega$ & $1$ \\
    $y^2=x^3+100$ & {900.d2} & $2^2\cdot 3^2\cdot 5^2$ & $\QQ(\sqrt{-3})$ & 1 & 1 & 1 & $\omega^2$ & $1$ \\
    $y^2=x^3-10000$ & {3600.b1} & $2^4\cdot 3^2\cdot 5^2$ & $\QQ(\sqrt{-3})$ & 1 & 1 & 1 & $-1$ & $1$ \\
    $y^2=x^3+3$ & {3888.i2} & $2^4\cdot 3^5$ & $\QQ(\sqrt{-3})$ & 1 & 1 & 1 & $-\omega$ & $1$ \\
    $y^2=x^3-4$ & {432.d1} & $2^4\cdot 3^3$ & $\QQ(\sqrt{-3})$ & 1 & 1 & 1 & $-\omega^2$ & $1$ \\
    $y^2=x^3-2$ & {1728.o3} & $2^6\cdot 3^3$ & $\QQ(\sqrt{-3})$ & 1 & 1 & 1 & $1$ & $-1$ \\
    $y^2=x^3-72$ & {15552.bg1} & $2^6\cdot 3^5$ & $\QQ(\sqrt{-3})$ & 1 & 1 & 1 & $\omega$ & $-1$ \\
    $y^2=x^3-288$ & {15552.bi1} & $2^6\cdot 3^5$ & $\QQ(\sqrt{-3})$ & 1 & 1 & 1 & $\omega^2$ & $-1$ \\
    $y^2=x^3+2$ & {1728.n4} & $2^6\cdot 3^3$ & $\QQ(\sqrt{-3})$ & 1 & 1 & 1 & $-1$ & $-1$ \\
    $y^2=x^3+8$ & {576.e4} & $2^6\cdot 3^2$ & $\QQ(\sqrt{-3})$ & 1 & 1 & 1 & $-\omega$ & $-1$ \\
    $y^2=x^3+288$ & {15552.y2} & $2^6\cdot 3^5$ & $\QQ(\sqrt{-3})$ & 1 & 1 & 1 & $-\omega^2$ & $-1$ \\
\end{longtable}\endgroup

We can see that the $2$-part of BSD holds for these $12$ curves.

\begin{remark}
    Note that the curve $X_8:y^2=x^3+8$ has $2$-torsion on $\QQ$; namely, $X_8(\QQ)\cong\ZZ\oplus\ZZ/2\ZZ$, while $\Sel_2(X_a/\QQ)\cong\ZZ/2\ZZ$. Thus, at $p=2$, inertness does not imply the irreducibility of $E[p]$, in contrast to Proposition \ref{prop:irreducible}.
\end{remark}

\begin{theorem}\label{thm:BSD 2,3 inert}
    Let $E/\QQ$ be an elliptic curve with complex multiplication by $\calK$ with analytic rank $1$. Let $p\in\{2,3\}$ be a prime inert in $\calK$. Then the $p$-BSD holds for $E/\QQ$.
\end{theorem}
\begin{proof}
    Since there is a curve $E^{\prime}$ in Tables \ref{tab:p_2 inert database} and \ref{tab:p_3 inert database} with $\psi_{E^{\prime},p}\cong\psi_{E,p}$, Theorem \ref{thm:main_transfer} shows that the $p$-BSD formula holds for $E/\QQ$.
\end{proof}

%% file: 107-BSD-ramified.tex
\section{Birch--Swinnerton-Dyer formula in the ramified case}\label{sec:ramified BSD}

In this section, we give the explicit classification of $\psi_p$ in the ramified case and search for a representative $E$ for each $\psi_p$ such that $\psi_{E,p}\cong\psi_p$. Recall that all the cases for $(p,\calK)$ are
\[
(2,\QQ(\sqrt{-1})),(2,\QQ(\sqrt{-2}))~\text{and}~(p,\QQ(\sqrt{-p})),
\]
with $p\in\{3,7,11,19,43,67,163\}$.

\subsection{The case $p>3$}

For $p>3$, we have $H_n:=\mathcal{O}_F^{\times}/U^n\cong C_{p-1}\oplus C_{p^{\lfloor n/2\rfloor}}\oplus C_{p^{\lfloor(n-1)/2\rfloor}}$ by the classification \eqref{eq:H_n classification}. Since $\mu(\QQ(\sqrt{-p}))=\{\pm1\}$, it follows that $\Hom(H_n,\mu(\calK))\cong C_2$ for all $n\geq1$. Hence $\psi_p$ is trivial on $U^1$. But $H_1$ is generated by a Teichmuller lift $g_p$ of a generator of $(\mathcal{O}_F/\pi\mathcal{O}_F)^{\times}\cong\FF_p^{\times}$, which actually lies in $\ZZ_p^{\times}$. Recall that $\psi_p\vert_{\ZZ_p^{\times}}=\eta_{F/\QQ_p}=\left(\frac{\cdot}{p}\right)$, the unique non-trivial character corresponding to $F/\QQ_p$. Thus $\psi_p(g_p)=-1$ since $p\neq2$ and $\bar{g}_p\in\FF_p^{\times}$ is not a square. Therefore, $\psi_p$ is uniquely determined by $\psi_p(\pi)\in\{\pm\pi\}$ and necessarily has conductor $v_F(\mathfrak{f}_{\psi})=1$, implying $v_p(N_E)=2$, consistent with Table \ref{table:CM_classification}.

The following table gives a set of realizations for each $p$, in which $\psi_{E,p}(\pi)$ and $S_{p}$ are calculated in Appendices \ref{sec:special values calculation} and \ref{sec:p-Sha calculation}. Note that the curves in the even rows are quadratic twists of the curves in the odd rows. We will denote the base curves by $E_p$ and the twisted curves by $E_p^d$, as explained in the appendices.

\begingroup\small\setlength{\tabcolsep}{2pt}\renewcommand{\arraystretch}{1.25}
\begin{longtable}{@{}cccccccc@{}}
\caption{Representatives for each $\psi_p$ with $p>3$ ramified in $\QQ(\sqrt{-p})$}\label{tab:p>3 ramified database}\\
\hline
Equation of $E$ & LMFDB label & $N_E$ & $\calK$ & $r_{\rm an}$ & $S_{p}$ & $S_{\rm an}$ & $\psi_{E,p}(\pi)$ \\\hline\endfirsthead
\multicolumn{8}{c}{\tablename\ \thetable\ (continued)}\\\hline
Equation of $E$ & LMFDB label & $N_E$ & $\calK$ & $r_{\rm an}$ & $S_{p}$ & $S_{\rm an}$ & $\psi_{E,p}(\pi)$ \\\hline\endhead
\hline\endfoot
    $y^2+xy+y=x^3-x^2-20x+46$ & 441.c4 & $3^{2}\cdot7^{2}$ & $\QQ(\sqrt{-7})$ & 1 & $1$ & $1$ & $\pi$\\
    $y^2=x^3-35x+98$ & 784.f4 & $2^{4}\cdot7^{2}$ & $\QQ(\sqrt{-7})$ & 1 & $1$ & $1$ & $-\pi$\\
    $y^2+y=x^3-x^2-7x+10$ & 121.b2 & $11^{2}$ & $\QQ(\sqrt{-11})$ & 1 & $1$ & $1$ & $\pi$\\
    $y^2=x^3+x^2-29x+53$ & 7744.bb2 & $2^{6}\cdot11^{2}$ & $\QQ(\sqrt{-11})$ & 1 & $1$ & $1$ & $-\pi$\\
    $y^2+y=x^3-38x+90$ & 361.a2 & $19^{2}$ & $\QQ(\sqrt{-19})$ & 1 & $1$ & $1$ & $\pi$\\
    $y^2=x^3-152x+722$ & 23104.bc2 & $2^{6}\cdot19^{2}$ & $\QQ(\sqrt{-19})$ & 1 & $1$ & $1$ & $-\pi$\\
    $y^2+y=x^3-860x+9707$ & 1849.b2 & $43^{2}$ & $\QQ(\sqrt{-43})$ & 1 & $1$ & $1$ & $\pi$\\
    $y^2+y=x^3-21500x+1213406$ & 46225.h2 & $5^{2}\cdot43^{2}$ & $\QQ(\sqrt{-43})$ & 1 & $1$ & $1$ & $-\pi$\\
    $y^2+y=x^3-7370x+243528$ & 4489.b2 & $67^{2}$ & $\QQ(\sqrt{-67})$ & 1 & $1$ & $1$ & $\pi$\\
    $y^2=x^3-29480x+1948226$ & 287296.h2 & $2^{6}\cdot67^{2}$ & $\QQ(\sqrt{-67})$ & 1 & $1$ & $1$ & $-\pi$\\
    $y^2+y=x^3-2174420x+1234136692$ & 26569.a2 & $163^{2}$ & $\QQ(\sqrt{-163})$ & 1 & $1$ & $1$ & $\pi$\\
    $y^2=x^3-8697680x+9873093538$ & \text{not located} & $2^{6}\cdot163^{2}$ & $\QQ(\sqrt{-163})$ & 1 & $1$ & $1$\footnote{This is computed by SageMath in appendix \ref{sec:p-Sha calculation}} & $-\pi$\\
\end{longtable}\endgroup

We can see that the $p$-BSD holds for these curves.

\subsection{The case $p=3$}

Now we take $p=3$ and $\calK=\QQ(\sqrt{-3})$. Hence $F=\QQ_3(\pi)$ and $\mathcal{O}_F=\ZZ_3[\pi]$, with $\pi:=\sqrt{-3}$. By classification \eqref{eq:H_n classification} we have $H_1\cong C_2$ and $H_n\cong C_2\oplus C_3\oplus C_{3^{\lceil n/2\rceil-1}}\oplus C_{3^{\lfloor n/2\rfloor-1}}$ for $n\geq2$. Thus
\[
\Hom(H_n,\mu_6)\cong
\begin{cases*}
    C_2 & \text{$n=1$}\\
    C_6 & \text{$n=2$}\\
    C_6\oplus C_3 & \text{$n=3$}\\
    C_6\oplus C_3\oplus C_3 & \text{$n\geq 4$}\\
\end{cases*}
\]
which stabilizes after $n=4$. Hence $\psi_3$ is trivial on $U^4$.

One can show that the images of the elements in the sequence $(-1,\omega,4,t:=1+3\pi)$ generate $H_4\cong C_2\oplus C_3\oplus C_3\oplus C_3$, while the images of the first $i$ terms generate $H_i$, $i=1,2,3$, viewed as subgroups via inflation. By normalization, we have $\psi_3(-1)=-1$ and $\psi_3(4)=1$. In particular, if $\psi_3$ is trivial on $U^3$, then $\psi_3$ is also trivial on $U^2$. Therefore, $v_F(\mathfrak{f}_{\psi})\in\{1,2,4\}$, corresponding to $v_3(N_E)\in\{2,3,5\}$, as shown in Table \ref{table:CM_classification}. Hence all $18$ possibilities are given by $\psi_3(\pi)\in\{\pm\pi\}$, $\psi_3(\omega)\in\mu_3$, and $\psi_3(t)\in\mu_3$, which can be computed as below, as explained in Appendix \ref{sec:special values calculation}.

\begin{proposition}\label{prop:calculation 3 ramified}
    For $X_a:y^2=x^3+a$, write $a=3^{v_3(a)}u$ with $u\in\ZZ_3^{\times}$. Then we have
    \[
    \psi_{X_a,3}(\pi)=(-1)^{v_3(a)}\left(\frac{u}{3}\right)\pi,\quad \psi_{X_a,3}(\omega)=\omega^{(u^2-1)/3},\quad \psi_{X_a,3}(t)=\omega^{v_3(a)}.
    \]
\end{proposition}

The following table realizes all $18$ possibilities.

\begingroup\small\setlength{\tabcolsep}{5pt}\renewcommand{\arraystretch}{1.25}
\begin{longtable}{@{}cccccccccc@{}}
\caption{Representatives for each $\psi_3$ with $3$ ramified in $\QQ(\sqrt{-3})$}\label{tab:p=3 ramified database}\\
\hline
Equation of $E$ & LMFDB label & $N_E$ & $\calK$ & $r_{\rm an}$ & $S_{3}$ & $S_{\rm an}$ & $\psi_{E,3}(\pi)$ & $\psi_{E,3}(\omega)$ & $\psi_{E,3}(t)$\\\hline\endfirsthead
\multicolumn{10}{c}{\tablename\ \thetable\ (continued)}\\\hline
Equation of $E$ & LMFDB label & $N_E$ & $\calK$ & $r_{\rm an}$ & $S_{3}$ & $S_{\rm an}$ & $\psi_{E,3}(\pi)$ & $\psi_{E,3}(\omega)$ & $\psi_{E,3}(t)$\\\hline\endhead
\hline\endfoot
$y^2+y=x^3+1$ & 225.c2 & $3^{2}\cdot5^{2}$ & $\QQ(\sqrt{-3})$ & 1 & $1$ & $1$ & $-\pi$ & $\omega^0$ & $\omega^0$\\
$y^2+y=x^3+11$ & 6075.v2 & $3^{5}\cdot5^{2}$ & $\QQ(\sqrt{-3})$ & 1 & $1$ & $1$ & $-\pi$ & $\omega^0$ & $\omega^{2}$\\
$y^2=x^3+3$ & 3888.i2 & $2^{4}\cdot3^{5}$ & $\QQ(\sqrt{-3})$ & 1 & $1$ & $1$ & $-\pi$ & $\omega^0$ & $\omega^1$\\
$y^2=x^3-4$ & 432.d1 & $2^{4}\cdot3^{3}$ & $\QQ(\sqrt{-3})$ & 1 & $1$ & $1$ & $-\pi$ & $\omega^{2}$ & $\omega^0$\\
$y^2+y=x^3+7031$ & 6075.s2 & $3^{5}\cdot5^{2}$ & $\QQ(\sqrt{-3})$ & 1 & $1$ & $1$ & $-\pi$ & $\omega^{2}$ & $\omega^{2}$\\
$y^2=x^3+12$ & 3888.l2 & $2^{4}\cdot3^{5}$ & $\QQ(\sqrt{-3})$ & 1 & $1$ & $1$ & $-\pi$ & $\omega^{2}$ & $\omega^1$\\
$y^2+y=x^3+31$ & 675.e4 & $3^{3}\cdot5^{2}$ & $\QQ(\sqrt{-3})$ & 1 & $1$ & $1$ & $-\pi$ & $\omega^1$ & $\omega^0$\\
$y^2+y=x^3+281$ & 6075.m2 & $3^{5}\cdot5^{2}$ & $\QQ(\sqrt{-3})$ & 1 & $1$ & $1$ & $-\pi$ & $\omega^1$ & $\omega^{2}$\\
$y^2=x^3+48$ & 3888.n2 & $2^{4}\cdot3^{5}$ & $\QQ(\sqrt{-3})$ & 1 & $1$ & $1$ & $-\pi$ & $\omega^1$ & $\omega^1$\\
$y^2+y=x^3+12$ & 441.d2 & $3^{2}\cdot7^{2}$ & $\QQ(\sqrt{-3})$ & 1 & $1$ & $1$ & $\pi$ & $\omega^0$ & $\omega^0$\\
$y^2=x^3+9$ & 972.b2 & $2^{2}\cdot3^{5}$ & $\QQ(\sqrt{-3})$ & 1 & $1$ & $1$ & $\pi$ & $\omega^0$ & $\omega^{2}$\\
$y^2+y=x^3-37$ & 11907.o1 & $3^{5}\cdot7^{2}$ & $\QQ(\sqrt{-3})$ & 1 & $1$ & $1$ & $\pi$ & $\omega^0$ & $\omega^1$\\
$y^2+y=x^3+600$ & 1323.k2 & $3^{3}\cdot7^{2}$ & $\QQ(\sqrt{-3})$ & 1 & $1$ & $1$ & $\pi$ & $\omega^{2}$ & $\omega^0$\\
$y^2=x^3+36$ & 972.a2 & $2^{2}\cdot3^{5}$ & $\QQ(\sqrt{-3})$ & 1 & $1$ & $1$ & $\pi$ & $\omega^{2}$ & $\omega^{2}$\\
$y^2+y=x^3-1801$ & 11907.p1 & $3^{5}\cdot7^{2}$ & $\QQ(\sqrt{-3})$ & 1 & $1$ & $1$ & $\pi$ & $\omega^{2}$ & $\omega^1$\\
$y^2=x^3-2$ & 1728.o3 & $2^{6}\cdot3^{3}$ & $\QQ(\sqrt{-3})$ & 1 & $1$ & $1$ & $\pi$ & $\omega^1$ & $\omega^0$\\
$y^2=x^3-18$ & 15552.x1 & $2^{6}\cdot3^{5}$ & $\QQ(\sqrt{-3})$ & 1 & $1$ & $1$ & $\pi$ & $\omega^1$ & $\omega^{2}$\\
$y^2+y=x^3-1$ & 243.a1 & $3^{5}$ & $\QQ(\sqrt{-3})$ & 1 & $1$ & $1$ & $\pi$ & $\omega^1$ & $\omega^1$\\
\end{longtable}\endgroup

\subsection{The case $p=2$}

First, we consider the case $p=2$ and $\calK=\QQ(\sqrt{-1})$, fixing the embedding sending $\sqrt{-1}$ to $i\in\CC$. By the classification in \eqref{eq:H_n classification} we have $H_1\cong C_1$, $H_2\cong C_2$, and $H_n\cong C_4\oplus C_{2^{\lceil(n-3)/2\rceil}}\oplus C_{2^{\lfloor(n-3)/2\rfloor}}$ for $n\geq 3$. Taking values in $\mu(\calK)=\mu_4$, one can show that $H_n^*:=\Hom(H_n,\mu_4)$ stabilizes after $n=7$, i.e., $\psi_2$ is trivial on $U^7$.

Choose $\pi:=1+i$ to be the uniformizer of $F=\QQ_2(i)$ and let $s:=1+\pi^3=-1+2i$. Then $\{i,s,5\}$ generates $H_7\cong C_4^{\oplus3}$. Under the normalization we have $\psi_2(-1)=-1$ and $\psi_2(5)=1$. So all $16$ possibilities are determined by $\psi_2(\pi)\in\{i^k\pi\mid k=0,1,2,3\}$, $\psi_2(i)\in\{\pm i\}$ and $\psi_2(s)\in\mu_4$. In particular, we have local conductor
\[
v_F(\mathfrak{f}_\psi)=
\begin{cases*}
    3 & \text{$\psi_2(s)=1$}\\
    4 & \text{$\psi_2(s)\neq1,~\psi_2(s^2)=1$}\\
    6 & \text{$\psi_2(s^2)\neq1$},
\end{cases*}
\]
corresponding to $v_2(N_E)=5,6,8$. Explicit calculations are as below, proved in Appendix \ref{sec:special values calculation}.

\begin{proposition}\label{prop:calculation 2 ramified D_K=-4}
    For the curve $Y_b:y^2=x^3+bx$, let $b=2^{v_2(b)}u$ with $u\in\ZZ_2^{\times}$. Then we have
    \[
    \psi_{Y_b,2}(\pi)=\pi i^{(u^2-1)/8},\qquad \psi_{Y_b,2}(i)=i\cdot(-1)^{(u^2-1)/8},\qquad \psi_{Y_b,2}(s)=(-1)\cdot(-i)^{v_2(b)-u+1}.
    \]
\end{proposition}

After checking LMFDB, representatives of each $\psi_2$ are given in the following table.

\begingroup\small\setlength{\tabcolsep}{5pt}\renewcommand{\arraystretch}{1.25}
\begin{longtable}{@{}cccccccccc@{}}
\caption{Representatives for each $\psi_2$ with $2$ ramified in $\calK=\QQ(\sqrt{-1})$}\label{tab:p=2 ramified database D_K=-4}\\
\hline
Equation of $E$ & LMFDB label & $N_E$ & $\calK$ & $r_{\rm an}$ & $S_{2}$ & $S_{\rm an}$ & $\psi_{E,2}(\pi)$ & $\psi_{E,2}(i)$ & $\psi_{E,2}(s)$\\\hline\endfirsthead
\multicolumn{10}{c}{\tablename\ \thetable\ (continued)}\\\hline
Equation of $E$ & LMFDB label & $N_E$ & $\calK$ & $r_{\rm an}$ & $S_{2}$ & $S_{\rm an}$ & $\psi_{E,2}(\pi)$ & $\psi_{E,2}(i)$ & $\psi_{E,2}(s)$\\\hline\endhead
\hline\endfoot
$y^2=x^3-5x$ & 800.e1 & $2^{5}\cdot5^{2}$ & $\QQ(i)$ & 1 & $1$ & $1$ & $-i\pi$ & $-i$ & $1$\\
$y^2=x^3+3x$ & 288.a2 & $2^{5}\cdot3^{2}$ & $\QQ(i)$ & 1 & $1$ & $1$ & $i\pi$ & $-i$ & $1$\\
$y^2=x^3-10x$ & 6400.k1 & $2^{8}\cdot5^{2}$ & $\QQ(i)$ & 1 & $1$ & $1$ & $-i\pi$ & $-i$ & $-i$\\
$y^2=x^3-250x$ & 6400.n1 & $2^{8}\cdot5^{2}$ & $\QQ(i)$ & 1 & $1$ & $1$ & $i\pi$ & $-i$ & $-i$\\
$y^2=x^3+5x$ & 1600.m2 & $2^{6}\cdot5^{2}$ & $\QQ(i)$ & 1 & $1$ & $1$ & $-i\pi$ & $-i$ & $-1$\\
$y^2=x^3+125x$ & 1600.l2 & $2^{6}\cdot5^{2}$ & $\QQ(i)$ & 1 & $1$ & $1$ & $i\pi$ & $-i$ & $-1$\\
$y^2=x^3-54x$ & 2304.m1 & $2^{8}\cdot3^{2}$ & $\QQ(i)$ & 1 & $1$ & $1$ & $-i\pi$ & $-i$ & $i$\\
$y^2=x^3-6x$ & 2304.c1 & $2^{8}\cdot3^{2}$ & $\QQ(i)$ & 1 & $1$ & $1$ & $i\pi$ & $-i$ & $i$\\
$y^2=x^3-49x$ & 1568.e3 & $2^{5}\cdot7^{2}$ & $\QQ(i)$ & 1 & $1$ & $1$ & $\pi$ & $i$ & $1$\\
$y^2=x^3-25x$ & 800.d3 & $2^{5}\cdot5^{2}$ & $\QQ(i)$ & 1 & $1$ & $1$ & $-\pi$ & $i$ & $1$\\
$y^2=x^3-2x$ & 256.b1 & $2^{8}$ & $\QQ(i)$ & 1 & $1$ & $1$ & $\pi$ & $i$ & $-i$\\
$y^2=x^3-50x$ & 6400.m1 & $2^{8}\cdot5^{2}$ & $\QQ(i)$ & 1 & $1$ & $1$ & $-\pi$ & $i$ & $-i$\\
$y^2=x^3+49x$ & 3136.m4 & $2^{6}\cdot7^{2}$ & $\QQ(i)$ & 1 & $1$ & $1$ & $\pi$ & $i$ & $-1$\\
$y^2=x^3+9x$ & 576.c4 & $2^{6}\cdot3^{2}$ & $\QQ(i)$ & 1 & $1$ & $1$ & $-\pi$ & $i$ & $-1$\\
$y^2=x^3+98x$ & 12544.f2 & $2^{8}\cdot7^{2}$ & $\QQ(i)$ & 1 & $1$ & $1$ & $\pi$ & $i$ & $i$\\
$y^2=x^3+18x$ & 2304.a2 & $2^{8}\cdot3^{2}$ & $\QQ(i)$ & 1 & $1$ & $1$ & $-\pi$ & $i$ & $i$\\
\end{longtable}\endgroup

Next, we consider the case $p=2$ and $\calK=\QQ(\sqrt{-2})$, fixing the embedding sending $\pi:=\sqrt{-2}$ to the upper half-plane of $\CC$. By the classification in \eqref{eq:H_n classification} we have $H_1\cong C_1$, $H_2\cong C_2$, $H_3\cong C_4$, and $H_n\cong C_2\oplus C_{2^{\lfloor n/2\rfloor}}\oplus C_{2^{\lceil n/2\rceil-2}}$ for $n\geq 4$. Taking values in $\mu(\calK)=\{\pm1\}$, it is easy to see that $H_n^*:=\Hom(H_n,\{\pm1\})$ stabilizes after $n=5$, i.e., $\psi_2$ is trivial on $U^5$.

Fix $\pi:=\sqrt{-2}$ to be the uniformizer, so $\mathcal{O}_F=\ZZ[\pi]$. One can show that the images of $\{-1,v:=1+\pi,5\}$ generate $H_5\cong C_2\oplus C_4\oplus C_2$. Under the normalization $\psi_2\vert_{\ZZ_2^{\times}}=\eta_{F/\QQ_2}$, we have $\psi_2(5)=-1$ and $\psi_2(-1)=-1$. Hence $\psi_2$ is non-trivial on $U^4$, giving $v_F(\mathfrak{f}_{\psi})=5$ and $v_2(N_E)=8$. Since the image of $v$ has order exactly $4$ in $H_5$ and $\mu(\calK)=\{\pm1\}$, all $4$ possibilities for $\psi_2$ are determined by $\psi_2(\pi)\in\{\pm\pi\}$ and $\psi_2(v)\in\{\pm1\}$, which can be realized as below, as explained in Appendices \ref{sec:special values calculation} and \ref{sec:p-Sha calculation}.

\begingroup\small\setlength{\tabcolsep}{5pt}\renewcommand{\arraystretch}{1.25}
\begin{longtable}{@{}ccccccccc@{}}
\caption{Representatives for each $\psi_2$ with $2$ ramified in $\calK=\QQ(\sqrt{-2})$}\label{tab:p=2 ramified database D_K=-8}\\
\hline
Equation of $E$ & LMFDB label & $N_E$ & $\calK$ & $r_{\rm an}$ & $S_{2}$ & $S_{\rm an}$ & $\psi_{E,2}(\pi)$ & $\psi_{E,2}(v)$\\\hline\endfirsthead
\multicolumn{9}{c}{\tablename\ \thetable\ (continued)}\\\hline
Equation of $E$ & LMFDB label & $N_E$ & $\calK$ & $r_{\rm an}$ & $S_{2}$ & $S_{\rm an}$ & $\psi_{E,2}(\pi)$ & $\psi_{E,2}(v)$\\\hline\endhead
\hline\endfoot
$y^2=x^3-x^2-83x+287$ & 6400.x2 & $2^{8}\cdot5^{2}$ & $\QQ(\sqrt{-2})$ & 1 & $1$ & $1$ & $-\pi$ & $-1$ \\
$y^2=x^3-30x+56$ & 2304.h2 & $2^{8}\cdot3^{2}$ & $\QQ(\sqrt{-2})$ & 1 & $1$ & $1$ & $-\pi$ & $1$ \\
$y^2=x^3+x^2-3x+1$ & 256.a2 & $2^{8}$ & $\QQ(\sqrt{-2})$ & 1 & $1$ & $1$ & $\pi$ & $-1$ \\
$y^2=x^3+x^2-163x+657$ & 12544.c2 & $2^{8}\cdot7^{2}$ & $\QQ(\sqrt{-2})$ & 1 & $1$ & $1$ & $\pi$ & $1$ \\
\end{longtable}\endgroup

We have seen that all the curves listed in the tables above satisfy $p$-BSD. Since they exhaust all possibilities, we obtain the following result from the transfer principle, Theorem \ref{thm:main_transfer}.

\begin{theorem}\label{thm:BSD ramified}
    Let $E/\QQ$ be an elliptic curve with complex multiplication by $\calK$, of analytic rank $1$. Let $p$ be a prime ramified in $\calK$. Then the $p$-BSD holds for $E/\QQ$.
\end{theorem}

%% file: 108-split.tex
\section{Birch--Swinnerton-Dyer formula in the split case} \label{sec:split}

In this section, we assume that $p$ splits in $\calK$. In this case, we prove directly that $\Delta(E/\calK) = 0$, without performing numerical experiments. The overall arguments proceed in the same manner as those from Section \ref{sec:plus-minus} to Section \ref{sec:magic}, and we often cite the results there. The advantages in the split case are twofold:
\begin{itemize}
    \item On the arithmetic side, the plus-minus theory in \cite{FanWan} is simpler, and more importantly, we have an explicit understanding of the regulator map, which we shall introduce in Section \ref{sec:regulator-split}.
    \item On the analytic side, the local toric integral at $p$ can be computed explicitly, as we shall see in Section \ref{sec:toric-split} and Proposition \ref{prop:U-sigma-k}.
\end{itemize}

\subsection{Arithmetic side: A local plus-minus theory in the split case} \label{sec:plus-minus-split}

In this subsection we adapt Section~\ref{sec:plus-minus}.

\subsubsection{Plus-minus subspaces in the split case}
We use \cite[Section 7]{FanWan} to deal with the case $p$ splits in $\calK$ individually. We put $p \calO_{\calK} = \frp_\sigma\barfrp_\sigma$. Then
\[
\rmH^{1}(\calK_{p}, T) = \rmH^{1}(\calK_{\frp}, T) \oplus \rmH^{1}(\calK_{\barfrp}, T).
\]
In this subsection, $\frp$ abbreviates $\frp_\sigma$ and $\barfrp$ its conjugate. Set $\rmH^1_+(\calK_p,V_\sigma)=\rmH^1(\calK_{\frp},V_\sigma)$ and $H^1_-(\calK_p,V_\sigma)=\rmH^1(\calK_{\barfrp},V_\sigma)$. Define their integral versions $\rmH^{1}_{\pm}(\calK_{p}, T_{\sigma})$ by their inverse images in $\rmH^1(\calK_p,T_\sigma)$. The following lemma is straightforward (for example, \eqref{eq:finite-line-split} follows in the same way as in Proposition \ref{prop:finite-lines}).

\begin{lemma} \label{lem:split-second-defect}
    For each $\sigma$, we have
    \begin{equation} \label{eq:finite-line-split}
    \dim_{L} \rmH^{1}(\calK_{\frp}, V_{\sigma}) = \dim_{L} \rmH^{1}(\calK_{\barfrp}, V_{\sigma}) = 1, \quad \rmH^{1}_{\rmf}(\calK_{p}, V_{\sigma}) = \rmH^{1}_{\rmf}(\calK_{\frp}, V_{\sigma}).
    \end{equation}
    Moreover,
    \[
    \rmH^{1}_{+}(\calK_p, T_{\sigma}) = \rmH^{1}(\calK_{\frp}, T_{\sigma}) \oplus \rmH^{1}(\calK_{\barfrp}, T_{\sigma})_{\tor}, \quad \rmH^{1}_{-}(\calK_p, T_{\sigma}) = \rmH^{1}(\calK_{\frp}, T_{\sigma})_{\tor} \oplus \rmH^{1}(\calK_{\barfrp}, T_{\sigma}).
    \]
    In particular, with the notations of \eqref{eq:delta2}, we have
    \[
    A_{\sigma} := \rmH^{1}(\calK_p, T_{\sigma})_{\tf} \cap \rmH^{1}_{\rmf}(\calK_p, V_{\sigma}) \simeq Q_{\sigma},
    \]
    and $\delta_{\sigma}^{\prime\prime} = (1)$.
\end{lemma}

\subsubsection{Control theorem in the split case}
In this subsection, we establish the counterparts of Theorems \ref{thm:fan_wan_control} and \ref{thm:fan_wan_tnc} in the split case. We keep the notation used there.

As in Theorem \ref{thm:classical_CM_main_conjecture}, we let $z_{\psi}$ be the elliptic-unit zeta morphism. We have the following result parallel to Theorem \ref{thm:fan_wan_control}, combined with Sections \ref{sec:appear-tam} and \ref{sec:appear-sha}.

\begin{theorem}[Control theorem: split case] \label{thm:fan_wan_control_split}
    We have
    \[
    \calO_{L} q_-\loc_{p}(z_{\psi}) = \dfrac{\abs{\Sel_{-}(\calK, W^{\dagger})}_{L}}{\abs{\rmH^{0}(\calK, W)}_{L} \abs{\rmH^{0}(\calK, W^{\dagger})}_{L}} \left( \prod_{v \in S_0} \abs{\rmH^{0}(\calK_v, W^{\dagger})}_{L} \right) Q_{\sigma}.
    \]
The first equality lies in $Q_\sigma\otimes L$. Moreover, in the tensor product of the two active lines we have
    \[
    \bigotimes_{\sigma = \iota, \iota^{c}} \calO_{L} \loc_{\frp_{\sigma}} z_{\psi, \sigma} = \left( \dfrac{\# \Sha(E/\calK)}{\# E(\calK)[p^\infty]^{2}}\prod_{v \nmid p} \Tam_{v}(E/\calK)\, \# E(\QQ_p)[p^{\infty}]^{2} \right) (\frl_{\iota} \frl_{\iota^{c}})^{2} \bigotimes_{\sigma = \iota, \iota^{c}} Q_{\sigma}.
    \]
\end{theorem}

\begin{proof}
    The proof of Theorem \ref{thm:fan_wan_tnc} carries over and proves the first equality. The last equality follows by taking the product of the equalities under the two embeddings $\iota$ and $\iota^{c}$ and applying Lemmas \ref{lem:torsion-appear}, \ref{lem:tamagawa-appear} and Propositions \ref{prop:length} and  \ref{prop:find_sha}. The second defect $\delta^{\prime\prime}_{\sigma}$ is the unit ideal for $\sigma \in \{\iota, \iota^{c} \}$ by Lemma \ref{lem:split-second-defect}.
\end{proof}

\subsubsection{Regulator maps} \label{sec:regulator-split}
We now establish the counterpart of Theorem \ref{thm:fan_wan_tnc}; a key improvement is that we have an explicit description of the regulator map appearing here.

For any Dirichlet character $\rho$ of conductor $p^{n}$ for $n>0$, we define the Gauss sum $G(\rho)$ by
\[
G(\rho) := \sum_{u \in (\ZZ/p^{n})^{\times}} \rho(u)^{-1} \zeta_{p^n}^{u}.
\]
Here $(\zeta_{p^{n}})$ is a compatible sequence of $p^{n}$-th roots of unity, and the element $u$ acts on $\zeta_{p^{n}}$ through the cyclotomic character, following \cite[\S 2.8]{LZ-ijnt}.

We write the local representation $V|_{G_{\calK_{\frp}}}$ as $V$ now for simplicity and identify $\calK_{\frp}$ with $\QQ_p$. Then we can write it as $V = U(1)\rho^{-1}$, where $U$ is an unramified rank-one representation of $G_{\QQ_p}$. Let $T_{U} \subseteq U$ be a $G_{\QQ_p}$-stable $\calO_{L}$-lattice. Choose a basis $\ttd_0$ of $\bfD_{\cris}(U)$ and write $\varphi(\ttd_0) = a \ttd_0$ for the crystalline Frobenius operator $\varphi$.

We choose an unramified $\ZZ_p$-extension $F_\infty/\QQ_p$, and put $G=\Gal(F_\infty(\mu_{p^\infty})/\QQ_p)$. For $K_{\infty} := F_\infty(\mu_{p^\infty})$, we define the Iwasawa cohomology
\[
\rmH^{1}_{\Iw}(K_{\infty}, T_{U}) := \varprojlim_{F \subseteq K_{\infty}} \rmH^{1}(F, T_{U}), \quad \rmH^{1}_{\Iw}(K_{\infty}, U) := \rmH^{1}_{\Iw}(K_{\infty}, T_{U}) \otimes_{\calO_{L}} L
\]
where $F$ ranges over the finite extensions of $\QQ_p$ contained in $K_{\infty}$ and the transition maps are corestriction. Let $\hatF_{\infty}$ denote the completion of $F_{\infty}$. The two-variable regulator map is given by
\[
\calL_{U}^{G}: \rmH^{1}_{\Iw}(K_{\infty}, U) \rightarrow \calH_{\hatF_{\infty}, L}(G) \otimes_{L} \bfD_{\cris}(U),
\]
where $\calH_{\hatF_{\infty}, L}(G) := \calH_{\hatF_{\infty}}(G) \hotimes_{\QQ_p} L$ with $\calH_{\hatF_{\infty}}(G)$ being the algebra of locally analytic distributions on $G$ with values in $\hatF_{\infty}$.

For any character $\eta: G \rightarrow L^{\times}$, evaluation of a distribution defines
\[
\calL_{U}^{G}(y)(\eta) := (\ev_{\eta} \otimes 1)\calL_{U}^{G}(y),
\]
for any fixed cohomology class $y \in \rmH^{1}_{\Iw}(K_{\infty}, U)$. Separately, twisting and specialization defines the specialization map
\[
\sp_{\eta}: \rmH^{1}_{\Iw}(K_{\infty}, U) \rightarrow \rmH^{1}(\QQ_p, U\eta^{-1}),
\]
and we write $y_{\eta, 0} := \sp_{\eta}(y)$.

We consider the special case where $\eta = \epsilon_{\cyc}^{-1}\rho$. Then we have
\[
\calL_{U}^{G}(y)(\epsilon_{\cyc}^{-1}\rho) \in \bfD_{\cris}(U).
\]

Let $\tte_j$ denote a basis of $L(j)$, with compatible tensor identifications $\tte_j\otimes\tte_l=\tte_{j+l}$, and let $\tte_\rho$ be an $L$-basis of the $L(\rho)$. Because $V \otimes_{L} L(1) \otimes_{L} L(\rho) = U$, we have an $L$-linear isomorphism
\[
\Theta_{\rho}: \bfD_{\dR}(V) \xrightarrow{\sim} \bfD_{\dR}(U) = \bfD_{\cris}(U)
\]
by $z \mapsto z \otimes t \tte_{1}$ when $\rho = 1$, and $z \mapsto z \otimes t \tte_{-1} \otimes G(\rho)^{-1} \tte_{\rho}$, when $\rho \neq 1$. 

We denote by $\log_{V}: \rmH^{1}_{\rmf}(\QQ_p,V) \rightarrow \bfD_{\dR}(V)$ the Bloch--Kato logarithm; note that $\Fil^{0} \bfD_{\dR}(V) = 0$. We then have the following theorem describing the regulator map $\calL_{U}^{G}(y)(\epsilon_{\cyc}^{-1}\rho)$.

\begin{theorem} \label{thm:LZV-regulator}
With the notation above, for any cohomology class $y \in \rmH^{1}_{\Iw}(K_{\infty}, U)$ whose specialization $y_{\epsilon_{\cyc}^{-1}\rho, 0} \in \rmH^{1}_{\rmf}(\QQ_p, V)$, we have
\[
\calL_{U}^{G}(y)(\epsilon_{\cyc}^{-1}\rho) = \begin{cases}
    \dfrac{1-a/p}{1-a^{-1}} \Theta_{1}(\log_{V}(y_{\epsilon_{\cyc}^{-1}, 0})), &\quad n=0, \\
    a^{n} \Theta_{\rho}(\log_{V}(y_{\epsilon_{\cyc}^{-1}\rho, 0}), & \quad n > 0
\end{cases} \quad  \in \bfD_{\cris}(U).
\]
\end{theorem}

\begin{proof}
    For odd $p$, we apply \cite[Theorem 4.15, Theorem B.5]{LZ-ijnt} with $j = -1$. In the unramified case, the regulator operator is $(1-p^{-1}\varphi)(1-\varphi^{-1})^{-1}$, and in the ramified case, the regulator operator is $G(\rho)^{-1} \varphi^{n}$. For $p=2$, use \cite[Theorem 1.3, Remark 1.4]{Venjakob-regulator}. Its extra hypotheses hold here: $U$ has weight zero, so $\Fil^{-1} \bfD_{\cris}(U^*(1))=\bfD_{\cris}(U^*(1))$; its Frobenius scalar is $p^{-1}a^{-1}$, which is neither $1$ nor $p^{-1}$. Indeed by the purity result, $a$ has complex absolute value $p^{1/2}$, and hence $a\ne1,p^{-1}$; also $a\ne p$.
\end{proof}

We write $\Reg_p(y) := \calL_{U}^{G}(y)(\epsilon_{\cyc}^{-1}\rho)$ for the expression appearing in Theorem \ref{thm:LZV-regulator}. We choose a Néron differential $\omega = \omega_{E}$ over $\ZZ_p$. Note that
\[
\bfD_{\dR}(V) = \bfD_{\dR}(V)/\Fil^{0}\bfD_{\dR}(V) \simeq \Lie(E)(\QQ_p) \otimes_{\QQ_p} L,
\]
we define $\omega^{\vee}$ to be the unique vector in $\bfD_{\dR}(V)$ satisfying $\omega(\omega^{\vee}) = 1$. For $y \in \rmH^{1}(\QQ_p,V)$, define $\log_{\omega}(x) := \omega(\log_{V}(y))$, which gives the identity
\[
\log_{V}(x) = \log_{\omega}(x) \omega^{\vee}.
\]

In the unramified case, i.e., $\rho = 1$, we have $V = U(1)$. Multiplication by $t \tte_{1}$ gives the map $\Theta_1$ as before. We choose explicitly $\ttd_{0} := \Theta_1(\omega^{\vee}) = \omega^{\vee} \otimes t \tte_{-1}$; then $\Theta_{1}(\log_{V}(y)) = \log_{\omega}(y)\ttd_{0}$. Theorem \ref{thm:LZV-regulator} gives
\begin{equation} \label{eq:regulator-unramified-ttd}
\Reg_{p}(y) = \dfrac{1-a/p}{1-a^{-1}} \log_{\omega}(y) \ttd_{0},
\end{equation}

Now suppose $\rho$ is ramified of conductor $p^{n}$. We put $\ttd_{1} := t^{-1} \tte_{1}$ and define $\ttd_{\rho} := G(\rho)\ttd_{0} \otimes \ttd_{1} \otimes \tte_{\rho^{-1}} \in \bfD_{\dR}(V)$, then one checks directly that $\Theta_{\rho}(\ttd_{\rho}) = \ttd_{0}$. Since both $\ttd_{\rho}$ and $\omega^{\vee}$ are bases of the same one-dimensional $L$-vector space, there exists a unique $b_{\rho} \in L^{\times}$ such that $\omega^{\vee} = b_{\rho}\ttd_{0}$. It follows that
\[
\Theta_{\rho}(\log_{V}(y)) = \Theta_{\rho}(\log_{\omega}(x)\omega^{\vee}) = b_{\rho} \log_{\omega}(x) \Theta_{\rho}(\ttd_{\rho}) = b_{\rho} \log_{\omega}(x) \ttd_{0}.
\]
Then Theorem \ref{thm:LZV-regulator} gives
\begin{equation} \label{eq:regulator-ramified-ttd}
\Reg_p(y) = a^{n} \Theta_{\rho}(\log_{V}(y)) = a^{n} b_{\rho} \log_{\omega}(y) \ttd_{0}.
\end{equation}
In either case, we define $\Reg_{p,0}(y) \in L$ as the unique scalar such that $\Reg_{p}(y) = \Reg_{p,0}(y) \ttd_{0}$. With the above preparation, especially equations \eqref{eq:regulator-unramified-ttd} and \eqref{eq:regulator-ramified-ttd} derived from Theorem \ref{thm:LZV-regulator}, we have the following corollary.

\begin{corollary} \label{cor:scalar-regulator}
    When $E$ has good reduction at $p$, let $\alpha_1$ be the $p$-adic unit root of the polynomial $X^2 - a_p(E) X + p$. Then
    \[
    \Reg_{p,0}(x) = \dfrac{1-\alpha_1/p}{1-\alpha_1^{-1}} \log_{\omega}(x), \quad \Reg_{p,0}(A_{\sigma}) = \dfrac{\# E(\QQ_p)[p^{\infty}]}{\Tam_{p}(E/\QQ) (1-\alpha_1^{-1})^{2}} \cdot \calO_{L},
    \]
    and when $E$ has additive reduction at $p$, then
    \[
    \Reg_{p,0}(x) = a^n b_{\rho} \log_{\omega}(x), \quad \Reg_{p,0}(A_{\sigma}) = a^{n} b_{\rho} \dfrac{\# E(\QQ_p)[p^{\infty}]}{\Tam_{p}(E/\QQ)} \cdot \calO_{L}.
    \]
\end{corollary}
\begin{proof}
    This is a combination of Lemma \ref{lem:tam-at-p-appear} and Theorem \ref{thm:LZV-regulator}. In fact, in the good-reduction case, denote the two roots of $X^{2}-a_p(E)X+p$ by $\alpha_1$ and $\beta_1$, where $\alpha_1$ is a $p$-adic unit and $\beta_1 = p \alpha_1^{-1}$. Then the expression for $\Reg_{p,0}(x)$ follows, and
    \begin{align*}
    \Reg_{p,0}(A_{\sigma}) = \dfrac{1-\alpha_1/p}{1-\alpha_1^{-1}} \log_{\omega_E}(A_{\sigma}) &= \dfrac{1-\alpha_1/p}{1-\alpha_1^{-1}} \dfrac{\# E(\QQ_p)[p^{\infty}]}{\Tam_{p}(E/\QQ) P_p(E)} \calO_{L} \\ &= \dfrac{\# E(\QQ_p)[p^{\infty}]}{\Tam_{p}(E/\QQ) (1-\alpha_1^{-1})^{2}} \calO_{L}.
    \end{align*}
    Here the second equality follows from Lemma \ref{lem:tam-at-p-appear}, and the third equality follows from $P_{p}(E) = (1-\alpha_1/p)(1-\alpha_1^{-1})$. The additive case follows similarly by using Lemma \ref{lem:tam-at-p-appear} and noting that $P_p(E) = 1$.
\end{proof}

Following the proof of Corollary \ref{cor:scalar-regulator}, we denote
\[
R_{\sigma} := \begin{cases}
    \dfrac{1-\alpha_1/p}{1-\alpha_1^{-1}} &\quad \text{when $E$ has good reduction at $p$}, \\ 
    a^{n} b_{\rho}, &\quad \text{when $E$ has additive reduction at $p$}.
\end{cases}
\]

\subsubsection{"Tamagawa number conjecture" in the split case}

We keep the notation of Section \ref{sec:plus-minus}. In Appendix \ref{app:elliptic-units}, we formulate the $\Lambda_{m}$-family of elliptic units $\bfz_{\bpsi,\sigma} \in \rmH^{1}_{\Iw}(\calK_{\infty}/\calK, V_{\bpsi,\sigma})$ in \eqref{eq:elliptic-unit-fan-wan-family}.

On the other hand, for the regulator, let $\chi_{U_{\sigma}}$ denote the Galois character of the fixed unramified representation $U_{\sigma}$. Then we define the local character family
\[
\bfeta_{\sigma}(g) := \chi_{U_{\sigma}}(g) \bpsi_{\sigma}(g)^{-1}, \quad\text{ for } g \in G_{\QQ_p}.
\]
At every specialization $\phi: \Lambda_{m} \rightarrow L_{\phi}$, we put $\eta_{\sigma,\psi} := \phi \circ \bfeta_{\sigma}$. Then
\[
U_{\sigma} \eta_{\sigma,\phi}^{-1} = V_{\bpsi,\sigma} \otimes_{\Lambda_m, \phi} L_{\phi}|_{G_{\frp_{\sigma}}},
\]
and at the origin, we have $\eta_{\sigma,0} = \epsilon_{\cyc}^{-1}\rho_{\sigma}$. Now, for any given cohomology class $y \in \rmH^{1}_{\Iw}(K_{\infty, U})$, evaluating the distribution along $\bfeta_{\sigma}$ gives the regulator
\[
\calL_{U}^{G}(y)(\bfeta_{\sigma}) := (\ev_{\bfeta_{\sigma}} \otimes 1) \calL_{U}^{G}(y) \in \Lambda_{m} \otimes_{L} \bfD_{\cris}(U).
\]
With the same choice of basis $\ttd_{0}$ as above, we define the $\Lambda_{m}$-adic regulator map $\bReg_{p,0}(y)$ to be the unique scalar in $\Lambda_m$ such that $\calL_{U}^{G}(y)(\bfeta_{\sigma}) = \bReg_{p,0}(y) \ttd_{0}$. At the origin, it specializes to the scalar regulator map $\Reg_{p,0}$ in the previous section.

We now define the \emph{elliptic-unit $p$-adic $L$-function} as
\begin{equation} \label{eq:elliptic-unit-L-function}
\scL_{\sigma}^{\elli} := \bReg_{p,0} (\log_{\omega}(\loc_{\frp_{\sigma}} \bfz_{\bpsi,\sigma})) \in \Lambda_{m}.
\end{equation}
Its value at the origin is exactly $R_{\sigma} \log_{\omega}(\loc_{\frp_{\sigma}} z_{\psi, \sigma})$.  Then the following theorem follows from Theorem \ref{thm:fan_wan_control_split}.

\begin{theorem}["Tamagawa number conjecture" in the split case] \label{thm:fan-wan-tnc-split}
    \[
    (\scL_{\sigma}^{\elli}(0)) = \dfrac{\abs{\Sel_{-}(\calK, W^{\dagger})}_{L}}{\abs{\rmH^{0}(\calK, W)}_{L} \abs{\rmH^{0}(\calK, W^{\dagger})}_{L}} \left( \prod_{v \in S_0} \abs{\rmH^{0}(\calK_v, W^{\dagger})}_{L} \right) \left( R_{\sigma} \dfrac{\# E(\QQ_p)[p^{\infty}]}{\Tam_{p}(E/\QQ) P_p(E)} \right).
    \]
\end{theorem}
In contrast with Theorem \ref{thm:fan_wan_tnc} in the nonsplit case, where $\delta(E,\scL_\sigma)^{-1}$ is nonexplicit, here, when $p$ splits in $\calK$, we have an explicit expression for this term thanks to our understanding of the regulator map. However, an additional comparison is needed between $\scL_{\sigma}^{\elli}$ and the (renormalized) anticyclotomic $p$-adic $L$-function $\scL_{\sigma}$ in the notation of Section \ref{sec:GZW-comparison}. This comparison is parallel to that in Section \ref{sec:transfer-principle}.

\subsection{Analytic side: The Gross--Zagier--Waldspurger comparison in the split case} \label{sec:GZW-split}

On the analytic side, the computations in Section \ref{sec:GZW-comparison} carry over directly to the case in which $p$ splits in $\calK$ upon putting $d_{p} = 0$ \footnote{See \cite[Proposition 3.10]{FanWan} which displays separately the split and nonsplit interpolation formulas.}, since neither the $p$-adic Waldspurger formula in \cite{FanWan} (as Theorems \ref{thm:fan-wan-p-adic-l-function} and \ref{thm:fan-wan-p-adic-waldspurger}) nor the Gross--Zagier formula (as Theorem \ref{thm:fan-wan-gross-zagier}) imposes any requirement on the decomposition type of $p$ in $\calK$.

We still keep the same renormalization of the $p$-adic $L$-function $\calL_{E}(\phi)$, as 
\begin{equation}\label{eq:renormalize-calL_E-split}
 \scL_E(\phi)=\frac{\Lambda(1,\eta)\Lambda(1,\omega_E^{-1})}
 {4\calE_E^{(p)}(\phi)}\calL_E(\phi).
\end{equation}
Then Theorem \ref{thm:final_comparison} carries over in this case, which we restate as follows.

\begin{theorem}[Gross--Zagier--Waldspurger comparison]
\label{thm:final_comparison_split}
With the notation of Section \ref{sec:GZW-comparison}, one has
\[
 \left(\frac{\scL_E(0)}{L'(\psi_E,1)}\right)^2
 =\left(\frac{r_p^{\tor}}{\hatr I_E}\right)^2
 \frac{k_{1}k_{2}}{k_{1}^ck_{2}^c}
 \frac{L_0^4}{R_0^2}.
\]
\end{theorem}

\subsection{The $p$-adic BSD defect and the conclusion}
In this section, we combine Sections \ref{sec:plus-minus-split} and \ref{sec:GZW-split} to obtain a $p$-adic BSD defect formula for $\Delta(E/\calK)$ parallel to Theorem \ref{thm:p-adic-defect}, and then identify the precise remaining condition for $\Delta_p(E/\calK)=0$.

\subsubsection{The $p$-adic BSD defect formula}
Recall that, in \eqref{eq:elliptic-unit-L-function}, we defined the elliptic-unit $p$-adic $L$-function $\scL_{\sigma}^{\elli}$. We define the constant $w_{\sigma} := \scL_{\sigma}(0)/\scL_{\sigma}^{\elli}(0)$. On the arithmetic side, we then have the following formula parallel to Proposition \ref{prop:embedding-arithmetic}.

\begin{proposition}[Arithmetic side identity: split case] \label{prop:embedding-arithmetic-split} 
We have the identity
\[
    (\scL_{E,\iota}(0) \cdot \scL_{E,\iota^{c}}(0)) = \left( \dfrac{\# \Sha(E/\calK)}{\# E(\calK)[p^{\infty}]^{2}} \prod_{v \nmid p} \Tam_{v}(E/\calK) \cdot \# E(\calK_{\frp})[p^{\infty}]^{2} \right) (\frl_{\iota} \frl_{\iota^{c}})^{2} \prod_{\sigma = \iota, \iota^{c}} \dfrac{w_{\sigma} R_{\sigma}}{C_{\sigma}}.
\]
\end{proposition}
\begin{proof}
    It follows from Theorem \ref{thm:fan-wan-tnc-split}, together with Lemma \ref{lem:tamagawa-appear} and Propositions \ref{prop:length} and \ref{prop:find_sha}.
\end{proof}

Next we focus on the analytic side. Following the same proof of Proposition \ref{prop:embedding-analytic}, we have the following.

\begin{proposition}[Analytic side identity: split case] \label{prop:embedding-analytic-split} 
    We have the following identity
    \[
    \left( \dfrac{L^{\ast}(E/\calK, 1)}{\Omega_{E/\calK} \Reg(E/\calK)} \right) (\frl_{\iota}\frl_{\iota^{c}})^{2} = (\scL_{E,\iota}(0) \cdot \scL_{E,\iota^{c}}(0)) \prod_{\sigma = \iota, \iota^{c}} \dfrac{C_{\sigma}^{2} \pi_{\sigma}}{r_{p,\sigma}^{\tor}},
    \]
    where $\pi_{\sigma} := \hatr_{\sigma} a_{E,\sigma}$ following the notation in Section \ref{sec:BSD-defect-formula}.
\end{proposition}

Combining Propositions \ref{prop:embedding-arithmetic-split} and \ref{prop:embedding-analytic-split}, we have
    \[
    \dfrac{L^{\ast}(E/\calK, 1)}{\Omega_{E/\calK} \Reg(E/\calK)}= \left( \dfrac{\# \Sha(E/\calK)}{\# E(\calK)[p^{\infty}]^{2}} \prod_{v \nmid p} \Tam_{v}(E/\calK) \cdot \# E(\calK_{\frp})[p^{\infty}]^{2} \right) \prod_{\sigma = \iota, \iota^{c}} \dfrac{w_{\sigma} R_{\sigma} C_{\sigma} \pi_{\sigma}}{r_{p,\sigma}^{\tor}}.
    \]
Now use the expression of $C_{\sigma}$ in Lemma \ref{lem:tam-at-p-appear} to obtain
\[
\dfrac{L^{\ast}(E/\calK, 1)}{\Omega_{E/\calK} \Reg(E/\calK)}= \left( \dfrac{\# \Sha(E/\calK)}{\# E(\calK)[p^{\infty}]^{2}} \prod_{v} \Tam_{v}(E/\calK) \right) \prod_{\sigma = \iota, \iota^{c}} \dfrac{w_{\sigma} R_{\sigma}P_p(E) \pi_{\sigma}}{r_{p,\sigma}^{\tor}}
\]
Here we use $\QQ_p \simeq \calK_{\frp} = \calK_{\barfrp}$, and hence $\Tam_{\frp}(E/\calK) = \Tam_{\barfrp}(E/\calK) = \Tam_{p}(E/\QQ)$. This gives the following $p$-adic BSD formula.

\begin{theorem}[$p$-adic BSD defect formula] \label{thm:p-adic-defect-split}
    We have
    \[
    \Delta_p(E/\calK) = v_p \left( P_p(E)^{2} \prod_{\sigma = \iota, \iota^{c}} \dfrac{w_{\sigma} R_{\sigma} \pi_{\sigma}}{r_{p,\sigma}^{\tor}} \right).
    \]
\end{theorem}
In contrast with the $p$-adic BSD defect formula in Theorem \ref{thm:p-adic-defect}, in this split case, the defects $\delta_{\sigma}$, $\delta_{\sigma}^{\prime}$, and $\delta_{\sigma}^{\prime\prime}$ disappear, and the more explicit regulator $R_{\sigma}$ appears. The next task is to show that the term on the right-hand side is actually zero, so $\Delta_{p}(E/K) = 0$.

\subsubsection{Computation of the finite toric integral} \label{sec:toric-split}
In this subsection, we aim to prove the following proposition.
\begin{proposition} \label{prop:r_p_tor-split}
    With the notation above, put $\varphi(\ttd_{0}) = a_{\sigma} \ttd_{0}$. Let $n_{\sigma}$ be the conductor of $\rho_{\sigma}$; then we have
    \[
    r_{p}^{\tor} =_p \begin{cases}
        (1-a_{\sigma}/p)^{2}, &\quad\text{when $n_{\sigma} = 0$}, \\ 
        (a_{\sigma}/p)^{n_{\sigma}} G(\rho_{\sigma}^{-1}), &\quad\text{when $n_{\sigma} > 0$}.
    \end{cases}
    \]
Here and henceforth, we use the notation "$=_p$" to indicate an equality up to a $p$-adic unit.
\end{proposition}

To prove this proposition, we first make the following four preparatory steps, aiming to obtain an automorphic expression for $r_{p}^{\tor}$ in \eqref{eq:automorphic-expression-rptor-split}.

\underline{(1) Setup.} To carry out the concrete computation, let $F = \QQ_p$ and adopt the convention that $\vol(\ZZ_p, \dif x) = \vol(\ZZ_p^{\times}, \dif^{\times} x) = 1$. Thus $\dif^{\times} x = (1-p^{-1})^{-1} \abs{x}^{-1} \dif x$. We choose and fix an additive character $\psi: F \rightarrow \CC^{\times}$ of conductor $\ZZ_p$ with $\psi(p^{-m}) = \zeta_{p^{m}}$. The \emph{Fourier transform} of a Schwartz function $F \in \scS(F)$ is given by
\[
\hatF(y) := \int_{F} F(x) \psi(xy) \dif x, \quad \widehat{\hatF}(x) = F(-x),
\]
and for any character $\theta: F^{\times} \rightarrow \barQQ_{p}^{\times}$, Tate's zeta integral is defined as
\[
Z(F, \theta,s) := \int_{F^{\times}} F(x) \theta(x) \abs{x}^{s} \dif^{\times} x,
\]
with the functional equation
\begin{equation} \label{eq:tate-functional-equation-split}
Z(\hatF, \theta^{-1}, 1-s) = \gamma(s, \theta, \psi_{p}) Z(F,\theta, s), \quad \gamma(s,\theta,\psi_p) = \epsilon(s,\theta,\psi_p) \dfrac{L(1-s,\theta^{-1})}{L(s,\theta)}.
\end{equation}

\underline{(2) Eisenstein sections.}
Following \cite[Section 3.3]{FanWan}, we consider the normalized principal series representation $I(1,1)$ satisfying
\[
I(1,1) := \left\{ h: \GL_2(F) \rightarrow \CC: h \left( \begin{bmatrix}
    a & b \\ 0 & d
\end{bmatrix} g \right) = \abs{\dfrac{a}{d}}^{1/2} h(g) \right\}.
\]
We define the \emph{Jacquet integral} of $h$ as
\[
J_{h}(g) := \int_{F} h(wn(x)g) \psi_{p}(x) \dif x, \quad w = \begin{bmatrix}
    0 & -1 \\ 1 & 0
\end{bmatrix}, \, n(x) = \begin{bmatrix}
    1 & x \\ 0 & 1
\end{bmatrix}.
\]
We put $F(x) = h(wn(x))$, then the Jacquet integral of $h$ simplifies to 
\[
J(h)(a(y)) = \abs{y}^{-1/2} \int_{F} F(x/y) \psi_p(x) \dif x = \abs{y}^{1/2} \int_{F} F(u) \psi_p(yu) \dif u = \abs{y}^{1/2} \hatF(y), \, a(y) := \begin{bmatrix}
    y & 0 \\ 0 & 1
\end{bmatrix}.
\]
Here we use the matrix identity $wn(x)a(y) = \begin{bmatrix}
    1 & 0 \\ 0 & y
\end{bmatrix} w n(x/y)$, and then a change of variables $x = yu$.

In particular, let $h^{\circ}$ be the unique spherical section with $h^{\circ}(1) = 1$. Its associated function $F^{\circ}(x) = h^{\circ}(wn(x))$ can be written as
\[
F^{\circ}(x) = \begin{cases}
    1, \quad & \abs{x} \leq 1, \\
    \abs{x}^{-1}, \quad & \abs{x} > 1,
\end{cases}
\]
and hence it is straightforward that $J_{h^{\circ}}(\bfone) = 1-p^{-1}$. The corresponding Whittaker functional sends $h^{\circ}$ to the normalized Whittaker function $W^{\circ}$ with $W^{\circ}(1) = 1$, and hence we have 
\begin{equation} \label{eq:whittaker-h}
W_{h}(a(y)) = \dfrac{\abs{y}^{1/2} \hatF(y)}{1-p^{-1}}.
\end{equation}

Now more generally, we consider $\mu = (\mu_1,\mu_2)$ with $\mu_1 \mu_2 = \abs{-}$. We choose $M$ sufficiently large so that $\mu_1$ is trivial on $H := 1+p^{M}\ZZ_p$. Then the normalized flat Whittaker function corresponding to the flat section $h_{\flat} \in I(\mu_1,\mu_2)$ is
\[
W_{\flat}(a(y)) = \dfrac{\bfone_{h}(y)}{\vol(H, \dif^{\times} y)}.
\]
The corresponding $F_{\flat} := F_{h_{\flat}}$ then satisfies $Z(\hatF_{\flat},\mu_1, 0) = 1-p^{-1}$ by \eqref{eq:whittaker-h}.

\underline{(3) Identifying the toric integral with Tate's zeta integral}.
We fix $g_{0} := \begin{bmatrix}
    1 & 0 \\ 1 & 1
\end{bmatrix}$ and set $h^{\tor} := \pi(g_0) h$. Then the toric integral is computed as follows:
\begin{align*}
T_{\mu_1}(h^{\tor}) = \int_{F^{\times}} h(g_0 a(t)) \mu_1(t) \abs{t}^{-1/2} \dif^{\times} t &= \int_{F^{\times}} F_{h}(t^{-1}) \mu_1(t) \abs{t}^{-1} \dif^{\times} t \\&= \int_{F^{\times}} F_h(x) \mu_1(x)^{-1} \abs{x} \dif^{\times} x = Z(F_h,\mu_1^{-1}, 1),
\end{align*}
where the last equality uses the change of variables $x = t^{-1}$. Applying the functional equation \eqref{eq:tate-functional-equation-split} with $\theta = \mu_1^{-1}$ and $s=1$, we have
\[
1 - p^{-1} = Z(\hatF_{\flat}, \mu_1, 0) = \gamma(1, \mu_1^{-1}, \psi_{p}) Z(F_{\flat}, \mu_1^{-1}, 1),
\]
which gives
\[
T_{\mu_1}(h_{\flat}^{\tor}) = \dfrac{1-p^{-1}}{\gamma(1,\mu_1^{-1},\psi_p)}.
\]
Now we are ready to compute the term $\lambda_{\mu}(h_{\flat}^{\tor})$: following its definition, we have
\[
\lambda_{\mu_1}(h_{\flat}^{\tor}) = \dfrac{L(1,\bfone)}{L(0,\mu_1) L(0, \mu_2)} T_{\mu_1}(h_{\flat}^{\tor}) = \dfrac{1}{L(0, \mu_1) L(0, \mu_2) \gamma(1,\mu_1^{-1}, \psi_p)},
\]
by noting that $L(1,\bfone) = (1-p^{-1})^{-1}$. Note that the relation $\mu_2 = \abs{-} \mu_1^{-1}$ gives $L(0,\mu_2) = L(1,\mu_1^{-1})$. Using also the epsilon factor given by \eqref{eq:tate-functional-equation-split}, we obtain
\begin{equation} \label{eq:lambda-flat-split}
\lambda_{\mu_1}(h_{\flat}^{\tor}) = \dfrac{1}{\epsilon(1,\mu_1^{-1}, \psi_p) L(0,\mu_1)^{2}} = \dfrac{\epsilon(0, \mu_1, \psi_p)}{\mu_1(-1) L(0, \mu_1)^{2}}.
\end{equation}
Here for the last equality, we note that $\epsilon(0, \mu_1, \psi_p) \epsilon(1,\mu_1^{-1},\psi_p) = \mu_1(-1)$.

\underline{(4) Compare the two Eisenstein integrals with the cuspidal integral.}
We write the auxiliary local representation $\pi_p$ as a normalised principal series $\pi_p = I(\nu_1, \nu_2)$. For the interpolating character $\phi_k$ of weight $k$, the corresponding toric character is $\chi_{k}^{\tor} = \chi_0 \phi_{k} \abs{-}_{\calK}$, and we write it as $\chi_k = (\chi_{1,k}, \chi_{2,k})$ with $\chi_{1,k} \chi_{2,k} = \omega_{\pi_p}^{-1}$. The two characters entering the Eisenstein calculation here are
\[
\theta_{E,k} := \nu_1 \chi_{2,k} \abs{-}^{1/2}, \quad \theta_{0,k} := \nu_2 \chi_{2,k} \abs{-}^{1/2},
\]
and the corresponding Eisenstein toric characters are $\mu_{j,k} = (\theta_{j,k}, \abs{-} \theta_{j,k}^{-1})$ for $j \in \{E,0\}$.

One computes that
\begin{equation} \label{eq:A_p-split-compute-gamma}
\calA_{p}(k) = \dfrac{\calC_p}{L_1(k) L_2(k)} \chi_{1,k}(-1) \gamma(1/2, \pi_p \otimes \chi_{2,k}, \psi_p) = \chi_{1,k}(-1) \calC_{p} \dfrac{\epsilon_p(k)}{L_2(k)^{2}},
\end{equation}
where 
\[
\calC_{p} = \dfrac{L_p(1,\eta) L_p(1,\pi,\ad)}{\zeta_p(2)}, \quad L_{2}(k) = L(1/2, \pi_p \otimes \chi_{2,k}), \quad L_{1}(k) = L(1/2, \pi_p^{\vee} \otimes \chi_{2,k}^{-1}),
\]
with $\epsilon(k) := \epsilon(1/2, \pi_p \otimes \chi_{2,k}, \psi_{p})$. The last equality in \eqref{eq:A_p-split-compute-gamma} follows from replacing the gamma factor by the epsilon factor using \eqref{eq:tate-functional-equation-split}. On the Eisenstein side, \eqref{eq:lambda-flat-split} gives
\[
\lambda_{j,p}(k) = \theta_{j,k}(-1)^{-1} \dfrac{\epsilon(0, \theta_{j,k}, \psi_p)}{L(0, \theta_{j,k})^{2}}, \quad j \in \{E, 0\}.
\]
Combining with \eqref{eq:A_p-split-compute-gamma}, we obtain
\[
\calV_{p}(k) = \dfrac{\lambda_{E,p}(k) \lambda_{0,p}(k)}{\calA_{p}(k)} = \calC_{p}^{-1} (\theta_{E,k}(-1) \theta_{0,k}(-1) \chi_{1,k}(-1))^{-1}.
\]
On a fixed branch, all three character values at $-1$ above are fixed, and the continuation to $k=0$ then gives
\begin{equation} \label{eq:final-vp0-split}
\calV_{p}(0) = \calC_p^{-1} (\theta_{E,0}(-1) \theta_{0,0}(-1) \chi_{1,0}(-1))^{-1}.
\end{equation}

\underline{(5) Remaining terms at $k=0$ and the "automorphic expression" for $r_{p}^{\tor}$}. The same computation as above gives
\begin{equation} \label{eq:final-ap0-split}
\calA_p^0 = \chi_{1,0}(-1)\mathcal C_p
\frac{ \epsilon(0,\theta_{E,0},\psi_p) \epsilon(0,\theta_{0,0},\psi_p)} {L(0,\theta_{E,0})^2L(0,\theta_{0,0})^2}.
\end{equation}
We also have
\begin{equation} \label{eq:final-lampda0p-split}
\lambda_{0,p}(0) = \theta_{0,0}(-1)^{-1} \frac{\epsilon(0,\theta_{0,0},\psi_p)}{L(0,\theta_{0,0})^2}.
\end{equation}
Combining \eqref{eq:final-vp0-split}, \eqref{eq:final-ap0-split} and \eqref{eq:final-lampda0p-split}, we have
\begin{equation} \label{eq:automorphic-expression-rptor-split}
r_p^{\tor} = \dfrac{\calA_{p}^{0}}{\lambda_{0,p}(0)} \calV_{p}(0) = \theta_{E,0}(-1)^{-1} \dfrac{\epsilon(0,\theta_{E,0},\psi_p)}{L(0,\theta_{E,0})^{2}}.
\end{equation}
This is the "automorphic expression" for $r_{p}^{\tor}$. We remark that, in this expression for $r_{p}^{\tor}$, the auxiliary part given by $\psi_0$ cancels, and the expression depends only on the local behavior of the elliptic curve $E$ at $p$.

\begin{proof}[Proof of Proposition \ref{prop:r_p_tor-split}]
    It remains to evaluate all the terms in \eqref{eq:automorphic-expression-rptor-split}. When $n = n_{\sigma} = 0$, the character $\theta_{E,0}$ is unramified, and hence $\epsilon(0, \theta_{E,0}, \psi_p) = 1$. Moreover, $L(0,\theta_{E,0})^{-1} = 1-a_{\sigma}/p$. This proves that $r_p^{\tor} =_p (1-a_{\sigma}/p)^{2}$. When $n = n_{\sigma} > 0$, the character $\theta_{E,0}$ is ramified, and hence we have $L(0,\theta_{E,0}) = 1$. Moreover, $\epsilon(0, \theta_{E,0}, \psi_p) = (a_{\sigma}/p)^{n}G(\rho^{-1})$. This proves that $r_p^{\tor} =_p (a_{\sigma}/p)^{n_{\sigma}} G(\rho_{\sigma}^{-1})$, as desired.
\end{proof}

\subsubsection{Auxiliary elliptic curve}
We choose and fix an elliptic curve $E_0/\QQ$ with the same $j$-invariant as $E$ and good ordinary reduction at $p$. If $E$ already has good reduction, take $E_0=E$. There is then an isomorphism $\varphi:E\to E_0$ over a sufficiently large finite extension of $\QQ$.

For fixed N\'eron differentials $\omega_{0}$ of $E_0/\QQ$ and $\omega_{E}$ of $E/\QQ$, put $\varphi^*\omega_0=d_\sigma\omega_E$ after fixing an embedding $\sigma: \calK \hookrightarrow \CC$. Then we have $a_{E,\sigma}=a_{0,\sigma}/d_\sigma$. We define $\pi_{0,\sigma}=\hatr_\sigma a_{0,\sigma}$. Thus $\pi_\sigma=_p\pi_{0,\sigma}/d_\sigma$. Since $\varphi_*(\omega_E^\vee)=d_\sigma\omega_0^\vee$ and $\ttd_\rho=G(\rho)\ttd_0\otimes t^{-1}\tte_1\otimes\tte_{\rho^{-1}}$, the scalar $b_{\rho}$ appearing in \eqref{eq:regulator-ramified-ttd} can be written as $b_\rho=d_\sigma/G(\rho)$.

\begin{proposition} \label{prop:last-p-adic-defect-split}
    We have 
    \[
    \pi_{0,\sigma} =_p \dfrac{P_p(E) R_{\sigma} \pi_{\sigma}}{r_{p,\sigma}^{\tor}}, \text{ and therefore } \Delta_{p}(E/\calK) = v_p \left( \prod_{\sigma = \iota, \iota^{c}} U_{\sigma} \pi_{0,\sigma} \right).
    \]
\end{proposition}
\begin{proof}
    We first consider the case in which $E$ has good ordinary reduction at $p$. Then we see that
    \[
    P_p(E) = (1-a_{\sigma}/p)(1-a_{\sigma}^{-1}), \quad R_\sigma = (1-a_{\sigma}/p)/(1-a_{\sigma}^{-1}), \quad r_{p,\sigma}^{\tor} =_p (1-a_{\sigma}/p)^{2},
    \]
    where the statement on $r_{p,\sigma}^{\tor}$ follows from Proposition \ref{prop:r_p_tor-split}. Thus $P_p(E)R_{\sigma}/r_{p,\sigma}^{\tor}=_p1$. Since $\pi_{\sigma} = \pi_{0, \sigma}$, this gives the result. We then consider the case in which $E$ is additive at $p$. In this case, we have
    \[
    P_p(E) = 1, \quad R_{\sigma} = a_{\sigma}^{n_{\sigma}}d_{\sigma}G(\rho_{\sigma})^{-1},\quad r_{p}^{\tor} =_p (a_{\sigma}/p)^{n_{\sigma}} G(\rho_{\sigma}^{-1}),
    \]
    since $b_{\rho} = d_{\sigma}/G(\rho_{\sigma})$ as we recalled at the beginning of this section, and the statement on $r_{p,\sigma}^{\tor}$ follows from Proposition \ref{prop:r_p_tor-split}. So the quotient
    \[
    \dfrac{P_p(E) R_{\sigma} \pi_{\sigma}}{r_{p,\sigma}^{\tor}} =_p \pi_{\sigma} d_{\sigma} \dfrac{p^{n_{\sigma}}}{G(\rho_{\sigma}) G(\rho_{\sigma}^{-1})} =_p \pi_{0,\sigma} 
    \]
    since $G(\rho_{\sigma}) G(\rho_{\sigma}^{-1}) = \rho_{\sigma}(-1)p^{n_{\sigma}}$. The last statement then follows from Theorem \ref{thm:p-adic-defect-split}.
\end{proof}

It remains to compare the anticyclotomic $p$-adic $L$-function $\scL_\sigma$ with the elliptic-unit $p$-adic $L$-function $\scL_{\sigma}^{\elli}$. The comparison is given in the following proposition, which is parallel to Proposition \ref{prop:compare-two-L-functions} in the nonsplit case.

\begin{proposition} \label{prop:U-sigma-k}
    There exists a positive integer $N$ such that, for $k \geq 1$ with $N \mid k$, if we put $\phi_{k} = (\psi_{E}^{c}/\psi_{E})^{k}$ and $m = 2k-1$, then we have
    \begin{align*}
        \scL_{\sigma}^{\elli}(k) &=_p (k-1)! \lambda_{E,p}(k) \dfrac{(2\pi i)^{k-1}}{a_{0,\sigma}^{m}} L(\psi_E \phi_{k}, 1), \\ 
        \scL_{\sigma}(k) &= \dfrac{i^{k-1}(k-1)!}{(4\pi)^{k}} r_{\sigma}^{m} L(\psi_{E}\phi_k, 1) \lambda_{E,p}(k).
    \end{align*}
Consequently, for $w_{\sigma}(k) := \scL_{\sigma}(k)/\scL_{\sigma}^{\elli}(k)$ we have
\begin{equation} \label{eq:expression-U-sigma}
w_{\sigma}(k) \pi_{0,\sigma} =_p i^{2k-1} (\pi_{0,\sigma}^{2}/2)^{k}.
\end{equation}
\end{proposition}
\begin{proof}
    We start from the formula for $\scL_{\sigma}(k)$. Equation \eqref{eq:scL-phi-k} also holds in the case in which $p$ splits in $\calK$, with $d_p = 0$, and gives
    \[
    \scL_{\sigma}(k) = \dfrac{1}{8} r_{\sigma}^{2k-1} \calC_{E}(k) \Gamma_{\CC}(k) L(\psi_{E}\phi_{k}, 1) \lambda_{E,p}(k).
    \]
    Note that $\calC_{E}(k) = i^{k-1}2^{2-k}$ (see Lemma \ref{lem:archimedean-computation-comparison}) and $\Gamma_{\CC}(k) = 2(2\pi)^{-k}(k-1)!$ gives
    \[
    \scL_{\sigma}(k) = \dfrac{i^{k-1}(k-1)!}{(4\pi)^{k}} r_{\sigma}^{2k-1}  L(\psi_{E}\phi_{k}, 1) \lambda_{E,p}(k).
    \]
For the expression for $\scL_{\sigma}^{\elli}(k)$, we use Theorem \ref{thm:eu-article-ERC}. For $V_{k} = U^{-m}(-j) \rho^{-1}$, where $m = 2k-1$ and $j=k-1$, we take $\ttd_{0}^{-m} := (\ttd_{0}^{\vee})^{\otimes m}$ as the basis of $\bfD_{\cris}(U^{-m})$; then $\varphi(\ttd_{0}^{-m}) = a^{-m} \ttd_{0}^{-m}$. We define accordingly $a_{k} := a^{-m} = a^{-(2k-1)}$. One computes that the removed local Euler factors, the regulator, and the toric factor in the renormalization factor of $\scL_{\sigma}$ are as follows:
\begin{equation}\label{eq:split-critical-local-table}
\begin{array}{c|c|c}
 &n=0&n>0\\ \hline
 \calE_{E,p}(k) &(1-p^{k-1}a_k)(1-p^{-k}a_k^{-1})&1\\ \hline
 R_{\sigma}(k) & \dfrac{1-p^{k-1}a_k}{1-p^{-k}a_k^{-1}}
 &\dfrac{p^{nk}a_k^n}{G(\rho)}\\ \hline
 \lambda_{E,p}(k) &(1-p^{k-1}a_k)^2 &p^{n(k-1)}a_k^nG(\rho^{-1}).
\end{array}
\end{equation}
Indeed,
\begin{itemize}
    \item The first row is the Euler polynomial at $p$ removed from the complex $L$-value; in the ramified case the local $L$-factor at $p$ is $1$.
    \item The second row is the regulator interpolation operator from Theorem~\ref{thm:LZV-regulator} at $j=k-1$.
    \item The third row follows from the same calculation as Proposition~\ref{prop:r_p_tor-split} with $\mu_1(p)=p^{k-1}a_k$.
\end{itemize}
It is straightforward to check that, in the unramified case, the product of the first two rows gives $(1-p^{k-1}a_k)^2$. In the ramified case, their product divided by the third row is
\[
 \frac{p^{nk}a_k^n/G(\rho)}
 {p^{n(k-1)}a_k^nG(\rho^{-1})}
 =\frac{p^n}{G(\rho)G(\rho^{-1})}=\rho(-1)^{-1},
\]
where the last equality uses the identity for Gauss sums $G(\rho)G(\rho^{-1}) = \rho(-1)^{-1}p^{n}$.
Also, $d_\calK$ appearing in Theorem \ref{thm:eu-article-ERC} is a $p$-adic unit because $p$ splits. This proves the expression for $\scL_{\sigma}^{\elli}(k)$.

To obtain the expression for $w_{\sigma}(k)$, using $r_\sigma a_{0,\sigma}=2\pi i\pi_{0,\sigma}$ gives
\[
 w_\sigma(k)=i^{2k-1}2^{-k} \pi_{0,\sigma}^{2k-1},
\]
which is \eqref{eq:expression-U-sigma}.
\end{proof}

We can now finish the split case by applying a continuity argument similar to that in Theorem \ref{thm:main_transfer}.

\begin{theorem} \label{thm:Nov-split-BSD}
    Let $E/\QQ$ be a CM elliptic curve with CM by $\calO_{\calK}$ and $\ord_{s=1}L(E/\QQ,s) = 1$. Suppose $p$ splits in $\calK$. Then we have $\Delta_{p}(E/\calK) = 0$. Equivalently, the $p$-part of BSD holds for $E/\calK$ and $E/\QQ$.
\end{theorem}
\begin{proof}
Note that $w_\sigma(k)=\scL_\sigma(k)/\scL_\sigma^{\elli}(k)$ is analytic and nonzero on a sufficiently small disc about zero. Proposition~\ref{prop:U-sigma-k} then gives
\[
 v_p(w_\sigma(k)\pi_{0,\sigma})  =k v_p(\pi_{0,\sigma}^2/2).
\]
For $k_a=Np^a$, the left side is eventually $v_p(w_\sigma(0)\pi_{0,\sigma})$ by continuity and the ultrametric inequality. Since $k_a$ is unbounded as a positive integer, the right side can be constant only if its slope is zero. It follows that
\[
 v_p(\pi_{0,\sigma}^2/2)=0,\qquad
 v_p(w_\sigma(0)\pi_{0,\sigma})=0.
\]
Proposition~\ref{prop:last-p-adic-defect-split} now gives
\[
 \Delta_p(E/\calK)
 =\sum_{\sigma=\iota,\iota^c} v_p\bigl(w_\sigma(0)\pi_{0,\sigma}\bigr)=0.
\]
Finally Lemma~\ref{lem:base-change} gives $\Delta_p(E/\calK)=2\Delta_p(E/\QQ)$. This proves the last assertion as well.
\end{proof}

%% file: 201-Setup.tex
\section{Setup for Part \ref{part:BDP}}
Throughout this part, we assume that $p \geq 5$. Let $E/\QQ$ be an elliptic curve with complex multiplication given by $\calO_{K}$, where $\calK$ is an imaginary quadratic field in which $p$ \emph{remains inert}. We write $N$ for the conductor of $E$.

\subsection{The biquadratic setting}
We fix another imaginary quadratic field $\calL$ such that
\begin{equation}\tag{spl-Heeg}
    \text{$p \calO_{\calL} = v_{0} \barv_{0}$, \quad every $\ell$ dividing $N$ splits in $\calL$, \quad $D_{\calL} < -4$ is odd.}
\end{equation}
We put $\calM = \calK\calL$ and let $\calF$ be its unique real quadratic subfield. It follows that $p \calO_{\calM} = v \barv$, and we set $v \mid v_{0}$. The prime of $\calF$ above $p$ is inert over $\QQ$, and
\[
\calL_{v_0} = \calL_{\barv_0} = \QQ_{p}, \quad \calM_{v} = \calM_{\barv} = \calF_{p} = \calK_{p} \simeq \QQ_{p^{2}}.
\]
Here $\QQ_{p^{2}}$ denotes the unramified quadratic extension of $\QQ_{p}$.

We choose a $p$-ordinary CM type of $\calL$ as $\Sigma_{0} := \{\sigma_{0}: \calL \hookrightarrow \CC\}$ which induces the $p$-adic place $v_{0}$ of $\calL$, and we choose the $p$-ordinary CM type of $\calM$ as 
\[
\Sigma := \{\sigma_1, \sigma_2\} := \{ \sigma: \calM \hookrightarrow \CC: \sigma|_{\calL} = \sigma_{0} \}.
\]
Both embeddings in $\Sigma$ induce the $p$-adic place $v$ of $\calM$, and their conjugates induce $\barv$. Thus, $\Sigma$ is indeed a $p$-ordinary CM type of $\calM/\calF$.

We keep the following assumption on $\calL$:
\begin{equation} \tag{rank-one} \label{eq:rank-one}
 \rank_{\ZZ} E(\calL) = 1, \quad  \# \Sha(E/\calL)[p^{\infty}] < \infty.    
\end{equation}
By the Gross--Zagier--Kolyvagin theorem, it is implied by the following more practical assumption:
\begin{equation} \tag{GZK} \label{eq:GZK}
    \ord_{s=1} L(E/\QQ, s) = 1, \quad L(E^{\calL}/\QQ, 1) \neq 0,
\end{equation}
where $E^{\calL}$ is the quadratic twist of $E$ attached to $\calL/\QQ$.

For any abelian number field $H$, we let $H_{\infty}$ be the compositum of all its $\ZZ_p$-extensions. Leopoldt's conjecture is valid for abelian number fields and states that $\Gamma_{H} := \Gal(H_{\infty}/H)$ is isomorphic to $\ZZ_p^{r_2(H)+1}$, where $r_{2}(H)$ is the number of pairs of complex places of $H$. Moreover, when $p \neq 2$, for an abelian CM number field $H/H^{+}$, the complex conjugation $c \in \Gal(H/H^{+})$ acts on $\Gamma_{H}$ by conjugation, and $\Gamma_{H}$ splits into eigenspaces $\Gamma_{H}^{\pm}$ on which the complex conjugation acts by $\pm 1$, with $\Gamma_{H}^{+} \simeq \ZZ_{p}^{[H^{+}:\QQ]}$ and $\Gamma_{H}^{-} \simeq \ZZ_{p}$. This gives
\begin{itemize}
    \item the \emph{cyclotomic $\ZZ_p$-extension} $H_{\infty}^{\cyc} := H_{\infty}^{\Gamma_{H}^{-}}$ over $H$, with $\Gamma_{H}^{\cyc} := \Gal(H_{\infty}^{\cyc}/H) \simeq \ZZ_{p}$, and
    \item the \emph{anticyclotomic $\ZZ_p$-extension} $H_{\infty}^{\ac} := H_{\infty}^{\Gamma_{H}^{+}}$ over $H$, with $\Gamma_{H}^{\ac} := \Gal(H_{\infty}^{\ac}/H) \simeq \ZZ_{p}^{[H^{+}:\QQ]}$.
\end{itemize}

Let $L/\QQ_p$ be a sufficiently large finite extension, and let $\calO$ be its ring of integers. Returning to our biquadratic setting, we define the following:
\begin{itemize}
    \item The $\ZZ_p^{3}$-extension $\calM_{\infty}/\calM$, with the Iwasawa algebra $\Lambda_{\calM} = \calO\lrbracket{\Gal(\calM_{\infty}/\calM)}$;
    \item The $\ZZ_p^{2}$-extension $\calM_{\infty}^{\calL} := \calM \calL_{\infty}/\calM$, with the Iwasawa algebra $\Lambda_{\calM}^{\calL} := \calO\lrbracket{\Gal(\calM_{\infty}^{\calL}/\calM)}$, isomorphic to $\Lambda_{\calL}$ via the restriction map $\Gal(\calM_{\infty}^{\calL}/\calM) \simeq \Gal(\calL_{\infty}/\calL)$;
    \item The $\ZZ_p$-extension $\calM_{\infty}^{\ac,\calL} := \calM \calL_{\infty}^{\ac}/\calM$, with the Iwasawa algebra $\Lambda_{\calM}^{\ac,\calL} = \calO\lrbracket{\Gal(\calM_{\infty}^{\ac, \calL}/\calM)}$, isomorphic to $\Lambda_{\calL}^{\ac}$ similarly.
\end{itemize}
Then restriction maps induce the quotient maps
\[
\Lambda_{M} \xrightarrow{\bmod (\gamma^{-,1}-1)} \Lambda_{M}^{\calL} \xrightarrow{\bmod (\gamma^{+} - 1)} \Lambda_{\calM}^{\ac,\calL}.
\]
where $\gamma^{-,1}$ and $\gamma^{+}$ are fixed topological generators of $\Gal(\calM_{\infty}/\calM_{\infty}^{\calL})$ and $\Gal(\calM_{\infty}^{\calL}/\calM_{\infty}^{\ac,\calL})$ (which are both isomorphic to $\ZZ_p$), respectively \footnote{Actually one may check that $\calM_{\infty} = \calM_{\infty}^{\calL} \calM_{\infty}^{\ac, \calK}$. Nevertheless, we are not using this relation in this article.}.

Let $\calO^{\ur}$ be the completed ring of integers of the maximal unramified extension of $L$. A superscript "$\ur$" means the completed tensor-product base change to $\calO^{\ur}$-coefficients.

\subsection{Characters and representations}
We fix the conventions in Part \ref{part:BDP}. 

\subsubsection{Hecke characters and their $p$-adic avatars}

For any algebraic Hecke character $\chi: \calH^{\times} \bs \AA_{\calH}^{\times} \rightarrow \CC^{\times}$ of a CM field $\calH/\calH^{+}$, we say that $\chi$ has infinite type $(a_\sigma, b_{\sigma})_{\sigma \in \Sigma}$ for a CM type $\Sigma$ of $\calH/\calH^{+}$ when the archimedean component of $\chi$ is
\[
\chi_{\infty}(z) = \prod_{\sigma \in \Sigma} \sigma(z)^{a_{\sigma}} \barsigma(z)^{b_{\sigma}}.
\]
Note that the associated ideal class character attached to $\chi$ satisfies the formula $\chi((\alpha)) = \alpha^{-a}\baralpha^{-b}$. 

We fix embeddings $\iota_{\infty}: \barQQ \hookrightarrow \CC$ and $\iota_p: \barQQ \hookrightarrow \CC_p$. For any embedding $\tau: \calH \hookrightarrow \barQQ$, we write $\tau_{\infty} = \iota_{\infty} \circ \tau$ and $\tau_{p} = \iota_{p} \circ \tau$ and let $\frP_{\tau}$ be the prime determined by $\tau_{p}$. For any Hecke character over $\calH$ of infinite type $(a_{\sigma}, b_{\sigma})_{\sigma \in \Sigma}$, we have
\[
\lambda_{\rmf}(\alpha) = \prod_{\sigma \in \Sigma} \sigma_{\infty}(\alpha)^{-a_{\sigma}} (c\sigma)(\alpha)^{-b_{\sigma}},
\]
since $\chi$ is trivial on $\calH^{\times}$. For a finite idèle $x = (x_w)_{w}$, we define
\[
\chi^{\dagger}(x_{p}) := \prod_{\sigma \in \Sigma} \sigma_{p}(x_{\frP_{\sigma}})^{a_{\sigma}} (c\sigma)_{p} (x_{\frP_{c\sigma}})^{b_{\sigma}}. 
\]
Throughout this part, we use the geometrically normalized reciprocity law, 
\[
\rec_{\calH}: \AA_{\calH}^{\times} \rightarrow \Gal_{\calH}^{\ab},
\]
    i.e., a local uniformizer is mapped to the geometric Frobenius. The \emph{$p$-adic avatar} of $\chi$ is the character $\hchi: \Gal_{\calH}^{\ab} \rightarrow \calO^{\times}$ defined by
\[
\hchi(\rec_{\calH}(x)) = \iota_{p} \iota_{\infty}^{-1}(\chi_{\rmf}(x)) \chi^{\dagger}(x_p).
\]

\subsubsection{Galois representations}
Let $\psi: \AA_{\calK}^{\times}/\calK^{\times} \rightarrow \CC^{\times}$ be the Hecke character over $\calK$ attached to $E$ in the sense that 
\[
\psi_{\infty}(z) = z^{-1}, \quad L(\psi, s) = L(E/\QQ, s),
\]
i.e. $\psi$ has infinite type $(-1,0)$.

If $\bfN$ is the norm character over $H$ with $\bfN(\frq) = \#(\calO_{H}/\frq)$ for any ideal $\frq$ of the ring of integers $\calO_{H}$ of $H$, then $\bfN$ has infinite type $(-1,-1)$ because, for any principal ideal $(\alpha)$, $\bfN((\alpha)) = \alpha \baralpha$. Moreover, $\widehat{\bfN^{-1}} = \epsilon_{\cyc}$ is the cyclotomic character.

We put $T_{E} := T_{p}(E) \otimes_{\ZZ_{p}} \calO$. Let $T_{\rho}$ be a one-dimensional $\Gal_{\calK}$-representation over $\calO$ with character $\rho := \hpsi^{-1}$. Then one checks that $T_{E} := \Ind_{\Gal_{\calK}}^{\Gal_{\QQ}} T_{\rho}$. More precisely, we enlarge $L$ so that $L$ contains $\calK_p$; then $\calO_{\calK_p} \otimes_{\ZZ_p} \calO \simeq \calO \times \calO$ gives the direct-sum decomposition
\[
T_{E} = T_{p}(E) \otimes_{\ZZ_p} \calO = e_{\iota} T_{E} \oplus e_{\iota^{c}} T_{E},
\]
where $e_{\iota}$ and $e_{\iota^{c}}$ are two idempotents, each giving a free rank-one $\calO$-summand of $T_{E}$. We denote $T_{\rho} = e_{\iota} T_{E}$, and the other summand is then $T_{\rho^{c}}$, carrying the action of the conjugate character $\rho^{c}$.

Throughout this article, we use the symbol "$\square$" to mean "no superscript" or "no subscript" at the corresponding places. For $\bullet \in \{\square; \calL; \ac,\calL\}$, we define
\[
(\Lambda_{\calM}^{\bullet})^{\vee} := \Hom_{\calO, \cts}(\Lambda_{\calM}^{\bullet}, L/\calO)
\]
to be the Pontryagin dual, where $L/\calO$ has the discrete topology, and this dual has the $\Lambda_{\calM}^{\bullet}$-action
\[
(a \cdot \phi)(b) = \phi(ab), \quad a,b\in \Lambda_{\calM}^{\bullet}.
\]
Let $q_{\bullet}: G_{\calM} \twoheadrightarrow \Gamma_{\calM}^{\bullet}$ be the natural quptient map, we define the \emph{(discrete) Iwasawa deformation}
\[
\bfT_{\rho}^{\bullet} := T_{\rho} \otimes_{\calO} (\Lambda_{\calM}^{\bullet})^{\vee},
\]
with the Iwasawa algebra acts on the second factor, and the $G_{\calM}$-action is given by
\[
g (t \otimes \phi) = \rho(g)t \otimes \bgamma_{\Lambda_{\calM}^{\bullet}}(q_{\bullet}(g))^{-1} \cdot \phi, \quad \text{for $g \in G_{\calM}$}.
\]
Likewise, for $\circ \in \{\square, \ac\}$, let $q_{\circ}: G_{\calL} \twoheadrightarrow \Gamma_{\calL}^{\circ}$, and define
\[
\bfT_{E}^{\circ} := T_{E} \otimes_{\calO} (\Lambda_{\calL}^{\circ})^{\vee},
\]
with the $G_{\calL}$-action given by
\[
g(t \otimes \phi) = (g \cdot t) \otimes \bgamma_{\Lambda_{\calL}^{\circ}}(q_{\circ}(g))^{-1} \cdot \phi, \quad\text{for $g \in G_{\calL}$},
\]
where $g \cdot t$ denotes the usual Galois action on $T_{E}$.

%% file: 202-Selmer.tex
\section{Arithmetic side: Selmer groups and control theorems} \label{sec:Selmer-group}
Let $S_{\calL}$ be a finite set of places of $\calL$ containing places over $pN$ and those ramified in $\calM/\calL$, and let $S_{\calM}$ be the places of $\calM$ lying above $S_{\calL}$. We define the following two types of Selmer groups.

\begin{definition}
Put $\bullet \in \{\square; \calL; \ac,\calL\}$ and $\circ \in \{\square, \ac\}$. We define the \emph{(primitive) Greenberg Selmer group} over $\calM$ as
    \[
    \Sel_{\Gr}(\bfT_{\rho}^{\bullet}/\calM) := \ker \left( \rmH^{1}(\Gal_{\calM, S_{\calM}}, \bfT_{\rho}^{\bullet}) \rightarrow \rmH^{1}(\calM_{v}, \bfT_{\rho}^{\bullet}) \times \prod_{w \in S_{\calM}, \, w \nmid p} \rmH^{1}(I_{w}, \bfT_{\rho}^{\bullet}) \right),
    \]
and the \emph{(primitive) Bertolini--Darmon--Prasanna (BDP) Selmer group} over $\calL$ as 
\[
    \Sel_{\BDP}(\bfT_{E}^{\circ}/\calL) := \ker \left( \rmH^{1}(\Gal_{\calL, S_{\calL}}, \bfT_{E}^{\circ}) \rightarrow \rmH^{1}(\calL_{v_{0}}, \bfT_{E}^{\circ}) \times \prod_{w \in S_{\calL}, \, w \nmid p} \rmH^{1}(I_{w}, \bfT_{E}^{\circ}) \right).
\]
We put $X_{\Gr}(\bfT_{\rho}^{\bullet}/\calM) := \Sel_{\Gr}(\bfT_{\rho}^{\bullet}/\calM)^{\vee}$ and $X_{\BDP}(\bfT_{E}^{\circ}/\calL) := \Sel_{\BDP}(\bfT_{E}^{\circ}/\calL)^{\vee}$.
\end{definition}

\subsection{Local invariant lemma}

The following local computations are crucial in proving the descent theorem (Proposition \ref{prop:descent-theorem}) and the control theorem for BDP Selmer groups (Proposition \ref{prop:control-theorem-BDP}).

\begin{lemma} \label{lem:invariant-vanishing}
Put $\bullet \in \{\square; \calL; \ac,\calL\}$ and $\circ \in \{\square, \ac\}$. The following assertions hold.
\begin{enumerate}[label = \rm (\arabic*)]
    \item $(\bfT_{\rho}^{\bullet})^{I_{v}} = (\bfT_{\rho}^{\bullet})^{I_{\barv}} = 0$, and $\rmH^{1}(\Gal_{\calM, S_{\calM}}, \bfT_{\rho}^{\bullet}) = 0$.
    \item If $w \nmid p$, then $(\bfT_{\rho}^{\bullet})^{I_{w}}$ is either zero or $\bfT_{\rho}^{\bullet}$. In particular, multiplication by $[\gamma^{\bullet^{\prime}}-1]$ on $(\bfT_{\rho}^{\bullet})^{I_{w}}$ is surjective for $\bullet^{\prime} \in \{ -,1; +\}$ accordingly.
    \item $(\bfT_{E}^{\circ})^{I_{v_{0}}} = (\bfT_{E}^{\bullet})^{I_{\barv_0}} = 0$, and $\rmH^{1}(\Gal_{\calL, S_{\calL}}, \bfT_{E}^{\circ}) = 0$.
    \item If $w \nmid p$, then multiplication by $[\gamma_{\ac}-1]$ on $(\bfT_{E}^{\ac})^{I_{w}}$ is surjective.
\end{enumerate}
\end{lemma}

\begin{proof}
    We first prove (1). We start by considering the action of the inertia group $I_{K_{p}}$ on $E[p]$, since $I_{K_p} = I_{v} = I_{\barv}$. Since $p$ remains inert in $\calK$, $E[p]$ is a one-dimensional vector space over $\calO_{\calK}/p = \FF_{p^{2}}$. The action of $I_{K_{p}}$ is therefore given by a character $\barrho_0: I_{K_p} \rightarrow \FF_{p^{2}}^{\times}$. If $\barrho_0$ were trivial, then $E[p]$ would be fixed by the inertia subgroup $I_{K_p}$. The nondegenerate Weil pairing $E[p] \times E[p] \rightarrow \mu_p$ would then imply that $K_{p}(\mu_p)/K_{p}$ is an unramified extension. This is impossible since $K_p/\QQ_p$ is unramified and $\QQ(\mu_p)/\QQ_p$ is totally ramified of degree $p-1$. Therefore, there exists $g \in I_{K_p}$ such that $a := \barrho_{0}(g) \in \FF_{p^{2}}^{\times}$ is not equal to $1$ and also $a^{p} \neq 1$.

    Note that $T_{E}/\varpi T_{E} \simeq E[p] \otimes_{\FF_p} \FF$, where $\FF = \calO/\varpi$. Since $\FF_{p^{2}} \otimes_{\FF_p} \FF \simeq \FF \times \FF$, the inertia actions on the two copies are given by either $\barrho_0(g)$ or $\barrho_0(g)^{p}$. Thus, $\rho(g)$ equals either $\barrho_0(g)$ or $\barrho_{0}(g)^{p}$, which is not equal to $1$. Consequently, $\rho(g) -1 \in \calO^{\times}$. Moreover, since $\bgamma_{\Lambda_{\calM}^{\bullet}}(g) \equiv 1 \pmod{\frm_{\Lambda_{\calM}^{\bullet}}}$, we see that $\rho(g)\bgamma_{\Lambda_{\calM}^{\bullet}}(g)^{-1} - 1 \in (\Lambda_{\calM}^{\bullet})^{\times}$. This proves the first part of (1). Since elements of $\bfT_{\rho}^{\bullet}$ invariant under $\Gal_{\calM, S_{\calM}}$ are also invariant under $I_{K_p}$, the second part of (1) follows.

    For (2), first note that when $w \nmid p$, the tautological character $\bgamma_{\Lambda_{\calM}^{\bullet}}$ is trivial on $I_{w}$. Since $E$ has CM, it has potentially good reduction at $w$. As $w \nmid p$, the Néron--Ogg--Shafarevich theorem implies that $I_{w}$ has finite image on $T_p(E)$. All CM endomorphisms are defined over $\calK \subseteq \calM$, and $T_p(E)$ is free of rank one over $\calO_{\calK_p}$. Hence this finite image lies in $\calO_{\calK_p}^{\times}$. Since $p$ is inert and $p \geq 5$, reduction modulo $p$ is injective on finite subgroups of $\calO_{\calK_p}^{\times}$: its kernel $1+p \calO_{\calK_p}$ is torsion free. The inertia image order therefore divides $\# \FF_{p^{2}}^{\times} = p^{2}-1$. In particular, the character $\rho|_{I_{w}}$ has finite order prime to $p$. If it is the trivial character, then $(\bfT_{\rho}^{\bullet})^{I_{w}}$ is $\bfT_{\rho}^{\bullet}$ itself. Otherwise, we can choose $g \in I_{w}$ such that $\zeta := \rho(g) \neq 1$ is a root of unity of prime-to-$p$ order, and so $\zeta - 1 \in \calO^{\times}$. Thus, the operator $g-1$ acts on $\bfT_{\rho}^{\bullet}$ by multiplication by this unit. Together with $\bgamma_{\Lambda_{\calM}^{\bullet}}(g) = 1$, this shows that $(\bfT_{\rho}^{\bullet})^{I_{w}} = 0$, as desired.

    The statement in (3) is proved as in the proof of (1), and the statement in (4) is proved as in the proof of (2) above: simply note that multiplication by $(\gamma_{\ac} - 1)$ is surjective on $(\Lambda_{\calL}^{\ac})^{\vee}$ and hence $(\gamma_{\ac}-1) \bfT_{E, \calL}^{\ac} = \bfT_{E, \calL}^{\ac}$.
\end{proof}

\subsection{Descent theorems}

\begin{proposition}[Descent theorem] \label{prop:descent-theorem}
    There are canonical isomorphisms
    \[
    \Sel_{\Gr}(\bfT_{\rho}^{\calL}/\calM) \simeq \Sel_{\Gr}(\bfT_{\rho}/\calM)[\gamma^{-,1}-1], \quad \Sel_{\Gr}(\bfT_{\rho}^{\ac, \calL}/\calM) \simeq \Sel_{\Gr}(\bfT_{\rho}^{\calL}/\calM)[\gamma^{+}-1].
    \]
Consequently, $X_{\Gr}(\bfT_{\rho}^{\ac, \calL}/\calM) \simeq X_{\Gr}(\bfT_{\rho}/\calM)/(\gamma^{-,1}-1, \gamma^{+}-1) X_{\Gr}(\bfT_{\rho}/\calM)$.
\end{proposition}

\begin{proof}
   We prove the first canonical isomorphism; the second follows in exactly the same way. Let $G$ be the profinite group equal to one of $\Gal_{\calM, S_{\calM}}, \calM_{v}$, or $I_{w}$ for $w \in S_{\calM}$ with $w \nmid p$. There is a short exact sequence of $G$-modules
   \[
   0 \rightarrow \bfT_{\rho}^{\calL} \rightarrow \bfT_{\rho} \xrightarrow{\times (\gamma^{-,1}-1)} \bfT_{\rho} \rightarrow 0. 
   \]
   Taking the $G$-cohomology gives the following short exact sequence 
   \[
   0 \rightarrow \bfT_{\rho}^{G}/(\gamma^{-,1}-1)\bfT_{\rho}^{G} \rightarrow \rmH^{1}(G, \bfT_{\rho}^{\calL}) \rightarrow \rmH^{1}(G, \bfT_{\rho})[\gamma^{-,1}-1].
   \]
   By Lemma \ref{lem:invariant-vanishing}, we see that the leftmost term above vanishes for $G$ equal to any of $\Gal_{\calM, S_{\calM}}, \calM_{v}$, or $I_{w}$ for $w \in S_{\calM}$ with $w \nmid p$. This gives the following commutative diagram with exact rows:
   \[
   \begin{tikzcd}
       0 & \Sel_{\Gr}(\bfT_{\rho}^{\calL}/\calM) & \rmH^{1}(\Gal_{\calM,S_{\calM}}, \bfT_{\rho}^{\calL}) & \rmH^{1}(\calM_{v}, \bfT_{\rho}^{\calL}) \times \prod_{w \in S_{\calM}, \, w \nmid p} \rmH^{1}(I_{w}, \bfT_{\rho}^{\calL}) \\ 
        0 & \ker(\res) & \rmH^{1}(\Gal_{\calM,S_{\calM}}, \bfT_{\rho})[z] & \rmH^{1}(\calM_{v}, \bfT_{\rho})[z] \times \prod_{w \in S_{\calM}, \, w \nmid p} \rmH^{1}(I_{w}, \bfT_{\rho})[z]
        \arrow[from = 1-1, to=1-2] \arrow[from = 1-2, to=1-3] \arrow[from = 1-3, to=1-4]
        \arrow[from = 2-1, to=2-2] \arrow[from = 2-2, to=2-3] \arrow[from = 2-3, to=2-4, "\res"]
        \arrow[from = 1-2, to = 2-2] \arrow[from = 1-3, to = 2-3, "\simeq"] \arrow[from = 1-4, to = 2-4, "\simeq"] 
   \end{tikzcd}
   \]
   where $z := \gamma^{-,1} - 1$ is a shorthand for space-saving. This gives an isomorphism
   \begin{equation}\label{eq:descent1}
   \Sel_{\Gr}(\bfT_{\rho}^{\calL}/\calM) \simeq \ker(\res)
   \end{equation}
   One verifies directly that an element of $\rmH^{1}(\calM_{v}, \bfT_{\rho})[z]$ has zero image via the restriction map $\res$ precisely when its image in $\rmH^{1}(\calM_{v}, \bfT_{\rho}) \times \prod_{w \in S_{\calM}, \, w \nmid p} \rmH^{1}(I_{w}, \bfT_{\rho})$ is zero, and hence 
   \begin{equation} \label{eq:descent2}
      \ker(\res) = \Sel_{\Gr}(\bfT_{\rho}/\calM)[\gamma^{-,1}-1].    
   \end{equation}
   Combining \eqref{eq:descent1} and \eqref{eq:descent2} gives the first canonical isomorphism. The second canonical isomorphism can be shown in the same way, by invoking Lemma \ref{lem:invariant-vanishing}. The last consequence follows easily. 
\end{proof}

\subsection{Shapiro's lemma for BDP Selmer groups}
In this section, we prove the following proposition, which is a version of Shapiro's lemma for BDP Selmer groups.

\begin{proposition} \label{prop:shapiro}
Put $\bullet \in \{\calL; \ac,\calL\}$ and $\circ \in \{\square, \ac\}$ accordingly. There is a canonical $\Lambda_{L}^{\ac}$-linear isomorphism
    \[
    \Sel_{\BDP}(\bfT_{E}^{\circ}/\calL) \simeq \Sel_{\Gr}(\bfT_{\rho}^{\bullet}/\calM)
    \]
    where we identify $\Lambda_{\calL}^{\circ} \simeq \Lambda_{\calM}^{\bullet}$.
\end{proposition}
In this article, we only use the case $\bullet = \ac, \calL$ and $\circ = \ac$. Nevertheless, proving the other pair is equally easy.

\begin{proof}
    This is essentially Shapiro's lemma since $T_{E} := \Ind_{\Gal_{\calK}}^{\Gal_{\QQ}} T_{\rho}$. Projection lemma gives $\bfT_{E}^{\circ} := \Ind_{\Gal_{\calK}}^{\Gal_{\QQ}} \bfT_{\rho}^{\bullet}$. Shapiro's lemma is applied three times.

    (1) Since $S_{\calL}$ contains the places of $\calL$ which are ramified in $\calM$, Shapiro's lemma applies to the global cohomology groups, i.e. there is a canonical isomorphism
    \[
    \rmH^{1}(\Gal_{\calM, S_{\calM}}, \bfT_{\rho}^{\bullet}) \xrightarrow{\sim} \rmH^{1}(\Gal_{\calL, S_{\calL}}, \bfT_{E}^{\circ}).
    \]

    (2) Second, we consider the local cohomology groups at $p$; in fact, only the comparison at the places $v$ is needed. For the whole decomposition group,
    \[
    \rmH^{1}(\calM_{v}, \bfT_{\rho}^{\bullet}) \xrightarrow{\sim} \rmH^{1}(\calL_{v_0}, \bfT_{\rho}^{\circ}),
    \]
    since $v$ is the unique place of $\calM$ lying over $v_{0}$.

    (3) Third, we consider the local cohomology groups away from $p$. In this case, we have the Mackey decomposition
    \[
    \Res_{I_{w}}^{G_{w}} \bfT_{E}^{\circ} = \Res_{I_{w}}^{G_{w}} \Ind_{G_{u}}^{G_{w}} \bfT_{\rho}^{\bullet} = \bigoplus_{r \in I_{w} \bs G_{w} / G_{u}} \Ind_{I_{w} \cap r G_{u} r^{-1}}^{I_{w}} {}^{r} \bfT_{\rho}^{\bullet}
    \]
    where $u$ is a place of $\calM$ lying over $w$, and here $I_{w} \cap r G_{u} r^{-1} = r I_{u} r^{-1}$. The $r=1$ summand is the restriction to $I_{u}$, while the other components are its conjugates. It follows that the restriction to $I_{w}$ vanishes if and only if the restriction to every $I_{u}$ vanishes.

    Combining (1),(2) and (3), the proposition is therefore proved.
\end{proof}

\subsection{BDP control theorems}

\subsubsection{Bloch--Kato Selmer groups and their invariants}
We take coefficients in $\ZZ_p$. We write $T := T_{p}(E)$, $V := T \otimes_{\ZZ_p} \QQ_p$, and $W := V/T := E[p^{\infty}]$. Over the base field $\calL$ and at any finite place $w$ of $\calL$, the Kummer short exact sequence gives natural inclusions
\[
E(\calL_{w})/p^{n} E(\calL_{w}) \hookrightarrow \rmH^{1}(\calL_{w}, E[p^{n}]).
\]
Taking inverse limits defines
\[
\rmH^{1}_{\rmf}(\calL_{w}, T) := \im\left( E(\calL_{w})_p^{\wedge} \rightarrow \rmH^{1}(\calL_{w}, T) \right), \quad E(\calL_{w})_{p}^{\wedge} := \varprojlim_{n} E(\calL_{w})/p^{n} E(\calL_{w}),
\]
and taking direct limit defines
\[
\rmH^{1}_{\rmf}(\calL_w, W) := \im \left( E(\calL_{w}) \otimes_{\ZZ} \QQ_p/\ZZ_p  \rightarrow \rmH^{1}(\calL_{w}, W) \right).
\]
For $M \in \{T, W\}$, we define the \emph{Bloch--Kato Selmer group} as
\[
\Sel_{\BK}(\calL, M) := \ker\left( \rmH^{1}(\Gal_{\calL, S_{\calL}}, M) \rightarrow \prod_{w \in S_{\calL}} \dfrac{\rmH^{1}(\calL_{w}, M)}{\rmH^{1}_{\rmf}(\calL_{w}, M)} \right)
\]

We also define the following two Selmer groups:
\[
\Sel_{\rel}(\calL, M) := \ker\left( \rmH^{1}(\Gal_{\calL, S_{\calL}}, M) \rightarrow \prod_{w \in S_{\calL}, w \neq \barv_{0}} \dfrac{\rmH^{1}(\calL_{w}, M)}{\rmH^{1}_{\rmf}(\calL_{w}, M)} \right),
\]
this imposes the relaxed Selmer condition at the place $\barv_{0}$ of $\calL$, so we call it the \emph{relaxed Selmer group}; we also define
\[
\Sel_{\BDP}(\calL, M) := \ker\left( \rmH^{1}(\Gal_{\calL, S_{\calL}}, M) \rightarrow  \rmH^{1}(\calL_{v_0}, M)  \times \prod_{w \in S_{\calL}, w \neq p} \dfrac{\rmH^{1}(\calL_{w}, M)}{\rmH^{1}_{\rmf}(\calL_{w}, M)} \right),
\]
This is a finite-level analogue of the BDP Selmer group defined previously.

\subsubsection{Length formula for BDP Selmer groups}
We gather local cohomological computations in the following lemma. These facts should be well-known.
\begin{lemma} \label{lem:gather-finite}
    For $w = v_{0}$ or $\barv_{0}$, we have 
    \begin{enumerate}[label = \rm (\arabic*)]
        \item $E(\calL_{w})[p^{\infty}] = 0$, $\rmH^{1}_{\rmf}(\calL_{w}, T) \simeq \ZZ_p$ and $\rmH^{1}_{\rmf}(\calL_{w}, W) \simeq \QQ_p/\ZZ_p$.
        \item We define 
        \[
        b_p(E) := \begin{cases}
            1, &\quad \text{$E$ has good reduction at $p$}, \\
            0, &\quad \text{$E$ has additive reduction at $p$}.
        \end{cases}
        \]
        The logarithm given by the Néron differential $\omega_{E}$ identifies $E(\calL_{w})_{p}^{\vee}$ with $p^{b_p(E)}\ZZ_p$.
    \end{enumerate}
For $w \nmid p$ any finite place of $\calL$, we have
    \begin{enumerate}
        \item[\rm (3)] $\rmH^{1}_{\ur}(\calL_{w}, W) = \rmH^{1}_{\rmf}(\calL_{w}, W) = 0$.
    \end{enumerate}
For all finite places $w$ of $\calL$, we have $p \nmid \Tam_{w}(E/\calL)$ and $E(\calL)[p] = 0$.
\end{lemma}

\begin{proof}
    As we have assumed that $p$ splits in $\calL$, each $\calL_{w}$ for $w \mid p$ is equal to $\QQ_p$.
    
    Next, we prove (2). We write $E_1(\QQ_p)$ for the formal group points and $E_{0}(\QQ_p)$ for points reducing to the identity components. The Néron differential then gives
    \[
    \log_{\omega_{E}}(E_1(\QQ_p)) = p\ZZ_p.
    \]

    Suppose $E$ has good reduction at $p$. Since $E$ is a CM elliptic curve, $E$ has good supersingular reduction at $p$. The requirement $p \geq 5$ gives $a_p(E) = 0$. This implies that $E(\QQ_p)/E_1(\QQ_p) \simeq \tilE(\FF_p)$ has order $p+1$, coprime to $p$. Therefore, $E(\QQ_p)_p^{\wedge} = E_1(\QQ_p)$, and hence the logarithm image is $p\ZZ_p$.

    Suppose $E$ has additive reduction at $p$. Then $E_{0}(\QQ_p)/E_1(\QQ_p) \simeq \FF_p^{+}$, of order $p$. Since we have seen that $E(\QQ_p)[p] = 0$, the logarithm extends injectively to $E_0(\QQ_p)$. Its image contains $p\ZZ_p$ with index $p$, and therefore its image is $\ZZ_p$. Since the component index $[E(\QQ_p):E_{0}(\QQ_p)]$ is coprime to $p$, we see that passing from $E_{0}(\QQ_p)$ to $E(\QQ_p)_{p}^{\wedge}$ gives the same result.

    Next we prove (3). The inflation-restriction exact sequence gives
    \[
    \rmH^{1}_{\ur}(\calL_{w}, W) := \ker(\rmH^{1}(\calL_{w}, W) \rightarrow \rmH^{1}(I_{w}, W)) \simeq W^{I_{w}}/(\Frob_{w}-1)W^{I_{w}}.
    \]
    We have seen in the proof of Lemma \ref{lem:invariant-vanishing} that the inertia $I_{w}$ has finite image of order prime to $p$, so $W^{I_{w}} = V^{I_{w}}/T^{I_{w}}$. By purity, $\Frob_{w}-1$ is invertible on $V^{I_{w}}$ and hence it is surjective on its divisible quotient $W^{I_{w}}$. This shows that $W^{I_{w}}/(\Frob_{w}-1)W^{I_{w}} = 0$. Separately, one sees directly that $E(\calL_{w}) \otimes \QQ_p/\ZZ_p = 0$, so $\rmH^{1}_{\rmf}(\calL_{w}, W) = 0$ as well. This shows (3).

    Lastly, it is known by Kodaira's computation that $\Tam_{w}(E/\QQ)$ is at most $4$. Since $p \geq 5$, we see that $p \nmid \Tam_{w}(E/\QQ)$. The claim that $E(\calL)[p] = 0$ follows from that $E(\QQ_p)[p] = 0$.
\end{proof}

We choose $Q \in E(\calL)$ to be a generator of $E(\calL)/E(\calL)_{\tor}$. Recall that this is part of the hypothesis \eqref{eq:rank-one}. We put
\begin{align*}
    r_{w} := \length_{\ZZ_p} \coker(E(\calL)\otimes\ZZ_p \rightarrow E(\calL_{w})_{p}^{\wedge}) = v_p(\log_{\omega_{E}}(Q)_{w}) - b_{p}(E),
\end{align*}
for $w \mid p$. Here the equality follows from Lemma \ref{lem:gather-finite}(2), and from the proof of the lemma, we see that $\log_{\omega}(Q)_{w} \neq 0$. Moreover, $r_{w} = r_{\barw}$. We denote this common value by $r$.

\begin{proposition} \label{prop:length-finite-BDP}
    Under the hypothesis \eqref{eq:rank-one}, the BDP Selmer group $\Sel_{\BDP}(E/W)$ is finite, and
    \[
    \length_{\ZZ_p}(\Sel_{\BDP}(E/W)) = v_p \left( \# \Sha(E/\calL)[p^{\infty}] \right) + 2r.
    \]
\end{proposition}
\begin{proof}
    First we consider the Bloch--Kato Selmer group. The Kummer short exact sequence gives
    \begin{equation} \label{eq:kummer-T}
    0 \rightarrow E(\calL) \otimes \QQ_p/\ZZ_p \rightarrow \Sel_{\BK}(\calL, W) \rightarrow \Sha(E/\calL)[p^{\infty}] \rightarrow 0.
    \end{equation}
    Moreover, taking inverse limit in the Kummer short exact sequence gives
    \[
    0 \rightarrow E(\calL)\otimes\ZZ_p \rightarrow \Sel_{\BK}(\calL, T) \rightarrow T_{p}(\Sha(E/\calL)) = 0.
    \]
    The last group is zero since we have assumed in \eqref{eq:rank-one} that $\Sha(E/\calL)[p^{\infty}]$ is finite. This implies that $\Sel_{\BK}(\calL, T) = \ZZ_p Q$ for the generator $Q \in E(\calL)$ of $E(\calL)/E(\calL)_{\tor}$ we have chosen.

    The Poitou--Tate exact sequence then gives
    \begin{equation} \label{eq:poitou-tate-bdp}
    0 \rightarrow \Sel_{\BK}(\calL, W) \rightarrow \Sel_{\rel}(\calL, W) \rightarrow \dfrac{\rmH^{1}(\calL_{\barv_0}, W)}{\rmH^{1}_{\rmf}(\calL_{\barv_0}, W)} \xrightarrow{\alpha} \Sel_{\BK}(\calL, T)^{\vee}.
    \end{equation}
    It follows from the definition that the map $\alpha$ is the multiplication by $u_{\barv_{0}} p^{r_{\barv_0}}$ on $\QQ_p/\ZZ_p$. This map is surjective and has a kernel of order $p^{r_{\barv_{0}}}$. The exact sequence \eqref{eq:poitou-tate-bdp} gives
    \[
    \#(\Sel_{\rel}(\calL, W) / \Sel_{\BK}(\calL, W)) = p^{r_{\barv_{0}}}.
    \]
    Applying \eqref{eq:kummer-T}, we see that 
    \[
    \#(\Sel_{\rel}(\calL, W) / E(\calL) \otimes \QQ_p/\ZZ_p) = \#(\Sha(E/\calL)[p^{\infty}]) p^{r_{\barv_{0}}}.
    \]
    The last step is to pass from $\Sel_{\rel}(\calL,W)$ to $\Sel_{\BDP}(\calL, W)$. It suffices to note that
    \[
    \Sel_{\BDP}(\calL, W) = \ker \left(\Sel_{\rel}(\calL,W) \rightarrow \rmH^{1}_{\rmf}(\calL_{v_0}, W) \right).
    \]
    The snake lemma provides the following short exact sequence
    \[
    0 \rightarrow \ker(E(\calL) \otimes \QQ_p/\ZZ_p \xrightarrow{\beta} \rmH^{1}(\calL_{v_0}, W)) \rightarrow \Sel_{\BDP}(\calL, W) \rightarrow \Sel_{\rel}(\calL, W)/(E(\calL) \otimes \QQ_p/\ZZ_p) \rightarrow 0.
    \]
    The map $\beta$ here is multiplication by $u_{v_0} p^{r_{v_0}}$, and hence has kernel of order $p^{r_{v_0}}$. Therefore, we conclude that
    \[
    \# \Sel_{\BDP}(\calL, W) = \#(\Sha(E/\calL)[p^{\infty}]) p^{r_{v_0} + r_{\barv_{0}}} = \#(\Sha(E/\calL)[p^{\infty}]) p^{2r},
    \]
    and the proposition is therefore proved.
\end{proof}

\subsubsection{Control theorem for BDP Selmer groups and the refined formula}
We are now ready to prove the BDP control theorem.
\begin{proposition}[Control theorem for BDP Selmer groups] \label{prop:control-theorem-BDP}
    Let $\gamma_{\ac}$ be a topological generator of $\Lambda_{\calL}^{\ac}$. Under \eqref{eq:rank-one}, we have the following results.
    \begin{enumerate}[label = \rm (\arabic*)]
        \item There is a canonical isomorphism of $\calO$-modules: 
        \[
        X_{\BDP}(\bfT_{E}^{\ac}/\calL)/(\gamma_{\ac}-1)X_{\BDP}(\bfT_{E}^{\ac}/\calL) \simeq (\Sel_{\BDP}(E/W) \otimes_{\ZZ_p} \calO)^{\vee},
        \]
        \item Moreover, $\length_{\calO}(X_{\BDP}(\bfT_{E}^{\ac}/\calL)/(\gamma_{\ac}-1)X_{\BDP}(\bfT_{E}^{\ac}/\calL)) = e \left( v_p \left( \# \Sha(E/\calL)[p^{\infty}] \right) + 2r \right)$, where $e:= v_{\varpi}(p)$ is the ramification index of $p$ in $\calO$.
    \end{enumerate}
    In particular, $X_{\BDP}(\bfT_{E}^{\ac}/\calL)$ is a finitely generated torsion $\Lambda_{\calL}^{\ac}$-module. 
\end{proposition}

\begin{proof}
    For (1), it suffices to prove the control theorem between $\Sel_{\BDP}(\bfT_{E}^{\ac}/\calL)$ and $\Sel_{\BDP}(E/W)$. Indeed, dualizing the injective map $\Lambda_{\calL}^{\ac} \rightarrow \Lambda_{\calL}^{\ac}$ given by multiplication by $\gamma_{\ac}-1$ gives an exact sequence 
    \[
    0 \rightarrow (\Lambda_{\calL}^{\ac}/(\gamma_{\ac}-1))^{\vee} \rightarrow (\Lambda_{\calL}^{\ac})^{\vee} \xrightarrow{\times (\gamma_{\ac}-1)} (\Lambda_{\calL}^{\ac})^{\vee} \rightarrow 0.
    \]
    Then $(\Lambda_{\calL}^{\ac})^{\vee}[\gamma_{\ac}-1] \simeq L/\calO$. Since $T_{E}$ is free over $\calO$, we have
    \[
    \bfT_{\calL}^{\ac}[\gamma_{\ac}-1] \simeq T_{E} \otimes_{\calO} (L/\calO) \simeq E[p^{\infty}]\otimes_{\ZZ_p} \calO = W \otimes_{\ZZ_p} \calO.
    \]
    Note that the tautological character becomes trivial modulo $\gamma_{\ac}-1$, so it is an isomorphism of Galois modules. The remaining proof follows in the same way as in the proof of Proposition \ref{prop:descent-theorem}, invoking Lemma \ref{lem:invariant-vanishing}(3)(4).

    The result in (2) follows from (1) and Proposition \ref{prop:length-finite-BDP}. Note that, for any finite abelian $p$-group with cyclic factors $\ZZ/p^{a_i}\ZZ$, its tensor product with $\calO$ over $\ZZ_p$ has factors $\calO/p^{a_i}\calO$, of $\calO$-length $ea_i$. Finally, the $\Lambda_{\calL}$-cotorsion property of $\Sel_{\BDP}(\bfT_{E}^{\ac}/\calL)$ follows from the finiteness of $\Sel_{\BDP}(E/W)$ by topological Nakayama's lemma together with the control theorem in (1). Note that the finiteness of $\Sel_{\BDP}(E/W)$ is already proved in Proposition \ref{prop:length-finite-BDP}.
\end{proof}

%% file: 203-analytic-comparison.tex
\section{Analytic side: $p$-adic $L$-functions and their comparison} \label{sec:comparison}
This is the most technical section in this part. We compare the interpolation formula for the projection of Katz's $p$-adic $L$-function over $\calM/\calF$ with the BDP $p$-adic $L$-function.

\subsection{Interpolation characters} \label{sec:interpolation-characters}
In this section, we introduce Katz's $p$-adic $L$-function and the BDP $p$-adic $L$-function. The latter is modified from the original article of Bertolini--Darmon--Prasanna to fit our setup.

For the setup, we write $\tpsi := \psi \circ \Norm_{\calM/\calK}$ and set $\lambda_{0} := \tpsi \bfN^{-1}$, so $\hlam_0 = \rho^{-1}\epsilon_{\cyc}$. Let $\nu$ be an anticyclotomic Hecke character of $\calL$ of infinite type $(n,-n)$, unramified at all finite places, whose $p$-adic avatar $\hnu$ factors through $\Gamma_{\calL}^{\ac}$. We call such Hecke characters $\nu$ over $\calL$ \emph{allowable} Hecke characters.

We show that there are sufficiently many allowable characters.

\begin{lemma} \label{lem:allowables}
    There is an integer $n_0 > 0$ divisible by $p-1$ and a Hecke character of $\calL$ of infinite type $(n_0, -n_0)$ such that $\nu_0$ is unramified at all finite places, and its $p$-adic avatar $\hnu_{0}$ is a nontrivial character of $\Gamma_{\calL}^{\ac}$. The characters $\hnu_{j} := \hnu_{0}^{p^{j}}$ are distinct and tend to $1$ as $j$ goes to infinity.
\end{lemma}
\begin{proof}
    Let $h$ be the class number of $\calL$. For every nonzero fractional ideal $\fra$, the ideal $\fra^{h}$ is principal, and we choose an element $\alpha_{\fra} \in \calL^{\times}$ such that $(\alpha_{\fra}) = \fra^{h}$. We define an ideal class character $\nu_{h}(\fra) := \overline{\alpha_{\fra}}/\alpha_{\fra}$. One checks directly that:
    \begin{itemize}
        \item The definition of $\nu_{h}(\fra)$ does not depend on the choice of the generator $\alpha_{\fra}$;
        \item The character $\nu_{h}$ such defined is multiplicative, i.e. $\nu_{h}(\fra\frb) = \nu_{h}(\fra) \nu(\frb)$. 
        \item The character $\nu_{h}$ is anticyclotomic, i.e. $\nu_{h}(\barfra) = \nu_{h}(\fra)^{-1}$.
    \end{itemize}
    The first two verifications use the assumption that $D_{\calL} < -4$, so $\calO_{\calL}^{\times} =\{ \pm 1\}$. For any principal ideal $(z)$, we can take $\alpha_{z} = z^{h}$. Therefore, $\nu_{h}((z)) = z^{-h}z^{h}$. This says that the corresponding Hecke character $\nu_{h}$ has infinite type $(h,-h)$, is anticyclotomic, and is unramified at all finite places since $\nu_{h}$ is defined on the group of all fractional ideals of $\calO_{\calL}$.

Let $x = (x_{w}) \in \AA_{\calL, \rmf}^{\times}$, and let $\fra(x) = \prod_{w} \frp_{w}^{v_{w}(x_w)}$. The $p$-adic avatar of $\nu_{h}$ is written as
\[
\hnu_{h}(x) = \iota_{p}(\nu_h(\fra(x))) x_{v_0}^{h} x_{\barv_0}^{-h},
\]
and one verifies directly that $\hnu_{h}(z) = 1$ for $z \in \calL^{\times}$. So $\hnu_{h}$ is trivial on $\calL^{\times}$, and hence defines a character of $\calL^{\times} \bs \AA_{\calL}^{\times}$. By global class field theory, we view it as a $p$-adic Galois character $\hnu_{h}: \Gal_{\calL}^{\ab} \rightarrow \barQQ_{p}^{\times}$. 

Next, we adjust the character $\nu_{h}$ so that $\hnu_{h}$ factors through $\Gamma_{\calL}^{\ac}$. We return to the construction of the anticyclotomic $\ZZ_p$-extension $\calL_{\infty}^{\ac}/\calL$. Let $\frX_{\calL}$ be the Galois group of the maximal abelian pro-$p$ extension of $\calL$ unramified outside $p$. Then global class field theory and Leopoldt's conjecture for the imaginary quadratic field $\calL$ imply that $\frX_{\calL, \tf} \simeq \ZZ_{p}^{\times}$. The complex conjugation $c \in \Gal(\calL/\QQ)$ decomposes the two-dimensional $\ZZ_p$-vector space $\frX_{\calL, \tf}$ into two rank-one eigenspaces, where the minus part $(\frX_{\calL,\tf})^{-} \simeq \ZZ_p$ gives the anticyclotomic $\ZZ_p$-extension, i.e.,
\[
\Gamma_{\calL}^{\ac} = (\frX_{\calL, \tf})^{-}.
\]
We choose a finite extension $\calO$ over $\ZZ_p$ containing the values of $\hnu_{h}$, with maximal ideal $\frm$ and residue field $k$. The image of $\hnu_{h}$ lies in $\calO^{\times}$, and the reduction of the image of $\hnu_{h}$ in $k^{\times}$ is finite, so let $e_{1}$ be its exponent. Then $\hnu_{h}^{e_1}$ takes values in $1+\frm$. Let $e_{2}$ be the exponent of the finite group $\frX_{\calL, \tor}$. Then $\hnu_{h}^{e_2}$ is trivial on $\frX_{\calL, \tor}$. We choose a positive integer $m$ divisible by $e_1$, $e_2$, and $p-1$, define $\nu_0 := \nu_{h}^{m}$, and put $n_0 := hm$. By construction, $\hnu_0$ is trivial on $\frX_{\calL,\tor}$. Moreover, it is trivial on $\Gamma_{\calL}^{\cyc} = (\frX_{\calL, \tf})^{+}$. Indeed, let $g \in \Gamma_{\calL}^{\cyc}$; then $cgc^{-1} = g$. Applying $\hnu_{0}$, we have
\[
\hnu_0(g) = \hnu_0(cgc^{-1}) = \hnu_0(g)^{-1},
\]
and hence $\hnu_0(g)^{2} = 1$. But $\hnu_0(g) \in 1+\frm$, and $1+\frm$ has no nontrivial element of order $2$ (note that $p \geq 5$ is an odd prime). Therefore, $\hnu_{0}(g) = 1$, i.e. $\hnu_0$ factors through $\Gamma_{\calL}^{\ac}$.

Take $u = 1+p$, place $u \in \calO_{\calL, v_0}^{\times}$ in the $v_0$-component of a finite idèle, and put $1$ in all other components. By abuse of notation, denote this idèle by $u$. Its associated fractional ideal is the unit ideal, and hence $\nu_{h}(\fra(u)) = 1$, so $\hnu_{h}(x) = (1+p)^{h}$. Since $(1+p)^{h}$ has infinite order (as $1+p\ZZ_p$ is torsion-free for odd $p$), the character $\hnu_{h}$ has infinite image, and therefore every $\hnu_j := \hnu_{0}^{p^j}$ is nontrivial.

Moreover, each $\hnu_j$ is the $p$-adic avatar of $\nu_j := \nu_0^{p^j}$ which has infinite type $(n_0 p^{j}, -n_0 p^{j})$. Choose a topological generator $\gamma_{\ac} \in \Gamma_{\calL}^{\ac}$ and put $b := \hnu_0(\gamma_{\ac}) \in 1+\frm$, then since $\hnu_0$ is nontrivial, $b \neq 1$ and has infinite order. Suppose $\hnu_{j} = \hnu_{k}$, then evaluation at $\gamma_{\ac}$ implies $b^{p^{k}(p^{j-k}-1} = 1$, which forces $j = k$. So the characters $\hnu_j$ are pairwise distinct.

Finally, since $b \in 1+\frm$, the $p$-adic logarithm of $b$ is well-defined and $b^{p^{j}} = \exp(p^{j} \log b)$. As $j \rightarrow \infty$, we see that $p^{j} \log b \rightarrow 0$, and hence $b^{p^{j}} \rightarrow 1$. More generally, every $g \in \Gamma_{\calL}^{\ac}$ can be written as $g = \gamma_{\ac}^{a}$ for $a \in \ZZ_p$, then
\[
\hnu_j(g) = b^{ap^{j}} = \exp(ap^{j} \log b) \rightarrow 1,
\]
uniformly in $a \in \ZZ_p$. Therefore, $\hnu_j \rightarrow 1$ in the $p$-adic character space of $\Gamma_{\calL}^{\ac}$.
\end{proof}

We label $\sigma_1$ so that its restriction to $\calK$ is the embedding at which $\psi$ has exponent $-1$. In the order $(\sigma_1,\sigma_2,\bar\sigma_1,\bar\sigma_2)$, one computes the infinite types of the following characters, where $\nu$ is an allowable Hecke character over $\calL$, as shown in Table \ref{tab:infinite-type}:
\begin{table}[H]
    \centering
    \begin{tabular}{c|rrrr}
 &$\sigma_1$&$\sigma_2$&$\bar\sigma_1$&$\bar\sigma_2$ \\ \hline
 $\tpsi$&$-1$&$0$&$0$&$-1$\\
 $\tnu$&$n$&$n$&$-n$&$-n$\\
$\bfN^{-1}$&$1$&$1$&$1$&$1$\\ \hline
 $\lambda_\nu:=\tpsi\tnu\bfN^{-1}$ &$n$&$n+1$&$1-n$&$-n$
    \end{tabular}
    \caption{Infinite types of characters over $\calM$}
    \label{tab:infinite-type}
\end{table}

\subsection{Two $p$-adic $L$-functions}
In this section, we introduce two $p$-adic $L$-functions: one is Katz's $p$-adic $L$-function for the CM field $\calM/\calF$, attached to the Hecke character $\lambda_{0}$ over $\calM$, and the other is the Bertolini--Darmon--Prasanna (BDP) $p$-adic $L$-function attached to $E$ over $\calL$. Our main goal is to compare the two $p$-adic $L$-functions, as established in Theorem \ref{thm:comparison}.

\subsubsection{Comparison of CM periods}
Before introducing the two $p$-adic $L$-functions, namely Katz's $p$-adic $L$-function and the BDP $p$-adic $L$-function, we introduce the complex and $p$-adic periods used in their interpolation formulas. These choices will be fixed throughout this article. We closely follow \cite[Section 4.5]{JSW} in this subsection.

\underline{(1) Complex and $p$-adic periods of CM elliptic curves: over $\calL$}

We have assumed that $p$ splits in $\calL$. Let $H/\calL$ be the Hilbert class field of $\calL$ and let $F$ be a finite extension of $H$. Let $E_{0}/F$ be an elliptic curve with complex multiplication by an order $\frD$ in the imaginary quadratic field $\calL$. We assume that $\frD \otimes \ZZ_{(p)}$ is a maximal $\ZZ_{(p)}$-order. We fix a complex uniformization $E_0(\CC) \simeq \CC/\fra$, where $\fra$ is a nonzero ideal of $\frD = \End(E_0)$. Enlarging $F$ when necessary, we may assume that $E_{0}$ has good reduction at all primes above $p$ and so it extends to an elliptic scheme $\calE_{0}/\calO_{F, (p)}$. The space of differentials $\rmH^{0}(\calE_{0}, \Omega^{1}_{\calE_{0}/\calO_{F, (p)}})$ is a free $\calO_{F, (p)}$-module of rank one. We can write
\[
\rmH^{0}(\calE_{0}, \Omega^{1}_{\calE_{0}/\calO_{F, (p)}}) = \calO_{F, (p)} \cdot \omega_{E_{0}}
\]
for some differential $\omega_{E_0} \in \rmH^{0}(\calE_{0}, \Omega^{1}_{\calE_{0}/\calO_{F, (p)}})$.

The \emph{complex period} of $E_{0}$ is defined as follows. We consider $\CC$ to be an $F$-algebra via the fixed inclusion $\iota_{\infty}: \barQQ \hookrightarrow \CC$ inducing the embedding $\sigma_0: \calL \hookrightarrow \CC$ in the CM type $\Sigma_{0}$ of $\calL$. Then
\[
\rmH^{0}(\calE_{0}, \Omega^{1}_{\calE_{0}/\calO_{F, (p)}}) \otimes_{\calO_{F, (p)}} \CC = \rmH^{0}(E_{0}, \Omega^{1}_{E_{0}/\CC}).
\]
As the invariant differential $\omega_{\CC} := \dif z$ on the complex torus $E_{0}(\CC) \simeq \CC/\fra$ gives a generator for $\rmH^{0}(E_{0}, \Omega^{1}_{E_{0}/\CC})$, there exists a scalar $\Omega_{\infty} \in \CC^{\times}$ such that 
\begin{equation} \label{eq:complex-period-L}
    \omega_{E_{0}} = \Omega_{\infty} \cdot (2 \pi i) \cdot \omega_{\CC}.
\end{equation}
The complex number $\Omega_{\infty}$ is called the \emph{complex period} of $\calL$ with respect to the CM type $\Sigma = \{\sigma_0\}$.

We define the \emph{$p$-adic period} of $E_{0}$ in much the same way. As $p$ splits in $\calL$ as $p \calO_{\calL} = \frp \barfrp$, we take the place $\frp$ to be that determined by the fixed embedding $\iota_{p}: \barQQ \hookrightarrow \barQQ_{p}$ and take $\iota_{\frp}: \barQQ \hookrightarrow \overline{\calL}_{\frp} = \barQQ_{p}$ to be $\iota_{p}$. Denote also by $\frp$ the place of $F$ determined by the embedding $\iota_{\frp}$, so $\iota_{\frp}$ identifies $F_{\frp}$ with a finite extension of $\calL_{\frp} = \QQ_p$ in $\barQQ_{p}$. Let $F_{\frp}^{\ur}$ be the maximal unramified extension of $F_{\frp}$ in $\barQQ_{p}$, and let $\calO^{\ur}$ be the ring of integers of $F_{\frp}^{\ur}$. Consider the good integral model $\calE_{0}/\calO^{\ur}$, and let $\widehat{\calE}_0$ be the formal completion over the identity section. As we are working over $\calO^{\ur}$, the formal completion $\widehat{\calE}_0$ is noncanonically isomorphic to $\widehat{\GG}_{\rmm}$. We fix such an isomorphism $\widehat{\calE}_{0} \simeq \widehat{\GG}_{\rmm}$. The latter is equivalent to fixing an isomorphism of $p$-divisible groups $\calE_{0}[p^{\infty}]^{0} = \calE_{0}[\frp^{\infty}] \simeq \mu_{p^{\infty}}$, which is uniquely determined up to the action of an element of $\ZZ_{p}^{\times}$. Here $\frp$ is the prime ideal of $\frD$ corresponding to $\frp$. The pullback of the differential $\dif t/t$ under the fixed isomorphism $\widehat{\calE}_{0} \simeq \widehat{\GG}_{\rmm}$ is an element
\[
\omega_{\can} \in \rmH^{0}(\calE_{0}, \Omega_{\calE_{0}/\calO^{\ur}}^{1}) = \rmH^{0}(\calE_{0}, \Omega_{\calE_{0}/\calO_{F, (p)}}^{1}) \otimes_{\calO_{F, (p)}} \calO^{\ur}.
\]
The \emph{$p$-adic period} $\Omega_{p} \in (\calO^{\ur})^{\times}$ is defined such that 
\begin{equation} \label{eq:padic-period-L}
\omega_{E_0} = \Omega_{p} \cdot \omega_{\can}.
\end{equation}

\underline{(2) Complex and $p$-adic periods of CM abelian surfaces: over $\calM$}

To relate the CM period for $\calL/\QQ$ and the CM period for $\calM/\calF$, we use Serre's tensor product construction. Regarding $\calO_{\calM}$ as a projective $\calO_{\calL}$-module of rank two, we form the Serre tensor $\calA_0$ as a CM abelian scheme over $\calO_{F, (p)}$ such that $\calA_{0}(R) = \calO_{\calM} \otimes_{\calO_{\calL}} \calE_{0}(R)$ for any $\calO_{F, (p)}$-algebra $R$. By the construction of Serre's tensor, we have
\[
\omega_{\calA_{0}} \simeq \calO_{M}^{\vee} \otimes_{\calO_{\calL}} \omega_{\calE_{0}}, \quad \calO_{M}^{\vee} := \Hom_{\calO_{\calL}}(\calO_{\calM}, \calO_{\calL}).
\]
For any $\calO_{\calL}$-algebra $S$ which carries a compatible $\CC$-module structure under the embedding $\sigma_{0}: \calL \rightarrow \CC$, we let
\[
\ell_{i, S} : \calO_{\calM} \otimes_{\calO_{\calL}} S \rightarrow S, \quad b \otimes x \mapsto \sigma_{i}(b)x
\]
for $i = 1,2$. Then we can define two differentials $\omega_{i} := \ell_{i,S} \otimes \omega_{E_0}$ where $\omega_{E_{0}}$ is the fixed differential in paragraph (1). 

We first consider the complex period of $\calA_{0}$. Under the embedding $\sigma_{0}: \calL \hookrightarrow \CC$ in the fixed CM type $\Sigma_{0}$ of $\calL$, we have
\[
\calA_{0}(\CC) \simeq \dfrac{\calO_{\calM} \otimes_{\calO_{\calL},\sigma_{0}}\CC}{\calO_{\calM} \otimes_{\calO_{\calL},\sigma_{0}}\fra} \simeq \dfrac{\CC^{2}}{\Lambda_{\calA_{0}}},
\]
where the last isomorphism is induced by the decomposition
\[
\Theta_{\infty}: \calM \otimes_{\calL, \sigma_0} \CC \xrightarrow{\sim} \CC^{2}, \quad b \otimes z \mapsto (\sigma_1(b)z, \sigma_2(b)z),
\]
where $\Sigma = \{ \sigma_1, \sigma_2\}$ is the fixed CM type of $\calM/\calF$ such that $\sigma_1|_{\calL} = \sigma_{2}|_{\calL} = \sigma_{0}$, and, under the map $\Theta_{\infty}$, the lattice $\Lambda_{\calA_{0}} := \{ (\sigma_1(b)z, \sigma_2(b)z): b \in \calO_{\calM}, z \in \fra \}$ is thereby defined. Under this decomposition, we fix the $\CC$-basis $\omega_{1,\CC}$ and $\omega_{2, \CC}$ of $\rmH^{0}(\calA_{0}(\CC), \Omega^{1}_{\calA_{0}(\CC)/\CC})$, i.e.,
\begin{equation} \label{eq:complex-period-M}
    \omega_{i,\CC} = \ell_{i,\CC} \otimes \omega_{\CC}, \quad\text{for } i=1,2.
\end{equation}
It then follows that
\[
\omega_{i} = \ell_{i,\CC} \otimes \omega_{E_0} \xlongequal{{\text{\eqref{eq:complex-period-L}}}} \Omega_{\infty} \cdot (\ell_{i,\CC} \otimes \omega_{\CC}) \xlongequal{{\text{\eqref{eq:complex-period-M}}}}  \Omega_{\infty} \cdot \omega_{i,\CC}
\]
for $i=1,2$, where $\Omega_{\infty} \in \CC^{\times}$ is the complex period in paragraph (1). By the definition of Katz's complex period, this shows that we can choose $\Omega_{\infty, \sigma_{i}} = \Omega_{\infty}$ for $i = 1, 2$.

The comparison of the $p$-adic period of $E_{0}$ and $\calA_{0}$ is similar. Recall we have put $v$ to be the prime of $\calM$ lying over $v_{0}$ selected by $\Sigma$. On the ordinary component, the Serre tensor construction gives
\begin{equation} \label{eq:serre-at-p}
\calA_0[p^{\infty}]^{0} \simeq \calO_{\calM, v} \otimes_{\ZZ_p} \calE_{0}[v_0^{\infty}]
\end{equation}
over $\calO^{\ur}$. The two $p$-adic avatars of $\sigma_1$ and $\sigma_2$ given by $\sigma_{i,p} = \iota_p \circ \sigma_{i}: \calM \hookrightarrow \CC_p$ for $i=1,2$ induce an integral decomposition
\[
\Theta_{p}: \calO_{M,v} \otimes_{\ZZ_p} \calO^{\ur} \xrightarrow{\sim} (\calO^{\ur})^{2}, \quad b \otimes x \mapsto (\sigma_{1,p}(b)x, \sigma_{2,p}(b)x).
\]
(Indeed, since $\calM_{v}/\QQ_p$ is an unramified quadratic extension, if $b_1, b_2$ is a $\ZZ_p$-basis of $\calO_{\calM,v}$ and consider the matrix $[\sigma_{i,p}(b_j)]_{i,j}$, then its determinant is $\disc(b_1, b_2) \in \ZZ_p^{\times}$, and hence it is indeed integral.) Applying these idempotents $\ell_{1,\calO^{\ur}}$ and $\ell_{2, \calO^{\ur}}$ to \eqref{eq:serre-at-p}, it follows that the canonical differentials of $\calA_{0}$ are
\begin{equation} \label{eq:padic-period-M}
\omega_{i, \can} = \ell_{i,\calO^{\ur}} \otimes \omega_{0,\can}, \quad\text{for } i=1,2,
\end{equation}
where $\omega_{0,\can}$ is the canonical differential of $E_{0}$ in paragraph (1). This yields that
\[
\omega_{i} = \ell_{i, \calO^{\ur}} \otimes \omega_{E_0} \xlongequal{{\text{\eqref{eq:padic-period-L}}}} \Omega_{p} \cdot (\ell_{i, \calO^{\ur}} \otimes \omega_{0, \can}) \xlongequal{{\text{\eqref{eq:padic-period-M}}}} \Omega_{p} \cdot \omega_{i,\can},
\]
where $\Omega_{p}$ is the $p$-adic period in paragraph (1). Consequently, the $p$-adic period $\Omega_{p, \sigma_{i}} = \Omega_{p}$ for $i=1,2$.

In conclusion, from this point onward, given the CM type $\Sigma$ of $\calM$, we fix compatible Katz periods $(\Omega_{\infty, \sigma}, \Omega_{p, \sigma})_{\sigma \in \Sigma}$ over $\calM$ such that
\begin{equation} \label{eq:comparison-period}
    \Omega_{\infty, \sigma_1} = \Omega_{\infty, \sigma_2} = \Omega_{\infty}, \quad \Omega_{p,\sigma_1} = \Omega_{p, \sigma_2} = \Omega_{p},    
\end{equation}
where $(\Omega_{\infty}, \Omega_{p})$ is the Katz's period over $\calL$ with CM type $\Sigma_{0}$.

\subsubsection{Katz's $p$-adic $L$-function}
Let $\frc := \cond^{(p)}(\lambda_{0}) = \cond^{(p)}(\tilpsi)$ be the prime-to-$p$ part of the conductor of $\lambda_{0}$. Let $\calM(\frc p^{\infty})$ be the union of the ray class fields of $\calM$ of conductor $\frc p^n$ for $n \geq 0$, and denote $G_{\frc} := \Gal(\calM(\frc p^{\infty})/\calM)$. Write $\Lambda_{\calO}^{\ur}(G_{\frc})= \calO^{\ur}\lrbracket{G_{\frc}}$ for the Iwasawa algebra attached to $G_{\frc}$.

Let $\Omega_{\infty, \sigma_1}, \Omega_{\infty, \sigma_2} \in \CC^{\times}$ (resp. $\Omega_{p, \sigma_1}, \Omega_{p,\sigma_2} \in \barZZ_p^{\times}$) be the \emph{complex CM period} (resp. \emph{$p$-adic CM period}) of $(\calM, \Sigma)$ fixed in \eqref{eq:comparison-period}.

\begin{theorem}[Katz $p$-adic $L$-function] \label{thm:katz}
There exists an element $\calL_{\Katz}(\calM, \frc, \Sigma) \in \Lambda_{\calO}^{\ur}(G_{\frc})$ characterized by the following interpolation property: for every Hecke character $\nu$ of $\calM$ of infinite type $k \Sigma + \kappa (1-\sfc)$, where 
\begin{enumerate}[label = \rm (\roman*)]
    \item either $k \geq 1$ and $\kappa \in \ZZ_{\geq 0}[\Sigma]$, or 
    \item $k \leq 1$ and $k \Sigma + \kappa \in \ZZ_{>0}[\Sigma]$,
\end{enumerate}
such that its $p$-adic avatar $\hnu$ factors through $G_{\frc}$,
we have
\begin{equation} \label{eq:interpolation-katz}
\calL_{\Katz}(\calM, \frc, \Sigma)(\hnu) = \left(\dfrac{\Omega_{p}}{\Omega_{\infty}}\right)^{k\Sigma+2\kappa} \dfrac{\pi^{\kappa} \Gamma_{\Sigma}(k\Sigma+\kappa)[\calO_{\calM}^{\times}:\calO_{\calF}^{\times}]}{\sqrt{\abs{D_{\calF}}_{\RR}} (\im \vartheta)^{\kappa}} t_{\calM} \Eul_{p}(\nu) \Eul_{\frc}(\nu) L^{(p\frc)}(\nu, 0).
\end{equation}
Here $L^{(p\frc)}(\nu, 0)$ denotes the Hecke $L$-function of $\nu$ with the Euler factors at the primes $w \mid p\frc$ removed, and
\[
\Eul_{p}(\nu) := \prod_{w \in \Sigma_{p}} \Eul_{w}(\nu), \quad \Eul_{\frc}(\nu) := \prod_{w \mid \frf} \Eul_{w}(\nu)
\]
are modified Euler factors where $\frc = \frf\frf^{\sfc}$ is a product of split primes in $\calM/\calF$ such that $(\frf, \frf^{\sfc}) = 1$, and
\[
\Eul_{w}(\nu) := \nu_{w}(2 \vartheta_{w}) \dfrac{L(\nu_{w}, 0)}{\epsilon(\nu_{w}, \psi, 0) L(1, \nu_{w}^{-1})},
\]
the number $\vartheta \in \calM$ is chosen explicitly in the proof, and $t_{\calM}$ is an explicit constant, namely a power of $2$ depending only on the number fields $\calF$ and $\calM$.
\end{theorem}

\begin{definition}[The $\lambda_0$-branch of the Katz $p$-adic $L$-function]
Recall that we fixed a Hecke character $\lambda_{0}$ over $\calM$ in Section \ref{sec:interpolation-characters}. We define \footnote{Here we omit $\frc$ and $\Sigma$ in the notation of $\calL_{\Katz}(\lambda_{0}/\calM)$ for simplicity, but it indeed depends on $\frc$ and $\Sigma$.} $\calL_{\Katz}(\lambda_{0}/\calM) \in \Lambda_{\calM}^{\ur}$ to be the image of $\calL_{\Katz}(\calM, \frc, \Sigma)$ given by Theorem \ref{thm:katz} under the composite map
\[
\Lambda_{\calO}^{\ur}(G_{\frc}) \xrightarrow{\Tw_{\hat{\lambda}_{0}}} \Lambda_{\calO}^{\ur}(G_{\frc}) \twoheadrightarrow \Lambda_{\calM}^{\ur},
\]
where $\Tw_{\hat{\lambda}_{0}}$ is the $\calO$-linear isomorphism given by $\gamma \mapsto \hlam_{0}(\gamma)\gamma$ for $\gamma \in G_{\frc}$, and the second arrow is the natural projection. Then the interpolation formula of $\calL_{\Katz}(\lambda_{0}/\calM)$ is merely replacing the appearance of $\nu$ in \eqref{eq:interpolation-katz} by $\lambda_{\nu}$ defined at the beginning of this section.     
\end{definition}

\begin{proof}[Proof of Theorem \ref{thm:katz}]
    This is Katz's $p$-adic $L$-function, originally constructed by Katz in \cite{katz78}. Here we use the version in \cite[Proposition 4.9]{HsiehMu}. More explicitly, we insert the following data into the setup in \cite{HsiehMu}.

    \underline{(1)} In \cite{HsiehMu}, the ideal $\frc = \frc^{+} \frc^{-}$, where $\frc^{+} = \frf \frf^{\sfc}$ is a product of split primes in $\calM/\calF$ such that $(\frf, \frf^{\sfc}) = 1$, and $\frc^{-}$ is a product of nonsplit primes in $\calM/\calF$. In our case, since every $\ell \mid N$ splits in $\calL$, it follows that every prime dividing $\frc$ is split in $\calM/\calF$, and hence $\frc^{-} = (1)$ and $\frc^{+} = \frc$.

    \underline{(2)} Let $\frd := p \frc^{2} D_{\calM/\calF}$. In \cite{HsiehMu}, Hsieh fixes a $\vartheta \in \calM$ such that 
    \begin{enumerate}[label = (\roman*)]
    \item $\vartheta^{\sfc} = -\vartheta$ and $\im(\sigma(\vartheta)) > 0$ for all $\sigma \in \Sigma$,
    \item $\frD_{\calF}^{-1}(2 \vartheta \frD_{\calM/\calF}^{-1})$ is coprime to $\frd$.
\end{enumerate}
    In our case, we take $\vartheta = a\sqrt{D_{\calL}}/2$, where $a \in \calF^{\times, +}$ is a totally positive element in $\calF$, and $a \in \calO_{\calF, \frp}^{\times}$ for every $\frp \mid p$. The condition (i) is satisfied immediately. We choose $a$ so that Hsieh's ideal $\frD_{\calF}^{-1}(2 \vartheta \frD_{\calM/\calF}^{-1})$ is coprime to $\frd$. Such an $a$ exists by the weak approximation theorem. 
\end{proof}

\subsubsection{BDP $p$-adic $L$-function}
We next recall the ``square-root'' $p$-adic $L$-function of Bertolini--Darmon--Prasanna \cite{BDP}. In this article, we use the construction of Hsieh in \cite{HsiehSpecial}. We write $N^{(p)} = Np^{-v_p(N)}$ and we choose an ideal $\frN^{(p)}$ of $\calL$ above $N^{(p)}$ such that $\calO_{\calL}/\frN^{(p)} \simeq \ZZ/N^{(p)}\ZZ$.

Let $\Pi$ be the unitary automorphic representation of $\GL_2(\AA_{\QQ})$ attached to $f_{E}$, so that $L(s, \Pi) = L(E, s+ \frac{1}{2})$, and moreover, we have the following comparison of $L$-values given by the automorphic induction
\begin{equation} \label{eq:compare-L-value}
L(\lambda_{\nu}, 0) = L(\tpsi\tnu, 1) = L(f_{E}/\calL, \nu, 1).
\end{equation}
We define
\[
\xi := \psi_p \abs{-}_{\calK_p}^{1/2}, \quad\text{then } \Pi_{p} := \AI_{\calK_{p}/\QQ_p}(\xi).
\]
We use the standard additive character $\Psi_{p}$ of conductor $\ZZ_p$ and $\Psi_{\calK_{p}} := \Psi_{p} \circ \Tr_{\calK_{p}/\QQ_{p}}$. Let $\nu$ be an \emph{allowable} anticyclotomic character of $\calL$ of infinite type $(n,-n)$ whose $p$-adic avatar factors through $\Gamma_{\calL}^{\ac}$, we define
\begin{equation} \label{eq:epnu-bdp}
e_{p}(\nu) := \dfrac{\epsilon(1/2, \Pi_{p} \otimes \nu_{\barv_0}, \Psi_{p})}{L(1/2, \Pi_{p} \otimes \nu_{\barv_0})^{2}}.
\end{equation}

We write $d := \abs{D_{\calL}}$ for simplicity. We are ready to state the following theorem.
\begin{theorem}[BDP $p$-adic $L$-function] \label{thm:BDP}
    There exists an element $\calL_{\BDP}(E/\calL) \in \Lambda_{\calL}^{\ac,\ur}$ characterised by the following interpolation property: for any anticyclotomic character $\nu$ of $\calL$ of infinite type $(n,-n)$ whose $p$-adic avatar $\hnu$ factors through $\Lambda_{\calL}^{\ac}$ and is allowable, we have
    \[
    \calL_{\BDP}(E/\calL)(\hnu) = \left( \dfrac{\Omega_{p}}{\Omega_{\infty}}\right)^{4n} e_{p}(\nu) \dfrac{(2\pi)^{2n-1}\Gamma(n) \Gamma(n+1) \nu(\frN^{(p)})^{-1}}{4d^{n-\frac{1}{2}}} L(f_{E}/\calL, \nu, 1).
    \]
\end{theorem}
\begin{proof}
   We begin with a harmless twist of the BDP $p$-adic $L$-function. We take and fix an anticyclotomic Hecke character $\eta$ over $\calL$ of infinite type $(1,-1)$ and unramified outside $p$, whose $p$-adic avatar $\heta$ factors through $\Gamma_{\calL}^{\ac}$. For the $p$-adic avatar $\heta$ as a character over $\Gamma_{\calL}^{\ac}$, we define an automorphism
\[
\Tw_{\heta^{-1}}: (\Lambda_{\calL}^{\ac})^{\ur} \rightarrow (\Lambda_{\calL}^{\ac})^{\ur}, \quad [g] \mapsto \heta^{-1}(g)[g]. 
\]
We denote the $p$-adic $L$-function $\scL_{\Sigma_{p}}(\pi, \lambda)$ in \cite[Theorem A]{HsiehSpecial} by $\Theta_{\ttH}(E/\calL) \in \Lambda_{\calL}^{\ac}$: we take the CM field $\calK/\calF$ there to be $\calL/\QQ$, the automorphic representation $\pi$ there to be $\Pi$ of weight $k=2$, and the branch character $\lambda$ there to be the character $\eta$. The conditions in \cite{HsiehSpecial} on the ideal $\frn^{-}$ are vacuous since we have assumed that every prime $\ell$ dividing $N$ splits in $\calL$. Thus, \cite[Theorem A]{HsiehSpecial} is applicable. We define
\[
\Theta_{\BDP}(E/\calL) := c_{\ttH}^{1/2} \Tw_{\heta^{-1}}(\Theta_{\ttH}(E/\calL)), \quad \calL_{\BDP}(E/\calL) := \Theta_{\BDP}(E/\calL)^{2}.
\]
where $c_{\ttH} := \sqrt{d}(2\eta(\frN^{(p)}))^{-1}$. We will review the constructions in Section \ref{sec:construct-BDP} because this is crucial for proving the $p$-adic Waldspurger formula (Theorem \ref{thm:waldspurger}); at this point, however, we only need to verify the interpolation formula stated in the theorem.

Applying \cite[Theorem A]{HsiehSpecial} to the character $\hphi = \hnu \heta^{-1}$ of infinite type $(n-1, 1-n)$, we have
\[
\Theta_{\ttH}(E/\calL)(\hphi)^{2} = \left( \dfrac{\Omega_{p}}{\Omega_{\calL}}\right)^{4n} e_{p}(\nu) \dfrac{\Gamma(n) \Gamma(n+1)}{(\sqrt{d}/2)^{2n}\cdot (4\pi)^{2n+1}} L(f_{E}/\calL, \nu, 1) \phi(\frN^{(p)})^{-1}.
\]
Here we have used $ [\calO_{\calL}^{\times}:\ZZ] = 1$, and we have chosen $\delta \in \calL$ so that $\im(\delta) = \sqrt{D_{\calL}}/2$. Using $\Omega_{\infty} = 2\pi i \Omega_{\calL}$, it can be rewritten as
\[
\Theta_{\ttH}(E/\calL)(\hphi)^{2} = \left( \dfrac{\Omega_{p}}{\Omega_{\infty}}\right)^{4n} e_{p}(\nu) \dfrac{(2\pi)^{2n-1}\Gamma(n) \Gamma(n+1)}{2d^{n}} L(f_{E}/\calL, \nu, 1) \phi(\frN^{(p)})^{-1}.
\]
Moreover, since $\phi(\frN^{(p)})^{-1} = \nu(\frN^{(p)})^{-1} \eta(\frN^{(p)}) = \eta(\frN^{(p)})$ as $\nu$ is allowable, multiplication by the term $c_{\ttH}$ now gives
\[
\calL_{\BDP}(E/\calL)(\hnu) = \left( \dfrac{\Omega_{p}}{\Omega_{\infty}}\right)^{4n} e_{p}(\nu) \dfrac{(2\pi)^{2n-1}\Gamma(n) \Gamma(n+1) \nu(\frN^{(p)})^{-1}}{4d^{n-\frac{1}{2}}} L(f_{E}/\calL, \nu, 1).
\]
This is exactly the interpolation formula in the theorem. Note that the term $c_{\ttH}$ is a $p$-adic unit in $\calO^{\ur}$.
\end{proof}

We are now ready to state the comparison theorem.
\begin{theorem}[Comparison theorem] \label{thm:comparison}
    Let $\pi_{\ac}: \Lambda_{\calM}^{\ur} \rightarrow \Lambda_{\calL}^{\ac,\ur}$ be the natural projection map. Then we have
    \[
    \pi_{\ac} \left( \calL_{\Katz}(\lambda_{0}/\calM)\right) = U \cdot \calL_{\BDP}(E/\calL)
    \]
    for an explicitly described unit $U \in (\Lambda_{\calL}^{\ac,\ur})^{\times}$.
\end{theorem}

The proof will be presented in the next section.

\subsection{Comparison of $p$-adic $L$-functions}
The proof of Theorem \ref{thm:comparison} proceeds by carefully comparing the interpolation formulas in Theorem \ref{thm:katz} and Theorem \ref{thm:BDP} term by term. Here we provide some details. Throughout this section, all characters $\nu$ that appear are \emph{allowable} Hecke characters over $\calL$.

\subsubsection{Comparison of Euler factors at $p$}
In this subsection, we prove the following lemma.
\begin{lemma} \label{lem:comparison-at-p}
We have the equality
    \[
    \Eul_p(\lambda_\nu) L_{p}(\lambda_{\nu}) = \dfrac{\lambda_{\nu,v}(2\vartheta_{v})}{\epsilon(0, \lambda_{\nu, v}) L(1, \lambda_{\nu,v}^{-1}) L(0, \lambda_{\nu, \barv})} = u_p e_{p}(\nu) 
    \]
for a unit $u_{p} = \psi_{p}(-2\vartheta_{v}) \in \calO^{\times}$, which is independent of $\nu$, and the factor $e_{p}(\nu)$ is defined in \eqref{eq:epnu-bdp}. 
\end{lemma}
\begin{proof}
    In the interpolation formula \eqref{eq:interpolation-katz}, removing the $p$-factors in $L(\lambda_{\nu}, 0)$ means dividing by $L(0, \lambda_{\nu, v}) L(0, \lambda_{\nu, \barv})$. Multiplying the definition of $\Eul_{p}(\lambda_{\nu})$ gives the first equality. We now prove the second. We write
    \[
    \alpha := \lambda_{\nu,v} = \xi (\nu_{v_0} \circ \Norm_{\calK_p/\QQ_p}) \abs{-}_{\calK_p}^{1/2}, \quad \beta := \lambda_{\nu,\barv} = \xi (\nu_{\barv_0} \circ \Norm_{\calK_p/\QQ_p}) \abs{-}_{\calK_p}^{1/2}.
    \]
    Since $\psi$ is the Hecke character attached to the elliptic curve $E$ with complex multiplication, we have $\xi^{c} = \xi^{-1}$. Since $\nu$ is anticyclotomic, we have $\nu_{\barv_0} = \nu_{v_0}^{-1}$ as characters of $\QQ_p^{\times}$. Consequently, we have $\beta^{c} = \alpha^{-1} \abs{-}_{\calK_{p}}$ and hence $L(0, \beta^{c}) = L(1, \alpha^{-1})$. Therefore,
    \begin{equation} \label{eq:euler-factor-p-L-factor}
    L(1,\alpha^{-1}) = L(0, \beta) = L(1/2, \Pi_p \otimes \nu_{\barv_0}).
    \end{equation}
    Here, we are using the observations that
    \begin{itemize}
        \item $L(0, \beta) = L(0, \beta^{c})$. Indeed, if $\beta$ is ramified, then both local factors are $1$. If $\beta$ is unramified, we pick the uniformizer $p \in \QQ_p \subseteq \calK_{p}$ which is fixed by $c$, so $\beta^{c}(p) = \beta(p)$.
        \item $L(1, \beta) = L(1/2, \Pi_p \otimes \nu_{\barv_0})$. This follows from the compatibility of the automorphic induction with twisting:
        \begin{align}
        \Pi_{p} \otimes \nu_{\barv_0} &= \AI_{\calK_{p}/\QQ_p}(\xi) \otimes \nu_{\barv_0} \notag \\ &= \AI_{\calK_{p}/\QQ_p}(\xi(\nu_{\barv_0} \circ \Norm_{\calK_p/\QQ_p})) = \AI_{\calK_{p}/\QQ_p}(\beta \cdot  \abs{-}_{\calK_p}^{-1/2}). \label{eq:induction-beta}
        \end{align}
    \end{itemize}
Note that the equality \eqref{eq:euler-factor-p-L-factor} holds trivially when one of the three $L$-factors is equal to $1$.

    We next deal with the $\epsilon$-factor. First recall that Tate's local functional equation gives $\epsilon(s, \alpha, \Psi_{\calK_{p}}) \epsilon(1-s, \alpha^{-1}, \Psi_{\calK_{p}}) = \alpha(-1)$. Taking $s=1$, it yields 
    \begin{equation} \label{eq:tate-functional-equation}
        \epsilon(0, \alpha, \Psi_{\calK_{p}}) \epsilon(1, \alpha^{-1}, \Psi_{\calK_{p}}) = \alpha(-1)
    \end{equation}
    Moreover, we have
    \begin{equation} \label{eq:epsilon-ai}
    \epsilon(1, \alpha^{-1}, \Psi_{\calK_{p}}) = \epsilon(0, \beta^{c}, \Psi_{\calK_p}) = \epsilon(0,\beta, \Psi_{\calK_p} \circ c^{-1}) = \epsilon(0,\beta, \Psi_{\calK_p}) = \epsilon(1/2, \Pi_{p} \otimes \nu_{\barv_0}, \Psi_{p}).
    \end{equation}
    Indeed,
    \begin{itemize}
        \item The second equality is given by the change of variables $x \mapsto c(x)$ in the Tate's local zeta integrals appearing in the functional equation involving $\epsilon(0, \beta^{c}, \Psi_{\calK_p})$.
        \item The third equality holds, since $\Psi_{\calK_p}$ is invariant under the complex conjugation.
        \item The last equation is obtained in the following way. By \eqref{eq:induction-beta}, applying the induction formula for epsilon factors, we have
        \[
        \epsilon(1/2, \Pi_p \otimes \nu_{\barv_0}, \Psi_p) = \lambda_{\calK_p/\QQ_p}(\Psi_p) \epsilon(0, \beta, \Psi_{\calK_p}),
        \]
        where $\lambda_{\calK_p/\QQ_p}(\Psi_p) := {\epsilon\left(s, \Ind_{\rmW_{\calK_p}}^{\rmW_{\QQ_p}} \bfone, \Psi_p \right)}/{\epsilon_{\calK_p}(s, \bfone, \psi_{\calK_p})}$ with $\rmW_{L}$ the Weil group of the local field $L$. The constant $\lambda_{\calK_p/\QQ_p}(\Psi_p)$ is independent of $s \in \CC$. Note that $\Ind_{\rmW_{\calK_p}}^{\rmW_{\QQ_p}} \bfone = \bfone \oplus \eta_{\calK_p/\QQ_p}$ where $\eta_{\calK_p/\QQ_p}$ is the unramified quadratic character attached to $\calK_p/\QQ_p$, we have
        \[
        \epsilon\left(s, \Ind_{\rmW_{\calK_p}}^{\rmW_{\QQ_p}} \bfone, \Psi_p \right) = \epsilon\left(s, \bfone, \Psi_p \right) \epsilon\left(s, \eta_{\calK_{p}/\QQ_p}, \Psi_p \right) = 1
        \]
        since both $\bfone$ and $\eta_{\calK_{p}/\QQ_p}$ are unramified, and our choice of $\Psi_{p}$ is such that it is unramified as well.
    \end{itemize}
Now combining \eqref{eq:tate-functional-equation} and \eqref{eq:epsilon-ai}, we have
    \begin{equation} \label{eq:euler-factor-p-eps-factor}
        \epsilon(0, \alpha, \Psi_{\calK_{p}}) = \alpha(-1) \epsilon(1/2, \Pi_{p} \otimes \nu_{\barv_0}, \Psi_{p})^{-1}.    
    \end{equation}
    Combining \eqref{eq:euler-factor-p-L-factor} and \eqref{eq:euler-factor-p-eps-factor}, we have
    \[
    \Eul_p(\lambda_\nu) L_{p}(\lambda_{\nu}) = \alpha(-2\vartheta_{v}) \dfrac{\epsilon(1/2, \Pi_{p} \otimes \nu_{\barv_0}, \Psi_{p})}{L(1/2, \Pi_p \otimes \nu_{\barv_0})^{2}}.
    \]

    Finally, the finite unramified character $\nu_{v_0}$ and the character $\abs{-}_{\calK_{p}}^{1/2}$ take the value $1$ on $-2\vartheta_{v}$ since it is a unit in $\calO$. Therefore, $\alpha(-2\vartheta_{v}) = \psi_{p}(-2\vartheta_{v})$. This proves the lemma.
\end{proof}

\subsubsection{Comparison of Euler factors away from $p$}
Now we compare the Euler factors away from $p$. We start with the finite places.
\begin{lemma} \label{lem:comparsion-away-from-p}
    We have the following equality
    \[
    \Eul_{\frc}(\lambda_{\nu})L_{\frc}(\lambda_{\nu}, 0) = \prod_{w \mid \frf} c_{w} \hnu(g_{w})
    \]
    where $c_{w} \in (\calO^{\ur})^{\times}$ is a fixed unit independent of $\nu$, and $g_{w}$ is defined in \eqref{eq:define-g_w} in the proof.
\end{lemma}

\begin{proof}
    We fix a place $w \mid \frf$ of $\calM$ and put $\alpha_{w} := \lambda_{0,w}$ for simplicity. Let $a_{w}$ be the conductor exponent of $\alpha_{w}$. Let $d_{w}$ be the exponent for which the additive character $\Psi_{w}$ is trivial on $\varpi_{w}^{-d_{w}}\calO_w$ but not on $\varpi_{w}^{-d_{w}-1}\calO_{w}$. We put $b_{w} := \varpi_{w}^{a_{w}+d_{w}}$. For any unramified character $\beta_{w}$, \cite[Equation (3.2.6.3)]{Tate} gives
    \[
    \epsilon(0,\alpha_{w}\beta_{w}, \Psi_{w}) = \beta_{w}(b_{w}) \epsilon(0, \alpha_{w}, \Psi_{w}).
    \]
    We now take $\alpha_{w} = \lambda_{w}$ and $\beta_{w} = \tnu_{w}$. Since $\alpha_{w}$ is ramified and $\beta_{w}$ is unramified, both $L$-factors $L(0, \lambda_{\nu, w})$ and $L(1, \lambda_{\nu,w}^{-1})$ in the definition of $\Eul_{w}(\lambda_{w})$ are $1$. Therefore
    \[
    \Eul_{w}(\lambda_{\nu}) = \dfrac{\alpha_w(2\vartheta_{w})}{\epsilon(0, \alpha_{w}, \Psi_{w})} \tnu_{w}(2\vartheta_{w}/b_{w}) =: c_{w} \hnu(g_{w}),
    \]
    where 
    \begin{equation} \label{eq:define-g_w}
    c_{w} := \dfrac{\alpha_w(2\vartheta_{w})}{\epsilon(0, \alpha_{w}, \Psi_{w})}, \quad g_{w} := \rec_{\calL}(\Norm_{\calM/\calL} j_{w}(2\vartheta_{w}/b_{w}))|_{\Gamma_{\calL}^{\ac}}
    \end{equation}
    where $j_{w}(z)$ denotes the finite idele of $\calM$ equal to $z$ at $w$ and $1$ elsewhere. Since it has trivial $p$-component, passing to $p$-adic avatar of $\nu$ will introduce no additional factors.

    We show that $c_{w} \in (\calO^{\ur})^{\times}$. Indeed, the numerator $\alpha_w(2\vartheta_{w})$ is a $p$-adic unit and is independent of $\nu$ since $\nu_{w}$ is unramified. For the $\epsilon$-factor appearing in the denominator, \cite[Equation (3.2.6.2)]{Tate} simplifies to
    \[
    \epsilon(0, \alpha_{w}, \Psi_{w}) = \alpha(\varpi_{w}^{b_{w}}) q_{w}^{d_{w}} \tau(\alpha_w^{-1}, \Psi_{w})
    \]
    where 
    \[
    \tau(\alpha_{w}^{-1}, \Psi_{w}) := \sum_{\baru \in \calO^{\times}/(1+\varpi^{b_{w}}\calO)} \alpha_{w}^{-1}(u) \Psi(u/\varpi_{w}^{b_{w}})
    \]
    is the local Gauss sum. Local Gauss sums satisfy the following product formula
    \[
    \tau(\alpha_{w}^{-1}, \Psi_{w}) \tau(\alpha_{w}, \Psi_{w}) = \alpha_{w}(-1)q_{w}^{a_{w}}.
    \]
    This implies that
    \[
    v_{p}(\tau(\alpha_{w}^{-1}, \Psi_{w})) + v_{p}(\tau(\alpha_{w}, \Psi_{w})) = v_{p}(\tau(\alpha_{w}^{-1}, \Psi_{w}) \tau(\alpha_{w}, \Psi_{w})) = v_p(\alpha_{w}(-1)q_{w}^{a_{w}}) = 0.
    \]
    Moreover, since both $\tau(\alpha_{w}^{-1}, \Psi_{w})$ and $\tau(\alpha_{w}, \Psi_{w})$ are algebraic integers, we have
    \[
    v_{p}(\tau(\alpha_{w}^{-1}, \Psi_{w})) \geq 0, \quad v_{p}(\tau(\alpha_{w}, \Psi_{w})) \geq 0.
    \]
    This implies that both $\tau(\alpha_{w}^{-1}, \Psi_{w})$ and $\tau(\alpha_{w}, \Psi_{w})$ are $p$-adic units. This shows that $c_{w} \in (\calO^{\ur})^{\times}$ as desired.
\end{proof}

We now come to the archimedean places. We define the following term at the archimedean places, informally, as the "Euler factor at infinity":
\[
\Eul_{\infty}(\lambda_{\nu}) := \left(\dfrac{\Omega_{p}}{\Omega_{\infty}}\right)^{k\Sigma+2\kappa} \dfrac{\pi^{\kappa} \Gamma_{\Sigma}(k\Sigma+\kappa)[\calO_{\calM}^{\times}:\calO_{\calF}^{\times}]}{\sqrt{\abs{D_{\calF}}_{\RR}} (\im \vartheta)^{\kappa}} t_{\calM}.
\]
With the parameter $\vartheta$ chosen in the proof of Theorem \ref{thm:katz}, we have the following simplification.
\begin{lemma} \label{lem:comparison-at-archemedian}
We have the following equality:
   \[
   \Eul_{\infty}(\lambda_{\nu}) =  \left(\dfrac{\Omega_{p}}{\Omega_{\infty}}\right)^{4n} \dfrac{\pi^{2n-1} \Gamma(n)\Gamma(n+1)}{\Norm_{\calF/\QQ}(a)^{n} \cdot d^{n-\frac{1}{2}}} t_{\calM}^{\prime\prime} = \left(\dfrac{\Omega_{p}}{\Omega_{\infty}}\right)^{4n} \dfrac{(2\pi)^{2n-1} \Gamma(n)\Gamma(n+1)}{d^{n-\frac{1}{2}}} t_{\calM}^{\prime\prime} \cdot \hnu(g_{v_0}^{-1}),
   \]
where $g_{v_0}$ is defined in \eqref{eq:define-g_v0}, and $t_{\calM}^{\prime\prime}$ is a $p$-adic unit in $\calO^{\times}$ independent of $\nu$.
\end{lemma}

\begin{proof}
We have computed the infinite type of $\lambda_{\nu}$ in the last row of Table \ref{tab:infinite-type}. Following the notation in Theorem \ref{thm:katz}, we have $k=1$ and $\kappa = (n-1)\sigma_{1}+n\sigma_{2}$. For all $n \geq 1$, it lies in the interpolation range.
Indeed,
\begin{itemize}
    \item $\Gamma_{\Sigma}(k\Sigma+\kappa) = \Gamma(n)\Gamma(n+1)$,
    \item $(\Omega_{p}/\Omega_{\infty})^{k\Sigma+2\kappa} = (\Omega_{p}/\Omega_{\infty})^{4n}$, under the choice of \eqref{eq:comparison-period}, and
    \item $\pi^{\kappa} = \pi^{2n-1}$,
\end{itemize}
Under the choice of $\vartheta$ in the proof of Theorem \ref{thm:BDP}, we have $\im \sigma_i(\vartheta) = a_i \sqrt{D_{\calL}}/2$. This gives
\begin{itemize}
    \item $(\im \vartheta)^{\kappa} = (a_1 \sqrt{D_{\calL}}/2)^{n-1} (a_2 \sqrt{D_{\calL}}/2)^{n}$.
\end{itemize}
We now consider $t_{\calM}[\calO_{\calM}^{\times}:\calO_{\calF}^{\times}]/\sqrt{D_{\calF}}$. We have seen that $t_{\calM}$ is a power of $2$. Since $p$ is unramified in $\calF$, the discriminant $\sqrt{D_{\calF}}$ is a $p$-adic unit. For $[\calO_{\calM}^{\times}:\calO_{\calF}^{\times}]$, one computes that
    \[
    [\calO_{\calM}^{\times}:\calO_{\calF}^{\times}] = [\calO_{\calM}^{\times}: \mu(\calO_{\calM}) \calO_{\calF}^{\times}] [\mu(\calO_{\calM}): \mu(\calO_{\calF})] = [\calO_{\calM}^{\times}: \mu(\calO_{\calM}) \calO_{\calF}^{\times}] \cdot \# \mu(\calO_{\calM})
    \]
    where $\mu(\calO)$ denotes the group of roots of unity in $\calO$. By \cite[Theorem 4.12]{washington}, we see that $[\calO_{\calM}^{\times}: \mu(\calO_{\calM}) \calO_{\calF}^{\times}] \in \{1, 2\}$, and the second equality follows from the fact that $\calF$ is a real quadratic field. Moreover, $\mu(\calO_{\calM})$ is generated by $\mu(\calO_{\calL})$ and $\mu(\calO_{\calK})$. Since we have assumed that $D_{\calL}$ is odd and smaller than $3$, this implies that $\mu(\calO_{\calL}) = \{\pm 1\}$ and hence $\# \mu(\calO_{\calM})$ is not divisible by any prime number $p \geq 5$. To sum up, $[\calO_{\calM}^{\times}:\calO_{\calF}^{\times}]$ is not divisible by $p$.
    
    To conclude, we have
    \[
    \Eul_{\infty}(\lambda_{\nu}) = \left(\dfrac{\Omega_{p}}{\Omega_{\infty}}\right)^{4n} \dfrac{\pi^{2n-1} \Gamma(n)\Gamma(n+1)}{(a_1 \sqrt{D_{\calL}}/2)^{n-1} (a_2 \sqrt{D_{\calL}}/2)^{n}} t_{\calM}^{\prime} = \left(\dfrac{\Omega_{p}}{\Omega_{\infty}}\right)^{4n} \dfrac{(2\pi)^{2n-1} \Gamma(n)\Gamma(n+1)}{\Norm_{\calF/\QQ}(a)^{n} \cdot d^{n-\frac{1}{2}}} t_{\calM}^{\prime\prime}
    \]
    where $t_{\calM}^{\prime}$ and $t_{\calM}^{\prime\prime} := a_1 t_{\calM}^{\prime}$ are fixed $p$-adic units independent of $\nu$. This proves the first equality of the lemma.

    Moreover, write $u := \Norm_{\calF/\QQ}(a) \in \QQ_{> 0} \cap \ZZ_p^{\times}$, and let $j_{v_0}(u)$ be the idèle of $\calL$ that is $u$ at $v_0$ and $1$ elsewhere. Define
    \begin{equation} \label{eq:define-g_v0}
    g_{v_0} := \rec_{\calL}(j_{v_0}(u))|_{\Gamma_{\calL}^{\ac}} = \rec_{\calL}(j_{v_0}(\Norm_{\calF/\QQ}(a)))|_{\Gamma_{\calL}^{\ac}}.
    \end{equation}
    It follows by the definition of the $p$-adic avatar $\hnu$ of $\nu$ that $\hnu(g_{v_0}) = \Norm_{\calF/\QQ}(a)^{n}$. This gives the final equality of the lemma.
\end{proof}

\subsubsection{The proof of the comparison theorem}
We are now ready to prove the comparison theorem for $p$-adic $L$-functions.

\begin{proof}[Proof of Theorem \ref{thm:comparison}]
Combining Lemmas \ref{lem:comparsion-away-from-p}, \ref{lem:comparison-at-p} and \ref{lem:comparison-at-archemedian} as well as \eqref{eq:compare-L-value}, we see that
\begin{multline} \label{eq:project-Hecke-formula}
\left( \pi_{\ac}\left( \calL_{\Katz}(\lambda_{0}/\calM) \right) \right)(\hnu) \\ = t_{\calM}^{\prime\prime} \left(\dfrac{\Omega_{p}}{\Omega_{\infty}}\right)^{4n} \dfrac{(2\pi)^{2n-1} \Gamma(n)\Gamma(n+1)}{d^{n-\frac{1}{2}}} \cdot \hnu \left(g_{v_0}^{-1} \prod_{w \mid \frf} g_{w}\right) \left( u_{p} \prod_{w \mid \frf} c_{w} \right) e_{p}(\nu)  L(f_E/\calL, \nu, 1).    
\end{multline}

We write $g_{\frN^{(p)}}$ for the image in $\Lambda_{\calL}^{\ac}$ under the reciprocity map $\rec_{\calL}$ of a finite idele of $\calL$ for $\frN^{(p)}$ with trivial $p$-components. Then $\nu((\frN^{(p)})^{-1}) = \hnu(g_{\frN^{(p)}})^{-1}$. We define the following element in $\Lambda_{\calL}^{\ac,\ur}$:
\begin{equation} \label{eq:define-U}
    U := 4 t_{\calM}^{\prime\prime} \left( u_{p} \prod_{w \mid \frf} c_{w} \right) \cdot [g_{v_0}^{-1}] \prod_{w \mid \frf} [g_{w}] [g_{\frN^{(p)}}], 
\end{equation}
where $[-]$ represents the corresponding element in $\Lambda_{\calL}^{\ac}$, all of them are units. The remaining factors are fixed $p$-adic units in $(\calO^{\ur})^{\times}$. This shows that $U \in (\Lambda_{\calL}^{\ac,\ur})^{\times}$. Comparing \eqref{eq:project-Hecke-formula} with the interpolation formula in Theorem \ref{thm:BDP}, we see that
\[
\left( \pi_{\ac}\left( \calL_{\Katz}(\lambda_{0}/\calM) \right) \right)(\hnu) = \hnu(U) \cdot \calL_{\BDP}(E/\calL)(\hnu).
\]
Lemma \ref{lem:allowables} provides infinitely many allowable Hecke characters $\nu$ over $\calL$; therefore, by the Weierstrass preparation theorem, Theorem \ref{thm:comparison} is proved.
\end{proof}

%% file: 204-main-conjecture.tex
\section{The Bertolini--Darmon--Prasanna main conjecture}
In this section, we prove the one-side divisibility of the BDP main conjecture for $E$ over $\calL$, as Theorem \ref{thm:BDP-one-divisibility}.

We recall the following CM Iwasawa main conjecture for $\calM/\calF$, which was independently established in the third-named author's doctoral thesis \cite[version 3]{Xu} and in \cite{ABL}. \footnote{Both works independently strengthen \cite[Theorem C]{hsieh2014eisenstein}. The authors of these two works plan to combine them into a single joint article, which will be made available.}

\begin{theorem}[CM main conjecture] \label{thm:CM-main-conjecture}
Let $p \geq 5$ be a prime. Then $X_{\Gr}(\bfT_{\rho}/\calM)$ is a finitely generated torsion $\Lambda_{\calM}$-module, and we have
\[
\Char_{\Lambda_{\calM}}(X_{\Gr}(\bfT_{\rho}/\calM)) \Lambda_{\calM}^{\ur} \subseteq (\calL_{\Katz}(\lambda_{0}/\calM))
\]
as ideals in $\Lambda_{\calM}^{\ur}$.
\end{theorem}

Next, by the control theorem in Section \ref{sec:Selmer-group} and the comparison of $p$-adic $L$-functions in Section \ref{sec:comparison}, we can deduce the corresponding one-side divisibility of the BDP main conjecture. The main algebraic tool is the following lemma, which is implicit in the proof of \cite[Corollary 3.4.2]{JSW}.

\begin{lemma}[Specialization lemma]\label{lem:specialization}
Let $\Lambda$ and $\Lambda^{\prime}$ be power-series rings over a complete discrete valuation ring, let $\pi: \Lambda \twoheadrightarrow \Lambda^{\prime}$ be a surjective homomorphism, and let $X$ be a finitely generated torsion $\Lambda$-module. Suppose that $Y := X \otimes_{\Lambda} \Lambda^{\prime}$ is torsion over $\Lambda^{\prime}$. If $\calL \in \Lambda$ satisfies $\Char_{\Lambda}(X) \subset (\calL)$, then $\Char_{\Lambda^{\prime}}(Y) \subset (\pi(\calL))$.
\end{lemma}
\begin{proof}
Write $\Fitt_{\Lambda}(X)$ for the zeroth fitting ideal. Since $X$ is a finitely generated torsion $\Lambda$-module, we have
\begin{equation}\label{eq:p3-fitting-inclusion}
 \Fitt_{\Lambda}(X)\subseteq\Char_{\Lambda}(X)\subseteq(\calL).
\end{equation}
Fitting ideals commute with arbitrary base change; hence, applying \eqref{eq:p3-fitting-inclusion}, we have
\begin{equation} \label{eq:fitting-inclusion}
 \Fitt_{\Lambda^{\prime}}(Y) = \Fitt_{\Lambda^{\prime}}(X \otimes_{\Lambda} \Lambda^{\prime}) = \Fitt_{\Lambda}(X) \Lambda^{\prime} = \pi(\Fitt_{\Lambda}(X)) \Lambda^{\prime} \subseteq(\pi(\calL)).
\end{equation}
We claim that $\pi(\calL) \neq 0$. Indeed, since $Y$ is a finitely generated torsion $\Lambda^{\prime}$-module, there exists a nonzero element $a \in \Lambda^{\prime}$ such that $aY = 0$. Choose a surjection $(\Lambda^{\prime})^{n} \twoheadrightarrow Y$ and a $\Lambda^{\prime}$-basis $\{\tte_{i}\}$ of $(\Lambda^{\prime})^{n}$. Because $aY=0$, each vector $a \tte_{i}$ belongs to the kernel of this surjection. Thus, one may include the columns of the diagonal matrix $a \bfone_{n}$ in a presentation matrix of $Y$. Its determinant $a^{n}$ is then one of its $(n \times n)$-minors. Thus, $0 \neq a^{n} \in \Fitt_{\Lambda^{\prime}}(Y)$, and therefore $\Fitt_{\Lambda^{\prime}}(Y) \neq 0$. It then follows from \eqref{eq:fitting-inclusion} that $\pi(\calL) \neq 0$.

Let $\wp$ be any height-one prime of $\Lambda^{\prime}$. Localizing \eqref{eq:fitting-inclusion} at $\wp$ gives 
\[
\Fitt_{\Lambda^{\prime}_{\wp}}(Y_{\wp}) = \left( \varpi_{\wp}^{\length_{\Lambda_{\wp}^{\prime}}(Y_{\wp})} \right) \subseteq \pi(\calL) \Lambda^{\prime}_{\wp},
\]
and hence $\length_{\Lambda_{\wp}^{\prime}}(Y_{\wp}) \geq v_{\wp}(\pi(\calL))$. If $f_{Y}$ is a generator of $\Char_{\Lambda^{\prime}}(Y)$, then $v_{\wp}(f_{Y}) = \length_{\Lambda_{\wp}^{\prime}}(Y_{\wp})$. Combining the two, since $\Lambda^{\prime}$ is a unique factorization domain, this implies that $\pi(\calL) \mid f_{Y}$, and hence $\Char_{\Lambda^{\prime}}(Y) \subset \pi(\calL)$, as desired.
\end{proof}

\begin{remark}
    The specialization lemma for the other side of the divisibility may not hold in general. For the special case in which $\Lambda^{\prime} := \Lambda/(T)$, where $T$ is one indeterminate of $\Lambda$ and $\pi: \Lambda \twoheadrightarrow \Lambda^{\prime}$ is the natural projection, continue to assume that $Y$ is a finitely generated torsion $\Lambda^{\prime}$-module. A sufficient condition to guarantee that
    \[
    \Char_{\Lambda}(X) \supset (\calL) \Rightarrow \Char_{\Lambda^{\prime}}(Y) \supset (\pi(\calL))
    \]
    is that $X$ has no nonzero pseudonull $\Lambda$-submodule.
\end{remark}

Now we are ready to prove the first main theorem of this part. We denote $X_{\BDP}(E/\calL) := X_{\BDP}(\bfT_{E}^{\ac}/\calL)$ for simplicity.
\begin{theorem}[The BDP main conjecture] \label{thm:BDP-one-divisibility}
    Let $p \geq 5$ be a prime and suppose that \eqref{eq:rank-one} is satisfied. Then we have the following divisibility
    \[
    \Char_{\Lambda_{\calL}^{\ac}}(X_{\BDP}(E/\calL)) \Lambda_{\calL}^{\ac} \subseteq (\calL_{\BDP}(E/\calL))
    \]
    as ideals in $\Lambda_{\calL}^{\ac,\ur}$.
\end{theorem}

\begin{proof}
    On the algebraic side, by the descent theorem (Proposition \ref{prop:descent-theorem}) and Shapiro's lemma (Proposition \ref{prop:shapiro}), we see that
    \[
    X_{\BDP}(E/\calL) \simeq X_{\Gr}(\bfT_{\rho}/\calM) \otimes_{\Lambda_{\calM}} \Lambda_{\calL}^{\ac}.
    \]
    Combining Theorem \ref{thm:comparison} and applying Lemma \ref{lem:specialization}, while noting that $X_{\BDP}(E/\calL)$ is a finitely generated torsion $\Lambda_{\calL}^{\ac}$-module by the control theorem (Proposition \ref{prop:control-theorem-BDP}), proves the theorem.
\end{proof}

We have the following corollary of the BDP main conjecture.

\begin{theorem} \label{thm:L_BDP_1}
    Let $p \geq 5$ be a prime and suppose that \eqref{eq:rank-one} holds. Then $\calL_{\BDP}(E/\calL)(1) \neq 0$ and 
    \[
    v_{p}(\calL_{\BDP}(E/\calL)(1)) \leq v_{p}(\# \Sha(E/\calL)[p^{\infty}]) + 2r.
    \]
\end{theorem}
\begin{proof}
    This is a combination of Theorem \ref{thm:BDP-one-divisibility} and Proposition \ref{prop:control-theorem-BDP}(2) using the specialization lemma (Lemma \ref{lem:specialization}). 

    More precisely, consider the quotient $Y := X_{\BDP}(\bfT_{E}^{\ac}/\calL)/(\gamma_{\ac}-1)X_{\BDP}(\bfT_{E}^{\ac}/\calL)$. Then
    \[
    \Fitt_{\calO^{\ur}}(Y \otimes_{\calO} \calO^{\ur}) \subseteq \Char_{\calO^{\ur}}(Y \otimes_{\calO} \calO^{\ur}) \subseteq (\calL_{\BDP}(E/\calL)(1)).
    \]
    Here the first inclusion is \eqref{eq:p3-fitting-inclusion}, and the second inclusion is obtained by applying Lemma \ref{lem:specialization} to Theorem \ref{thm:BDP-one-divisibility} in the case $\Lambda := \Lambda_{\calL}^{\ur}$ and $\Lambda^{\prime} = \calO^{\ur}$, with the projection map given by quotienting out $\gamma_{\ac}-1$. Using the same argument as in the proof of Lemma \ref{lem:specialization}, the finiteness of $Y$ implies that $\calL_{\BDP}(E/\calL)(1) \neq 0$. The inequality in the theorem then follows from Proposition \ref{prop:control-theorem-BDP}(2) by noting that, for any finitely generated torsion $\calO$-module $Y$, we have $\Fitt_{\calO}(Y) = (\varpi^{\length_{\calO}(Y)})$.
\end{proof}

%% file: 205-Waldspurger.tex
\section{A $p$-adic Waldspurger formula}
In the following two sections, we finish the proof of the $p$-part of the BSD formula for $E/\QQ$, with $p \geq 5$ inert in $\calK$. Starting from Theorem \ref{thm:L_BDP_1}, the missing components are the $p$-adic Waldspurger formula and the Gross--Zagier formula. We establish the $p$-adic Waldspurger formula in this section.

\subsection{Heegner points} \label{sec:Heegner-point}

We start with some basic setup for Heegner points. Let $f_{E}$ be a \emph{normalized} cuspidal newform attached to $E$. We choose and fix a modular parameterization $\phi: X_{0}(N) \rightarrow E$ over $\QQ$ with $\phi(\infty) = 0$. For a fixed Néron differential $\omega_{E}$ of $E$, we write
\[
\phi^{\ast} \omega_{E} = c_{\phi} \omega_{f}
\]
for $\omega_{f} = 2 \pi i f(z) \dif z$ and some nonzero rational number $c_{\phi} \in \QQ$, called the \emph{Manin constant} for the modular parametrization $\phi: X_{0}(N) \rightarrow E$. Let $H_{\calL}$ be the Hilbert class field of $\calL$ and put $\frN = \frN^{(p)}(\varpi_{v_0})^{v_p(N)}$. We let $\sfx \in X_{0}(N)(H_{\calL})$ be the CM point represented by $\CC/\calO_{\calL} \rightarrow \CC/\frN^{-1}$. Define the \emph{Heegner point} of $E$ over $\calL$ by
\[
P := \Tr_{H_{\calL}/\calL} \phi(\sfx) \in E(\calL),
\]
and write $m := [E(\calL)/E(\calL)_{\tor}: \ZZ \barP]$. Under the stronger assumption \eqref{eq:GZK}, it follows from the Gross--Zagier formula that $P$ is not a torsion point in $E(\calL)$, and hence $m > 0$ is a positive integer.

For any prime-to-$pN$ ideal $\fra$ of $\calL$, we write
\[
\sfx_{\fra}: \CC/\fra^{-1} \rightarrow \CC/(\fra^{-1}\frN^{-1})
\]
for the associated point on $X_{0}(N)$. Its kernel is the cyclic group $\fra^{-1} \frN^{-1}/\fra^{-1}$ of order $N$. The CM reciprocity law then gives
\[
P = \sum_{[\frb] \in \Cl(\calL)} \phi(\sfx_{\frb}) \in E(\calL),
\]
where $\Cl(\calL)$ is the class group of $\calL$. Moreover, let $\Cl_{\ram}(\calL)$ be the subgroup of $\Cl(\calL)$ generated by the ideal classes ramified in $\calL$. Denote $s := \# \Cl_{\ram}(\calL)$. Let $D$ be a complete set of representatives of $\Cl(\calL)/\Cl_{\ram}(\calL)$. Then we can write
\[
P = \sum_{[\frb] \in \Cl(\calL)} \phi(\sfx_{\frb}) = \sum_{(\fra, c) \in D \times \Cl_{\ram}(\calL)} \phi(\sfx_{\fra c}),
\]
noting that the map $D \times \Cl_{\ram}(\calL) \rightarrow \Cl(\calL)$ given by $(\fra, c) = [\fra]c$ is bijective.

\subsection{$p$-adic Waldspurger formula}
We have the following result.

\begin{theorem}[$p$-adic Waldspurger formula] \label{thm:waldspurger}
    We put 
    \[
    \scE_{p} = \begin{cases}
        1 - \frac{a_p}{p} + \frac{1}{p} = \frac{p+1}{p}, \quad& \text{$E$ has good reduction at $p$, i.e. $p \nmid N$}, \\
        1, \quad& \text{$E$ has bad reduction at $p$, i.e. $p \mid N$}. \\
    \end{cases}
    \]
    Then
    \begin{equation} \label{eq:bdp-formula}
    \calL_{\BDP}(E/\calL)(1) = u_{\rmW} \left( \dfrac{\scE_{p}}{c_{\phi}} \log_{\omega_{E}}(P)\right)^{2} 
    \end{equation}
    for some $u_{\rmW} \in (\calO^{\ur})^{\times}$.
\end{theorem}

The proof deals with the two cases separately.

\begin{proof}[Proof of Theorem \ref{thm:waldspurger}: the case $p \nmid N$]
    The classical BDP formula \cite[Theorem 5.13]{BDP} (see also the normalization in \cite[Formula (BDP)]{casetal25}) gives
    \begin{equation} \label{eq:BDP-good-case}
    \calL_{\BDP}(E/\calL)(1) = \dfrac{1}{u_{\calL}^{2} c_{\phi}^{2}} \left( \dfrac{1-a_p(E)+p}{p} \right)^{2} \log_{\omega_{E}}(P)^{2}.
    \end{equation}
    In our case, $E$ has good supersingular reduction at $p$ and $p \geq 5$. This implies that $a_{p}(E) = 0$, and $u_{\calL} = \frac{1}{2} \# \calO_{\calL}^{\times} = 1$ since we have required $D_{\calL} < -3$. This gives the good-reduction formula \eqref{eq:bdp-formula} with $u_{\rmW} = 1$. Careful readers may wonder whether the BDP $p$-adic $L$-functions in \cite{BDP, casetal25} are identical to the one normalized in Theorem \ref{thm:BDP}. This indeed holds in the good-reduction case by a direct comparison of interpolation formulas.
\end{proof}

The remainder of this section deals with the bad-reduction case, i.e., the case $p \mid N$. To prepare for the proof, we review the setup in \cite{LZZ}. Following \cite[Notation 3.1.2]{LZZ}, let $\bPi$ be the space of modular parameterizations
\[
\bPi = \varinjlim_{U} \Hom_{\xi_{U}}(X_{U}, E), 
\]
where
\begin{itemize}
    \item the colimit is taken over all open compact subgroups $U$ of $\GL_2(\AA_{\rmf})$. Note that $X_{0}(N)$ corresponds to the level
    \[
    U_{0}(N) := \left\{ g \in \GL_2(\hatZZ): \begin{bmatrix}
        a & b \\ c & d
    \end{bmatrix}: c \equiv 0 \pmod N  \right\}.
    \]
    \item $X_{U}$ is the compact modular curve with level $U$,
    \item $\xi_{U}$ is the normalized Hodge class on $X_{U}$,
    \item $\Hom_{\xi_U}(X_U, E)$ denotes the $\QQ$-vector space of modular parameterizations, i.e., the abelian group of morphisms from $X_U$ to $E$ that send $\xi_{U}$ to a torsion point, tensored with $\QQ$.
\end{itemize}

\begin{remark} \label{rem:identification-LZZ}
In Section \ref{sec:Heegner-point}, we have fixed a modular parametrization $\phi: X_{0}(N) \rightarrow E$. Pulling back the normalized invariant differential $\omega_{E}^{\nor} := c_{\phi}^{-1} \omega_{E}$ along $\phi$, we obtain $\omega_{f_{E}} = \phi^{\ast} \omega_{E}^{\nor}$, where $\omega_{f_{E}} = 2 \pi i f(z)\dif z$. Fixing the normalized invariant differential $\omega_{E}^{\nor}$, we consider the map 
\begin{equation} \label{eq:identification-LZZ}
\calD_{\omega_{E}^{\nor}}: \bPi \rightarrow \varinjlim_{U} \rmH^{0}(X_{U}, \Omega_{X(U)}^{1}), \quad [h] \mapsto h^{\ast}\omega_{E}^{\nor}.
\end{equation}
Our normalization gives $\calD_{\omega_{E}^{\nor}}([\phi]) = \omega_{f_{E}}$. This map identifies $\bPi$ with the finite part of the automorphic representation $\Pi_{\rmf}$ attached to $f_{E}$ (upon the normalization that its first Fourier coefficient $a_{1}(f_{E}) = 1$).

Moreover, let $N$ be the unipotent subgroup of $\GL_2$ as a $\ZZ$-scheme. We also fix an additive character $\psi: \QQ_p \rightarrow \CC_p^{\times}$ of conductor $\ZZ_p$. Then we can consider the Kirillov model $\sfK(\Pi_{p}, \psi)$ with respect to $(N, \psi)$, which is a subspace of smooth $\CC$-valued functions on $\QQ_p^{\times}$.  Then we can realize both $\bPi_{p}$ and $\Pi_{p}$ inside $\sfK(\Pi_{p}, \psi)$, and they share the same image inside $\sfK(\Pi_{p}, \psi)$ via the identification $\calD_{\omega_{E}^{\nor}}$.
\end{remark}

Following \cite[Definition 4.1.1]{LZZ}, we define the stable vectors in $\bPi$.

\begin{definition}
    We say that the element $g \in \bPi \otimes_{\QQ} L$ or $\bPi_{p} \otimes_{\QQ} L$ is a \emph{stable vector} if $g$ is fixed by $N(\ZZ_p)$ and it satisfies the relation
    \[
    \sum_{h \in N(p^{-1})/N(\calO_{p})} \bPi_{p}(h) g = 0.
    \]
    Equivalently, if we realize $\bPi_{p}$ in the Kirillov model $\sfK(\Pi_{p}, \psi)$ with respect to the pair $(N, \psi)$, then $g$ is a stable vector if and only if $g_{p}$ is supported on $\ZZ_p^{\times}$. We denote by $\bPi_{L}^{\hsuit}$ (resp. $\bPi_{L, p}^{\hsuit}$) the subset of $\bPi \otimes_{\QQ} L$ (resp. $\bPi_{L, p}^{\hsuit}$) consisting of stable vectors.
\end{definition}

Next, we review the Mellin transform constructed in \cite[Theorem 2.3.17]{LZZ}. Following \cite[Definition 2.3.8]{LZZ}, we denote the space of \emph{$L$-convergent modular forms of weight $w$ and $p$-Iwahori level $m$} by $\calM^{w}(m,L)$, which is naturally a complete $L$-vector space. Here $m \in \NN \cup \{\infty\}$ and $w \in \ZZ$, and $\frX(m)$ is the formal ordinary model at $p$-Iwahoric level $m$ (see \cite[Section 2.3]{LZZ}). A convergent modular form of weight $0$ is simply called a \emph{convergent modular form}. In \cite[Definition 2.3.10]{LZZ}, the authors also define when a convergent modular form is called \emph{stable}, and denote the space of stable convergent modular forms by $\calM^{w}(m,L)^{\hsuit}$.

For $w \in \NN$, there is an \emph{Atkin--Serre operator} 
\[
\Theta_{\ord}: \calM^{w}(m,L) \rightarrow \calM^{w+1}(m,L)
\]
given in \cite[Equation (2.18)]{LZZ}. Let $\calD(\calO_{p}^{\anti}, \QQ_p)$ be the algebra of locally analytic distributions on $\calO_{p}^{\anti} \simeq \ZZ_p^{\times}$ (since $p$ splits in $\calL$).

\begin{theorem}[{Mellin transform, \cite[Theorem 2.3.17]{LZZ}}] \label{thm:mellin}
    Let $f \in \calM^{w}(m,L)^{\hsuit}$ be a stable convergent modular form for some $w, m \in \NN$. Then there is a unique element
    \[
    \bfM(f) \in \calM^{0}(\infty, L) \hotimes_{\QQ_p} \calD(\calO_{p}^{\anti}, \QQ_p),
    \]
    such that for every $k \in \NN$, 
    \[
    \bfM(f)(\lrangle{w+k}) = (\Theta_{\ord}^{k} f) \omega_{\nu}^{-w-k},
    \]
    where $\lrangle{j}$ denotes the power character $x \mapsto x^{j}$, and $\omega_{\nu}$ is the global Lubin--Tate differential defined in \cite[Definition 2.3.6]{LZZ} with respect to the homomorphism $\nu: \mathrm{LT} \rightarrow \widehat{\GG}_{m}$ specified in \cite[Notation 1.8.3]{LZZ}. Moreover, if $w \geq 1$, then we have
    \begin{equation} \label{eq:theta-relation}
    \Theta_{\ord}(\bfM(f)(\lrangle{w-1})\omega_{\nu}^{w-1}) = f.
    \end{equation}
\end{theorem}
A particular case of Theorem \ref{thm:mellin} occurs when $w = 1$ in \eqref{eq:theta-relation}, which gives $\Theta_{\ord}(\bfM(f)(1)) = f$, where $1$ means the specialization $\lrangle{0}$ at the trivial character.

Now we are ready to state the $p$-adic Waldspurger hidden in \cite[Equation (4.6)]{LZZ} in the proof of \cite[Proposition 4.3.6]{LZZ}.

\begin{theorem}[Liu--Zhang--Zhang's $p$-adic Waldspurger formula] \label{thm:LZZ}
    Let $g$ be a stable vector inside $\bPi_{L}^{\hsuit}$ which can be written as a $L$-linear combination of modular parametrizations $X_{U} \rightarrow E$ over sufficiently small tame level $U^{p}$. Let $\omega$ be an invariant differential on $E$. Then we have
    \[
    \bfM((g^{\ast}\omega)_{\ord})(1) = (\log_{\omega} \circ g)_{\ord}.
    \]
    This identity is restricted to the ordinary locus $X_{U}^{\ord}$ (see \cite[Notation 2.4.8]{LZZ}), and here $1$ means the specialization $\lrangle{0}$ at the trivial character.
\end{theorem}
\begin{proof}
    Note that the original \cite[Proposition 4.3.6]{LZZ} required $g$ to be $n$-admissible, which is stronger than being stable. Nevertheless, one checks in the proof that the deduction of \cite[Equation (4.6)]{LZZ} carries over when $g$ is only a stable vector.
\end{proof}
We remark that for $g = \sum_{i} b_i [g_i] \in \bPi_{L}$, the notation means
\[
g^{\ast} \omega := \sum_{i} b_i g_{i}^{\ast}\omega, \quad \log_{\omega} \circ g = \sum_{i} b_i (\log_{\omega} \circ g_{i}).
\]

We shall see in the remainder of this section that the term $\bfM((g^{\ast}\omega)_{\ord})(1)$ corresponds to the BDP $p$-adic $L$-function evaluated at the trivial character $1$, and the right-hand side $(\log_{\omega} \circ g)_{\ord}$ corresponds to the logarithm of the Heegner point $P$.

\subsection{Review: the construction of BDP $p$-adic $L$-function} \label{sec:construct-BDP}
We need to review the construction of the BDP $p$-adic $L$-function $\calL_{\BDP}(E/\calL)$ in Theorem \ref{thm:BDP}, which is due to \cite{HsiehSpecial}. We will use the same notation and conventions as in the proof of Theorem \ref{thm:BDP}.

\subsubsection{The choice of toric modular forms}
Let $\chi = \eta \phi$ be a Hecke character over $\calL$, where $\eta$ is the auxiliary Hecke character over $\calL$ of infinite type $(1,-1)$ and $\phi$ is the character whose $p$-adic avatar $\hphi$ lies in the critical range $\frX^{\mathrm{crit}}$ in \cite{HsiehSpecial}. In particular, the Hecke character $\chi$ is unramified outside $p$.

In \cite[Definition 3.12]{HsiehSpecial}, Hsieh constructed explicitly a \emph{$p$-modified toric Whittaker function} $W_{\chi} \in \sfW(\Pi, \psi)$ as a product of local Whittaker functions
\[
W_{\chi, u} = W_{\chi, \infty} W_{\chi, f}^{(p)} W_{\chi, f, u} \in \sfW(\Pi, \psi), \quad W_{\chi} = W_{\chi, \infty} W_{\chi, f}^{(p)} W_{\chi, f}^{\flat} \in \sfW(\Pi, \psi),
\]
where $u \in \ZZ_p^{\times}$. They are defined as follows in \cite[Section 3.6]{HsiehSpecial}:

\underline{(1)} At the archimedean place, $W_{\chi, \infty}$ is defined in \cite[Section 3.6.1]{HsiehSpecial}, which we shall not recall.

For the local Whittaker functions at finite places, Hsieh first chooses a function $\bfa_{\chi,\ell}: \QQ_{\ell}^{\times} \rightarrow \CC^{\times}$ in the Kirillov model $\sfK(\Pi_{\ell}, \psi)$ and takes its corresponding Whittaker function $W_{\chi, \ell}$ under the identification map $\sfW(\Pi_p, \psi) \simeq \sfK(\Pi_p,\psi)$ introduced in \cite[Equation (1.1)]{HsiehSpecial}.

\underline{(2)} For the local Whittaker function at $p$, the Kirillov function is given by
    \[
    \bfa_{\chi,p, u}(a) := \bfone_{u(1+p\ZZ_p)}(a) \chi_{\barv_0}(a^{-1}).
    \]
for some $u \in \ZZ_p^{\times}$. This is called the \emph{$u$-component} of the local Kirillov function at $p$. We also define 
\[
    \bfa_{\chi,p}^{\flat}(a) := \bfone_{\ZZ_p^{\times}}(a) \chi_{\barv_0}(a^{-1}),
\]
then $W_{\chi, p}^{\flat} := W_{\bfa_{\chi,p}^{\flat}} = \sum_{u \in \calU_{p}} W_{\chi, p,u}$, where $\calU_{p}$ is the torsion subgroup of $\ZZ_p^{\times}$.

\underline{(3)} For the local Whittaker function at $\ell \neq p$, suppose $\ell$ splits in $\calL$ as $\frl \bar{\frl}$, then the Kirillov function chosen in \cite[Section 3.6.2]{HsiehSpecial} depends on the local behavior of $\Pi_{\ell}$: 
    \begin{enumerate}[label = \rm (\roman*)]
        \item If $L(s, \Pi_{\ell} \otimes \chi_{\frl}) = 1$, Hsieh takes
        \[
        \bfa_{\chi,\ell} = \bfone_{\ZZ_{\ell}^{\times}}(a) \chi_{\frl}(a^{-1}) = \bfone_{\ZZ_{\ell}^{\times}}(a) 
        \]
        since $\chi_{\frl}$ is trivial on units.
        \item If $\Pi_{\ell} = \DS(\mu,\nu)$ is unramified, then
        \begin{align*}
        \bfa_{\chi,\ell} &= \bfone_{\ZZ_{\ell}^{\times}}(a) \chi_{\frl}(a^{-1}) \abs{a}_{\ell}^{1/2} \sum_{i+j=r, i,j \geq 0} (\mu\chi)_{\bar{\frl}}(\omega_{\bar{\frl}}^{i}) \cdot (\nu\chi)_{\bar{\frl}}(\omega_{\bar{\frl}}^{j}) \\
        &= \bfone_{\ZZ_{\ell}}(a) \abs{a}_{\ell}^{1/2} \sum_{i+j=r, i, j \geq 0} \mu(\varpi_{\bar{\frl}})^{i} \nu(\varpi_{\bar{\frl}})^{j}
        \end{align*}
        where $r = v_{\ell}(a)$.
        \item If, in $\Pi_{\ell} = \DS(\mu, \nu)$, one character $\mu_{i}$ is ramified, then he takes
        \[
        \bfa_{\chi,\ell} = \bfone_{\ZZ_{\ell}}(a)\mu_{i}(a) \abs{a}^{1/2},
        \]
        \item If $\Pi_{\ell}$ is special, he takes
        \[
        \bfa_{\chi,\ell} = \bfone_{\ZZ_{\ell}}(a)\mu(a) \abs{a}^{1/2}
        \]
    \end{enumerate}
These are exactly the respective Kirillov functions of the normalized newvector of $\Pi_{\ell}$.

\underline{(4)} For the local Whittaker function at $\ell \neq p$ with $\ell$ being inert or ramified in $\calL$. We first introduce the averaging operator introduced in \cite[Section 3.6.3]{HsiehSpecial}. Consider the torus embedding $\iota: \calL \hookrightarrow \rmM_{2}(\QQ)$ given by \cite[Equation (2.1)]{HsiehSpecial}, and let $\varsigma_{\ell} \in \GL_2(\QQ_{\ell})$ be the matrices defined in \cite[Equation (3.3)]{HsiehSpecial} for all rational primes $\ell$. At each $\ell \mid d$, we choose a local uniformizer $\varpi_{\calL,\ell} \in \calO_{\calL,\ell}$ and set the averaging operators 
\[
\calP_{\varsigma} := \prod_{\ell \mid d} \dfrac{1}{2}(1+\sfR(t_{\ell})), \quad\text{ for } t_{\ell} := \varsigma_{\ell}^{-1} \iota(\varpi_{\calL,\ell}) \varsigma_{\ell}, 
\]
where $\sfR(t)$ means the right-translation action by $t \in \GL_2(\AA_{\QQ,\rmf})$, whenever this action makes sense. This is exactly the operator $\calP_{\chi, \varsigma}$ defined in \cite[Section 3.6.3]{HsiehSpecial} since $\chi_{\ell}(\varpi_{\calL, \ell}) = 1$. Moreover, by our Heegner assumption, $\Pi_{\ell}$ is unramified. Thus, the term $\calR_{\ell}^{m_{\ell}(\chi, \Pi_{\ell})}$ in \textit{op.cit.} is $1$.

Then the local Whittaker function $W_{\chi, \ell}$ is taken to be
\[
W_{\chi, \ell} := \calP_{\varsigma} W_{\ell}^{\circ},
\]
where $W_{\ell}^{\circ}$ is the normalized spherical Whittaker function in $\sfW(\Pi_{\ell}, \psi)$, which is available since $\Pi_{\ell}$ is unramified.

After choosing the local toric Whittaker functions, Hsieh attached an automorphic form $\varphi_{\chi, u} \in \calA(\Pi)$ to $W_{\chi, u}$ in \cite[Equation (3.24)]{HsiehSpecial}. Moreover, Hsieh attached to $\varphi_{\chi, u}$, in \cite[Definition 5.1]{HsiehSpecial}, a complex Hilbert modular form $\bff_{\chi, u}^{\star}: \frH \times \GL_2(\AA_{\rmf}) \rightarrow \CC$ by
\begin{equation} \label{eq:complex-modular-form-bff-star}
\bff_{\chi, u}^{\star}(\tau, g_{\rmf}) := \abs{\det \varsigma_{\rmf}}_{\rmf}^{-1} \varphi_{\chi, u}(g) J(g_{\infty}, i)^{2} (\det g_{\infty})^{-2} \abs{\det g}_{\AA},
\end{equation}
where $g_{\infty} \in \GL_2(\RR)$ satisfies $g_{\infty} i = \tau$ for $g = (g_{\infty}, g_{\rmf}) \in \GL_{2}(\AA_{\infty}) \times \GL_{2}(\AA_{\rmf}) = \GL_2(\AA)$. The same definitions are also made after removing $u$ in the subscript. Putting $\chi = 1$, we have the complex modular form $\bff_{1,u}^{\star}$.

Now we fix an integer $c$ that is coprime to $pN$. \footnote{Note that the subscript $c$ here has nothing to do with the previous subscript $u \in \ZZ_p^{\times}$.} In \cite[Proposition 5.5]{HsiehSpecial}, by interpolating their Fourier coefficients, Hsieh constructed a $p$-adic family of modular forms $\calF_{\eta, c}: \calC(\Gamma_{\calL}^{\ac}) \rightarrow V(c, L, \barZZ_p)$ such that, for every $\hphi \in \frX_{p}^{\crit}$ of infinite type $(m,-m)$ for $m \geq 0$, we have
\begin{equation} \label{eq:theta-property-plus}
\calF_{\eta, c}(\hphi) = \Theta^{m} \cdot \widehat{\bff}_{\eta\phi, c}^{\star}.
\end{equation}
where $\Theta$ is the Dwork--Katz $p$-adic differential operator in \cite[Corollary (2.6.25)]{katz78}, and $\widehat{\bff}_{\eta\phi,c}^{\star}$ is the $p$-adic avatar of $\bff_{\eta\phi,c}^{\ast}$ defined in \cite[Equation (4.4)]{HsiehSpecial}.

With the $p$-adic family of modular forms $\calF_{\eta,c}$, we specialize it at the character $\heta^{-1}$ to obtain a $p$-adic modular form $\calF_{\eta,c}(\heta^{-1})$. Note that since $\eta^{-1}$ has infinite type $(-1,1)$, we are specializing $\calF_{\eta,c}$ outside its interpolation range. Nevertheless, using the same proof of \cite[Proposition 5.5]{HsiehSpecial}, we have the following proposition.

\begin{proposition} \label{prop:theta-property-minus}
With the notation above, we have
\[
\Theta \cdot \calF_{\eta, c}(\heta^{-1}) = \widehat{\bff}_{1, c}^{\star}.
\]   
\end{proposition}

\begin{proof}
We remark that this is \emph{not} a consequence of \eqref{eq:theta-property-plus} since $\eta^{-1}$ has infinite type $(-1,1)$, which is outside the interpolation range $\frX_{p}^{\crit}$. Nevertheless, the proof is the same by comparison of Fourier coefficients on both sides: the Fourier coefficient expression
\[
\calF_{\eta,c}(q)(\phi) = \sum_{\beta > 0 } \bfa_{\beta}(\bff_{\eta}^{\star}, c) \phi(\rec_{p}(\beta)) q^{\beta}.
\]
appearing in the proof of \cite[Proposition 5.5]{HsiehSpecial} is valid for all $\phi \in \calC(\Gamma_{\calL}^{-})$ and it is actually the defining property of the $p$-adic family $\calF_{\eta,c}$. It is then possible to apply $\phi = \eta^{-1}$ and compare the Fourier coefficients of both sides of the identity $\Theta \cdot \calF_{\eta, c}(\heta^{-1}) = \widehat{\bff}_{1, c}^{\star}$ that we are proving.
\end{proof}

\subsubsection{CM points}
Following \cite[Section 4.6]{HsiehSpecial}, we define 
\[
\sfx: \AA_{\calL}^{\times} \rightarrow \frH \times \GL_2(\AA_{\rmf}), \quad a = (a_{\infty}, a_{\rmf}) \mapsto \sfx(a) = (\vartheta, \iota(a_{\rmf})\varsigma_{\rmf}). 
\]
For $a \in (\AA_{\calL,\rmf}^{(pN)})^{\times}$, we have a $\CC$-quintuple $(\ulA(a), j(a))$ associated with $\sf\sfx(a)$ following \cite[Section 4.4]{HsiehSpecial}. It gives rise to a complex point $[\sf\sfx(a)]$ in the Igusa tower $I_{L}(c(a))(\CC)$, where $c(a) = c(\calO_{L}) \Norm(\fra)^{-1}$ with $\fra = \il_{\calK}(a) := a \widehat{\calO}_{K} \cap K$, and the point $[\sfa(a)]$ descends to $I_{K}(\ZZ_p^{\ur})$. These points are called the \emph{CM points} on the modular curve $X_{0}(N)$.

\subsubsection{The construction of BDP $p$-adic $L$-function}
We are now ready to state the construction of $\Theta_{\ttH}(E/\calL)$ in \cite{HsiehSpecial}, given in \cite[Equation (5.9)]{HsiehSpecial}. For any anticyclotomic character $\kappa$ over $\calL$ factoring through $\Gamma_{\ac}^{-}$,
\begin{align*}
\Theta_{\ttH}(E/\calL)(\kappa) &= C_{\Pi}(\eta)^{-1/2} \sum_{a \in D} (\eta \abs{-}_{\AA_{\calL}}^{-1})(a) \calF_{\eta, c(a)}(\kappa|[a])(\sfx(a)) \\ 
&= C_{\Pi}(\eta)^{-1/2} \sum_{a \in D} (\eta \abs{-}_{\AA_{\calL}}^{-1})(a) \kappa(a) \calF_{\eta, c(a)}(\kappa)(\sfx(a)),
\end{align*}
where $C_{\Pi}(\eta)$ is the constant defined in \cite[Proposition 5.2]{HsiehSpecial}.

We now compute its specialization. Substituting $\kappa = \heta^{-1}$ and using $a_p = 1$ so that $\heta(a) = \eta(a)$, we see that 
\[
\Theta_{\ttH}(E/\calL)(\heta^{-1}) = C_{\Pi}(\eta)^{-1/2} \sum_{a \in D} (\eta \abs{-}_{\AA_{\calL}}^{-1})(a) \heta(a)^{-1} \calF_{\eta, c(a)}(\heta^{-1})(\sfx(a)).
\]
Since $(\eta \abs{-}_{\AA_{\calL}}^{-1})(a) \heta(a)^{-1} = \eta(a)\abs{a}_{\AA_{\calL}}^{-1} \eta(a)^{-1} = \abs{a}_{\AA_{\calL}}^{-1} = \Norm(\fra)$, it further simplifies to 
\[
\Theta_{\ttH}(E/\calL)(\heta^{-1}) = C_{\Pi}(\eta)^{-1/2} \sum_{a \in D} \Norm(\fra) \cdot \calF_{\eta, c(a)}(\heta^{-1})(\sfx(a)),
\]
and hence by the twist in the definition of $\Theta_{\BDP}$ in the proof of Theorem \ref{thm:BDP}, we have
\begin{equation} \label{eq:construction-BDP}
\Theta_{\BDP}(E/\calL)(1) = \left(\dfrac{c_{\ttH}}{C_{\Pi}(\eta)}\right)^{1/2} \sum_{a \in D} \Norm(\fra) \cdot \calF_{\eta, c(a)}(\heta^{-1})(\sfx(a)).
\end{equation}

\subsection{Liu--Zhang--Zhang's formula: the $p$-adic $L$-function side}
Returning to Liu--Zhang--Zhang's formula, i.e., Theorem \ref{thm:LZZ}, we now focus on its left-hand side.

Given the fixed modular parameterization $\phi: X_{0}(N) \rightarrow E$, we define its average $[\phi^{\av}] := \calP_{\varsigma} [\phi] \in \bPi_{K}$. We start by the following observation.

\begin{lemma}
    With the notation above, retain the assumption that $E$ has additive reduction at $p$ (in particular, $p \mid N$). Then $g := [\phi^{\av}]$ is a stable vector in $\bPi_{L}$.
\end{lemma}

\begin{proof}
This follows from the fact that, in our setting, $f_{E}$ is the normalized cuspidal newform attached to the CM elliptic curve $E$ and $p \mid N$. We explain this below.

We start with $[\phi] \in \bPi_{L}$. Under the identification $\calD_{\omega_{E}^{\nor}}$ in Remark \ref{rem:identification-LZZ}, we have seen that 
\[
\calD_{\omega_{E}^{\nor}}([\phi]) = \phi^{\ast} \omega_{E}^{\nor} = 2 \pi i f_{E}(z) \dif z.
\]
To check that $[\phi] \in \bPi_{L}$ is a stable vector, it suffices to check that the automorphic form $\varphi_{f}$ attached to $f_{E}$ inside $\Pi$ is "stable" in the sense that, in the abstract tensor-product factorization $\varphi_{f} = \otimes_{\ell} \varphi_{f,\ell}$, the local vector $\varphi_{f,p} \in \Pi_{p}$ corresponds to a Kirillov function in the Kirillov model $\sfK(\Pi_{p}, \psi)$ supported on $\ZZ_p^{\times}$.

In our setting, $E$ is a CM elliptic curve with complex multiplication given by $\calK$ where $p$ remains inert. Now that $p \mid N$, this forces
\[
L_p(f_{E}, s) = \sum_{r \geq 0} a_{p^{r}} p^{-rs} = 1,
\]
which implies that $a_{p^{r}} = 0$ for all $r \geq 1$. The Kirillov function $\kappa_{p}$ correspond to $\varphi_{f,p}$ is given by 
\[
a_{p^{r}}(f) = p^{r} \kappa(p^{r}), \quad \text{for $r \geq 0$}.
\]
This implies that $\kappa_{p}$ is supported on $\ZZ_p^{\times}$. The averaging operator $\calP_{\varsigma}$ only modifies the local Kirillov functions at $\ell$ not dividing $p$, so it leaves $\varphi_{f,p}$ unchanged. This verifies that $[\phi^{\av}]$ is a stable vector inside $\bPi_{L}^{\hsuit}$.
\end{proof}

\begin{remark}[Good reduction case] \label{rem:good-reduction-remark}
    This is the crucial point at which the proof of Theorem \ref{thm:waldspurger} uses the assumption that $p \mid N$ (and also that $E$ has CM, so that $E$ has additive reduction at $p$). In the case $p \nmid N$, to continue the proof, one needs to apply the \emph{$p$-depletion operator} $\sfD_{p} := 1 - p^{-1} \sum_{a \in \ZZ/p} \sfR(n(a/p))$ to $\phi$, i.e., to take
    \[
    \phi^{[p]} := \sfD_{p} \phi, \quad [\phi^{\av}] := \calP_{\varsigma} [\phi^{[p]}]
    \]
    in place of $[\phi^{\av}]$ in our current proof. This is the source of the factor $(1-a_p(E)+p)/p$ in the BDP formula \eqref{eq:BDP-good-case} for the case $p \nmid N$. After this adaptation, we can obtain a unified proof of Theorem \ref{thm:waldspurger} for both reduction types using the method of this section. We do not pursue this here.
\end{remark}

We are now ready to prove the following proposition.
\begin{proposition}[$p$-adic $L$-function side identity] \label{prop:p-adic-L-LZZ}
    With the notation above, we have
\[
\bfM(((\phi^{\av})^{\ast} \omega_{E}^{\nor})_{\ord})(1)(\sfx(a)) = \Norm(\fra) \cdot \calF_{\eta, c(a)}(\heta^{-1})(\sfx(a)).
\]
As a result, we have
\begin{equation} \label{eq:LZZ-BDP-final}
\bfM(((\phi^{\av})^{\ast} \omega_{E}^{\nor})_{\ord})(1)(\sfx(a)) = \left(\dfrac{c_{\ttH}}{C_{\Pi}(\eta)}\right)^{-1/2} \Theta_{\BDP}(E/\calL)(1).
\end{equation}
\end{proposition}

\begin{proof}
We fix an element $g_{\rmf} \in \GL_2(\AA_{\rmf})$; it gives a map
\[
\frH \rightarrow X_{U}(\CC) = \GL_2(\QQ)^{+} \bs \frH \times \GL_{2}(\AA_{\rmf})/U, \quad \tau \mapsto [\tau, g_{\rmf}], 
\]
and hence fixes a \emph{chart} $\sfC_{g_{\rmf}}$ of $X_{U}(\CC)$. We shall work with each individual chart.

We put $h_{g_{\rmf}}(\tau) := \varphi_{0}(g_{\infty}, g_{\rmf}) J(g_{\infty}, i)^{2} (\det g_{\infty})$, where $\varphi_{0}$ is the adelization of $\calD_{\omega_{E}^{\nor}}([\phi^{\av}])$, i.e. $(\phi^{\av})^{\ast} \omega_{E}^{\nor} = 2 \pi i h_{g_{\rmf}}(\tau) \dif \tau$ over the chart $\sfC_{g_{\rmf}}$. Specializing the definition of $\bff_{\chi}^{\star}$ in \eqref{eq:complex-modular-form-bff-star} at $\chi = 1$, we see that 
\[
\bff_{1}^{\star}(\tau, g_{\rmf}) = \abs{\det \varsigma_{\rmf}}_{\rmf}^{-1} \varphi_{0}(g) J(g_{\infty}, i)^{2} (\det g_{\infty})^{-2} \abs{\det g}_{\AA} = D(g_{\rmf}) \cdot h_{g_{f}}(\tau)
\]
where $D(g_{\rmf}) = \abs{\det \varsigma_{\rmf}}_{\rmf}^{-1} \abs{\det g_{\rmf}}$. After considering a cusp on $\PP^{1}(\QQ)$ in this chart and moving it to the infinite cusp, we denote the local parameter by $q_{w} = \exp(2 \pi i \tau/w)$. Thus, we can write
\[
h_{g_{\rmf}}(\tau) = \sum_{n > 0} b_{n} q_{w}^{n}.
\]
Rewrite $2 \pi i \dif \tau = w q_{w}^{-1} \dif q_{w}$, we can express
\[
(\phi^{\av})^{\ast} \omega_{E}^{\nor} = w \sum_{n > 0} b_{n} q_{w}^{n} \dfrac{\dif q_{w}}{q_{w}}.
\]
Using the relation \eqref{eq:theta-relation} in Theorem \ref{thm:mellin}, together with the description of the Atkin--Serre operator at this cusp, we have
\[
\bfM(((\phi^{\av})^{\ast} \omega_{E}^{\nor})_{\ord})(1) = w \sum_{n > 0} \dfrac{b_{n}}{n} q_{w}^{n}.
\]
On the other hand, we can express
\[
\bff_{1,c}^{\star}(\tau, g_{\rmf}) = D(g_{\rmf}) \sum_{n>0} b_{n}q_{w}^{n}.
\]
By Proposition \ref{prop:theta-property-minus}, together with the description of the Dwork--Serre operator at this cusp \footnote{Note that, under the normalisation of \cite{HsiehSpecial}, an extra factor $w$ appears: $\Theta(q_{w}^{n}) = \frac{n}{w} q_{w}^{n}$, i.e., $\Theta = w^{-1} q_{w} \dfrac{\rmd}{\dif q_{w}}$.}
we see that
\[
\calF_{\eta, c}(\heta^{-1})(g_{\rmf}) = D(g_{\rmf}) \cdot w \sum_{n > 0} \dfrac{b_{n}}{n} q_{w}^{n}.
\]
Thus, by the $q$-expansion principle, we obtain on this chart that
\[
\bfM(((\phi^{\av})^{\ast} \omega_{E}^{\nor})_{\ord})(1) = D(g_{\rmf}) \cdot \calF_{\eta, c}(\heta^{-1}). 
\]
In particular, we take the chart given by the CM point $\sfx(a)$, i.e., we take $g_{\rmf}$ to be $\iota(a)\varsigma_{\rmf}$. Then
\[
D(\iota(a)\varsigma_{\rmf}) = \abs{\det(\iota(a))}_{\rmf} = \abs{\Norm_{\calL/\QQ}(a)}_{\rmf} = (\Norm \fra)^{-1}.
\]
This shows the first equation in the proposition. The final equality \eqref{eq:LZZ-BDP-final} is then a direct consequence of \eqref{eq:construction-BDP}.
\end{proof}

\subsection{Liu--Zhang--Zhang's formula: the Heegner points side}
In this section, we establish the following identity.

\begin{proposition}[Heegner point side identity] \label{prop:Heegner-point-LZZ}
    With the notation above, we have
    \[
    (\log_{\omega} \circ \phi^{\av})_{\ord} = \dfrac{\log_{\omega_E} P}{ s c_{\phi}},
    \]
    where $P$ is the Heegner point and $s = \#\Cl_{\ram}(\calL)$, both defined in Section \ref{sec:Heegner-point}.
\end{proposition}

\begin{proof}
We define functions 
\[
F(z) := \dfrac{1}{c_{\phi}} \log_{\omega_{E}} \circ \phi, \quad F^{\av}:= \calP_{\varsigma} F = \log_{\omega_{E}^{\nor}} \circ \phi^{\av}.
\]
as functions on the modular curve $X_{0}(N)$. In particular, we shall apply the function $F^{\av}$ to CM points $\sfx(a) = (\vartheta, \iota(a_{\rmf})\varsigma_{\rmf})$. For $t \in \AA_{\calL,\rmf}^{\times}$, one checks that the right translation $\sfR(t_{\ell})$ acts on $\iota(a_{\ell})\varsigma_{\ell}$ by $\iota(a_{\ell})\varsigma_{\ell} \varsigma_{\ell}^{-1} \iota(t) \varsigma(\ell) = \iota(a_{\ell}t_{\ell}) \varsigma_{\ell}$. When $t = \varpi_{\calL,\ell}$ is the uniformizer at a prime $\frl \mid d$, the above equality shows that the CM point $\frx(a)$ is translated to $\sfx(a \frl)$. Let $r_{d}$ be the number of primes of $\QQ$ ramified in $\calL$. Expanding the definition of $\calP_{\varsigma}$ gives
\[
F^{\av}(\sfx(a)) = \dfrac{1}{2^{r_{d}}} \sum_{\epsilon \in \{0,1\}^{r_{d}}} F \left(\sfx_{\fra \prod_{\ell \mid d} \frl^{\epsilon}} \right).
\]
Since $\frl^{2} = (\ell)$ is a principal ideal, each generator $[\ell]$ has order at most two in the class group. The map
\[
(\ZZ/2)^{r_d} \rightarrow \Cl_{\ram}(\calL), \quad (\epsilon_{\ell})_{\ell} \mapsto \prod_{\ell \mid d} [\frl]^{\epsilon_{\ell}}
\]
is surjective by the definition of $\Cl_{\ram}(\calL)$, with kernel of cardinality $2^{r_d}/s$, and each fiber is a coset of the kernel. Regrouping the preceding finite sum gives
\[
F^{\av}(\sfx(a)) = \dfrac{1}{2^{r_{d}}} \dfrac{2^{r_{d}}}{s} \sum_{c \in \Cl_{\ram}(\calL)} F(\sfx_{\fra c}) = \dfrac{1}{s} \sum_{c \in \Cl_{\ram}(\calL)} F(\sfx_{\fra c}).
\]
Using the definition of $F$ and summing over $D$, we have
\begin{align} \label{eq:LZZ-Heegner-final}
\sum_{a \in D} F^{\av}(\sfx(a)) = \dfrac{1}{s} \sum_{[\frb] \in \Cl(\calL)} F(\sfx_{\frb}) &= \dfrac{1}{s c_{\phi}} \sum_{[\frb] \in \Cl(\calL)} \left( \log_{\omega_{E}} \phi(\sfx_{\frb}) \right) \notag \\ &= \dfrac{1}{s c_{\phi}} \log_{\omega_{E}} \left(  \sum_{[\frb] \in \Cl(\calL)} \phi(\sfx_{\frb}) \right) =  \dfrac{\log_{\omega_{E}} P}{s c_{\phi}},
\end{align}
where the first equality follows from the fact that the map $D \times \Cl_{\ram}(\calL) \rightarrow \Cl(\calL)$ given by $(\fra, c) = [\fra]c$ is bijective. The proposition is therefore proved.
\end{proof}

\begin{remark}
    We make the following two remarks.
    \begin{itemize}
        \item From this short reduction in this section, one can see that Hsieh's averaging operator is really indispensable in the proof: otherwise we are unable to obtain the summation over all CM points indexed by the class group, and thus cannot obtain the Heegner point on this side.
        \item As a continuation of Remark \ref{rem:good-reduction-remark}, on the Heegner point side, we also need to replace $F$ by the $p$-depleted version $\sfD_{p}(F)$. Then it turns out that
    \[
    \log_{\omega_{E}}(\sfD_{p}[\phi]) = F - \dfrac{a_{p}(f_{E})}{p} V_{p} F + \dfrac{1}{p} V_{p}^{2} F
    \]
    when $p \nmid N$. This corresponds to the extra factor in \eqref{eq:bdp-formula}.
    \end{itemize}
\end{remark}

\subsection{Proof of the $p$-adic Waldspurger formula}

We are now ready to prove Theorem \ref{thm:waldspurger}.

\begin{proof}[Proof of Theorem \ref{thm:waldspurger}]
Applying Liu--Zhang--Zhang's $p$-adic Waldspurger formula in Theorem \ref{thm:LZZ}, together with the $p$-adic $L$-function identity in Proposition \ref{prop:p-adic-L-LZZ} and the Heegner side identity in Proposition \ref{prop:Heegner-point-LZZ}, we have
\[
\Theta_{\BDP}(E/\calL)(1) = \left(\dfrac{c_{\ttH}}{C_{\Pi}(\eta)}\right)^{1/2} \dfrac{\log_{\omega_{E}} P}{s c_{\phi}}.
\]
To deduce Proposition \ref{thm:waldspurger} it remains to write down explicitly the coefficient $c_{\ttH}/C_{\Pi}(\eta)$. By the definition of the constant $C_{\Pi}(\eta)$ in \cite[Proposition 5.2]{HsiehSpecial}, we see that
\[
C_{\Pi}(\eta) = C^{\prime}_{\Pi}(\eta) \cdot \dfrac{\sqrt{d}}{4s^{2}},
\]
where the constant $C^{\prime}_{\Pi}(\eta)$ is defined in \cite[Equation (3.27)]{HsiehSpecial}, which simplifies to 
\begin{align*}
C^{\prime}_{\Pi}(\eta) = 8 \prod_{w \mid \frN^{(p)}} \epsilon(1/2, \Pi_{\ell} \otimes \eta_{\barw}, \psi_{\ell}) &= 8 \prod_{w \mid \frN^{(p)}} \eta_{\barw}(\ell)^{v_{\ell}(N)} \epsilon(1/2, \Pi_{\ell}, \psi_{\ell}) \\
&= 8 \eta(\frN^{(p)})^{-1} \prod_{w \mid \frN^{(p)}} \epsilon(1/2, \Pi_{\ell}, \psi_{\ell})
\end{align*}
where $\ell$ is the rational prime lying below $w$, and we use the fact that $\eta_{\barw}$ is unramified. Combining these expressions gives
\[
\dfrac{c_{\ttH}}{C_{\Pi}(\eta)} = \dfrac{\sqrt{d}}{2 \eta(\frN^{(p)})} \dfrac{1}{8 \eta(\frN^{(p)})^{-1} \prod_{w \mid \frN^{(p)}} \epsilon(1/2, \Pi_{\ell}, \psi_{\ell})} \dfrac{4s^{2}}{\sqrt{d}} = \dfrac{s^{2}}{4 \epsilon^{\frN^{(p)}}},
\]
where $\epsilon^{\frN^{(p)}} := \prod_{w \mid \frN^{(p)}} \epsilon(1/2, \Pi_{\ell}, \psi_{\ell})$. Therefore,
\[
\calL_{\BDP}(E/\calL)(1) = (\Theta_{\BDP}(E/\calL)(1))^{2} = \dfrac{1}{4\epsilon^{\frN^{(p)}}} \left( \dfrac{\log_{\omega_E} P}{c_{\phi}} \right)^{2}.
\]
Each local representation $\Pi_{\ell}$ is self-dual with trivial central character, and hence the functional equation for $\epsilon$-factors as well as the change of additive character $\psi_{\ell}$ gives
\[
\epsilon(1/2, \Pi_{\ell}, \psi_{\ell}) \epsilon(1/2, \Pi_{\ell}^{\vee}, \psi_{\ell}^{-1}) =  \epsilon(1/2, \Pi_{\ell}, \psi_{\ell}) = 1
\]
and 
\[
\epsilon(1/2, \Pi_{\ell}, \psi_{\ell}^{-1}) = \omega_{\Pi_{\ell}}(-1) \epsilon(1/2, \Pi_{\ell}, \psi_{\ell}) = \epsilon(1/2, \Pi_{\ell}, \psi_{\ell}).
\]
Putting these together with $\Pi_{\ell} \simeq \Pi_{\ell}^{\vee}$, we see that each of these $\epsilon$-factors has square $1$, and hence $\epsilon^{\frN^{(p)}} \in \{\pm 1\}$. This shows that $u_{\rmW} := (4 \epsilon^{\frN^{(p)}})^{-1}$ is a $p$-adic unit when $p \geq 5$. Theorem \ref{thm:waldspurger} is therefore proved. 
\end{proof}

We will use the following corollary of the $p$-adic Waldspurger formula (Theorem \ref{thm:waldspurger}).

\begin{corollary}[$p$-adic Waldspurger formula] \label{cor:length-p-adic-waldspurger}
    With the notation above, we have
    \[
    v_{p}(\calL_{\BDP}(E/\calL)(1)) = 2 v_p(m) - 2 v_p(c_{\phi}) + 2r.
    \]
\end{corollary}
\begin{proof}
Choose the sign of the generator $Q$ of $E(\calL)/E(\calL)_{\tor}$ so that $P = m Q + Q^{\prime}$, where $Q^{\prime}$ is a torsion point in $E(\calL)_{\tor}$. Then $\log_{\omega_E} P = m \log_{\omega_E} Q$. Hence
\[
v_p(\log_{\omega_E} P) = v_p(m) + v_p(\log_{\omega_{E}} Q) = v_p(m) + r + b_{p}(E)
\]
where the last equality uses the notation of Lemma \ref{lem:gather-finite} together with the definition of $r$. One verifies directly from the definitions that $v_{p}(\scE_{p}) = -b_p(E)$, which gives
\[
v_p(\log_{\omega_E} P) = v_p(m) + v_p(\log_{\omega_{E}} Q) = v_p(m) + r - v_{p}(\scE_{p}).
\]
Taking the $p$-adic valuation of both sides of \eqref{eq:bdp-formula}, we obtain
\begin{align*}
    v_p(\calL_{\BDP}(E/\calL)(1)) &= 2 v_p(\scE_p) - 2 v_p(c_{\phi}) + 2 v_p(\log_{\omega_E}(P)) \\
    &= 2 v_p(\scE_p) - 2 v_p(c_{\phi}) + 2 v_p(m) + 2 r - 2 v_{p}(\scE_{p}) \\ 
    &= 2 v_p(m) - 2 v_p(c_{\phi}) + 2r.
\end{align*}
This proves the corollary, as desired.
\end{proof}

%% file: 206-final.tex
\section{Birch--Swinnerton-Dyer formula in the general inert case: lower bound}

In this section, we finish the proof of the one-side divisibility of the $p$-part of the BSD formula in the case in which $p \geq 5$ remains inert in $\calK$. The notation and conventions are inherited from the previous section.

\subsection{Gross--Zagier formula}
Let $D = D_{\calL}$, and let $E^{D}$ be the quadratic twist of $E$ by $\calL$. Then we have the following reformulation of the Gross--Zagier formula:
\begin{proposition}[Gross--Zagier formula] \label{prop:Gross-Zagier}
    Under the assumption \eqref{eq:GZK}, we have the following equality
    \[
    v_p \left( \dfrac{L^{\prime}(E/\QQ, 1)}{\Omega_{E} \Reg(E/\QQ)} \dfrac{L(E^{D}/\QQ, 1)}{\Omega_{E^{D}}} \right) = 2 v_p(m/c_{\phi}),
    \]
where both sides are nonzero rational numbers.
\end{proposition}
\begin{proof}
    This follows from the explicit Gross--Zagier formula in \cite[Theorem 1.1]{CST}. One verifies that, in our setup, the hypotheses there are satisfied, and the formula gives
    \begin{equation} \label{eq:explicit-GZ}
    L^{\prime}(E/\calL, 1) = 2^{-\mu(N,D)} \dfrac{8 \pi^{2} \lrangle{f,f}}{u_{\calL}^{2} \sqrt{d} \deg(\phi)} \hath_{\calL}(P).
    \end{equation}
    Here the terms are as follows:
    \begin{itemize}
        \item $\mu(N,D)$ is the number of primes dividing both $N$ and $D$. Here since every prime dividing $N$ splits in $\calL$, we have $\mu(N,D) = 0$,
        \item $u_{\calL} := [\calO_{\calL}^{\times}:\{\pm 1\}]$. Since we have taken $D_{\calL} < -4$ to be an odd number, $u_{\calL} = 1$;
        \item $\hath_{\calL}$ is the Néron--Tate height on $E$ over $\calL$,
        \item $\deg \phi$ is the degree of the modular parametrization $\phi: X_{0}(N) \rightarrow E$,
        \item $\lrangle{f,f}$ is the Petersson norm
        \[
        \lrangle{f,f} := \int_{\Gamma_{0}(N) \bs \frH} \abs{f(z)}^{2} \dif x \dif y.
        \]
    \end{itemize}
    We deduce the formula in the proposition with the following steps.

    \underline{(1) We rewrite the quotient $\lrangle{f,f}/\deg(\phi)$}. We put 
    \[
    \calA_{E} := \dfrac{i}{2} \int_{E(\CC)} \omega_{E} \wedge \baromega_{E}.
    \]
    Pulling back along $\phi$, we have
    \begin{align}
    (\deg \phi) \calA_{E} &= \dfrac{i}{2} \int_{X_0(N)(\CC)} \phi^{\ast}\omega_{E} \wedge \overline{\phi^{\ast}\omega_{E}} \notag \\ 
    &= \dfrac{i}{2} 4\pi^{2} c_{\phi}^{2} \int_{\Gamma_{0}(N)\bs\frH} \abs{f(z)}^{2} (-2i) \dif x \dif y = 4\pi^{2}c_{\phi}^{2} \lrangle{f,f}. \label{eq:calA_E}
    \end{align}
    Here we have used that $\dif z \wedge \dif \barz = -2i \dif x \wedge \dif y$.

    \underline{(2) Next we obtain the periods $\Omega_{E}$ and $\Omega_{E^{D}}$ from $\calA_{E}$}. Consider the complex uniformization
    \[
    E(\CC) \simeq \CC/\Lambda, \quad \Lambda = \left\{ \int_{\gamma} \omega_{E}: \gamma \in \rmH_{1}(E(\CC),\ZZ) \right\}.
    \]
    Consider the complex-conjugation action on $\rmH_{1}(E(\CC), \ZZ)$, and let $\rmH_{1}^{+}(E(\CC), \ZZ)$ and $\rmH_{1}^{-}(E(\CC), \ZZ)$ denote the two eigenspaces of complex conjugation with eigenvalues $\pm 1$. One checks that the index $\kappa := [\rmH_{1}(E(\CC), \ZZ): \rmH_{1}^{+}(E(\CC), \ZZ) \oplus \rmH_{1}^{-}(E(\CC), \ZZ)]$ is either $1$ or $2$. Both $\rmH_{1}^{+}(E(\CC), \ZZ)$ and $\rmH_{1}^{-}(E(\CC), \ZZ)$ have rank one; we therefore choose their generators $\gamma_{+}$ and $\gamma_{-}$ and set
    \[
    \lambda_{+} := \abs{\int_{\gamma_{+}} \omega_{E}}, \quad \lambda_{-} := \abs{\dfrac{1}{i} \int_{\gamma^{-}} \omega_{E}}.
    \]
    The rectangle spanned by $\lambda_{+}$ and $i \lambda_{-}$ in $\CC$ therefore has area $\kappa$ times the area of the fundamental parallelogram of $\Lambda$, with the latter being exactly $\calA_{E}$. This gives
    \[
    \calA_{E} = \lambda_{+} \lambda_{-}/\kappa.
    \]
    We choose the twisting isomorphism $\tau: E_{\calL}^{D} \rightarrow E_{\calL}$ and write $\tau^{\ast} \omega_{E} = a_{D} \sqrt{D} \omega_{E^{D}}$ for some nonzero rational number $a_{D}$; moreover, one checks that $a_{D}$ is a $p$-adic unit. Since $p$ is unramified in $\calL$, the factor $\sqrt{D}$ is also a $p$-adic unit. The map $\tau_{\ast}$ exchanges the positive lattice $\rmH_1^{+}(E^{D}(\CC),\ZZ)$ with $\rmH_1^{-}(E(\CC),\ZZ)$. Thus, the positive period of $E^{D}$ is $\lambda_{-}/(a_{D} \sqrt{d})$. We denote $n_{C} := \pi_0(C(\RR)) \in \{1,2\}$ for $C = E$ or $E^{D}$. Integration then gives
    \begin{equation} \label{eq:period-appear}
    \Omega_{E/\QQ} = n_{E} \lambda_{+}, \quad \Omega_{E^{D}/\QQ} = \dfrac{n_{E^{D}} \lambda_{-}}{\abs{a_{D}}\sqrt{d}}.
    \end{equation}
    Combining \eqref{eq:calA_E} and the relation in \eqref{eq:period-appear} gives
    \begin{equation} \label{eq:GZ-step-2}
    L^{\prime}(E/\calL, 1) = \dfrac{2}{c_{\phi}^{2}} \cdot \dfrac{\abs{a_{D}}}{\kappa n_{E} n_{E^{D}}} \Omega_{E/\QQ} \Omega_{E^{D}/\QQ} \cdot \hath_{\calL}(P).
    \end{equation}

    \underline{(3) Finally we obtain the regulator from the Néron--Tate height $\hath_{\calL}(P)$}. Let $Q_{\QQ}$ be a generator of $E(\QQ)/E(\QQ)_{\tor}$; then its image in $E(\calL)/E(\calL)_{\tor}$ is $j\QQ$ for $j = 1$ or $2$. For the absolute quadratic Néron--Tate height, we have
    \[
    \hath(P) = m^{2} \hath(Q), \quad \Reg(E/\QQ) = j^{2} \hath(Q),
    \]
    and hence $\hath(P) = m^{2} j^{-2} \Reg(E/\QQ)$. We write $\hath_{\calL}(P) = \kappa_{\Ht} \hath(P)$ where $\kappa_{\Ht}$ is a power of $2$ when comparing the absolute quadratic Néron--Tate height with the canonical Néron--Tate pairing. So this gives
    \[
    \hath_{\calL}(P) = \kappa_{\Ht} m^{2} j^{-2} \Reg(E/\QQ).
    \]
    Substituting this into \eqref{eq:GZ-step-2}, and noting that $L^{\prime}(E/\calL, 1) = L^{\prime}(E/\QQ, 1) L(E^{D}/\QQ, 1)$, we have
    \[
    \dfrac{L^{\prime}(E/\QQ, 1) L(E^{D}/\QQ, 1)}{\Omega_{E/\QQ} \Omega_{E^{D}/\QQ} \Reg(E/\QQ) } = \dfrac{2 \abs{a_{D}} \kappa_{\Ht}}{\kappa n_{E} n_{E^{D}} j^{2}} \dfrac{m^{2}}{c_{\phi}^{2}}.
    \]
    For $p \geq 5$, every factor in the first fraction on the right-hand side is a $p$-adic unit. Taking valuations proves the proposition, as desired.
\end{proof}

\subsection{The final proof}
We are now ready to conclude the proof.

\begin{theorem}\label{thm:Sha-lower-bound}
    Let $E/\QQ$ have complex multiplication by $\calO_{K}$ and $\ord_{s=1} L(E/\QQ, s) = 1$. Let $p \geq 5$ be a prime number which remains inert in $\calK$. Then
    \[
    v_{p} \left( \dfrac{L^{\prime}(E/\QQ, 1)}{\Omega_{E/\QQ} \Reg(E/\QQ)} \right) \leq v_{p} \left( \dfrac{\# \Sha(E/\QQ) \prod_{\ell} \Tam_{\ell}(E/\QQ)}{\# E(\QQ)_{\tor}^{2}} \right).
    \]
\end{theorem}
\begin{proof}
    The existence of the imaginary quadratic field $\calL$ satisfying \eqref{eq:rank-one} is well-known. Combining Proposition \ref{prop:Gross-Zagier}, Theorem \ref{thm:waldspurger} (or more precisely Corollary \ref{cor:length-p-adic-waldspurger}) and Theorem \ref{thm:L_BDP_1}, we have
    \begin{align*}
    v_p \left( \dfrac{L^{\prime}(E/\QQ, 1)}{\Omega_{E/\QQ} \Reg(E/\QQ)} \right) + v_p \left( \dfrac{L(E^{D}/\QQ, 1)}{\Omega_{E^{D}}} \right) &= 2v_p(m/c_{\phi}) \\ 
    &\leq v_{p}(\# \Sha(E/\calL)[p^{\infty}])  \\
    &= v_{p}(\# \Sha(E/\QQ)[p^{\infty}]) + v_{p}(\# \Sha(E^{D}/\QQ)[p^{\infty}]),
    \end{align*}
    where the last equality follows from the decomposition of the Tate--Shafarevich group
    \[
    \Sha(E/\calL)[p^{\infty}] \simeq \Sha(E/\QQ)[p^{\infty}] \oplus \Sha(E^{D}/\QQ)[p^{\infty}].
    \]
    We have seen that $E^{D}(\QQ)[p] = 0$ and all Tamagawa numbers $\Tam_{\ell}(E^{D}/\QQ)$ are coprime to $p$. The rank-zero $p$-part of the BSD formula for $E^{D}$ is already proved in \cite{BurungaleFlach}. This simplifies the above inequality to
    \[
    v_p \left( \dfrac{L^{\prime}(E/\QQ, 1)}{\Omega_{E/\QQ} \Reg(E/\QQ)} \right) \leq v_{p}(\# \Sha(E/\QQ)[p^{\infty}]).
    \]
    Since $E(\QQ)[p] = 0$ and all Tamagawa numbers $\Tam_{\ell}(E/\QQ)$ are coprime to $p$, the theorem is therefore proved.
\end{proof}

%% file: 207-Kolyvagin.tex
\section{Birch--Swinnerton-Dyer formula in the general inert case: upper bound}

In this section, suppose $p \geq 5$. Let $E/\QQ$ still be an elliptic curve with CM by the field $\calK$ and analytic rank $1$. We assume that $p$ is inert in $\calK$. We use the Kolyvagin system method. The aim of this section is to give the divisibility in the other direction:

\begin{equation}\label{eq:Sha-upper-bound}
    v_{p} \left( \dfrac{L^{\prime}(E/\QQ, 1)}{\Omega_{E/\QQ} \Reg(E/\QQ)} \right) \geq v_{p} \left( \dfrac{\# \Sha(E/\QQ) \prod_{\ell} \Tam_{\ell}(E/\QQ)}{\# E(\QQ)_{\tor}^{2}} \right).
\end{equation}
Together with the lower bound, Theorem \ref{thm:Sha-lower-bound}, this proves the $p$-BSD formula for $E/\QQ$ (Theorem \ref{thm:BSD inert}).

\begin{remark}\label{rmk:Sha-upper-bound}
    As in the proof of Theorem \ref{thm:Sha-lower-bound}, we only need to show, for some Heegner field $\calL$, an upper bound for the Tate--Shafarevich group:
    \[
    v_p(\#\Sha(E/\calL)[p^{\infty}])\leq 2v_p(m/c_{\phi}).
    \]
\end{remark}

We recall the notation from above. The elliptic curve $E/\QQ$ has conductor $N_E$ and CM field $\calK$. Let $\calL$ be a Heegner field in which all the primes dividing $N_E$ split. We may assume that $E$ is a strong Weil curve with an optimal modular parametrization $\phi:X_0(N_E)\rightarrow E$. We then obtain the Heegner point $P\in E(\calL)$ and Heegner index $m:=[E(\calL)/E(\calL)_{\tor}:\ZZ \bar{P}]\in\ZZ_+$.

\subsection{The Kolyvagin system method}

By the method of Kolyvagin systems, one can obtain an exact order for $\Sha(E/\calL)[p^{\infty}]$ in terms of divisibility indices. But here a rough upper bound suffices.

\begin{theorem}[Kolyvagin, special case of {\cite[Corollary 13]{kolyvagin2007euler}}, see also {\cite[Theorem 0.3]{matar2019kolyvagin}} for the formulation here]\label{thm:Kolyvagin-bound}
    Assume that $D\neq-3,-4$ and that
    \begin{itemize}
        \item $\forall n_1,n_2>0$, $\rmH^1(\calL(E[p^{n_1+n_2}])/\calL, E[p^{n_1}])=0$.
        \item the $\pmod{p}$ Galois representation $\bar{\rho}_{E,p}:\Gal_{\QQ}\rightarrow\Aut_{\mathbb{F}_p}(E[p])\cong\GL_2(\mathbb{F}_p)$ is irreducible.
        \item $E(\calL)[p]=0$.
    \end{itemize}
    If $P\in E(\calL)$ is non-torsion, then $E(\calL)/\ZZ P$ is finite and
    \[
    v_p(\#\Sha(E/\calL)[p^{\infty}])\leq 2v_p(m).
    \]
\end{theorem}

\begin{remark}
    The irreducibility of $\bar{\rho}_{E,p}$ implies $E(\calL)[p]=0$ for $p\neq2$, by an easy check in linear algebra. Moreover, \cite{lawson2015vanishing} gives a complete calculation, which implies that the cohomology groups above also vanish when $\bar{\rho}_{E,p}$ is irreducible and $p\neq2$. See also \cite[Page 457]{matar2019kolyvagin} for an explanation.
\end{remark}

It is in fact easy to show that, in our inert case, the residual representation is irreducible.

\begin{proposition}\label{prop:irreducible}
    Assume that $E/\mathbb{Q}$ has CM by $\calK$ and that the prime $p\geq3$ is inert in $\calK$. Then the residual representation $\bar{\rho}_{E,p}$ is irreducible.
\end{proposition}
\begin{proof}
    It suffices to show that $\bar{\rho}_{E,p}$ is irreducible after restricting to $\Gal_\calK$.

    Let $\mathcal{O}:=\End_{\bar{\QQ}}(E)$ denote the CM order in $\mathcal{O}_\calK$, with conductor $f$. By a well-known classification of $(E,\mathcal{O},\calK)$, we have $p\nmid f$ and $\mathcal{O}/p\mathcal{O}\cong\FF_{p^2}$. Since the action of $\mathcal{O}$ on $E[p]$ is $\FF_p$-linear, there is a natural map
    \[
    \mathcal{O}/p\mathcal{O}\rightarrow\End_{\FF_p}(E[p]).
    \]
    It is actually injective. Indeed, if a nonzero isogeny $\alpha\in\mathcal{O}$ acts trivially on $E[p]$, then $E[p]\subseteq\ker(\alpha)$, which implies $p^2=\#E[p]\mid\#\ker(\alpha)=\Norm^\calK_{\QQ}(\alpha)$. Since $p$ is inert in $\calK$ and $p\nmid f$, this forces $\alpha\in p\mathcal{O}_\calK\cap\mathcal{O}=p\mathcal{O}$. Hence, we obtain an injection
    \[
    \theta:\FF_{p^2}^{\times}\cong(\mathcal{O}/p\mathcal{O})^{\times}\hookrightarrow\Aut_{\FF_p}(E[p])\cong\GL_2(\FF_p),
    \]
    which sends the subgroup $\FF_p^{\times}\cong(\ZZ/p\ZZ)^{\times}\subset(\mathcal{O}/p\mathcal{O})^{\times}$ onto the diagonals.

    Moreover, since all endomorphisms in $\mathcal{O}$ are defined over $\calK$, the action of $\sigma\in\Gal_\calK$ on $E[p]$ is $\FF_{p^2}\cong\mathcal{O}/p\mathcal{O}$-linear. But, as an $\FF_{p^2}$-vector space, $E[p]$ has dimension $1$, by comparing orders. It follows that the action of $\sigma\in\Gal_\calK$ on $E[p]$ is multiplication by some $a_{\sigma}\in\mathcal{O}/p\mathcal{O}$ and thus factors through $\theta$, i.e., there is a character
    \[
    \psi:\Gal_\calK\rightarrow\FF_{p^2}^{\times}\cong(\mathcal{O}/p\mathcal{O})^{\times},
    \]
    satisfying $\bar{\rho}_{E,p}\vert_{\Gal_\calK}=\theta\circ\psi$.

    Recall that by Weil pairing, the determinant of $\bar{\rho}_{E,p}$ is given by the residue cyclotomic character $\bar{\chi}_p:\Gal_{\QQ}\rightarrow\FF_p^{\times}$, which is surjective. We claim that $\det(\bar{\rho}_{E,p}\vert_{\Gal_\calK})$ is still full $\FF_p^{\times}$. This is the case since $\bar{\rho}_{E,p}$ factors through $\Gal(\QQ(E[p])/\QQ)$, while $p$ is ramified in $\QQ(E[p])$ but remains inert in $\calK$. It follows that $\Gal_\calK$ surjects onto $\Gal(\calK\QQ(E[p])/\calK)\cong\Gal(\QQ(E[p])/\QQ)$ and has the same image as $\Gal_{\QQ}$ under $\bar{\rho}_{E,p}$.

    Now assume that $E[p]$ has a nontrivial $\Gal_\calK$-subspace $L$, generated by $v\in E[p]$ as an $\FF_p$-subspace. Then, for any $\sigma\in\Gal_\calK$, there is some $b_{\sigma}\in\FF_p^{\times}$ such that
    \[
    a_\sigma\cdot v=\sigma(v)=b_{\sigma}\cdot v.
    \]
    Since $v\neq0$ and $\FF_{p^2}$ acts $1$-transitively on $E[p]$, we deduce that $a_{\sigma}=b_{\sigma}\in\FF_p$, i.e. $\psi$ has image contained in diagonals $\FF_p\subset\FF_{p^2}$. But then $\det(\bar{\rho}_{E,p}(\Gal_\calK))\subseteq(\FF_p^{\times})^2$, a contradiction to the surjectivity of $\bar{\chi}_{p}\vert_{\Gal_\calK}$, provided $p\neq2$. Therefore $\bar{\rho}_{E,p}$ is irreducible on $\Gal_\calK$ and thus irreducible on $\Gal_{\QQ}$.    
\end{proof}

So by Theorem \ref{thm:Kolyvagin-bound} and Proposition \ref{prop:irreducible} we obtain an upper bound, without the Manin constant term.

\begin{proposition}\label{prop:Sha-upper-bound-no-Manin}
    Notations as in Theorem \ref{thm:Kolyvagin-bound}, we have
    \[
    v_p(\#\Sha(E/\calL)[p^{\infty}])\leq 2v_p(m).
    \]
\end{proposition}

Therefore, the only remaining issue is a version of Manin's conjecture\footnote{This was recently proved using AI; see \cite{chi2026maninconstant}. We will not use this result, however, and will instead rely on our own transfer method.}. We will reduce it, by the transfer method, to finitely many cases.

\subsection{The case $p\geq11$}

\begin{theorem}\label{thm:Sha-upper-bound p>=11}
    For optimal $E/\QQ$ with CM by $\calK$ as above, assume $p\geq 11$ is inert in $\calK$. Then the Manin constant $c_{\phi}$ is a $p$-adic unit. Thus the desired inequality \eqref{eq:Sha-upper-bound} follows.
\end{theorem}
\begin{proof}
    By \cite[Theorem 3]{edixhoven1991manin} we know that $v_p(c_{\phi})=0$, except for some cases when $E$ has potentially ordinary reduction at $p$. But by the classification in Table \ref{table:CM_classification}, $E$ is potentially supersingular at $p$, since $p$ is inert in $\calK$. Together with Proposition \ref{prop:Sha-upper-bound-no-Manin} we obtain the inequality \eqref{eq:Sha-upper-bound}.
\end{proof}

\subsection{The case $p=5,7$}

Here we use the transfer method and check the Manin constants case by case. Consequently, there are only $10$ cases to check.

We still denote by $F:=\calK_p$ the unique unramified quadratic extension of $\QQ_p$, and put $U^n:=1+p^n\mathcal{O}_F$. First, the classification in \eqref{eq:G_n classification} shows that $G_n:=\mathcal{O}_F^{\times}/(\ZZ_p^{\times}U^n)\cong C_{p+1}\oplus C_{p^{n-1}}$. Note that $\calK$ contains $\mu_2$, $\mu_4$, or $\mu_6$ as its full group of roots of unity, and $\psi_p$ must factor through some $G_n$. Hence, to exhaust all possibilities, we may assume that $\psi_5$ takes values in $\mu_6$ and $\psi_7$ takes values in $\mu_4$, after checking that $5$ is inert in $\QQ(\sqrt{-3})$ and $7$ is inert in $\QQ(\sqrt{-1})$. Moreover, we have
\[
\lvert\Hom(C_6\oplus C_{5^{n-1}},\mu_6)\rvert=\lvert\Hom(C_6,\mu_6)\rvert=6
~\text{and}~
\lvert\Hom(C_8\oplus C_{7^{n-1}},\mu_4)\rvert=\lvert\Hom(C_8,\mu_4)\rvert=4
\]
for all $n\geq1$. Thus, $\psi_5$ and $\psi_7$ must factor through $G_1$, i.e., the target curves must satisfy $v_p(\mathfrak{f}_{\psi})=0$ (good supersingular) or $1$ (additive potentially good supersingular), or equivalently $v_p(N_E)=0$ or $2$, which is consistent with Table $\ref{table:CM_classification}$.

As in the setup of Proposition \ref{prop:calculation 3,5,7 inert}, we set $g_5:=\Teich(2+\omega)$ and $g_7:=\Teich(2+i)$. Note that $\psi_p$ is now determined by its value on $g_p$. We then searched LMFDB to realize all such $\psi_p$. The following two tables display families of representatives for $\psi_5$ and $\psi_7$, respectively.

\begingroup\small\setlength{\tabcolsep}{5pt}\renewcommand{\arraystretch}{1.25}
\begin{longtable}{@{}ccccccc@{}}
\caption{Representatives for each $\psi_5$ with $5$ inert}\label{tab:p_5 inert database}\\
\hline
Equation of $E$ & LMFDB label & $N_E$ & $\calK$ & $r_{\rm an}$ & Manin constant (all optimal) & $\psi_{E,5}(g_5)$ \\\hline\endfirsthead
\multicolumn{7}{c}{\tablename\ \thetable\ (continued)}\\\hline
Equation of $E$ & LMFDB label & $N_E$ & $\calK$ & $r_{\rm an}$ & Manin constant (all optimal) & $\psi_{E,5}(g_5)$ \\\hline\endhead
\hline\endfoot
$y^2+y=x^3+2$ & 1728.n4 & $2^{6}\cdot 3^{3}$ & $\QQ(\sqrt{-3})$ & 1 & 1 & $1$ \\
$y^2+y=x^3+11$ & 6075.v2 & $3^{5}\cdot 5^{2}$ & $\QQ(\sqrt{-3})$ & 1 & 1 & $\zeta_6^5$ \\
$y^2+y=x^3+56$ & 6075.q2 & $3^{5}\cdot 5^{2}$ & $\QQ(\sqrt{-3})$ & 1 & 1 & $\zeta_6^4$ \\
$y^2+y=x^3+281$ & 6075.m2 & $3^{5}\cdot 5^{2}$ & $\QQ(\sqrt{-3})$ & 1 & 1 & $\zeta_6^3$ \\
$y^2+y=x^3+1406$ & 6075.o2 & $3^{5}\cdot 5^{2}$ & $\QQ(\sqrt{-3})$ & 1 & 1 & $\zeta_6^2$ \\
$y^2+y=x^3+7031$ & 6075.s2 & $3^{5}\cdot 5^{2}$ & $\QQ(\sqrt{-3})$ & 1 & 1 & $\zeta_6^1$ \\
\end{longtable}\endgroup

\begingroup\small\setlength{\tabcolsep}{5pt}\renewcommand{\arraystretch}{1.25}
\begin{longtable}{@{}ccccccc@{}}
\caption{Representatives for each $\psi_7$ with $7$ inert}\label{tab:p_7 inert database}\\
\hline
Equation of $E$ & LMFDB label & $N_E$ & $\calK$ & $r_{\rm an}$ & Manin constant (all optimal) & $\psi_{E,7}(g_7)$ \\\hline\endfirsthead
\multicolumn{7}{c}{\tablename\ \thetable\ (continued)}\\\hline
Equation of $E$ & LMFDB label & $N_E$ & $\calK$ & $r_{\rm an}$ & Manin constant (all optimal) & $\psi_{E,7}(g_7)$ \\\hline\endhead
\hline\endfoot
$y^2=x^3-2x$ & 256.b1 & $2^{8}$ & $\QQ(\sqrt{-1})$ & 1 & 1 & $1$ \\
$y^2=x^3-7x$ & 3136.l1 & $2^{6}\cdot 7^{2}$ & $\QQ(\sqrt{-1})$ & 1 & 1 & $i$ \\
$y^2=x^3+49x$ & 3136.m4 & $2^{6}\cdot 7^{2}$ & $\QQ(\sqrt{-1})$ & 1 & 1 & $-1$ \\
$y^2=x^3-343x$ & 3136.r1 & $2^{6}\cdot 7^{2}$ & $\QQ(\sqrt{-1})$ & 1 & 1 & $-i$ \\
\end{longtable}\endgroup

\subsection{The final proof}

Note that every possible $\psi_p$ is represented by a curve with Manin constant $1$.

\begin{theorem}\label{thm:BSD inert}
    Let $E/\QQ$ be an elliptic curve with complex multiplication by $\calK$ with analytic rank $1$. Let $p$ be a prime inert in $\calK$. Then the $p$-BSD holds for $E/\QQ$.
\end{theorem}
\begin{proof}
    For $p=2,3$, this is proved in Theorem \ref{thm:BSD 2,3 inert}. For $p=5,7$, all the curves in Tables \ref{tab:p_5 inert database} and \ref{tab:p_7 inert database} have Manin constant $1$. Thus, Remark \ref{rmk:Sha-upper-bound} and Proposition \ref{prop:Sha-upper-bound-no-Manin} show that \eqref{eq:Sha-upper-bound} holds for those curves. Together with Theorem \ref{thm:Sha-lower-bound}, the $p$-BSD formula holds. Since any such curve $E$ must have a local character $\psi_{E,p}$ isomorphic to that of some curve $E^{\prime}$ in the tables, it follows that the $p$-BSD formula holds for $E$ by the transfer theorem, Theorem \ref{thm:main_transfer}. For $p\geq11$, this follows directly from Theorems \ref{thm:Sha-lower-bound} and \ref{thm:Sha-upper-bound p>=11}.
\end{proof}

We are now able to conclude our main theorem.
\begin{proof}[Proof of Theorem \ref{thm:intro-cm-bsd}]
    This is a combination of Theorem \ref{thm:BSD ramified} (when $p$ ramifies in $\calK$), Theorem \ref{thm:Nov-split-BSD} (when $p$ splits in $\calK$) and Theorem \ref{thm:BSD inert} (when $p$ remains inert in $\calK$).
\end{proof}

%% file: AppA-Elliptic-units.tex
\section{Elliptic units}\label{app:elliptic-units}

In this appendix, we provide the necessary background on elliptic units used in the main text. The primary reference is \cite[\S 15]{kato2004padic}.

\subsection{Basic setup}\label{subsec:eu-conventions}

We fix the imaginary quadratic field $\calK$, the elliptic curve $E/\QQ$ with complex multiplication by $\calO_{\calK}$, and the embeddings into $\CC$ and $\overline{\QQ}_p$ as in the article. Let $M/\calK$ be a sufficiently large finite extension and $L/\QQ_p$ be a sufficiently large finite extension containing the image of $M$. Put $O=\calO_L$. For an integral ideal $\mathfrak m\ne0$, write $\calK(\mathfrak m)$ for its ray class field and $N\mathfrak m=\#(\calO_{\calK}/\mathfrak m)$. For $H/\calK$ finite, $S_H$ denotes the places above a specified set $S$.

More explicitly, let $I(\mathfrak m)$ be the group of fractional ideals prime to $\mathfrak m$, and let $P_1(\mathfrak m)$ be the subgroup of principal ideals $(a)$ with $a\equiv1\pmod{\mathfrak m}$ at its prime divisors. Class field theory specifies the ray class field by
\[
 I(\mathfrak m)/P_1(\mathfrak m)
 \xrightarrow{\ \sim\ }\operatorname{Gal}(\calK(\mathfrak m)/\calK).
\]
We put $\calO_{H,S_H} = \calO_{H}[1/S_{H}]$ as the ring obtained by inverting the finite primes in $S_H$. An \emph{$S_H$-unit} is an element of $H^\times$ of valuation zero at every finite prime outside $S_H$.

We use the arithmetic convention for the Artin symbol in class field theory: at an unramified prime $\mathfrak l$, $\sigma_{\mathfrak l}$ induces the arithmetic Frobenius $x\mapsto x^{N\mathfrak l}$ on residue fields. To be consistent with Part \ref{part:BDP}, we say that a Hecke character $\vartheta: \calK^{\times} \bs \AA_{\calK}^{\times} \rightarrow \CC^{\times}$ has infinity type $(a,b)$ if its archimedean component is $z^a\bar z^b$, and
\begin{equation}\label{eq:eu-Hecke-L}
 L(\vartheta,s)=\prod_{\mathfrak l\nmid\mathfrak f_\vartheta}
 \left(1-\vartheta(\mathfrak l)(N\mathfrak l)^{-s}\right)^{-1}.
\end{equation}
These are also consistent with Kato's conventions in \cite[\S15.7]{kato2004padic}. 

We put $\epsilon=\chi_{\mathrm{cyc}}$ as the cyclotomic character. A compatible system $\varepsilon=(\zeta_{p^n})_{n\geq1}$, i.e. $\zeta_{p^{n+1}}^p=\zeta_{p^n}$ for any $n \geq 1$, defines a basis $e_1$ of $\ZZ_p(1)$ with $g(e_1)=\epsilon(g)e_1$. Put $e_j=e_1^{\otimes j}$, using dual bases for negative $j$. We have
\begin{equation}\label{eq:eu-native-T}
 R_p=\calO_{\calK}\otimes\ZZ_p,\quad T_E=T_pE,\quad
 T_\sigma=T_E\otimes_{R_p,\sigma}O,\quad V_\sigma=T_\sigma\otimes_O L.
\end{equation}
Its Galois character is $\rho_\sigma$. Thus $g(t_\sigma)=\rho_\sigma(g)t_\sigma$ for a basis $t_\sigma$, and away from $p$ and the conductor, $\rho_\sigma(\sigma_{\mathfrak l})=\sigma(\psi_E(\mathfrak l))$ with the corresponding CM label. The character of $T_\sigma^\dagger=\operatorname{Hom}_O(T_\sigma,O)(1)$ is $\epsilon\rho_\sigma^{-1}$.

\subsection{The construction of elliptic units}\label{subsec:eu-theta}

A \emph{CM pair} of modulus $\mathfrak m$ over $H\supset\calK$ consists of an elliptic curve $\mathscr E/H$ with $\calO_{\calK}$-action and a point $\alpha$ whose annihilator is $\mathfrak m$. The action on $\operatorname{Lie}(\mathscr E/H)$ must be the given embedding $\calK\hookrightarrow H$. Assume
\begin{equation}\label{eq:eu-fine-modulus}
 \calO_{\calK}^{\times}\longrightarrow(\calO_{\calK}/\mathfrak m)^\times
 \quad\text{is injective}.
\end{equation}
An automorphism fixing $\alpha$ is a CM unit congruent to $1$ modulo $\mathfrak m$, hence is $1$.

The main theorem of complex multiplication supplies a canonical pair $(\mathscr E_{\mathfrak m},\alpha_{\mathfrak m})$ over $\calK(\mathfrak m)$ with analytic description $(\CC/\mathfrak m,1\bmod\mathfrak m)$. For $\mathfrak a$ prime to $\mathfrak m$, it gives
\[
 (\mathscr E_{\mathfrak m}/\mathscr E_{\mathfrak m}[\mathfrak a],
 \alpha_{\mathfrak m}\bmod\mathscr E_{\mathfrak m}[\mathfrak a])
 \simeq(\mathscr E_{\mathfrak m}^{\sigma_{\mathfrak a}},
 \sigma_{\mathfrak a}\alpha_{\mathfrak m}).
\]
We use \cite[(15.3.1)--(15.3.3)]{kato2004padic} for these existence and reciprocity statements.

For any Hecke character $\vartheta$ over $\calK$ of infinite type $(-r,0)$, Kato defined a one-dimensional $M$ space
\[
V_{M}(\vartheta) = \rmH^{1}_{\rmB}(\scE(\CC), \QQ)^{\otimes_{\calK} r} \otimes_{\calK} M,
\]
where $\scE = \scE_{\frm}$ is the CM elliptic curve appearing in the canonical CM-pair; see \cite[\S15.8]{kato2004padic}. In particular, for $\vartheta = \psi_{E}$, we have
\begin{equation}\label{eq:eu-cohom-homology}
 V_M(\psi_E)\otimes_{M,\sigma}L\simeq V_\sigma^\vee,\qquad
 V_M(\psi_E)(1)\otimes_{M,\sigma}L\simeq V_\sigma^\dagger.
\end{equation}

Let $\mathscr E/H$ be a CM elliptic curve, and let $\mathfrak a$ be a proper nonzero ideal prime to $6$. Set
\[
 D_{\mathfrak a}=N\mathfrak a[0]-\sum_{P\in\mathscr E[\mathfrak a]}[P].
\]
The following proposition, as stated in \cite[\S 15.4]{kato2004padic}, follows by the same argument as in the proof of \cite[Proposition 1.3]{kato2004padic}.

\begin{proposition}\label{prop:eu-theta}
There is a unique rational function $ {}_{\mathfrak a}\theta_{\mathscr E}$ over $H$ satisfying
\begin{equation}\label{eq:eu-theta-normalization}
 \operatorname{div}({}_{\mathfrak a}\theta_{\mathscr E})=D_{\mathfrak a},
 \qquad N_{[m]}({}_{\mathfrak a}\theta_{\mathscr E})={}_{\mathfrak a}\theta_{\mathscr E}
 \quad\text{if }(m,\mathfrak a)=1.
\end{equation}
Moreover, for any isogeny $q:\mathscr E\to\mathscr E'$ of degree prime to $N\mathfrak a$, such that $\End(\scE) \simeq \calO_{\calK}$ and $\End(\scE^{\prime}) \simeq \calO_{\calK}$, we have
\begin{equation}\label{eq:eu-theta-isogeny}
 N_q({}_{\mathfrak a}\theta_{\mathscr E})={}_{\mathfrak a}\theta_{\mathscr E'}.
\end{equation}
\end{proposition}

For $(\mathfrak a,6\mathfrak m)=1$ define
\begin{equation}\label{eq:eu-unit-definition}
 {}_{\mathfrak a}z_{\mathfrak m}
 ={}_{\mathfrak a}\theta_{\mathscr E_{\mathfrak m}}(\alpha_{\mathfrak m})^{-1}
 \in\calK(\mathfrak m)^\times.
\end{equation}
The point is outside the divisor's support, so the value is finite and nonzero. The integrality theorem gives $ {}_{\mathfrak a}z_{\mathfrak m}\in\calO_{\calK(\mathfrak m)}[1/\mathfrak m]^\times$; it is a global unit if $\mathfrak m$ has at least two distinct prime ideal divisors (see \cite[\S15.5]{kato2004padic} or \cite[\S3.2, equation (25)]{BurungaleFlach}). The following relation is recorded as \cite[(15.4.3)--(15.4.4)]{kato2004padic}.
\begin{lemma}\label{lem:eu-smoothing}
For $(\mathfrak a\mathfrak b,6\mathfrak m)=1$, we have
\begin{equation}\label{eq:eu-smoothing}
 (N\mathfrak b-\sigma_{\mathfrak b})\,{}_{\mathfrak a}z_{\mathfrak m}
 =(N\mathfrak a-\sigma_{\mathfrak a})\,{}_{\mathfrak b}z_{\mathfrak m}.
\end{equation}
\end{lemma}

\subsection{Elliptic units over the ray class tower}\label{subsec:eu-kummer}

Let $H$ be a number field and let $S_H$ contain the primes above $p$. The Kummer sequence defines
\[
\kappa_{H,n}:\calO_{H,S_H}^{\times}/(\calO_{H,S_H}^{\times})^{p^n} \rightarrow \rmH^1(\calO_{H,S_H},\mu_{p^n}).
\]
Write $B^\wedge_p=\varprojlim B/p^nB$ for completion of an abelian group, using $p^n$-th powers in multiplicative notation.  We then have the isomorphism
\begin{equation}\label{eq:eu-Kummer-iso}
 \kappa_{H}: (\calO_{H,S_H}^{\times})^\wedge_p
 \xrightarrow{\sim} \rmH^1(\calO_{H,S_H},\ZZ_p(1)).
\end{equation}
Moreover, the following lemma is easily checked and gives the compatibility of the Kummer maps $\kappa_{H}$ under corestriction maps.

\begin{lemma}\label{lem:eu-norm-cor}
For $H'/H$ finite, with compatible sets of places, we have $\cor_{H^{\prime}/H}\kappa_{H^{\prime}}(u)=\kappa_H(N_{H^{\prime}/H}u)$.
\end{lemma}

Fix $\frf \ne 0$ and set
\[
 \calK_\infty=\bigcup_{n\geq0}\calK(p^n\mathfrak f),\qquad
 G=\operatorname{Gal}(\calK_\infty/\calK),\qquad \Lambda=\ZZ_p \lrbracket{G}.
\]
Here $p^n\mathfrak f=(p^n)\mathfrak f$. For sufficiently large $n$, this modulus satisfies \eqref{eq:eu-fine-modulus}. Fix a proper ideal $\mathfrak a$ prime to $6p\mathfrak f$. Then elliptic units satisfy the following distribution property:
\begin{equation}\label{eq:eu-tower-norm}
 N_{\calK(p^{n+1}\mathfrak f)/\calK(p^n\mathfrak f)}
 ({}_{\mathfrak a}z_{p^{n+1}\mathfrak f})
 ={}_{\mathfrak a}z_{p^n\mathfrak f}.
\end{equation}
See \cite[\S15.5]{kato2004padic} or \cite[\S3.2]{BurungaleFlach}.

For $S$ containing $p\mathfrak f$ and infinity, define
\[
 \rmH^1_{\mathrm{Iw}}(\calK_\infty/\calK,\ZZ_p(1))=\varprojlim_{H\subset\calK_\infty} \rmH^1(\calO_{H,S_H},\ZZ_p(1)),
\]
with corestriction transition maps. Equation \eqref{eq:eu-tower-norm} and Lemma~\ref{lem:eu-norm-cor} give
\[
 {}_{\mathfrak a}\mathbf z_{p^\infty\mathfrak f}
 =\bigl(\kappa_{\calK(p^n\mathfrak f)}
 ({}_{\mathfrak a}z_{p^n\mathfrak f})\bigr)_n
 \in \rmH^1_{\mathrm{Iw}}(\calK_\infty/\calK,\ZZ_p(1)).
\]
The action of $G$ is the usual action on field elements and their Kummer classes.

Let $Q(\Lambda)$ be the total quotient ring. Following \cite[\S 15.6]{kato2004padic}, Kato proves that $N\mathfrak a-\sigma_{\mathfrak a}$ is a non-zero-divisor for the proper ideal $\mathfrak a$ above, and defines
\begin{equation}\label{eq:eu-normalized-z}
 \mathbf z_{p^\infty\mathfrak f}
 =(N\mathfrak a-\sigma_{\mathfrak a})^{-1}
 {}_{\mathfrak a}\mathbf z_{p^\infty\mathfrak f}
 \in \rmH^1_{\mathrm{Iw}}(\calK_\infty/\calK,\ZZ_p(1))
 \otimes_\Lambda Q(\Lambda).
\end{equation}
Independence of $\mathfrak a$ follows by multiplying two such expressions by the product of their denominators and applying \eqref{eq:eu-smoothing}. This is the normalization of \cite[\S15.6]{kato2004padic} and \cite[\S3.2]{BurungaleFlach}.

\subsection{Two variations of elliptic units}\label{subsec:eu-twisting}

\subsubsection{With $\calO_{L}$-coefficients}
Let $\Lambda^\#$ be free of rank one over $\Lambda$, with $G_{\calK}$-action $g(\lambda)=[\bar g]^{-1}\lambda$, where $\bar g\in G$. Shapiro's lemma gives the identification
\begin{equation}\label{eq:eu-Shapiro}
 \rmH^{1}_{\Iw}(\calK_{\infty}/\calK, \ZZ_p(1)) \simeq  \rmH^1(\calO_{\calK,S},\Lambda^\#(1)).
\end{equation}
Let $T=O \ttt$ be a free $O$-module of rank one, with an $O$-basis $\ttt$, and with a $\Gal_{\calK}$-action given by a character $\rho$ factoring through $G$. Define
\begin{equation}\label{eq:eu-evaluation}
 \operatorname{ev}_{\rho}:\Lambda\longrightarrow O,\qquad
 [g]\longmapsto\epsilon(g)\rho(g)^{-1}.
\end{equation}
It gives the coefficient map
\begin{equation}\label{eq:eu-coefficient-map}
 \Phi_{\ttu}: \Lambda^\#(1)\longrightarrow T,\qquad
 \lambda\otimes \tte_1\longmapsto\operatorname{ev}_{\rho}(\lambda)\ttt.
\end{equation}
Here the letter $\ttu$ represents the choice of basis $\tte_1$ and $\ttt$. The map $\Phi_u$ is Galois equivariant, since
\[
 \operatorname{ev}_{\rho}([g]^{-1})\epsilon(g)
 =(\epsilon(g)\rho(g)^{-1})^{-1}\epsilon(g)=\rho(g).
\]
Combining these maps, we have a composition
\[
\rmH^{1}_{\Iw}(\calK_{\infty}/\calK, \ZZ_p(1)) \simeq  \rmH^1(\calO_{\calK,S},\Lambda^\#(1)) \xrightarrow{(\Phi_{\ttu})_{\ast}} \rmH^{1}(\calO_{\calK,S}, T) \rightarrow \rmH^{1}(\calO_{\calK,S}, V).
\]
Applying this map to the elliptic unit $\bfz_{p^{\infty}\frf}$, we obtain
\begin{equation}\label{eq:eu-indexed-class}
 z_{T,S}(\ttu)\in H^1(\calO_{\calK,S},V),\qquad V=T\otimes_O L.
\end{equation}

\subsubsection{Over the anticyclotomic family in \cite{FanWan}}

We keep the notation of Section \ref{sec:plus-minus}. We now consider $V_{\bpsi} = \Lambda_{m} \tte_{\bpsi,\sigma}$, a free $\Lambda_{m}$-module of rank one, with a $\Lambda_{m}$-basis $\tte_{\bpsi, \sigma}$ and a $\Gal_{\calK}$-action given by $g(\tte_{\bpsi,\sigma}) = \bpsi_{\sigma}(g) \tte_{\bpsi,\sigma}$. Proceeding as in Section \ref{subsec:eu-twisting}, we have a coefficient map
\[
\Psi_{\ttu}: \Lambda^{\#}(1) \rightarrow V_{\bpsi,\sigma}, \quad \lambda \otimes \tte_{1} \mapsto \ev_{\bpsi,\sigma}(\lambda) \tte_{\bpsi,\sigma},
\]
where $\ev_{\bpsi,\sigma}: \Lambda \rightarrow \Lambda_{m}$ is the evaluation map $[g] \mapsto \epsilon(g)\bpsi_{\sigma}(g)^{-1}$. Again, the letter $\ttu$ represents the choice of the bases $\tte_{1}$ and $\tte_{\bpsi, \sigma}$. Combining these maps, we have a composition
\[
\rmH^{1}_{\Iw}(\calK_{\infty}/\calK, \ZZ_p(1)) \simeq  \rmH^1(\calO_{\calK,S},\Lambda^\#(1)) \xrightarrow{(\Phi_{\ttu})_{\ast}} \rmH^{1}(\calO_{\calK,S}, V_{\bpsi,\sigma}).
\]
Applying this map to the elliptic unit $\bfz_{p^{\infty}\frf}$, we obtain an element 
\begin{equation} \label{eq:elliptic-unit-fan-wan-family}
\bfz_{\bpsi,\sigma} \in \rmH^{1}(\calO_{\calK,S}, V_{\bpsi,\sigma}).
\end{equation}
This is \emph{the family of elliptic units} parameterized by $\Lambda_{m}$ used in \cite{FanWan} and also in the main text.

\subsection{Explicit reciprocity laws of elliptic units}\label{subsec:eu-reciprocity}

\subsubsection{The cohomological realization}

For a local field $F/\QQ_p$ and an $L$-representation $V$, its de Rham module and filtration are
\[
 \bfD_{\mathrm{dR}}(F,V) =(\bfB_{\mathrm{dR}}\otimes_{\QQ_p}V)^{G_F},
 \quad \operatorname{Fil}^j \bfD_{\mathrm{dR}}(F,V)
 =(\operatorname{Fil}^j \bfB_{\mathrm{dR}}\otimes_{\QQ_p}V)^{G_F}.
\]
Here $\bfB_{\mathrm{dR}}$ is Fontaine's de Rham period ring. The representation $V$ is de Rham when $\bfD_{\dR}(F,V)$ has dimension $\dim_L V$ over $F\otimes_{\QQ_p}L$. The related background results are summarized in \cite[\S\S9.1--9.2]{kato2004padic}.

Let $\vartheta$ be a Hecke character over $\calK$ with infinity type $(-r,0)$, with $r\geq1$. Kato constructs a one-dimensional vector space over $M$:
\[
V_M(\vartheta) = \rmH^{1}(\scE(\CC), \QQ)^{\otimes_{\calK} r} \otimes_{\calK} M,
\]
as the Betti realization, and a one-dimensional $M$-vector space $S(\vartheta)$ as the de Rham realization (see \cite[page 257]{kato2004padic} for a precise definition), and a complex period map
\begin{equation}\label{eq:eu-period-maps}
 \operatorname{per}_\tau:S(\vartheta)\otimes_{M,\tau}\CC
 \xrightarrow{\ \sim\ }V_M(\vartheta)\otimes_{M,\tau}\CC.
\end{equation}
For $V_L(\vartheta)=V_M(\vartheta)\otimes_{M,\sigma}L$, with its specified Galois action, the de Rham realization satisfies
\[
 \per_{p}^{-1}: \operatorname{Fil}^1 \bfD_{\mathrm{dR}}(\calK_p,V_L(\vartheta))
 =\operatorname{Fil}^r \bfD_{\mathrm{dR}}(\calK_p,V_L(\vartheta))
 \simeq S(\vartheta)\otimes_{M,\sigma}L.
\]
Here $\calK_p=\prod_{\mathfrak p\mid p}\calK_{\mathfrak p}$ and the de Rham space is the direct sum over these primes, as in the article. The left-hand side has dimension one over $L$. The constructions are in \cite[\S15.8]{kato2004padic}.

To specify the Tate coordinate, let $t_{\mathrm{dR}}=\log([\varepsilon])\in \bfB_{\mathrm{dR}}$ and $\ttd_1^{\mathrm{Tate}}=t_{\mathrm{dR}}^{-1}\tte_1$. This is a basis of $\bfD_{\mathrm{dR}}(\QQ_p(1))$. The isomorphism $\bfD_{\dR}(V_{L}(\vartheta)(1)) = \bfD_{\dR}(V_{L}(\vartheta)) \otimes \bfD_{\dR}(\QQ_p(1))$ identifies
\[
 \operatorname{Fil}^0 \bfD_{\mathrm{dR}}(V_L(\vartheta)(1))
 =\operatorname{Fil}^1 \bfD_{\mathrm{dR}}(V_L(\vartheta))
 \otimes_{\calK_p} \calK_p \cdot \ttd_1^{\mathrm{Tate}}.
\]
Then we have the dual exponential map
\[
\exp_{\vartheta}^{\ast}:= \per_{p}^{-1} \circ (\id{} \otimes (\ttd_1^{\Tate})^{\vee}) \circ \exp_{\BK, V_{L}(\vartheta)(1)}: \rmH^{1}(\calK_p, V_{L}(\vartheta)(1)) \rightarrow S(\vartheta) \otimes_{M, \sigma} L,
\]
where $(\ttd_{1}^{\Tate})^{\vee}: \bfD_{\dR}(V_{L}(\vartheta)(1)) \rightarrow \calK_p$ is the linear functional sending $\ttd_{1}^{\Tate}$ to $1$, and $\exp_{\BK, V_{L}(\vartheta)(1)}$ is the Bloch--Kato dual exponential map.

\subsubsection{Kato's explicit reciprocity law}

We choose and fix $\gamma\in V_M(\vartheta)$. Then $\ttu := \gamma \otimes e_{1}$ gives a basis for $V_{L}(\vartheta)(1)$. Following \eqref{eq:eu-indexed-class}, we denote the resulting element $z_{V_{L}(\vartheta(1)), \frf}(\ttu)$ by $z_{\vartheta,\mathfrak f}(\gamma)$. The following theorem appears as \cite[Proposition 15.9]{kato2004padic} with $\calK'=\calK$. \footnote{Its case $r=1$ is restated in \cite[Proposition 3.1]{BurungaleFlach}.}

\begin{theorem}[Kato's explicit reciprocity law]\label{thm:eu-Kato-ERC}
The vector $\exp_\vartheta^*(\operatorname{loc}_p z_{\vartheta,\mathfrak f}(\gamma))$ belongs to $S(\vartheta) \otimes_{M} L$. For each compatible complex coefficient embedding $\tau$,
\begin{equation}\label{eq:eu-vector-ERC}
 \operatorname{per}_\tau\!
 \left(\exp_\vartheta^*(\operatorname{loc}_p z_{\vartheta,\mathfrak f}(\gamma))\right)
 =L_{p\mathfrak f}(\overline{\vartheta_\tau},r)\gamma.
\end{equation}
\end{theorem}

Choose a basis $\ttb$ of $S(\vartheta)$ and define $\Omega_\tau(\ttb,\gamma)$ by
\begin{equation}\label{eq:eu-period-definition}
 \per_\tau(\ttb)=\Omega_\tau(\ttb,\gamma)\gamma.
\end{equation}
For nonzero $\gamma$, \eqref{eq:eu-vector-ERC} is equivalent to
\begin{equation}\label{eq:eu-scalar-ERC}
 \exp_\vartheta^*(\operatorname{loc}_p z_{\vartheta,\mathfrak f}(\gamma))
 =\frac{L_{p\mathfrak f}(\overline{\vartheta_\tau},r)}
 {\Omega_\tau(\ttb,\gamma)}\,\ttb \in L.
\end{equation}

\begin{remark}
If we change basis as $\gamma'=a\gamma$ and $\ttb'=c\ttb$, then $z_{\vartheta, \frf}(\gamma')=a z_{\vartheta, \frf}$, and
\[
\per_{\tau}(\ttb^{\prime}) = \tau(c) \Omega_{\tau}(\ttb,\gamma)\gamma = \tau(c/a) \Omega_{\tau}(\ttb,\gamma)\gamma^{\prime},
\]
and hence $\Omega_\tau(\ttb',\gamma')=\tau(c/a)\Omega_\tau(\ttb,\gamma)$.
\end{remark}

\subsubsection{Explicit reciprocity law for mixed infinity types in \cite{FanWan}}
More generally, consider a Hecke character $\Psi$ over $\calK$ of infinite type $(-r_1, r_2)$, with $r_1\geq1$ and $r_2\geq0$. We put $\Theta=\Psi|\cdot|_{\calK}^{-r_2}$; then its infinite type is $(-r_1-r_2,0)$. We put
\begin{equation}\label{eq:eu-mixed-twist}
 V_M(\Psi) :=V_M(\Theta)(r_2), \quad S(\Psi) := S(\Theta)(r_2).
\end{equation}
Then \cite[Proposition 5.9 and Remark 5.10]{FanWan} gives the corresponding explicit reciprocity law as a generalization of Theorem \ref{thm:eu-Kato-ERC}, with $L_{p\mathfrak f}(\overline\Psi,r_1-r_2)$ replacing $L_{p\frf}(\overline{\vartheta_\tau}, r)$. More precisely, we have the following theorem.

\begin{theorem}[Mixed explicit reciprocity law]\label{thm:eu-fan-wan-ERC}
The vector $\exp_\Psi^*(\operatorname{loc}_p z_{\Psi,\mathfrak f}(\gamma))$ belongs to $S(\Psi) \otimes_{M} L$. For each compatible complex coefficient embedding $\tau$,
\begin{equation}\label{eq:eu-vector-ERC-mixed}
 \operatorname{per}_\tau\!
 \left(\exp_\Psi^*(\operatorname{loc}_p z_{\Psi,\mathfrak f}(\gamma))\right)
 =L_{p\mathfrak f}(\overline{\Psi_\tau},r_1 - r_2)\gamma.
\end{equation}
\end{theorem}
In our application, we use the range $r_1=k$ and $r_2=k-1$ for $k\geq1$. We shall make this explicit in the following arguments, with the goal of proving Theorem \ref{thm:eu-article-ERC}.

\subsubsection{The cohomological realization, (bis)}
For our CM elliptic curve $E/\QQ$ with CM given by $\calO_{\calK}$, the Betti homology $\rmH_{1}(E(\CC),\ZZ)$ is a rank-one free $\calO_{\calK}$-module. We choose an $\calO_{\calK}$-basis $\ttt_{E}$ so that $\rmH_{1}(E(\CC),\ZZ) = \calO_{\calK} \ttt_{E}$.

Let $\mathfrak d_{\calK}^{-1}=\{x\in\calK:\operatorname{Tr}_{\calK/\QQ}(x\calO_{\calK})\subset\ZZ\}$ be the codifferent of $\calK/\QQ$ and let $d_{\calK}$ generate the different. Then we have
\begin{equation}\label{eq:eu-codifferent}
 \rmH^1(E(\CC),\ZZ)=\operatorname{Hom}_{\ZZ}(\rmH_{1}(E(\CC),\ZZ),\ZZ) =\mathfrak d_{\calK}^{-1} \ttt_E^\vee,
\end{equation}
where $\ttt_{E}^{\vee}$ is the linear functional sending $x \ttt_{E} \in \rmH_{1}(E(\CC),\ZZ)$ to $x$. An $\calO_{\calK}$-basis for $\rmH^1(E(\CC),\ZZ)$ is therefore
\begin{equation}\label{eq:eu-gamma}
 \gamma_E=d_{\calK}^{-1} \ttt_E^\vee,\quad \text{with }
 \gamma_E(x \ttt_E)=\operatorname{Tr}_{\calK/\QQ}(x/d_{\calK}).
\end{equation}

Let $\omega_E$ be an invariant differential of $E/\QQ$ fixed in the main text, and let $a_E=\int_{\ttt_E}\omega_E$. Given the embedding $\iota: \calK \hookrightarrow \CC$, integration over $x \ttt_E$ gives $a_E x$, whereas $\gamma_E$ has component $x/d_{\calK}$. Hence
\begin{equation}\label{eq:eu-elliptic-period}
 \operatorname{per}_{E}(\omega_E)=a_Ed_{\calK}\gamma_E,
\end{equation}

\subsubsection{The explicit reciprocity law with explicit bases} \label{sec:ERC-basis}
Now let $\phi_k=(\psi_E^c/\psi_E)^k$ be a Hecke character over $\calK$ of infinite type $(k,-k)$. With $m=2k-1$ and $\Psi_k=(\psi_E\phi_k)^c$, we have $\Psi_k=\psi_E^{\,2k-1}|\cdot|_{\calK}^{k-1}$ with infinite type $(-k,k-1)$: indeed, substituting $\psi_E^c=|\cdot|_{\calK}^{-1}\psi_E^{-1}$ into $\psi_E^c(\psi_E/\psi_E^c)^k$ gives the result.

Consider the Betti cohomology $\rmH^{1}(E(\CC),\QQ)$, which is a one-dimensional $\calK$-vector space. For the fixed complex embedding $\iota: \calK \hookrightarrow \CC$ used to integrate $\omega_{E}$ above, we put
\[
\rmH^{1}_{E, \iota} := \rmH^{1}(E(\CC),\QQ) \otimes_{\calK,\iota} \CC, \quad \rmH^{1}_{E, p,\sigma} := \rmH^{1}_{\et}(E_{\bar{\calK}}, \QQ_p) \otimes_{\calK_p,\sigma} L
\]
after choosing an embedding $\sigma: \calK_p := \calK \otimes \QQ_p \hookrightarrow L$. We put
\begin{equation}\label{eq:eu-WV-k}
 W_k := (\rmH^{1}_{E,p,\sigma})^{\otimes_{\calK}(2k-1)}(k-1),\qquad
 V_k=W_k(1)=(\rmH^{1}_{E,p,\sigma})^{\otimes_{\calK}(2k-1)}(k).
\end{equation}

We start with the multiplicative group $\GG_{\rmm} := \Spec \QQ[T,T^{-1}]$. We let $\eta \in \rmH^{1}(\CC^{\times}, \QQ)$ with $\eta(T) = 1$, and on the de Rham side, we put $\omega_{1} := \dif T / T$. Integration gives
\[
\per_{\GG_{\rmm}}(\omega) = \left(\int_{T} \dfrac{\dif T}{T} \right) \eta = 2 \pi i \eta.
\]
We then define the bases for the dual $\QQ(1)$: 
\begin{itemize}
    \item let the Betti basis $\tte_{1}^{\rmB}$ be the element of $\rmH^{1}(\CC^{\times}, \QQ)^{\vee}$ satisfying $\tte_{1}^{\rmB}(\eta) = 1$;
    \item let the de Rham basis $\ttb_{1}^{\Tate}$ be the element of $\rmH^{1}_{\dR}(\GG_{\rmm}/\QQ)^{\vee}$ satisfying $\ttb_{1}^{\Tate}(\dif T/T)$ = 1.
\end{itemize}
Then we have $\per_{\QQ(1)}(\ttb_{1}^{\Tate}) = (2 \pi i)^{-1} \tte_{1}^{\rmB}$. For $j \geq 0$, we put $\tte_{j}^{\rmB} := (\tte_{1}^{\rmB})^{\otimes j}$ and $\ttb_{j}^{\Tate} = (\ttb_{1}^{\Tate})^{\otimes j}$. Then we have 
\[
 \operatorname{per}_{\QQ(j)}(\ttb_j^{\mathrm{Tate}})=(2\pi i)^{-j}\tte_j^B.
\]
Their $p$-adic counterparts are the fixed $\tte_j$ and $\ttd_j^{\mathrm{Tate}}:= t_{\mathrm{dR}}^{-j} \tte_j$. We define
\begin{equation}\label{eq:eu-tensor-bases}
 \gamma_{E,k}= \gamma_E^{\otimes(2k-1)}\otimes \tte_{k-1}^{\rmB} \in (\rmH^{1}_{E,\iota})^{\otimes (2k-1)}(k-1),\quad 
 \ttb_{E,k}=\omega_E^{\otimes(2k-1)}\otimes \ttb_{k-1}^{\mathrm{Tate}} \in W_{k}.
\end{equation}
Tensoring \eqref{eq:eu-elliptic-period} gives
\begin{equation}\label{eq:eu-tensor-period}
 \operatorname{per}_{E}(\ttb_{E,k})
 =\frac{(a_Ed_{\calK})^{2k-1}}{(2\pi i)^{k-1}}\gamma_{E,k}.
\end{equation}

\begin{theorem}[Explicit reciprocity law with explicit bases]\label{thm:eu-article-ERC}
For $k\geq1$, let $z_{\Psi_{k},\frf}(\gamma_{E,k})$ be the class indexed by $\gamma_{E,k}$, with modulus $\mathfrak f$ containing the conductor of $\Psi_k$. Then
\begin{equation}\label{eq:eu-article-ERC}
 \exp_{\Psi_{k}}^*(\operatorname{loc}_p z_{\Psi_{k},\frf}(\gamma_{E,k}))
 =L_{p\mathfrak f}(\psi_E\phi_k,1)
 \frac{(2\pi i)^{k-1}}{(a_Ed_{\calK})^{2k-1}} \ttb_{E,k}.
\end{equation}
\end{theorem}
\begin{proof}
Applying the explicit reciprocity law to $\Psi_k$ in Theorem \ref{thm:eu-fan-wan-ERC}, with $r_1=k,r_2=k-1$, and noting that $\overline{\Psi_{k}} = \psi_{E}\phi_k$, we have
\[
 \operatorname{per}_{E}\bigl(\exp_k^*(\operatorname{loc}_p z_{\Psi_{k},\frf}(v_{E,k}))\bigr)
 =L_{p\mathfrak f}(\psi_E\phi_k,1)\gamma_{E,k}.
\]
Write the vector inside $\operatorname{per}_{E}$ as $x_k \ttb_{E,k}$ and use \eqref{eq:eu-tensor-period}. It gives
\[
 x_k\frac{(a_Ed_{\calK})^{2k-1}}{(2\pi i)^{k-1}}
 =L_{p\mathfrak f}(\psi_E\phi_k,1).
\]
Solving $x_{k}$ proves the explicit reciprocity law in the theorem.
\end{proof}

%% file: AppBC-calculations.tex
\section{Computations on special values of $\psi_{E,p}$}\label{sec:special values calculation}

\subsection{Special CM families}

Recall that we studied the families $X_a:y^2=x^3+a$ and $Y_b:y^2=x^3+bx$ with $a,b\in\QQ^{\times}$ in Section \ref{sst:practice}. Here we make a slight generalization.

For all $9$ possible CM fields $\calK$, or more precisely all $13$ possible CM orders $\mathcal{O}_D=\ZZ+f\mathcal{O}_{\calK}$ with conductor $f$ and discriminant $D=f^2D_{\calK}$, we define an integral representative
\[
E_D:y^2=x^3+A_Dx+B_D
\]
with specified constants $A_D$ and $B_D$ in the following table, such that each $E_D$ has CM by $\mathcal{O}_D$ and the integral pair $(A_D,B_D)$ is minimal under the coordinate-wise partial order\footnote{$(A,B)\leq(A^{\prime},B^{\prime})\iff\abs{A}\leq\abs{A^{\prime}},\abs{B}\leq\abs{B^{\prime}}$.}. Here $n_D:=\#\Aut_{\bar{\QQ}}(E_D)=\#\mathcal{O}_D^{\times}$. 

\begingroup\small\setlength{\tabcolsep}{10pt}\renewcommand{\arraystretch}{1.25}
\begin{longtable}{@{\hspace{5pt}}ccccccc@{\hspace{5pt}}}
\caption{Models of all $13$ $j$-invariants}\label{tab:j-models}\\
\hline
${\calK}$ & $f$ & $D$ & $j_D$ & $A_D$ & $B_D$ & $n_D$\\\hline\endfirsthead
\multicolumn{7}{c}{\tablename\ \thetable\ (continued)}\\\hline
${\calK}$ & $f$ & $D$ & $j_D$ & $A_D$ & $B_D$ & $n_D$\\\hline\endhead
\hline\endfoot
$\mathbb Q(i)$ & $1$ & $-4$ & $1728$ & $1$ & $0$ & $4$\\
$\mathbb Q(i)$ & $2$ & $-16$ & $287496$ & $-11$ & $14$ & $2$\\
$\mathbb Q(\sqrt{-2})$ & $1$ & $-8$ & $8000$ & $-30$ & $56$ & $2$\\
$\mathbb Q(\sqrt{-3})$ & $1$ & $-3$ & $0$ & $0$ & $1$ & $6$\\
$\mathbb Q(\sqrt{-3})$ & $2$ & $-12$ & $54000$ & $-15$ & $22$ & $2$\\
$\mathbb Q(\sqrt{-3})$ & $3$ & $-27$ & $-12288000$ & $-120$ & $506$ & $2$\\
$\mathbb Q(\sqrt{-7})$ & $1$ & $-7$ & $-3375$ & $-35$ & $-98$ & $2$\\
$\mathbb Q(\sqrt{-7})$ & $2$ & $-28$ & $16581375$ & $-595$ & $5586$ & $2$\\
$\mathbb Q(\sqrt{-11})$ & $1$ & $-11$ & $-32768$ & $-264$ & $1694$ & $2$\\
$\mathbb Q(\sqrt{-19})$ & $1$ & $-19$ & $-884736$ & $-152$ & $722$ & $2$\\
$\mathbb Q(\sqrt{-43})$ & $1$ & $-43$ & $-884736000$ & $-3440$ & $77658$ & $2$\\
$\mathbb Q(\sqrt{-67})$ & $1$ & $-67$ & $-147197952000$ & $-29480$ & $1948226$ & $2$\\
$\mathbb Q(\sqrt{-163})$ & $1$ & $-163$ & $-262537412640768000$ & $-8697680$ & $9873093538$ & $2$\\
\end{longtable}
\endgroup

By the theory of elliptic curves over $\CC$, each $E/\QQ$ with CM by $\mathcal{O}_D$ is $\CC$-isomorphic to $E_D$. In fact, these isomorphisms can be taken to be algebraic. For $d\in\QQ$, define the twisted curves
\[
E_D^d:
\begin{cases}
    y^2=x^3+d, & D=-3\\
    y^2=x^3+dx, & D=-4\\
    y^2=x^3+A_Dd^2x+B_Dd^3, & \text{otherwise.}
\end{cases}
\]
In particular, $E_D=E_D^1$. See also \cite[Chapter X, \S5, Proposition 5.4]{Silverman}

\begin{proposition}
    Let $E/\QQ$ be an elliptic curve with CM by $\mathcal{O}_D$, then $E\cong_{\QQ}E_D^d$ for some $d\in\QQ^{\times}$.
\end{proposition}
\begin{proof}
    Write $E:y^2=x^3+Ax+B$ after some $\QQ$-isomorphism. Since $E$ has CM by $\mathcal{O}_D$, it follows that $j(E)=j_D=j(E_D)$. Recall
    \[
    j(E)=1728\frac{4A^3}{4A^3+27B^2}.
    \]
    If $j(E)=0$ or $1728$, i.e. $A=0$ or $B=0$, the result follows directly. Otherwise we have $B_D^2/B^2=A_D^3/A^3$. Take $d=\frac{A_DB}{AB_D}$ and thus $d^2=\frac{A}{A_D}$, $d^3=\frac{B}{B_D}$. We conclude that $E\cong_{\QQ}E_D^d$.
\end{proof}

For each $d\in\QQ^{\times}$, choose an $n_D$-th root of $d$, say $t$, and define
\begin{align*}
    \phi_D^t:E_D=E_D^1 & \rightarrow E_D^d\\
    (x,y) & \mapsto(t^2x,t^3y).
\end{align*}
It is an isomorphism defined over $\QQ(t)$. In particular,
\[
d/d^{\prime}\in(\QQ^{\times})^{n_D}\iff E_D^d\cong_{\QQ}E_D^{d^{\prime}}.
\]
Thus, $d\in\QQ^{\times}/(\QQ^{\times})^{n_D}$ parametrizes all $E/\QQ$ with CM by $\mathcal{O}_D$.

\subsection{Norm residue symbols and CM Hecke characters}

Recall the global and local reciprocity maps
\[
\rec_{\calK}:\AA_{\calK}^{\times}/{\calK}^{\times}\rightarrow\Gal({\calK}^{\ab}/{\calK}),\quad\rec_{{\calK}_{\mathfrak{p}}}:{\calK}_{\mathfrak{p}}^{\times}\rightarrow\Gal({\calK}_{\mathfrak{p}}^{\ab}/{\calK}_{\mathfrak{p}})
\]
from class field theory. Assume $\mu_n\subset\calK$. Then reciprocity gives the $n$-th norm residue symbol\footnote{Some authors use different conventions, for example $(y,b)_{n,\calK_{\mathfrak{p}}}:=\frac{\rec_{\calK_{\mathfrak{p}}}(b)(\sqrt[n]{y})}{\sqrt[n]{y}}$.}
\[
H_n(x,a):=\frac{\rec_{\calK}(x)(\sqrt[n]{a})}{\sqrt[n]{a}}\in\mu_n,\quad (y,b)_{n,\calK_{\mathfrak{p}}}:=\frac{\rec_{\calK_{\mathfrak{p}}}(y)(\sqrt[n]{b})}{\sqrt[n]{b}}\in\mu_n,
\]
for any $x\in\AA_{\calK}^{\times}$, $a\in\calK^{\times}$ and $y,b\in\calK_{\mathfrak{p}}^{\times}$.

\begin{theorem}[CM comparison]\label{thm:CM comparison}
    For any $x\in\AA_{\calK}^{\times}$ and $y\in\calK_{\mathfrak{p}}^{\times}$, one has
    \[
    \psi_{E_D^d}(x)=H_{n_D}^{-1}(x,d)\psi_{E_D^1}(x),\quad \psi_{E_D^d,\mathfrak{p}}(y)=(y,d)_{n_D,\calK_{\mathfrak{p}}}^{-1}\psi_{E_D^1,\mathfrak{p}}(y).
    \]
\end{theorem}
\begin{proof}
    By \cite[Chapter II, Theorem 9.1]{silverman2013advanced}, for any $x\in\AA^{\times}_\calK$, the finite part of $\psi_{E}(x)$, say $\alpha_{E}(x)\in \calK^{\times}$, is characterized as the unique element satisfying
    \begin{enumerate}
        \item $\alpha_{E}(x)\mathcal{O}_\calK=(x)$
        \item for any fractional ideal $\mathfrak{a}\subset \calK$ and any analytic isomorphism\footnote{Here by an 'analytic isomorphism', in \cite{silverman2013advanced}, the author only meant a bi-holomorphic map. However it is equivalent to requiring $u$ to be an analytic group isomorphism, since any analytic isomorphism is a composition of an analytic group isomorphism and a translation, while a translation commutes with Galois action. So hereafter we understand 'analytic isomorphism' as an analytic group isomorphism.} $u:\CC/\mathfrak{a}\rightarrow E(\CC)$, there is a commutative diagram
        \begin{equation}\label{eq:CM psi characterization}
            \begin{tikzcd}
                \calK/\mathfrak{a} \arrow[r, "\alpha_{E}(x)x^{-1}"] \arrow[d, "u"] & \calK/\mathfrak{a} \arrow[d, "u"] \\
                E(\calK^{\ab}) \arrow[r, "\rec_{\calK}(x)"]                                & E(\calK^{\ab}).
            \end{tikzcd}
        \end{equation}
    \end{enumerate}

    Assume we already know $\alpha_{E_D^1}(x)$, we show that $H_{n_D}^{-1}(x,d)\alpha_{E_D^d}(x)$ satisfies the requirements above. For any analytic isomorphism $u:\CC/\mathfrak{a}\rightarrow E_D^d(\CC)$, we consider $(\phi_D^t)^{-1}\circ u:\CC/\alpha\rightarrow E_D^1(\CC)$, which is also analytic isomorphic. By definition of $\alpha_{E_D^1}(x)$, one has
    \begin{equation}\label{eq:universal property psi}
        (\phi_D^t)^{-1}(u (\alpha_{E_D^1}(x)x^{-1}z))=\rec_{\calK}(x)((\phi_D^t)^{-1}(u(z))).
    \end{equation}
    
    It is easy to verify that $u(\zeta z)=[\zeta](u(z))=\phi_D^{\zeta^{-1}}(u(z))$ for any $\zeta\in\mu_{n_D}\cong\Aut(E_D^d)$, and that
    \begin{equation}\label{eq:rec phi commutes}
        \rec_{\calK}(x)\circ\phi_D^t=\phi_D^{\rec_{\calK}(x)(t)}\circ\rec_{\calK}(x):E_D^1(\calK^{\ab})\rightarrow E_D^d(\calK^{\ab}).
    \end{equation}
    
    Since $H_{n_D}(x,d)\in\mu_{n_D}$, it follows that
    \begin{align*}
        u(H_{n_D}^{-1}(x,d)\alpha_{E_D^1}(x)x^{-1}z)&=\phi_D^{H_{n_D}(x,d)}(u(\alpha_{E_D^1}(x)x^{-1}z))\\
        &=\phi_D^{H_{n_D}(x,d)}(\phi_D^t(\rec_{\calK}(x)((\phi_D^t)^{-1}(u(z)))))\tag{by \eqref{eq:universal property psi}}\\
        &=\phi_D^{H_{n_D}(x,d)}(\phi_D^t((\phi_D^{\rec_{\calK}(x)(t)})^{-1}(\rec_{\calK}(x)(u(z)))))\tag{by \eqref{eq:rec phi commutes}}\\
        &=\phi_D^{H_{n_D}(x,d)}(\phi_D^{H_{n_D}^{-1}(x,d)}(\rec_{\calK}(x)(u(z))))\\
        &=\rec_{\calK}(x)(u(z)).
    \end{align*}
    Hence $H_n^{-1}(x,d)\alpha_{E_D^1}(x)$ fits into Diagram \eqref{eq:CM psi characterization} for $E=E_D^d$. By uniqueness we have $\alpha_{E_D^d}(x)=H_n^{-1}(x,d)\alpha_{E_D^1}(x)$. Since $\psi_{E_D^d,\infty}(x)=x^{-1}=\psi_{E_D^1,\infty}(x)$, the first conclusion follows. The local version follows from the local-global compatibility.
\end{proof}

\begin{remark}
    Note that this result applies for $E_D^d$ with \emph{any} rank, without adding any rank $1$ hypothesis!
\end{remark}

The following are a few lemmas computing norm residue symbols, whose proofs we omit.

\begin{lemma}[Tame symbol, {\cite[Chapter XIV, \S 3, Proposition 8 and Corollary]{serre2013local}}]\label{lem:calculation local norm residue symbol}
    Let $F$ be a finite extension of $\QQ_p$ with residue field $k$ of order $q$. Assume that $\chari(k)\nmid n$ and $F$ contains $\mu_n$, equivalently $n\mid q-1$. Then for $s,t\in F^{\times}$, one has
    \[
    (s,t)_{n,F}=\Teich\left((-1)^{v_F(s)v_F(t)}\frac{s^{v_F(t)}}{t^{v_F(s)}}\right)^{-(q-1)/n}\in\mu_n.
    \]
\end{lemma}

\begin{lemma}[Second Hilbert symbol, {\cite[Chapter XIV, \S 4]{serre2013local}}]\label{lem:calculation local norm residue symbol n=p=2}
    When $F=\QQ_2$ and $n=2$, one has
    \[
    (2,2)_{2,\QQ_2}=1,\quad (c,d)_{2,\QQ_2}=(-1)^{(a-1)(b-1)/4},\quad (2,d)_{2,\QQ_2}=(-1)^{(d^2-1)/8},
    \]
    for any $c,d\in\ZZ_2^{\times}$.
\end{lemma}

\begin{lemma}[$p$-th Hilbert symbol with $p$ odd, see for example {\cite[Exercise 2.13]{cassels1967algebraic}}]\label{lem:calculation Hilbert symbol}
    For an odd prime $p$, let $F=\QQ_p(\zeta_p)$ with uniformizer $\lambda:=1-\zeta_p$. Define $\eta_i:=1-\lambda^i\in U^i$ for $i\geq1$. One has
    \begin{enumerate}
        \item if $x\in U^{p+1}$, then $x\in (F^{\times})^{p}$, hence $(x,y)_{p,F}=1$ for any $y\in F^{\times}$.
        \item if $i+j\geq p+1$, then $(a,b)_{p,F}=1$ for any $a\in U^i$ and $b\in U^j$.
        \item $(\eta_i,\eta_j)_{p,F}=(\eta_i,\eta_{i+j})_{p,F}(\eta_{i+j},\eta_j)_{p,F}(\eta_{i+j},\lambda)_{p,F}^{-j}.$
        \item $(\eta_i,\lambda)_{p,F}=
        \begin{cases}
            1,& 1\leq i\leq p-1\\
            \zeta_p^{-1},& i=p.
        \end{cases}
        $
    \end{enumerate}
\end{lemma}

\subsection{The inert case}

Here we compute the special values of $\psi_{E,p}$ in Tables \ref{tab:p_3 inert database}, \ref{tab:p_2 inert database}, \ref{tab:p_5 inert database} and \ref{tab:p_7 inert database}. We focus on the family $X_a=E_{-3}^a$ and $Y_b=E_{-4}^b$.

\begin{proposition}[Proposition \ref{prop:calculation 3,5,7 inert}]
    If $p$ is inert in $\calK$, let $g_p$ denote the Teichmuller lift of a generator of $(\mathcal{O}_F/p\mathcal{O}_F)^{\times}\cong\FF_{p^2}^{\times}$ whose image generates $G_1$.
    \begin{enumerate}
        \item In the case $(p,\calK)=(5,\QQ(\sqrt{-3}))$, choose $g_5:=\Teich(2+\bar{\omega})$, then $\psi_{X_a,5}(g_5)=\zeta_6^{-v_5(a)}$.
        \item In the case $(p,\calK)=(7,\QQ(\sqrt{-1}))$, choose $g_7:=\Teich(2+\bar{i})$, then $\psi_{Y_b,7}(g_7)=i^{v_7(b)}$.
        \item In the case $(p,\calK)=(3,\QQ(\sqrt{-1}))$, choose $g_3:=\Teich(1-\bar{i})$, then $\psi_{Y_b,3}(g_3)=i^{v_3(b)}$.
    \end{enumerate}
\end{proposition}
\begin{proof}
    In the case $(p,\calK)=(5,\QQ(\sqrt{-3}))$, $F=\QQ_{5^2}$, we first consider $X_1:y^2=x^3+1$ with conductor $2^2\cdot3^2$. Since $p=5$ is good for $X_1$, it follows that $\psi_{X_1,5}$ is trivial on $\calO_F^{\times}=\ZZ_{5^2}^{\times}$. For a general $X_a$, we use the comparison Theorem \ref{thm:CM comparison} and deduce that $\psi_{X_a,5}(g_5)=(g_5,a)_{6,F}^{-1}$. By Lemma \ref{lem:calculation local norm residue symbol} we have
    \[
    (g_5,a)_{6,F}=g_5^{-4v_5(a)}.
    \]
    Since $g_5^4=\Teich((1+\bar{\omega})^4)=\Teich(-\bar{\omega})=-\omega=\zeta_6^{-1}$, it follows that $\psi_{X_a,5}(g_5)=\zeta_6^{-v_5(a)}$.

    For $Y_b$, we consider first $Y_1:y^2=x^3+x$ with conductor $2^6$. So $\psi_{Y_1,p}$ is trivial on $\mathcal{O}_F^{\times}$ for $p=3,7$, and thus by CM comparison Theorem \ref{thm:CM comparison},
    \[
    \psi_{Y_b,p}(g_p)=(g_p,b)_{4,F}^{-1}=g_p^{v_p(b)(p^2-1)/4}.
    \]
    Since $g_7^{12}=\Teich((2+\bar{i})^{12})=\Teich(\bar{i})=i$ and $g_3^2=\Teich((1-\bar{i})^2)=\Teich(\bar{i})=i$, the desired results follow.
\end{proof}

\begin{proposition}[Proposition \ref{prop:calculation 2 inert}]
    For $X_a:y^2=x^3+a$ with $a\in\QQ^{\times}$, write $a=2^m u$ with $u\in\ZZ_2^{\times}$ and
    $e=(u-1)/2\pmod2$. Recall that in $F=\QQ_2(\sqrt{-3})$ we defined $f_2:=-\omega-2$ and $h_2:=1+4\omega$. Under the identifications as in Section \ref{sst:calculation}, the special values are
    \[
    \psi_{X_a,2}(f_2)=\omega^{-m+1}(-1)^{e+m},\qquad \psi_{X_a,2}(h_2)=(-1)^m.
    \]
\end{proposition}
\begin{proof}
    We consider the curve $X_{1/4}:y^2=x^3+\frac{1}{4}$, with conductor $3^3$. Combining the equations for $X_a$ and $X_{1/4}$ in Theorem \ref{thm:CM comparison}, we obtain the relative version $\psi_{X_a,2}(x)=(x,4a)_{6,F}^{-1}\psi_{X_{1/4},2}(x)$. Since $X_{1/4}$ has good reduction at $2$, it follows that $\psi_{X_a,2}(x)=(x,4a)^{-1}_{6,F}$. To determine $z:=(x,4a)_{6,F}$, we consider $z^2=(x,4a)_{3,F}$ and $z^3=(x,4a)_{2,F}$.

    On the one hand, by Lemma \ref{lem:calculation local norm residue symbol} we have $(f_2,4a)_{3,F}=\Teich(\bar{\omega})^{-v_2(4a)}=\omega^{-m-2}$ and $(h_2,4a)_{3,F}=\Teich(\bar{1})^{-v_2(4a)}=1$. On the other hand, $(x,4a)_{2,F}=(x,a)_{2,F}=(\Norm^F_{\QQ_2}(x),a)_{2,\QQ_2}$ by compatibility of local reciprocity. Thus
    \[
    (f_2,4a)_{2,F}=(3,2^mu)_{2,\QQ_2}=(-1)^m\cdot(3,u)_{2,\QQ_2}=(-1)^{m+(u-1)/2}=(-1)^{m+e}
    \]
    and
    \[
    (h_2,4a)_{2,F}=(13,2^mu)_{2,\QQ_2}=(-1)^m\cdot(13,u)_{2,\QQ_2}=(-1)^{m+3(u-1)}=(-1)^m.
    \]
    Therefore $(f_2,4a)_{6,F}=(f_2,4a)_{2,F}/(f_2,4a)_{3,F}=\omega^{m+2}(-1)^{e+m}$ and $(h_2,4a)_{6,F}=(-1)^m$. We conclude by $\psi_{X_a,2}(x)=(x,4a)^{-1}_{6,F}$ that
    \[
    \psi_{X_a,2}(f_2)=\omega^{-m+1}(-1)^{e+m},\qquad
    \psi_{X_a,2}(h_2)=(-1)^m,
    \]
    as desired.
\end{proof}

\subsection{The ramified case}

Now we compute the special values of $\psi_{E,p}$ in the ramified case, namely in Tables \ref{tab:p>3 ramified database}, \ref{tab:p=3 ramified database}, \ref{tab:p=2 ramified database D_K=-4} and \ref{tab:p=2 ramified database D_K=-8}. We again denote $\calK_{\mathfrak{p}}$ by $F$, where $\mathfrak{p}$ is the unique prime above $(p)$. 

First we deal with $p>3$ ramified in $\calK=\QQ(\sqrt{-p})$, i.e. the curves $E_p$\footnote{By a slight abuse of notation, we denote $E_{-p}^1$ by $E_p$ and $E_{-p}^d$ by $E_p^d$ for the $p>3$ ramified case.} and their twists $E_p^d$ in Table \ref{tab:p>3 ramified database}.

\begingroup\small\setlength{\tabcolsep}{10pt}\renewcommand{\arraystretch}{1.25}
\begin{longtable}{@{\hspace{5pt}}cccc@{\hspace{5pt}}}
\caption{$E_p$ and their quadratic twists $E_p^d$}\label{tab:p>3 ramified twisted databse}\\
\hline
$p$ & Equation of $E_p$ & LMFDB label & $d$\\\hline\endfirsthead
\multicolumn{4}{c}{\tablename\ \thetable\ (continued)}\\\hline
$p$ & Equation of $E_p$ & LMFDB label & $d$\\\hline\endhead
\hline\endfoot
    $7$ & $y^2+xy+y=x^3-x^2-20x+46$ & 441.c4 & $4/3$\\
    $11$ & $y^2+y=x^3-x^2-7x+10$ & 121.b2 & $2$\\
    $19$ & $y^2+y=x^3-38x+90$ & 361.a2 & $2$\\
    $43$ & $y^2+y=x^3-860x+9707$ & 1849.b2 & $5$\\
    $67$ & $y^2+y=x^3-7370x+243528$ & 4489.b2 & $2$\\
    $163$ & $y^2+y=x^3-2174420x+1234136692$ & 26569.a2 & $2$\\
\end{longtable}\endgroup

\begin{theorem}
    For all $p=7,11,19,43,67,163$, we have $\psi_{E_p,p}(\pi)=\pi$ and $\psi_{E_p^d,p}(\pi)=-\pi$.
\end{theorem}
\begin{proof}
    We compute $(\pi,d)_{2,F}=(d,\pi)_{2,F}$, where $\pi=\sqrt{-p}$ and $d=4/3,2,2,5,2,2$ for $p=7,11,19,43,67,163$ respectively. Since $d$ is a $p$-adic unit, by Lemma \ref{lem:calculation local norm residue symbol} it follows that
    \[
    (\pi,d)_{2,F}=\Teich(\bar{d})^{-(p-1)/2}=\Teich(\bar{d})^{(p-1)/2}=\Teich(\bar{d}^{(p-1)/2}).
    \]
    But one may check that $\bar{d}^{(p-1)/2}\equiv-1\mod{p}$ in all cases. Therefore we obtain that $\psi_{E_p^d,p}(\pi)=-\psi_{E_p,p}(\pi)$. It remains only to verify that $\psi_{E_p,p}(\pi)=\pi$.
    
    Recall that we have normalization $\psi_{E_p}(x)=1$ for $x\in \calK^{\times}\subset\AA_{\calK}^{\times}$. Since $E_p$ is good outside $p$ for $p=11,19,43,67,163$ and $\pi$ is a $q$-adic unit for $q\neq p$, it follows that $\psi_{E_p,p}(\pi)=\psi_{E_p,\infty}(\pi)^{-1}=\pi$. For $p=7$, since $E_7$ has conductor $3^2\cdot7^2$, we only have $\psi_{E_7,3}(\pi)\cdot\psi_{E_7,7}(\pi)=\pi$.
    
    Next we compute $\psi_{E_7,3}(\pi)$. Note that $3$ is inert in $\calK=\QQ(\sqrt{-7})$, so $\calK_3$ is the unique unramified extension of $\QQ_3$ of degree $2$. Moreover $\pi\equiv (1-\pi)^2\in\mathcal{O}_{\calK_3}/{3\mathcal{O}_{\calK_3}}\cong\FF_9$ and $2(1-\bar{\pi})\neq0\in\FF_9$ since $\pi^2\equiv-1\in\FF_9$, it follows by Hensel's lemma that $\pi$ is a square in $\mathcal{O}_{\calK_3}^{\times}$. Because $\psi_{E_7,3}(\mathcal{O}_{\calK_3}^{\times})$ takes value in $\mathcal{O}_{\calK}^{\times}=\{\pm1\}$, this implies $\psi_{E_7,3}(\pi)=1$. So we deduce that $\psi_{E_7,7}(\pi)=\pi$, and the theorem follows.
\end{proof}

Next we consider the case $p=3$ ramified in $\calK=\QQ(\sqrt{-3})$, i.e. cases in Table \ref{tab:p=3 ramified database}. Recall that we took $\pi:=\sqrt{-3}$ to be a uniformizer of $F=\calK_3$, embedded into the upper half plane, and $t:=1+3\pi$.

\begin{proposition}[Proposition \ref{prop:calculation 3 ramified}]
    For $X_a:y^2=x^3+a$, write $a=3^{v_3(a)}u$ with $u\in\ZZ_3^{\times}$. Then we have
    \[
    \psi_{X_a,3}(\pi)=(-1)^{v_3(a)}\left(\frac{u}{3}\right)\pi,\quad \psi_{X_a,3}(\omega)=\omega^{(u^2-1)/3},\quad \psi_{X_a,3}(t)=\omega^{v_3(a)}.
    \]
\end{proposition}
\begin{proof}
    We first compute $(x,a)_{6,F}=(x,a)_{2,F}/(x,a)_{3,F}$ for $x=\pi,\omega,t$ and $a\in\QQ^{\times}$. Note that $v_F(a)=2v_3(a)$ since $F/\QQ_3$ is ramified of degree $2$. By Lemma \ref{lem:calculation local norm residue symbol} we obtain
    \[
    (\pi,a)_{2,F}=\Teich((-1)^{v_3(a)}u)=(-1)^{v_3(a)}\left(\frac{u}{3}\right),\quad (\omega,a)_{2,F}=1,\quad (t,a)_{2,F}=1,
    \]
    where we use $\Teich(\bar{\omega})=\Teich(\bar{1})=\Teich(\bar{t})=1$ and the fact that $v_F(a)$ is even since $a\in\QQ^{\times}$.

    To calculate $(x,a)_{3,F}$, we need Lemma \ref{lem:calculation Hilbert symbol} on the $p$-th Hilbert symbol. In our notation, $\eta_1=\omega$ and $\pi=\lambda/\eta_1^2$. We have $(\eta_1,\pi)_{3,F}=(\eta_1,\lambda)_{3,F}/(\eta_1,\eta_1)_{3,F}^2=1$ by Lemma \ref{lem:calculation Hilbert symbol}. Here we used $(x,x)_{3,F}=1$ since $(x,x)^2_{3,F}=1$ and $(x,x)_{3,F}\in\mu_3$. Additionally, $(\eta_1,\eta_2)=\omega^2$.
    
    Write $a=3^{v_3(a)}u=(-1)^{v_3(a)}\pi^{2v_3(a)}u$ with $u\in\QQ^{\times}\cap\ZZ_3^{\times}$. Note that $u^2\equiv1\mod{3}$, hence $u^2\in U^2$. We thus compute by
    \[
    (x,a)_{3,F}=(x,a)_{3,F}^{-2}=(x,\pi)^{-v_3(a)}_{3,F}(x,u^2)_{3,F}^{-1}.
    \]
    Therefore, the target values can be decomposed and simplified as
    \begin{align*}
        (\pi,a)_{3,F}&=(\lambda,u^2)_{3,F}^{-1}(\eta_1,u^2)_{3,F}^{-1}\\
        (\omega,a)_{3,F}&=(\eta_1,u^2)_{3,F}^{-1}\\
        (t,a)_{3,F}&=(t,\lambda)^{-v_3(a)}_{3,F}(t,\eta_1)^{-v_3(a)}_{3,F}(t,u^2)_{3,F}^{-1}=(t,\lambda)_{3,F}^{-v_3(a)},
    \end{align*}
    where the last equation used the second assertion in Lemma \ref{lem:calculation Hilbert symbol}.

    For the terms containing $x=u^2$ or $t$, we use $3$-adic logarithm to decompose $x$ as $x=\eta_2^s\eta_3^tv$ with $v\in U^4$. Recall that $\log(x)\equiv x-1 \mod{\lambda^4}$ for $x\in U^2$. In particular, one has $\log(\eta_2)\equiv-\lambda^2\mod{\lambda^4}$ and $\log(\eta_3)\equiv-\lambda^3\mod{\lambda^4}$. Since
    \[
    \log(u^2)\equiv u^2-1\equiv-\frac{u^2-1}{3}\lambda^2-\frac{u^2-1}{3}\lambda^3\mod{\lambda^4},\quad\log(t)\equiv3\pi\equiv-\lambda^3\mod{\lambda^4},
    \]
    it follows that $u^2\in\eta_2^{\frac{u^2-1}{3}}\eta_3^{\frac{u^2-1}{3}}U^4$ and $t\in\eta_3 U^4$. Substituting these expressions into the simplifications above and using Lemma \ref{lem:calculation Hilbert symbol} again, we obtain
    \[
    (\pi,a)_{3,F}=1,\quad (\omega,a)_{3,F}=\omega^{(u^2-1)/3},\quad (t,a)_{3,F}=\omega^{v_3(a)}.
    \]
    Together with the values $(x,a)_{2,F}$, we conclude by $(x,a)_{6,F}=(x,a)_{2,F}/(x,a)_{3,F}$ that
    \[
    (\pi,a)_{6,F}=(-1)^{v_3(a)}\left(\frac{u}{3}\right),\quad (\omega,a)_{6,F}=\omega^{-(u^2-1)/3},\quad (t,a)_{6,F}=\omega^{-v_3(a)}.
    \]
    
    By the local CM comparison Theorem \ref{thm:CM comparison}, we need only verify the formulas for one curve. We take $X_{1/4}:y^2=x^3+1/4$ as the base curve, since $X_{1/4}$ is only bad at $3$. Uniformization $\psi_{X_a(x)}=x^{-1}\prod_{\mathfrak{q}}\psi_{X_a,\mathfrak{q}}(x)=1$ for $x\in\calK^{\times}$ gives $\psi_{X_{1/4},3}(\pi)=\pi$ and $\psi_{X_{1/4},3}(\omega)=\omega$, as expected.
    
    For $t=1+3\sqrt{-3}$, notice that $t=2\rho$ with $\rho=2+3\omega$ of norm $7$ and $7\mathcal{O}_{\calK}=(\rho)(\bar{\rho})$. By uniformization we have
    \[
    \psi_{X_{1/4},3}(2)\psi_{X_{1/4},2}(2)=2,\quad \psi_{X_{1/4},3}(\rho)\psi_{X_{1/4},(\rho)}(\rho)=\rho.
    \]
    Since $2$ is inert in $\calK$, it follows that $2$ is a uniformizer, so $\psi_{E_{1/4},2}(2)=-2$ by normalization $\psi_{E,p}\vert_{\QQ_p^{\times}}=\eta_{F/\QQ_p}\cdot\abs{\cdot}_{p}^{-1}$ for inert primes $p$. Computation gives $a_7(E_{1/4})=-1$ so $\psi_{X_{1/4},(\rho)}(\rho)$ is a root of $x^2+x+7=0$, namely $-\rho$ or $-\bar{\rho}$. But $-\bar{\rho}$ is not a generator of $(\rho)$. It follows that $\psi_{E_{1/4},(\rho)}(\rho)=-\rho$. Combining these facts gives
    \[
    \psi_{X_{1/4},3}(t)=t/((-2)(-\rho))=1,
    \]
    as expected, completing the proof.
\end{proof}

For the case $p=2$ ramified in $\calK=\QQ(i)$, i.e. the curves in Table \ref{tab:p=2 ramified database D_K=-4}, recall that we defined $\pi:=1+i$ as a uniformizer of $F=\QQ_2(i)$ and $s:=1+\pi^3=-1+2i$.

\begin{proposition}[Proposition \ref{prop:calculation 2 ramified D_K=-4}]
    For the curve $Y_b:y^2=x^3+bx$, let $b=2^{v_2(b)}u$ with $u\in\ZZ_2^{\times}$. Then we have
    \[
    \psi_{Y_b,2}(\pi)=\pi i^{(u^2-1)/8},\qquad \psi_{Y_b,2}(i)=i\cdot(-1)^{(u^2-1)/8},\qquad \psi_{Y_b,2}(s)=(-1)\cdot(-i)^{v_2(b)-u+1}.
    \]
\end{proposition}
\begin{proof}
    We again first consider $(x,b)_{4,F}$ for $x=\pi,i,s$. By an explicit formula for the fourth Hilbert symbol, see for example \cite[Chapter VII, \S 4]{fesenko2002local}, one has
    \[
    (\pi,b)_{4,F}=i^{-(u^2-1)/8},\quad (i,b)_{4,F}=(-1)^{-(u^2-1)/8}.\quad (s,b)_{4,F}=i^{v_2(b)-u+1}.
    \]
    Thus, we need only verify the formula for a single curve, again by the local CM comparison theorem (Theorem \ref{thm:CM comparison}).
    
    Consider the curve $Y_{-1}:y^2=x^3-x$ with conductor $2^5$. Since $\pi$ has norm $2$ and $\omega$ is a unit, we easily obtain $\psi_{Y_{-1},2}(\pi)=\pi$ and $\psi_{Y_{-1},2}(i)=i$ by uniformization. For $s=-1+2i$ with norm $5$, while $5=s\bar{s}$ splits in $\calK=\QQ(i)$, we have $\psi_{Y_{-1},2}(s)\psi_{Y_{-1},(s)}(s)=s$. A computation or check in LMFDB database gives $a_5(Y_{-1})=-2$, so $\psi_{Y_{-1},(s)}(s)$ is a root of $x^2+2x+5=0$, namely $s$ or $\bar{s}$. But $\bar{s}$ does not generate $(s)$, so $\psi_{Y_{-1},(s)}(s)=s$. Hence $\psi_{Y_{-1},2}(s)=1$, the same as expected in the formula, taking $b=-1$. Thus the theorem follows.
\end{proof}

The final case is $p=2$ ramified in $\calK=\QQ(\sqrt{-2})$, i.e. the curves in Table \ref{tab:p=2 ramified database D_K=-8}. Here we choose $\pi:=\sqrt{-2}$ to be a uniformizer of $F=\QQ_2(\sqrt{-2})$ and $v:=1+\pi$. After some translation, each curve is of the form $E_{-8}^d$. We choose a square-free integer representative for each $d\in\QQ^{\times}/(\QQ^{\times})^2$.

\begingroup\small\setlength{\tabcolsep}{5pt}\renewcommand{\arraystretch}{1.25}
\begin{longtable}{@{}ccccc@{}}
\caption{Curves in Table \ref{tab:p=2 ramified database D_K=-8} viewed as $E_{-8}^d$}\label{tab:translated p=2 ramified database D_K=-8}\\
\hline
Equation of $E$ & LMFDB label & $d$ & $\psi_{E,2}(\pi)$ & $\psi_{E,2}(v)$\\\hline\endfirsthead
\multicolumn{5}{c}{\tablename\ \thetable\ (continued)}\\\hline
Equation of $E$ & LMFDB label & $d$ & $\psi_{E,2}(\pi)$ & $\psi_{E,2}(v)$\\\hline\endhead
\hline\endfoot
$y^2=x^3-x^2-83x+287$ & 6400.x2 & $5$ & $-\pi$ & $-1$ \\
$y^2=x^3-30x+56$ & 2304.h2 & $3$ & $-\pi$ & $1$ \\
$y^2=x^3+x^2-3x+1$ & 256.a2 & $1$ & $\pi$ & $-1$ \\
$y^2=x^3+x^2-163x+657$ & 12544.c2 & $7$ & $\pi$ & $1$ \\
\end{longtable}\endgroup

By local CM comparison theorem (Theorem \ref{thm:CM comparison}), one has
\[
\psi_{E_{-8}^d,2}(x)=(x,d)_{2,F}^{-1}\psi_{E_{-8}^1,2}(x)
\]
for any $x\in F^{\times}$. But by compatibility of $\rec_F$ in local class field theory, $(x,d)_{2,F}=(\Norm^F_{\QQ_2}(x),d)_{2,\QQ_2}$. Hence Lemma \ref{lem:calculation local norm residue symbol n=p=2} gives
\[
(\pi,d)_{2,F}=(2,d)_{2,\QQ_2}=(-1)^{(u^2-1)/8},\quad (v,d)_{2,F}=(3,d)_{2,\QQ_2}=(-1)^{v_2(d)+(u-1)/2},
\]
where $u=d/2^{v_2(d)}$.

We consider $E_{-8}^1$, which is isomorphic to 256.a2 over $\QQ$, with conductor $2^8$. Since $\Norm^{\calK}_{\QQ}(\pi)=2$, it follows that $\psi_{E_{-8}^1,2}(\pi)=\pi$ by uniformization. For $v=1+\pi$ with norm $3=v\bar{v}$, which splits, we have $\psi_{E_{-8}^1,2}(v)\psi_{E_{-8}^1,(v)}(v)=v$. Since $a_3(E_{-8}^1)=-2$, $\psi_{E_{-8}^1,(v)}(v)$ is a root of $x^2+2x+3=0$, namely $-v$ or $-\bar{v}$. But only $-v$ generates $(v)$, so $\psi_{E_{-8}^1,(v)}(v)=-v$ and $\psi_{E,2}(v)=-1$. Moreover, we obtain the following theorem.

\begin{proposition}
    For $E_{-8}^d$ with $d\in\QQ^{\times}$, one has
    \[
    \psi_{E_{-8}^d,2}(\pi)=\pi\cdot(-1)^{(u^2-1)/8},\qquad \psi_{E_{-8}^d,2}(v)=(-1)^{v_2(d)+(u+1)/2},
    \]
    where $u:=d/2^{v_2(d)}$.
\end{proposition}

Table \ref{tab:p=2 ramified database D_K=-8} and the equivalent Table \ref{tab:translated p=2 ramified database D_K=-8} are thereby verified.

\section{Computations on $p$-Shafarevich-Tate groups}\label{sec:p-Sha calculation}

We first verify $S_{\rm{an}}=1$ for the twisted curve $E_{163}^{2}:y^2=x^3-8697680x+9873093538$ outside the LMFDB database. Then we show that $\Sha(E/\QQ)[p]=0$ for all curves $E$ displayed in the previous sections, Sections \ref{sec:inert BSD 2,3} and \ref{sec:ramified BSD}.

\subsection{Analytic $\Sha$ for $E_{163}^{2}$}

We use \href{https://sagecell.sagemath.org/}{SageMath Cell Server} to do the computation.\\

Codes (SageMath):\\
\begin{codebox}
    E = EllipticCurve(QQ, [0,0,0,-8697680,9873093538])
    
    print('conductor:', E.conductor().factor())
    print('root number:', E.root_number())
    print('analytic rank:', E.analytic_rank(algorithm='pari'))
    l, err = E.lseries().deriv_at1(12000,128)
    print('L\'(E,1):', l)
    reg = E.regulator()
    print('regulator:', reg)
    omega = E.period_lattice().omega(bsd_normalise=True)
    print('period:', omega)
    t = E.torsion_order()
    print('torsion order', t)
    local_Tamagawa=[(p, E.tamagawa_number(p)) for p in E.conductor().prime_divisors()]
    Tamagawa = prod(c for p,c in local_Tamagawa)
    print('local Tamagawa numbers', local_Tamagawa)
    
    S_an = (l*t**2)/(reg*omega*Tamagawa)
    print('analytic Sha:', S_an)
\end{codebox}

Outputs:\\
\begin{codebox}
    conductor: 2^6 * 163^2
    root number: -1
    analytic rank: 1
    L'(E,1): 6.7330958470029679499518962105133984991
    regulator: 5.99890930721786
    period: 0.561193335503651
    torsion order 1
    local Tamagawa numbers [(2, 1), (163, 2)]
    analytic Sha: 1.00000000000000
\end{codebox}

This shows that $S_{\rm an}$ for $E_{163}^{2}:y^2=x^3-8697680x+9873093538$ is indeed $1$.

\subsection{$2,3$-descent via Magma}

In the following, we show that $\Sha(E/\QQ)[p]=0$ for the displayed curves $E$. When $p=2,3$, i.e. for the curves in Tables \ref{tab:p_3 inert database}, \ref{tab:p_2 inert database}, \ref{tab:p=3 ramified database}, \ref{tab:p=2 ramified database D_K=-4} and \ref{tab:p=2 ramified database D_K=-8}, we use the \href{https://magma.maths.usyd.edu.au/calc/}{Magma calculator}\footnote{The algorithm of $2,3$-descent in SageMath is essentially from Magma.} and the exact sequence
\[
0\rightarrow E(\QQ)/pE(\QQ)\rightarrow\Sel_p(E/\QQ)\rightarrow\Sha(E/\QQ)[p]\rightarrow0.
\]
For instance, the code and output below show $S_3=1$ for $y^2=x^3-2x$ in Table \ref{tab:p_3 inert database}, and $S_2=1$ for $y^2=x^3-x^2-83x+287$ in Table \ref{tab:p=2 ramified database D_K=-8}.\\

Codes (Magma):\\
\begin{codebox}
    E_1 := EllipticCurve([ Rationals() | 0,0,0,-2,0]);
    print "E_1:", E_1;
    
    M := MordellWeilGroup(E_1);
    print "Mordell-Weil group E_1(Q):", M, "";
    
    M_3 := quo<M | 3*M>;
    print "E_1(Q)/3E_1(Q) has order", #M_3, "\n";
    
    Sel_3 := ThreeSelmerGroup(E_1);
    print "3-Selmer group of E_1:", Sel_3, "";
    
    S_3 := #Sel_3/#M_3;
    print "#Sha(E_1/Q)[3] =", S_3, "\n";
    
    E_2 := EllipticCurve([ Rationals() | 0,-1,0,-83,287]);
    print "E_2:", E_2;
    
    N := MordellWeilGroup(E_2);
    print "Mordell-Weil group E_2(Q):", N, "";
    
    N_2 := quo<N | 2*N>;
    print "E_2(Q)/2E_2(Q) has order", #N_2, "\n";
    
    Sel_2 := TwoSelmerGroup(E_2);
    print "2-Selmer group of E_2:", Sel_2, "";
    
    S_2 := #Sel_2/#N_2;
    print "#Sha(E_2/Q)[2] =", S_2;
\end{codebox}

Outputs:\\
\begin{codebox}
    E_1: Elliptic Curve defined by y^2 = x^3 - 2*x over Rational Field
    Mordell-Weil group E_1(Q):
    Abelian Group isomorphic to Z/2 + Z
    Defined on 2 generators
    Relations:
        2*M.1 = 0
    
    E_1(Q)/3E_1(Q) has order 3
    
    3-Selmer group of E_1:
    Abelian Group isomorphic to Z/3
    Defined on 1 generator
    Relations:
        3*Sel_3.1 = 0
    
    #Sha(E_1/Q)[3] = 1
    
    E_2: Elliptic Curve defined by y^2 = x^3 - x^2 - 83*x + 287 over Rational Field
    Mordell-Weil group E_2(Q):
    Abelian Group isomorphic to Z/2 + Z
    Defined on 2 generators
    Relations:
        2*N.1 = 0
    
    E_2(Q)/2E_2(Q) has order 4
    
    2-Selmer group of E_2:
    Abelian Group isomorphic to Z/2 + Z/2
    Defined on 2 generators
    Relations:
        2*Sel_2.1 = 0
        2*Sel_2.2 = 0
    
    #Sha(E_2/Q)[2] = 1
\end{codebox}

Consequently, we checked all the curves in Tables \ref{tab:p_3 inert database}, \ref{tab:p_2 inert database}, \ref{tab:p=3 ramified database}, \ref{tab:p=2 ramified database D_K=-4} and \ref{tab:p=2 ramified database D_K=-8}, and concluded that $S_2,S_3=1$ in the tables above. The remaining cases are exactly those with $p>3$ ramified in the CM field $\calK$, namely the curves in Table \ref{tab:p>3 ramified database}.

\subsection{$p$-isogeny descent for ramified $p>3$}

For $p>3$ ramified in $\calK=\QQ(\sqrt{-p})$ and the curves in Table \ref{tab:p>3 ramified database}, \cite[Theorem 1.1]{creutz2012second} first shows that the full BSD conjecture holds for elliptic curves $E/\QQ$ of conductor $N_E<5000$ with analytic rank $\leq1$; thus, $6$ such curves in Table \ref{tab:p>3 ramified database} are verified. For each remaining $p=11,19,43,67,163$, we denote the curve with smaller conductor by $E_p$. These curves have CM by the maximal order $\mathcal{O}_\calK$. The curves with larger conductor are quadratic twists of $E_p$, which we denote by $E_p^d$, with $d=2,2,5,2,2$, respectively.

\begingroup\small\setlength{\tabcolsep}{10pt}\renewcommand{\arraystretch}{1.25}
\begin{longtable}{@{\hspace{5pt}}cccc@{\hspace{5pt}}}
\caption{Remaining curves}\label{tab:p>3 ramified remaining databse}\\
\hline
$p$ & Equation of $E_p$ & LMFDB label & $d$\\\hline\endfirsthead
\multicolumn{4}{c}{\tablename\ \thetable\ (continued)}\\\hline
$p$ & Equation of $E_p$ & LMFDB label & $d$\\\hline\endhead
\hline\endfoot
    $11$ & $y^2+y=x^3-x^2-7x+10$ & 121.b2 & $2$\\
    $19$ & $y^2+y=x^3-38x+90$ & 361.a2 & $2$\\
    $43$ & $y^2+y=x^3-860x+9707$ & 1849.b2 & $5$\\
    $67$ & $y^2+y=x^3-7370x+243528$ & 4489.b2 & $2$\\
    $163$ & $y^2+y=x^3-2174420x+1234136692$ & 26569.a2 & $2$\\
\end{longtable}\endgroup

\begin{theorem}[{\cite[Corollary 7.3]{miller2013explicit}}]\label{thm:ramified Sha upper bound}
    Let $E/\QQ$ be an elliptic curve with a rational isogeny $\phi$ of prime degree $p>3$, and assume that $E[\phi]\subseteq E(\RR)$. Let $M:=\QQ(E[\phi])$ be the splitting field of $E[\phi]$. If $E$ does not have split multiplicative reduction anywhere, and $M$ has class number $h_M$ prime to $p$, then
    \[
    \rank(E/\QQ)+\dim_{\FF_p}\Sha(E/\QQ)[p]\leq 2(1-\dim_{\FF_p}E(\QQ)[\phi])-w,
    \]
    where $w\in\{0,1\}$ is an integer related to $\phi$ and the invariant differential.
\end{theorem}

As an example, the author proved that $\Sha(E_{163}/\QQ)[163]=0$ in \cite[Example 7.4]{miller2013explicit}. So the remaining curves are essentially the twisted ones. But here we may explain more. Let $\alpha=\sqrt{-p}$ and $\mathfrak{p}:=(\alpha)$ be the only prime above $p$. Then $[\alpha]\in\End_{\bar{\QQ}}(E)$ factors through a $p$-isogeny $\phi_{\alpha}:E\rightarrow E^{\prime}:=E/E[\mathfrak{p}]$, which induces $[\bar{\alpha}]:E'\rightarrow E$, making $[\alpha]=[\bar{\alpha}]\circ\phi_{\alpha}$.

\vspace{0.5em}
\begin{center}
    \begin{tikzcd}[row sep=large]
        E \arrow[rr, "{[\alpha]}"] \arrow[rd, "\phi_{\alpha}"'] & & E \\
        & {E^{\prime}:=E/E[\mathfrak{p}]} \arrow[ru, "{[\bar{\alpha}]}"'] &  
    \end{tikzcd}
\end{center}
\vspace{0.5em}

Since $[\alpha]$ is defined over $\calK$, any $\sigma\in\Gal_{\QQ}$ acts on $[\alpha]$ by $\pm1$, depending on the image in $\Gal(\calK/\QQ)$. It follows that
\[
\sigma E[\mathfrak{p}]=\sigma\ker([\alpha])=\ker(\sigma[\alpha])=\ker([\pm\alpha])=\ker[\alpha]=E[\mathfrak{p}],
\]
for any $\sigma\in\Gal_{\QQ}$, i.e. $E[\mathfrak{p}]$ is $\Gal_{\QQ}$-stable. Hence $\phi_p$ and $E^{\prime}$ are both defined over $\QQ$. Note that $[\bar{\alpha}]:E^{\prime}\rightarrow E$ is only a $\calK$-isomorphism.

\begin{proposition}\label{prop:real kernel}
    For curves $E_p$ and $E_p^d$ in the Table \ref{tab:p>3 ramified remaining databse}, we have $E_p[\mathfrak{p}]\subset E_p(\RR)$ and $E_p^d[\mathfrak{p}]\subset E_p^d(\RR)$.
\end{proposition}
\begin{proof}
    First we need the following lemma.
    \begin{lemma}\label{lem:real kernel equiv}
        Let $E$ be an elliptic curve given by $y^2=4x^3-g_2x-g_3$ with $g_2,g_3\in\RR$. Let $\Lambda$ be the period lattice satisfying $g_2(\Lambda)=g_2$ and $g_3(\Lambda)=g_3$, with complex uniformization
        \[
        u:\CC/\Lambda\rightarrow E(\CC)
        \]
        sending $z$ to $(\wp_{\Lambda}(z),\wp^{\prime}_{\Lambda}(z))$. If $\gamma\in i\RR^{\times}$ satisfies $\gamma\Lambda\subset\Lambda$, then
        \[
        \ker([\gamma])\subset E(\RR)\iff w+\bar{w}\in\gamma\Lambda,~\forall w\in\Lambda.
        \]
    \end{lemma}
    \begin{proof}
        Since $g_2(\Lambda),g_3(\Lambda)\in\RR$, we have $\bar{\Lambda}=\Lambda$ and thus $u(\bar{z})=\overline{u(z)}$. Hence $u(z)\in E(\RR)\iff z-\bar{z}\in\Lambda$. Note that $\ker([\gamma])=u(\gamma^{-1}\Lambda/\Lambda)$. It follows that
        \begin{align*}
            \ker([\gamma])\subset E(\RR)&\iff u(z)\in E(\RR),\forall z\in\gamma^{-1}\Lambda\\
            &\iff u(\gamma^{-1}w)\in E(\RR),\forall w\in\Lambda\\
            &\iff w+\bar{w}\in\gamma\Lambda,\forall w\in\Lambda.
        \end{align*}
    \end{proof}
    For those CM curves in Table \ref{tab:p>3 ramified remaining databse}, $E_p$ can be transformed over $\RR$ into a model as in Lemma \ref{lem:real kernel equiv}, with $g_2=c_4/12$ and $g_3=c_6/216$. The corresponding $\Lambda$ satisfies $\Lambda=\lambda\mathcal{O}_\calK$ for some $\lambda\in\CC^{\times}$, since $h_{\calK}=1$. We claim that $\lambda\in i\RR^{\times}$.

    \begin{lemma}
        Viewing $\mathcal{O}_\calK$ as a lattice, we have $g_2(\mathcal{O}_\calK),g_3(\mathcal{O}_\calK)\in\RR_+$ for $\calK=\QQ(\sqrt{-p})$ with $p\in\{11,19,43,67,163\}$.
    \end{lemma}
    \begin{proof}
        Let $\tau:=\frac{1+\sqrt{-p}}{2}$, so that $\mathcal{O}_\calK=\ZZ+\ZZ(\tau)$. Put $q=e^{2\pi i\tau}=-e^{-\pi\sqrt{p}}\in\RR$. The Eisenstein series gives
        \[
        g_2(\mathcal{O}_\calK)=\frac{4\pi^4}{3}E_4(\tau),\quad E_4(\tau)=1+240\sum_{n\geq1}\sigma_3(n)q^n
        \]
        and
        \[
        g_3(\mathcal{O}_\calK)=\frac{8\pi^6}{27}E_6(\tau),\quad E_6(\tau)=1-504\sum_{n\geq1}\sigma_5(n)q^n.
        \]
        Using $\sigma_k(n)\leq n^{k+1}\leq 2^{(n-1)(k+1)}$, we obtain that
        \[
        \abs{E_4(\tau)-1}\leq 240\abs{q}\sum_{n\geq1}\abs{16q}^{n-1}=\frac{240\abs{q}}{1-16\abs{q}}
        \]
        and
        \[
        \abs{E_6(\tau)-1}\leq 504\abs{q}\sum_{n\geq1}\abs{64q}^{n-1}=\frac{504\abs{q}}{1-64\abs{q}}.
        \]
        Since $\abs{q}<\frac{1}{10000}$ for $p\geq11$, we conclude that $\abs{E_4(\tau)-1}<1$ and that $\abs{E_6(\tau)-1}<1$. Therefore $g_2(\mathcal{O}_\calK),g_3(\mathcal{O}_\calK)\in\RR_+$.
    \end{proof}

    Recall that $c_4/12=g_2(\Lambda)=\lambda^{-4}g_2(\mathcal{O}_\calK)$ and that $c_6/216=g_3(\Lambda)=\lambda^{-6}g_3(\mathcal{O}_\calK)$.
    Since the curves $E_p$ with the Weierstrass equations in Table \ref{tab:p>3 ramified remaining databse} all satisfy $c_4\in\RR_+$ and $c_6\in\RR_-$, it follows that $\lambda^{-4}\in\RR_+$ and $\lambda^{-6}\in\RR_-$. This forces $\lambda\in i\RR^{\times}$.

    Then we can use the criterion in Lemma \ref{lem:real kernel equiv}. For any $w\in\Lambda=\lambda\mathcal{O}_\calK$, write $w=\lambda(a+b\tau)$ with $a,b\in\ZZ$. Then
    \[
    w+\bar{w}=\lambda(a+b\tau)-\lambda(a+b\bar{\tau})=\lambda b(\tau-\bar{\tau})=\alpha\lambda b\in\alpha\Lambda.
    \]
    Therefore $E_p[\mathfrak{p}]\subset E(\RR)$ for $E_p$ in Table \ref{tab:p>3 ramified remaining databse}.

    Finally for twists $E_p^d$ with $d>0$, there is an $\RR$-isomorphism $E_p\rightarrow E_p^d$ sending $(x,y)$ to $(x/d,y/d^{3/2})$. Thus, the same result follows as for $E_p$. It can also be checked directly via $c_4(E_p^d)=d^2c_4(E_p)$ and $c_6(E_p^d)=d^3c_6(E_p)$. This completes the proof.
\end{proof}

By Proposition \ref{prop:real kernel} and after putting $\phi=\phi_{\alpha}$, it suffices to show that $p\nmid h_M$ for $E=E_p$ or $E_p^d$; then
\[
\dim_{\FF_p}\Sha(E/\QQ)[p]\leq1-w,
\]
since $E/\QQ$ has rank $1$. Recall that $\dim_{\FF_p}\Sha(E/\QQ)[p]$ should be even by the parity lemma. Therefore we can conclude that $\Sha(E/\QQ)[p]=0$ as expected. In the remainder, we prove that $p\nmid h_M$ for $M=M(E):=\QQ(E[\phi_{\alpha}])$ with $E=E_p$ or $E_p^d$. Let $\QQ(\zeta_p)^+$ denote the maximal totally real subfield of $\QQ(\zeta_p)$, namely $\QQ(t_p)$ with $t_p=\zeta_p+\zeta_p^{-1}$.

\begin{proposition}
    We have $M(E_p)\subseteq\QQ(\zeta_p)^+$ and $M(E_p^d)\subseteq\QQ(\zeta_p)^+(\sqrt{d})$\footnote{In fact, one can show that both equalities hold.}.
\end{proposition}
\begin{proof}
    Since $E_p$ has conductor $p^2$, it follows that $M(E_p)/\QQ$ is only ramified at $p$. Note that $\Gal(M(E_p)/\QQ)$ injects into $\Aut(E[\mathfrak{p}])\cong\FF_p^{\times}$, so $[M(E_p):\QQ]$ divides $p-1$. By Kronecker-Weber we know $M(E_p)\subseteq\QQ(\zeta_p)$. Moreover $M(E_p)$ is totally real since $E_p[\mathfrak{p}]\subset E(\RR)$. Therefore $M(E_p)\subseteq\QQ(\zeta_p)^+$.

    Let $\lambda:\Gal_{\QQ}\rightarrow \Aut(E_p[\mathfrak{p}])\cong\FF_p^{\times}$ be the character on $E_p[\mathfrak{p}]$. Then the corresponding action on $E_p^d$ is just $\lambda\chi_d$, where $\chi_d$ is the quadratic character associated to $\QQ(\sqrt{d})/\QQ$, with kernel $\Gal_{\QQ(\sqrt{d})}$. The assertion above shows that $\Gal_{\QQ(\zeta_p)^+}\subseteq\ker(\lambda)$. Since $p\nmid 2d$ in our cases, it follows that $\QQ(\zeta_p)^+$ and $\QQ(\sqrt{d})$ are linearly disjoint. Therefore $\Gal_{\QQ(\zeta_p)^+(\sqrt{d})}=\Gal_{\QQ(\zeta_p)^+}\cap\Gal_{\QQ(\sqrt{d})}\subseteq\ker(\lambda\chi_d)$. We conclude that $M(E_p^d)=\bar{\QQ}^{\ker(\lambda\chi_d)}\subseteq\QQ(\zeta_p)^+(\sqrt{d})$.
\end{proof}

\begin{lemma}\label{lem:class groups injection}
    Let $E/F$ be any Galois extension of number fields of degree $n$. If $p\nmid n$, then the extension map $\Cl_F[p^{\infty}]\rightarrow\Cl_E[p^{\infty}]$ is injective.
\end{lemma}
\begin{proof}
    Let $j:\Cl_F\rightarrow\Cl_E$ be the extension map and $\Norm^E_F:\Cl_E\rightarrow\Cl_F$ be the norm map. It is easy to see that $\Norm^E_F\circ j$ is $[n]:\Cl_F\rightarrow\Cl_F$. Since $p\nmid n$, it follows that $[n]$ is an automorphism on $\Cl_F[p^{\infty}]$. Therefore $j$ is injective on $p$-primary parts.
\end{proof}

Let $M_p:=\QQ(\zeta_p)^+$ and $M_p^d:=\QQ(\zeta_p)^+(\sqrt{d})$. Since $p\nmid [M_p:M(E_p)]$ and $p\nmid [M_p^d:M(E_p^d)]$, in order to show $p\nmid h_{M(E_p)}$ and $p\nmid h_{M(E_p^d)}$, by Lemma \ref{lem:class groups injection}, it suffices to show that $p\nmid h_{M_p}$ and that $p\nmid h_{M_p^d}$.

\begin{theorem}\label{thm:class number indivisibility}
    Both $p\nmid h_{M_p}$ and $p\nmid h_{M_p^d}$ hold for $(p,d)$ in Table \ref{tab:p>3 ramified remaining databse}.
\end{theorem}
\begin{proof}
    By the verification of Kummer-Vandiver conjecture for $p\leq163577856$ in \cite{buhler2011irregular}, we see that $p\nmid h_{M_p}$ holds for all $p$ in Table \ref{tab:p>3 ramified remaining databse}. Note that $M_p^2$ is contained in $\QQ(\zeta_{8p})^+$ when $p\neq2$. By \cite[Theorem 1.1]{miller2014class} we have $h_{\QQ(\zeta_{8p})^+}=1$ for $p=11,19$. Since $[\QQ(\zeta_{8p})^+:\QQ]=\phi(8p)/2=2(p-1)$ is prime to $p$, we deduce that $p\nmid h_{M_{11}^2}$ and $p\nmid h_{M_{19}^2}$ by Lemma \ref{lem:class groups injection}.

    The remaining cases are $(p,d)=(43,5),(67,2),(163,2)$. We consider the field $M_n:=\QQ(\zeta_n)^+$. By \cite[Proposition 2.9, Equation 2.7]{hakkarainen2009computation} one has the decomposition
    \begin{equation}\label{eq:class group decomposition}
        h_{F,p}=\prod_{\tilde{\chi}}h_{\chi,p}
    \end{equation}
    for any real abelian number field $F$ with degree prime to $p$, where $\tilde{\chi}$ is taken over all rational orbits of complex characters of $\Gal(F/\QQ)$, and $h_{\chi,p}$ means the $p$-part of the $\chi$-class number $h_{\chi}$. In particular, $h_{\chi}$ only depends on $\chi$ as a Dirichlet character, via some embedding $F\hookrightarrow \QQ(\zeta_{f_F})$ into a cyclotomic field, not on the field $F$ itself. In \cite[Chapter 8, Table 1]{hakkarainen2009computation}, the author computed all possible prime divisors $p<10000$ for $h_{\chi}$ with $\chi$ of conductor $f_{\chi}\leq2000$.

    Put $D=5$ when $d=5$ and $D=8$ when $d=2$, so that $\QQ(\sqrt{d})\subset\QQ(\zeta_D)$. Recall that $\Gal(M_p^d/\QQ)\cong\Gal(\QQ(\zeta_p)^+/\QQ)\times\Gal(\QQ(\sqrt{d})/\QQ)\cong((\ZZ/p\ZZ)^{\times}/\{\pm1\})\times C_2$. Let $g$ be a generator of $(\ZZ/p\ZZ)^{\times}/\{\pm1\}$, $\theta$ a character sending $g$ to $\zeta_n$, where $n:=[\QQ(\zeta_p)^+:\QQ]=\frac{p-1}{2}$, and $\eta$ the non-trivial character of $\Gal(\QQ(\sqrt{d})/\QQ)$. Thus, all characters $\chi$ of $\Gal(M_p^d/\QQ)$ are of the form $\theta^a\eta^b$ with $0\leq a\leq n$ and $b\in\{0,1\}$. We view $\chi$ as an even Dirichlet character via $M_p^d\subset\QQ(\zeta_{pD})^+$.
    
    All possibilities for $f_{\chi}$ are $1,p,D,pD$, by an easy check. From \cite[Chapter 8, Table 1]{hakkarainen2009computation} we see that $p\nmid h_{\chi}$ when $f_{\chi}=1,pD$. When $f_{\chi}=D=5,8$, one can show that $\chi$ also appears in the decomposition \eqref{eq:class group decomposition} for $F=\QQ(\sqrt{D})$. But the Minkowski bound gives $h_{\QQ(\sqrt{D})}=1$, so such $h_{\chi}$ must equal $1$. When $f_{\chi}=p$, recall that verification of Vandiver's conjecture \cite{buhler2011irregular} again gives $h_{\QQ(\zeta_p)^+,p}=1$. It follows that such $h_{\chi,p}$ must also be $1$, using \eqref{eq:class group decomposition} for $F=\QQ(\zeta_p)^+$. Therefore $h_{M_p^d,p}=1$ for the pairs $(p,d)=(43,5),(67,2),(163,2)$, and all the desired indivisibilities follow.
\end{proof}

Combining Proposition \ref{prop:real kernel}
and Theorem \ref{thm:class number indivisibility}, we can apply the main input, Theorem \ref{thm:ramified Sha upper bound}, to obtain the upper bound $\dim_{\FF_p}\Sha(E/\QQ)[p]\leq1$. By the parity lemma, we conclude that $\Sha(E/\QQ)[p]=0$ for $E$ in Table \ref{tab:p>3 ramified database}.